\documentclass[11pt,a4paper]{amsart}
\usepackage[colorlinks =  true, 
            citecolor  =  red, 
            linkcolor  =  blue, 
            urlcolor   =  blue]{hyperref}

\usepackage[top    =  20mm, 
            bottom =  20mm, 
            right  =  20mm, 
            left   =  20mm]{geometry}
            
\usepackage[normalem]{ulem}
\usepackage{amsmath, amsthm, amssymb, bm, color}
\usepackage[capitalize,noabbrev]{cleveref}
\usepackage{float}
\usepackage{tikz,tikz-cd,tikz-3dplot}
\usepackage{dsfont}
\usepackage{mathrsfs}
\usepackage{pgfplots}
\usepackage{booktabs}
\usepackage{lineno}
\usepackage{comment}

\usetikzlibrary{patterns}
\usepgfplotslibrary{colormaps}
\pgfplotsset{compat=1.18} 
\pgfdeclarelayer{pre main}

\numberwithin{equation}{section}

\newtheorem{theorem}{Theorem}[section]
\newtheorem*{theorem*}{Theorem}
\newtheorem{proposition}[theorem]{Proposition}
\newtheorem{lemma}[theorem]{Lemma}
\newtheorem{corollary}[theorem]{Corollary}

\theoremstyle{definition}
\newtheorem{definition}[theorem]{Definition}
\newtheorem{example}[theorem]{Example}
\newtheorem{remark}[theorem]{Remark}
\newtheorem{notation}[theorem]{Notation}

\newtheorem{maintheorem}{Theorem}[section]

\definecolor{darkorange}{rgb}{1.0, 0.55, 0.0}
\definecolor{darkblue}{rgb}{0.0, 0.0, 0.55}
\definecolor{darkgreen}{rgb}{0.0, 0.2, 0.13}
\definecolor{darkred}{rgb}{0.75, 0.0, 0.0}

\newcommand{\CC}{\mathbb{C}}
\newcommand{\RR}{\mathbb{R}}
\newcommand{\RRi}{\overline{\mathbb{R}}}
\newcommand{\QQ}{\mathbb{Q}}
\newcommand{\ZZ}{\mathbb{Z}}

\newcommand{\PP}{\mathbb{P}}

\newcommand{\suchthat}{\;\ifnum\currentgrouptype=16 \middle\fi|\;}

\newcommand{\an}{{\mathrm{an}}}

\newcommand{\Pscr}{\mathscr{P}}

\newcommand{\cone}{\mathrm{cone}}
\newcommand{\relint}{\mathrm{relint}}

\newcommand{\sh}{\mathrm{sh}}

\newcommand{\ang}{\mathrm{ang}}
\renewcommand{\star}{\mathrm{star}}
\renewcommand{\tilde}{\widetilde}

\newcommand{\T}{\mathbb{T}}

\DeclareMathOperator{\Mscr}{\mathscr{M}}

\DeclareMathOperator{\conv}{conv}

\DeclareMathOperator{\rec}{rec}
\DeclareMathOperator{\trop}{trop}

\DeclareMathOperator{\val}{val}
\DeclareMathOperator{\Spec}{Spec}

\newcommand{\norm}[1]{\left\lVert#1\right\rVert}

\newcommand{\Ocal}{\mathcal{O}}

\newcommand{\scrP}{\mathscr{P}}

\tikzset{nvert/.style={draw,circle,fill=black,minimum size=3pt,inner sep=0pt}}

\title[Polar decompositions of complex varieties and tropical homotopy theory]
{Polar decompositions of complex varieties \\ and tropical homotopy theory}%

\author{Yassine {Elmaazouz}}
\address{Yassine {Elmaazouz} (Caltech)}
\email{maazouz@caltech.edu}

\author{Paul Alexander Helminck}
\address{Paul Alexander Helminck (IMPAN Warsaw, Instytut Matematyczny Polskiej Akademii Nauk)}
\email{phelminck@impan.pl}


\date{\today}

\begin{document}
\begin{abstract}
We take the first steps towards developing a tropical homotopy theory, which we use to construct explicit finite CW complex structures on smooth complex varieties.
Our approach decomposes a variety into a set of radial and angular parts, a construction we call \emph{polar decomposition}.
We derive a combinatorial formula for the $f$-vector of the resulting CW complex in terms of a weighted sum of local angle $f$-vectors, and we determine these $f$-vectors explicitly for toric varieties and complete intersections in terms of Stirling numbers.

To construct these decompositions, we introduce the notion of a \emph{rational tropical manifold}, which we view as a tropical manifold with a non-reduced structure. 
The local topology of these is dictated by Kummer coverings of affine linear spaces in a complex torus, as opposed to the affine linear spaces in the theory of tropical manifolds. 
Examples of such manifolds include the Fermat and Dwork hypersurfaces.

Using Kato-Nakayama spaces, we associate a diagram of topological spaces to any rational tropical manifold. We then construct a  cofibrant replacement of this diagram using \emph{sealed Kato-Nakayama spaces} to show that the homotopy colimit is homotopy equivalent to a general nearby fiber.
To construct CW complex structures, we then show that the spaces in the diagram can be replaced by their angle sets.
More precisely, for any Kummer covering of a very affine linear space, we prove that the angle map is a homotopy equivalence, extending earlier work of Salvetti and Bj\"{o}rner--Ziegler. 
We construct compatible cellulations for these angle spaces and show that these give rise to an explicit CW complex structure on the nearby fibers.
For complete intersections, we make these spaces and all the CW complex structures completely explicit using \emph{refined Coxeter arrangements}.
\end{abstract}

\maketitle

\setcounter{tocdepth}{1}
\tableofcontents

\section{Introduction}

Tropical geometry can be seen as the study of one-parameter logarithmic limit sets of algebraic varieties. In this process, we discard the local angular data associated to the degeneration. For many applications in algebraic geometry, however, this angular data is reintroduced in various guises.
For instance, in tropical intersection theory and the study of tropical Gromov-Witten invariants, one introduces local intersection multiplicities as in \cite{Mikhalkin05,HuberSturmfels1995,MS15}\footnote{The local intersection multiplicities essentially arise from the index of a sum of integral tangent spaces inside the ambient lattice. In other words, we are calculating the determinant of an integral matrix. This determinant exactly gives the cardinality of the local angle space. See Section~\ref{sec:ClassicalExamplesD01} for more details.}, in tropical homology one attaches Orlik-Solomon algebras to the cells of the tropicalization as in \cite{IKMZ19,GS23}, and for tropical K3 surfaces, one introduces the notion of a structure of an affine integral manifold with singularities as in \cite{KontsevichSoibelman,OO21,Yamamoto22}. We will see in this paper that these are all forms of purely angular data. We can thus view the above tropical results as the outcome of a \emph{successful separation of the radial and angular parts of a degeneration}. Next to finding the individual building blocks, this can be seen as one of the main goals of tropical geometry.

 In this paper we will prove a homotopical separation of the radial and angular parts for the one-parameter degeneration associated to a \emph{rational tropical manifold}. 
 The latter generalizes the notion of a tropical manifold used in \cite{IKMZ19}. For instance, it applies to Fermat and Dwork hypersurfaces, see Section~\ref{sec:FermatAndDwork} for more details. 
 The main goal of this paper will be to make the following slogan precise:
 \begin{center}
\begin{figure}[ht]
    \begin{displaymath}
        \begin{tikzcd}
            \boxed{\text{Complex variety } X} \hspace{0.8cm} \simeq & \boxed{\text{Angular part }\Theta_{\sigma}}\hspace{0cm}
              \hspace{0.3cm} \times \hspace{0.3cm}
            \boxed{\text{Radial part }{\Sigma}} \quad.
            \end{tikzcd}
    \end{displaymath}
\end{figure}
\end{center}
\vspace{-0.3cm}
See Theorem~\ref{thm:MainTheoremPolarDecomposition} for a precise statement. 
We call this a \emph{polar decomposition}, in analogy with the standard polar-coordinate decomposition $\CC^{\times}=S^{1}\times \RR_{>0}$. 
 Our construction of this decomposition gives rise to an explicit finite CW complex that is homotopy equivalent to a nearby fiber of an associated semistable degeneration. We give explicit formulas for the $f$-vectors arising from this construction, which allow us to see the topological complexity from a combinatorial point of view. This CW complex structure also allows us to define tropical fundamental groups and integral tropical homology groups. Overall, this paper can be seen as a first step towards constructing a tropical homotopy theory.
 
\begin{figure}[ht]
\centering
\includegraphics[scale=0.44]{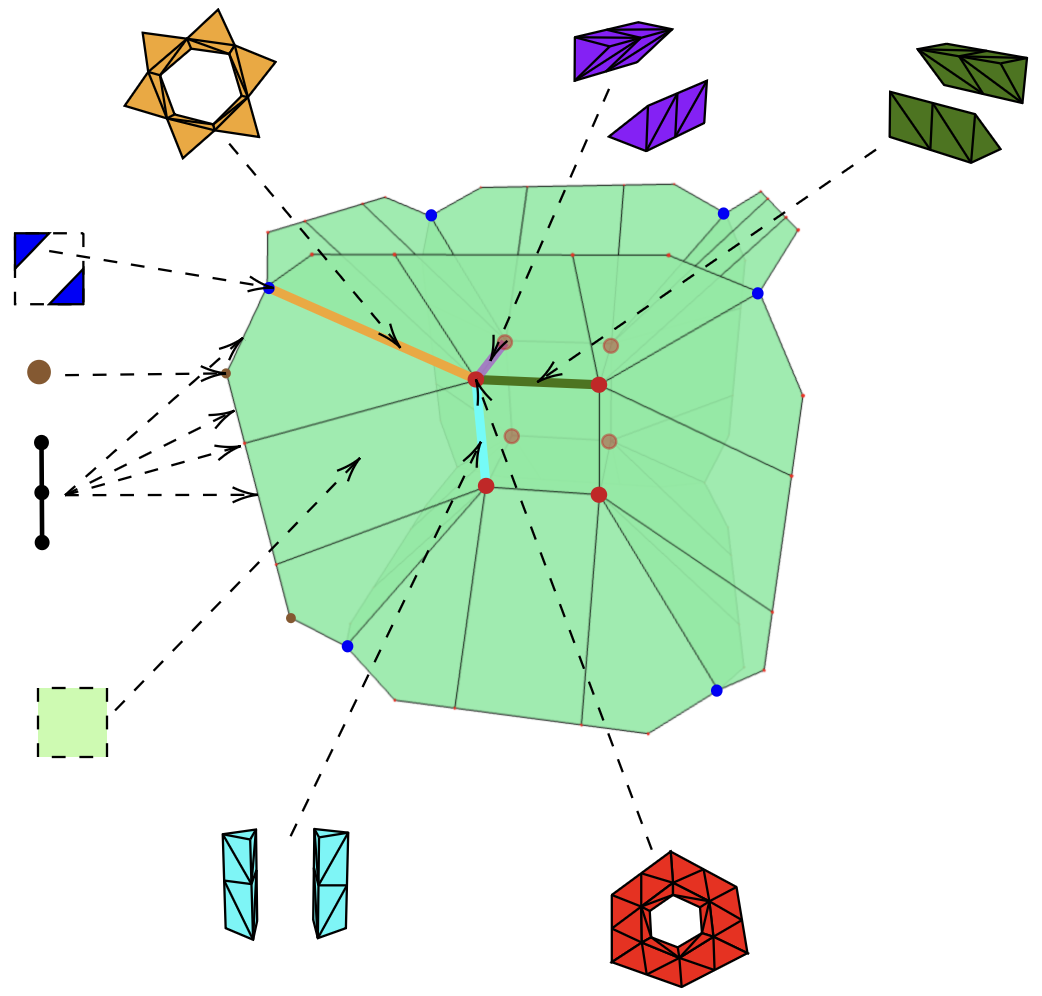}
\caption{\label{fig:AperyK33} The diagram of angle spaces associated to the Ap\'{e}ry-Fermi pencil. By gluing these angle spaces together over the tropicalization, we obtain an explicit CW complex structure on the general fiber of the associated compactified family, which is a K3 surface. 
The angle spaces here are drawn in a fundamental domain of $\RR^d / 2\ZZ^d$ with $d = 0, 1, 2, 3$.
}
\end{figure}

 \begin{example}
     The construction for a K3 surface is summarized in Figure
    ~\ref{fig:AperyK33}. This arises from the Apéry-Fermi pencil given by the local toric equation $f=1+t(x+y+z+x^{-1}+y^{-1}+z^{-1})$. In this case, our $f$-vector formula yields
 \begin{equation*}
\mathfrak{f}(\Delta(\mathcal{F}^{\triangle})) = \big(656, 3960, 9136, 8976, 3168\big),
 \end{equation*}
 with Euler characteristic $\chi = 24$. 
 The cells in this CW complex structure arise from gluing the local
cellulations of the angle spaces over the tropicalization. For this example and various connections to the SYZ program, see Section~\ref{sec:CYMonodromy}.\footnote{Overall, we will use our techniques to obtain three different CW complex structures on K3 surfaces. These arise from the degree-$4$ Fermat and Dwork hypersurfaces, and the Ap\'{e}ry-Fermi pencil. See Examples~\ref{exa:FermatK3Example} and \ref{exa:DworkFamilyK3}, and Section~\ref{sec:CYMonodromy}. } 

 \end{example}

\subsection{Summary of the main results}

We now summarize the main results in this paper. These will be discussed in more detail in Sections~\ref{sec:IntroRTM}--\ref{sec:IntroFVector}.

\begin{enumerate}
    \item In Theorem~\ref{thm:MainTheoremHCKNFunctor}, we prove a homotopy colimit theorem that describes the topology of a compact variety in terms of local Kato-Nakayama spaces. This generalizes results by Danilov and Welker--Ziegler--Živaljevic in the case of toric varieties \cite{Danilov78,WZZ99}.
    \item In Theorem~\ref{thm:MainTheoremKummerCoveringAngleHE}, we prove that the angle map associated to a Kummer covering of a very affine linear space is a homotopy equivalence onto its image. This extends results of Salvetti \cite{Salvetti87} and Björner--Ziegler \cite{BjornerZiegler92}. 
    Combined with Theorem~\ref{thm:MainTheoremHCKNFunctor}, this gives rise to a polar decomposition of a smooth complex variety, see Theorem~\ref{thm:MainTheoremPolarDecomposition}.
    \item We construct CW complex structures on the homotopy colimit of the angle functor and we prove a formula for its $f$-vector in terms of local angle $f$-vectors. We make these explicit for complete intersections in terms of refined Coxeter arrangements. See Theorem~\ref{mainthm:Fermat_FVector} for the case of Fermat hypersurfaces.  
\end{enumerate}

In all of these theorems, our main examples are associated to \emph{rational tropical manifolds}, which are central to our polar decompositions. We explain this concept in more detail in the next section.

\subsection{Rational tropical manifolds}\label{sec:IntroRTM} 
The main tropical object we introduce in this paper is a rational tropical manifold $(X,\Sigma,D)$, which we shall now explain. Let $M$ be a lattice and $N={\rm Hom}(M,\ZZ)$ its dual. We write $X = V(I)\subset \Spec(K[M])$ for a closed subscheme cut out by an ideal $I$. Here $K=\CC(t)\subset \CC(\!(t)\!)$ with the $t$-adic valuation. We write $\Sigma$ for a choice of a polyhedral structure on $\trop(X) \subset N_{\RR}$, and $D$ for a choice of boundary components in the toric variety $Y(\Delta)$ corresponding to the recession fan $\Delta$ of $\Sigma$. 
  
 A tropical manifold is locally the Bergman fan of a matroid, see \cite{IKMZ19}. In the embedded realizable case, this means that the initial degenerations are very affine linear spaces $Z\subset \Spec(\CC[M])$. In other words, there is a basis $m_{i}$ of $M$ such that $Z$ is cut out by affine linear equations $\sum_{i=1}^{n}c_{i}\chi^{m_{i}}$ with respect to this basis, see Definition~\ref{def:TropicalSmoothness}. We generalize this notion here by allowing Kummer coverings of a very affine linear space $Z$. To define this, consider a lattice $\Lambda$ of finite index in $M$. We then have a map of tori
 \[
    \kappa_{M,\Lambda}:\Spec(\CC[M])\to \Spec(\CC[\Lambda]).
 \]
 We say that $X\subset \Spec(\CC[M])$ is a Kummer covering of a very affine linear space $Z\subset \Spec(\CC[\Lambda])$ if $\kappa^{-1}_{M,\Lambda}(Z)=X$.
 For instance, the Fermat hypersurfaces $X = V(f)\subset \Spec(\CC[x^{\pm},y^{\pm},z^{\pm}])$ for 
 \begin{equation*}
     f=1+x^{n_{1}}+y^{n_{2}}+z^{n_{3}}
 \end{equation*}
 are all Kummer coverings of the very affine linear space $Z=V(1+x+y+z)$.\footnote{Note that $Z$ is isomorphic to the essential projective hyperplane complement $\mathbb{P}^{2}\backslash\,V\big(XYZ(-X-Y-Z)\big)$. We are thus interested in certain natural \'{e}tale coverings of these complements. } We then say that $(X,\Sigma,D)$ is a rational tropical manifold if the initial degenerations of $X$ are all Kummer coverings of very affine linear spaces. For more details, see Section~\ref{sec:RTM}.

 These rational tropical manifolds satisfy several properties that make them similar to ordinary tropical manifolds, but also very amenable to explicit methods.
 
 \begin{center}
 \boxed{\textbf{Key features of rational tropical manifolds}}
 \end{center}
 \begin{enumerate}
 \item The underlying \emph{set} of the rational tropical manifold can be identified with the underlying {set} of a tropical manifold. However, the \emph{multiplicities} can be nontrivial.
 \item A stable intersection of rational tropical manifolds is again a rational tropical manifold, see Proposition~\ref{pro:RationalManifoldStableIntersection}. This can be seen as a tropical analogue of Bertini or Kleiman's theorem.
 \item Let $X=V(f)$ be a hypersurface such that the induced subdivision of the Newton polytope of $f$ is a triangulation. We do not assume that this triangulation is unimodular. Then $(X,\Sigma,D)$ is a rational tropical manifold for any choice of a polyhedral structure $\Sigma$, and for any boundary divisor $D$. For instance, Fermat and Dwork hypersurfaces give rational tropical manifolds.   
 \end{enumerate} 

The first quickly follows from the fact that a vector $v$ is in a cone if and only if $(1/m)\cdot v$ is in the cone. This property allows us to see rational tropical manifolds as certain \emph{tropical nilpotent thickenings}, as in scheme theory. Combinatorially, we can summarize this in terms of the matroids over rings studied in \cite{FM16}: a rational tropical manifold is locally given by a \emph{matroid over $\ZZ$}.  

The additional flexibility provided by rational tropical manifolds allows us to construct a larger pool of examples to which our methods can be applied. 
For instance, in the hypersurface case, we obtain a tropical manifold only when the induced subdivision of the Newton polytope is a regular unimodular triangulation. 
On the other hand, if we have any regular triangulation of a polytope, then we automatically obtain a rational tropical manifold. 
By using the fact that stable intersections of rational tropical manifolds are rational tropical manifolds, we  can then take stable intersections of hypersurfaces to obtain rational complete intersections. 
Note that a stable intersection of two tropical manifolds is rarely a tropical manifold, a fact well known in homotopy continuation, see \cite{HuberSturmfels1995}. 
By working with this more general class, we are thus able to naturally perform generic tropical intersections by taking stable intersections. 
We study this class in great detail in Section~\ref{sec:CompleteIntersections}.
 
\subsection{Tropical varieties as categories}

 Rational tropical manifolds $(X,\Sigma,D)$ are our radial objects. 
 We now wish to enhance these with angular data. 
 To that end, we first interpret these combinatorial objects from tropical geometry as categories.
 
 We start by using the boundary divisor $D$ to partially compactify $\Sigma$. 
 This gives a combinatorial object $\Sigma_{D}$ that is locally a polyhedral complex.
 To make this categorical, we define a poset $\Pscr(\Sigma_{D})$ that captures the different cells in $\Sigma_{D}$. 
 If $\Sigma_{D} = \overline{\Sigma}$ is the full compactification, then $\Pscr(\overline{\Sigma})$ is the poset of the natural regular CW complex structure on $\overline{\Sigma}$. 
 Explicitly, elements of this tropical poset are given by pairs $(\tau,\sigma)$, where $\tau\in \Sigma$, $\sigma\in \Delta$, and ${\rm rec}(\tau)\supset \sigma$. The inclusions between the cells and cones then give rise to a natural partial order, see Definition~\ref{def:TropicalPoset}. We view the induced poset $\Pscr(\Sigma_{D})$ as a category by assigning an arrow $z\to w$ whenever we have a relation $z\geq w$, and we will use the same notation for this category. We emphasize here that this poset only contains the {\it{radial data}} of our polar decompositions. 

As in homotopy theory, we now consider various diagrams of spaces over these categories. That is, we consider functors  
 \begin{equation*}
\mathcal{F}:\Pscr(\Sigma_{D})\to {\rm Top}.
 \end{equation*}
 These can be seen as ways to attach angular data.
 A nontrivial example can be found in Figure~\ref{fig:AperyK33}, which arises from the Ap\'{e}ry-Fermi pencil $$f=1+t(x+y+z+x^{-1}+y^{-1}+z^{-1}).$$
 This defines a hypersurface $X=V(f)\subset \Spec(K[M])$, where $K=\CC(t)\subset \CC(\!(t)\!)$. 
 After choosing a polyhedral structure $\Sigma$ on $\trop(X)$, we can take the compactification $\overline{X}$ of $X$ in the toric variety $Y(\Delta)$ corresponding to the recession fan of $\Sigma$. This gives a K3 surface over $\CC(t)$ in this case.
 The polyhedral structure also gives rise to a strictly semistable model $\mathcal{X}\to U\subset \Spec(\CC[t])$, where $U$ is some Zariski open containing $t=0$. See \cite{Gubler2013} for more on this toric construction.  
The periods of $\overline{X}$ are related to the sequences used in Ap\'ery's proof of the irrationality of $\zeta(3)$, see \cite{PS89}.  
 The tropicalization in this example is homotopy equivalent to a $2$-sphere. By adding the spaces indicated in Figure~\ref{fig:AperyK33} and gluing in the homotopy category (that is, we take a homotopy colimit), we obtain an explicit CW complex structure that faithfully represents the homotopy type of a nearby K3 surface. This can be seen as a combined result of our upcoming Theorems~\ref{thm:MainTheoremHCKNFunctor},~\ref{thm:MainTheoremKummerCoveringAngleHE} and \ref{thm:MainTheoremPolarDecomposition}. 

\subsection{Constructing suitable functors}

We now discuss how to construct suitable functors 
\[
\mathcal{F}:\Pscr(\Sigma_{D})\to {\rm Top}
\] 
that recover the topology of a nearby fiber of an associated strictly semistable model $\mathcal{X}\to U$. Let $(\tau,\sigma)\in \Pscr(\Sigma_{D})$. We want to assign a topological space to this cell. Algebraically, we can assign the initial degeneration of $X=V(I)$ at $w\in {\rm relint}(\tau)$ in the corresponding ambient torus determined by $\sigma$. We note however that this is not functorial, since there is no canonical map between initial degenerations. We use the theory of Kato-Nakayama spaces to amend this. We first associate a natural space $Z^{\sigma,\log}_{\tau}$ to every element $z=(\tau,\sigma)\in \Pscr(\Sigma_{D})$. For elements $(\tau_{1},\sigma)\geq (\tau_{2},\sigma)$, we obtain a natural inclusion map
 \begin{equation*}
 Z^{\sigma,\log}_{\tau_{1}}\to Z^{\sigma,\log}_{\tau_{2}}.
 \end{equation*}        
 If we change the cone $\sigma$ however, then we obtain a quotient map of tori. Overall, these spaces give us our  
 first diagram, which we call the Kato-Nakayama functor $\mathcal{F}_{{\rm KN}}$. We prove the following:
 \begin{maintheorem}\label{thm:MainTheoremHCKNFunctor}
 Let $(X,\Sigma,D)$ be a rational tropical manifold with associated strictly semistable model $\mathcal{X}\to U$ and Kato-Nakayama functor $\mathcal{F}_{\rm KN}$. 
 The homotopy colimit of the Kato-Nakayama functor $\mathcal{F}_{\rm KN}$ is homotopy equivalent to a nearby fiber $\mathcal{X}_{t}\backslash \mathcal{D}_{t}$ of the strictly semistable model $\mathcal{X}\to U\subset \CC$ attached to $(X,\Sigma,D)$. 
 \end{maintheorem}
  
 If the poset only consists of cells in the interior, then this theorem follows from a \v{C}ech-covering argument. At the boundary however, these maps are no longer cofibrations, so we cannot use the standard theory. To fix this, we define a natural cofibrant replacement (which can be viewed as a sealed Kato-Nakayama space) and show that this is homotopy equivalent to the general fiber.
 
 When applied to smooth toric varieties, this yields Danilov's spectral sequence \cite[Theorem 12.2]{Danilov78} for the homology, which was reproved in \cite{WZZ99}. Our result recovers this statement, and it extends it to arbitrary snc-compactifications of hyperplane arrangements, and many more compactifications of smooth algebraic varieties. 

\subsection{Angle spaces}

Next, we show how to replace the Kato-Nakayama spaces by suitable angle spaces.  This reduction to angle spaces will 
allow us to construct explicit CW complex structures on algebraic varieties.  
Consider a Kummer covering $Z_{m}\to X\to Z$ associated to a very affine linear space $Z\subset \Spec(\CC[M])=:T$ for some lattice $M$. We can identify the complex points of $Z$ and $Z_{m}$ with subsets of $\mathrm{Hom}(M,\CC^{\times})=T(\CC)$. We now consider the angle map
\begin{equation*}
    \CC^{\times}\to S^{1}, \quad z\mapsto \dfrac{z}{|z|}.
\end{equation*}
This also gives an angle map
\[
T(\CC)=\mathrm{Hom}(M,\CC^{\times})\to \mathrm{Hom}(M,S^{1})
\]
by composition. Note that after choosing a basis of $M$, we can identify this map with 
\begin{equation*}
    (z_{1},...,z_{n})\mapsto \left(\dfrac{z_{1}}{|z_{1}|},...,\dfrac{z_{n}}{|z_{n}|}\right)\in (S^{1})^{n}.
\end{equation*}
The angle set of $X\subset T=\Spec(\CC[M])$ is the image of $X(\CC)\subset T(\CC)$ under the angle map.
For instance, the angle set of $X=V(1+x^5+y^5)$ is depicted in Figure~\ref{fig:Examplen=2FermatHypersurface}.\footnote{Technically speaking, this is the closed angle set of the natural Kato-Nakayama space associated to $X$. See Section~\ref{sec:AngleSetsCI} for the distinction. We note that these two angle sets have the same homotopy type. }  
We denote this {\it{angle set}} or {\it{coamoeba}} by $\Theta$. 

        \begin{figure}[ht]
        \centering    
            \begin{tikzpicture}[scale=0.4, line width = 0.5, clip=true]
                \def\mincoord{-5}
                \def\maxcoord{5}

                \foreach \i in {-2,...,2} {
                    \foreach \j in {-2,...,2} {
                        \pgfmathsetmacro\xshift{2*\i}
                        \pgfmathsetmacro\yshift{2*\j}
            
                        \filldraw[fill=orange!70, draw=black] 
                            ({-1 + \xshift}, {0 + \yshift}) --
                            ({0 + \xshift}, {1 + \yshift}) --
                            ({-1 + \xshift}, {1 + \yshift}) -- cycle;
            
                        \filldraw[fill=orange!70, draw=black] 
                            ({0 + \xshift}, {-1 + \yshift}) --
                            ({1 + \xshift}, {-1 + \yshift}) --
                            ({1 + \xshift}, {0 + \yshift}) -- cycle;
                    }
                }

                \foreach \m in {-2,...,2} {
                    \draw[black, line width = 0.2] (2*\m, \mincoord) -- (2*\m, \maxcoord);
                }

                \foreach \m in {-3,...,2} {
                    \draw[blue, line width = 2] (2*\m + 1, \mincoord) -- (2*\m + 1, \maxcoord);
                }

                \foreach \m in {-2,...,2} {
                    \draw[black, line width = 0.2] (\mincoord, 2*\m) -- (\maxcoord, 2*\m);
                }
                \foreach \m in {-3,...,2} {
                    \draw[green, line width = 2] (\mincoord, 2*\m + 1) -- (\maxcoord, 2*\m + 1);
                }

                \foreach \m in {0,...,5} {
                    \draw[black, line width = 0.2] 
                        (\maxcoord - 2*\m, \mincoord) -- 
                        (\maxcoord, \mincoord + 2*\m);
                }
                \foreach \m in {0,...,5} {
                    \draw[black, line width = 0.2] 
                        (\mincoord, \mincoord + 2*\m) -- 
                        (\maxcoord - 2*\m, \maxcoord);
                }

                \foreach \m in {0,...,4} {
                    \draw[red, line width = 2] 
                        (\mincoord, \mincoord + 2*\m + 1) -- 
                        (\maxcoord - 2*\m - 1, \maxcoord);
                }
                \foreach \m in {0,...,4} {
                    \draw[red, line width = 2] 
                        (\maxcoord - 2*\m-1, \mincoord) -- 
                        (\maxcoord, \mincoord + 2*\m+1);
                }

                \draw[->, black, line width = 1] (-5, 0) -- (5.5, 0) node[right] {$x$};
                \draw[->, black, line width = 1] (0, -5) -- (0, 5.5) node[above] {$y$};

                \draw[black, line width = 3] (1,1) -- (-1,1) -- (-1,-1) -- (1,-1) -- cycle;

                 \foreach \m in {-5,...,5} {
                        \foreach \n in {-5,...,5} {
                                \draw [fill=black] (\m, \n) circle (2pt);
                    }
                }
                
            \end{tikzpicture}
            \caption{The lifting of the angle set of $V(1+x+y)\subset (\CC^{\times})^{2}$ to $\RR^{2}$. The current fundamental domain determines the topology of the open Fermat curve $V(1+x^5+y^5)\subset (\CC^{\times})^{2}$.}
        \label{fig:Examplen=2FermatHypersurface}
    \end{figure} 
We show that these angle maps give rise to homotopy equivalences for all Kummer coverings of very affine linear spaces. 
This is one of the main results in this paper.

\begin{maintheorem}\label{thm:MainTheoremKummerCoveringAngleHE} 
    Consider a Kummer covering $X \to Z$ of a very affine linear space $Z\subset \Spec(\CC[\Lambda])$. Then the angle map 
\begin{equation*}
    \ang \colon X(\CC)\to \Theta
\end{equation*}
is a homotopy equivalence. 
 \end{maintheorem}

 \begin{remark}
Note that the tropical variety of a very affine linear variety $Z$ is a polyhedral fan, and thus contractible. 
In other words, the topological information carried by the radii of the coordinates is trivial. 
Moreover, the contractibility remains true even when the coefficients of the linear equations depend on a parameter, see \cite[Theorem 2.9]{Speyer08}. 
It is thus natural to expect that the angle set contains all the topological information, which is indeed the case as stated in Theorem~\ref{thm:MainTheoremKummerCoveringAngleHE}. 
Note that the underlying ordinary matroid in general is not enough to capture the topology either, see \cite{Rybnikov11}.
In the language of matroids over ordered blueprints and hyperfields as in \cite{BL21} and \cite{BB19} however, we can then interpret our result as follows: the associated matroid over the \emph{phase hyperfield} recovers the topology of the complement. Moreover, we will see that the phase hyperfield matroid is functorial with respect to taking initial ideals. See \cite{AD19} for more on these types of matroids.
\end{remark}

To prove this theorem, we first make an initial observation. 
The fibers of the angle map are polyhedral and thus convex. In particular, they are contractible. One would now like to use the Smale-Vietoris theorem \cite{SmaleVietoris57} to conclude that the angle map $\mathrm{ang}(\cdot)$ is a homotopy equivalence. 
Unfortunately, the angle map is not proper, so we cannot use this result. Moreover, a modification of this theorem in general does not work for non-proper maps of semialgebraic sets. Indeed, consider the map 
\[
    Z = V\big(x(xy-1)\big) \subset \mathbb{R}^{2}\to \mathbb{R}, \quad (x,y)\mapsto x.
\]
Its fibers are polyhedra and thus contractible. However, $Z$ has three components and $\RR$ only one, so these spaces are not homotopy equivalent.

For the angle map associated to very affine linear spaces, a somewhat stronger property holds. Namely, we have a homotopy equivalence over arbitrarily small open neighborhoods of the angle set. 
To formalize this, we work with slope coordinates on $S^{1}$. Using these coordinates, we show how to construct arbitrarily small contractible neighborhoods whose inverse images under the angle map are again contractible.
These neighborhoods give a \emph{basis-like open covering} of the angle set $\Theta$, so we can conclude using a result by McCord \cite[Theorem 6]{McCord66}.

We then combine Theorems~\ref{thm:MainTheoremHCKNFunctor} and \ref{thm:MainTheoremKummerCoveringAngleHE} to obtain our polar decompositions. The idea is first to replace 
$Z^{\sigma,\log}_{\tau}$ by its open variant $Z^{\sigma,\log,\circ}_{\tau}$, consisting of points that intersect the torus corresponding to $\sigma$. This space is homotopy equivalent to the original $Z^{\sigma,\log}_{\tau}$, since it gives the interior of a manifold with boundary. We note that this requires fixing a homotopy inverse. The space $Z^{\sigma,\log,\circ}_{\tau}$ is homotopy equivalent to its angle set $\Theta_{Z^{\sigma,\log,\circ}_{\tau}}$ by Theorem~\ref{thm:MainTheoremKummerCoveringAngleHE}. 
We then set $\mathcal{F}_{\Theta}(\tau,\sigma) := \Theta_{Z^{\sigma,\log,\circ}_{\tau}}$. 
The morphisms are induced from the maps of Kato-Nakayama spaces.  This gives a natural transformation 
\begin{equation*}
    \mathcal{F}_{KN}\to \mathcal{F}_{\Theta}
\end{equation*}
that is a homotopy equivalence on objects. The following polar decomposition then directly follows.  

\begin{maintheorem}\label{thm:MainTheoremPolarDecomposition}
     Let $\mathcal{F}_{\Theta}$ be the angle functor associated to a rational tropical manifold $(X,\Sigma,D)$ with strictly semistable  model $\mathcal{X}\to U$. There is then a homotopy equivalence 
 \begin{equation*}
     {\rm hocol}(\mathcal{F}_{\Theta})\to \mathcal{X}^{\log,\zeta}_{s}.
 \end{equation*}
 In particular, we find that ${\rm hocol}(\mathcal{F}_{\Theta})$ is homotopy equivalent to $\mathcal{X}_{t}\backslash \mathcal{D}_{t}$ for general $t\in U(\CC)$.
 \end{maintheorem}  
 
 The last part of Theorem~\ref{thm:MainTheoremPolarDecomposition} follows from a theorem of Usui \cite[Equation (0.2)]{Usui01}. 
 This result can be seen as a generalization of various constructions in algebraic geometry. For instance, in the zero-dimensional case, this is Bernstein's theorem on solutions to generic polynomial systems with fixed monomial supports, and in the one-dimensional case this gives a pair-of-pants decomposition of a compact Riemann surface. If the variety $X$ over $\CC(t)$ is a Calabi-Yau variety, then this gives rise to a homotopical analogue of the torus fibration. We investigate these in more detail in Section~\ref{sec:CYMonodromy}. 

 \subsection{Constructing CW complex structures}
 
 We now wish to use Theorem~\ref{thm:MainTheoremPolarDecomposition} to construct explicit CW complex structures on various varieties. We show that if we take suitable compatible cellulations on the individual angle spaces, then this gives an explicit finite CW complex structure on the homotopy colimit. To explain this in more detail, we introduce some notation. We let $\mathcal{F}:\Pscr(\Sigma_{D})\to {\rm Top}$ be a fixed diagram of spaces associated to the poset $\Pscr(\Sigma_{D})$. We assume that all of the spaces are homotopy equivalent to $Z^{\sigma,\log}_{\tau}$. We call this a tropical diagram of spaces, see Definition~\ref{def:TropicalDiagram}. We then fix a set of CW complex structures on the spaces $\mathcal{F}(\tau,\sigma)$. Moreover, we require the maps between the different spaces to be cellular. In other words, we have a functor 
 \begin{equation*}
     \mathcal{F}^{\triangle}:\Pscr(\Sigma_{D})\to {\rm Finite-CW}
 \end{equation*}
 to the category of finite CW complexes. For each cell $z=(\tau,\sigma)\in \Pscr(\Sigma_{D})$, we have a CW complex $\mathcal{F}^{\triangle}(z)$ and we write $|\mathcal{F}^{\triangle}(z)_{m}|$ for the number of $m$-cells. We define the function
    \begin{equation*}
       \mathfrak{g}_{k}:\ZZ \to \ZZ, \quad m \mapsto \sum_{(z_{0}<...<z_{k})}|\mathcal{F}^{\triangle}(z_{k})_{m}|. 
    \end{equation*}
    Here the sum runs over all chains $(z_{0}<...<z_{k})$ of length $k$ in $\Pscr(\Sigma_{D})$. These are exactly the $k$-dimensional cells of the order complex $\Delta(\Pscr(\Sigma_{D}))$ of $\Pscr(\Sigma_{D})$.   Note that $\mathfrak{g}_{k}(m)=0$ for $m<0$ and $m$ large. For any integer $\ell\geq{0}$, we write  $\sh^{\ell}(\mathfrak{g}_{k})$ for the shifted function defined by $\sh^{\ell}(\mathfrak{g}_{k})=\mathfrak{g}_{k}(m-\ell)$. We define $$\mathfrak{f}_{k}=\sh^{k}(\mathfrak{g}_{k}).$$
    We call $\mathfrak{f}_{k}$ the shifted $f$-vector for chains of length $k$. Using this, we obtain the following formula for the $f$-vector of the induced finite CW complex.

\begin{corollary}
Let $(X,\Sigma,D)$ be a rational tropical manifold and let $\mathcal{F}^{\triangle}$ be a cellulation of a tropical diagram associated to $(X,\Sigma,D)$. Its geometric realization is homotopy equivalent to a nearby fiber $\mathcal{X}_{t}\backslash \mathcal{D}_{t}$, and its $f$-vector is given by  
    \begin{equation*}
        \mathfrak{f}(\Delta(\mathcal{F}^{\triangle}))=\sum_{k=0}^{n}\mathfrak{f}_{k}.
    \end{equation*}
\end{corollary}

We note that this $\infty$-categorical formula involving the order complex of the poset is necessary,
as ordinary colimits give the wrong $f$-vector, already in the case of one-dimensional tropicalizations. See Example~\ref{exa:PairOfPantsg2} for an instructive example using angle sets.

\subsection{Explicit \texorpdfstring{$f$}{f}-vectors}\label{sec:IntroFVector}

We make the local $f$-vectors explicit for complete intersections using refined Coxeter arrangements. As an example, we state the corresponding $f$-vector formula here for the complete Fermat hypersurface. We first define the $f$-vector for the local angular spaces. Locally, every space is a Kummer covering of $V(x_{0}+...+x_{s-1})/\CC^{\times}$, viewed as sitting inside an $n$-dimensional torus. We show in Theorem~\ref{thm:FVectorCI} that if the degree of the Kummer covering is $m$, then the number of cells of dimension $d$ in the Coxeter cellulation of the angle set is given by 
 $$f(\Delta_{n,s,m})(d)=m \ d! \ g_{n,s}(d)\quad \text{for all } 0 \leq  d \leq n
,$$ where 
   \[
 g_{n,s}(d) 
 =  2^n
     \left\{ \begin{matrix} n+1 \\ d+1 \end{matrix} \right\} - 2^{n+1-s} \sum_{\ell = 1}^{n + 1 - d} \left[ \binom{n+1}{\ell} - \binom{n+1-s}{\ell} \right] \left\{ \begin{matrix} n + 1 - \ell \\ d \end{matrix} \right\}.
 \]
 Here $\left\{ \begin{matrix} n  \\ d \end{matrix} \right\}$ denotes the Stirling number of the second kind, counting partitions of $n$ elements into $d$ blocks. 
 This gives the $f$-vectors for the angular part. We then need to know the corresponding vectors for the radial part. That is, we need to know how many cells of which type there are in the order complex of $\Pscr(\overline{\Sigma})$. This is given by 
     \begin{equation*}
        N_n(s,d,k) = \frac{(n+1)!}{(n+1-d)! \ (d-s)! \ s!}
\sum_{j=0}^{k}(-1)^{k-j}\binom{k}{j}(2j+1)^{d-s}, \quad \text{for} \quad \begin{matrix}
    0 \leq  s \leq d \leq n-1 \\ 
    0 \leq k \leq n - 1.
\end{matrix}
    \end{equation*}
To explain this formula, consider a chain  
$(z_{0}<...<z_{k})$ of length $k$ in $\Pscr(\overline{\Sigma})$. We write $z_{k}=(\tau,\sigma)$. We then have $s={\rm dim}(\sigma)$ and $d={\rm dim}(\tau)$. The function $N_{n}(s,d,k)$ then counts the different pairs $z_{k}=(\tau,\sigma)$ arising from a chain of length $k$. Combining this, we obtain the following.  

\begin{maintheorem}\label{mainthm:Fermat_FVector}
     Consider the cellulation of the angle functor of  the complete Fermat hypersurface $X_{m}=V(\chi^{me_{0}}+...+\chi^{me_{n}})\subset \mathbb{P}^{n}$ of degree $m$ induced by the Coxeter arrangement. 
     Its geometric realization is homotopy equivalent to $X_{m}$.
     Then \begin{equation*}
    \mathfrak{f}_{k}=\mathrm{sh}^{k} \left(\sum_{0 \leq s \leq d \leq n-1} N_{n}(s,d,k)f(\Delta_{n-s,n+1-d,m^{n-s}})\right)
    \end{equation*} 
    and 
    \begin{equation*}
        \mathfrak{f}=\sum_{k=0}^{n-1} \mathfrak{f}_{k}.
    \end{equation*}
\end{maintheorem}
 For example, for Fermat surfaces, the entries are polynomials in $m^3$, see 
 Corollary~\ref{cor:FermatSurfaceFVector}.
 
 \begin{corollary}
    The $f$-vector of the $m$-th tropical Fermat complex grows as $O(m^{n})$. The same holds in the toric case, where we intersect these with the dense open torus in $\mathbb{P}^{n}(\CC)$.  
 \end{corollary}    
 
More details can be found in Section~\ref{sec:Examples}. There, we also 
treat a wide variety of other examples:
from the classical zero-dimensional and one-dimensional cases, to toric varieties, Fermat and Dwork hypersurfaces. We also show how to obtain three different CW complex structures on K3 surfaces,
and we show how to topologically interpret the monodromy of integral affine structures, as it arises naturally from our angle diagram. This then also allows us to see this monodromy $p$-adically.

\begin{remark}
    We remark here that our constructions also allow us to construct explicit triangulations on smooth algebraic varieties over $\CC$, although we do not pursue this here. To do this, one has to triangulate spaces of the form
    $P\times \Delta^{k}$,
    where $P$ is a convex polytope and $\Delta^{k}$ is the standard $k$-simplex. For this, one can for instance use the staircase triangulation.  
\end{remark}

\subsection{Context and related work}
In this section we provide more context for our results and discuss related work. 
 We start by discussing connections between our polar decompositions, patchworking methods and Morse theory in Section~\ref{sec:RelatedWorkPolarPatchworkingMorse}. In Section~\ref{sec:RelatedWorkAngleSets}, we consider results in the literature related to hyperplane complements and angle sets. Finally, in Section~\ref{sec:RelatedWorkHyperGeometricHodgeMirror}, we discuss connections to Hodge theory and mirror symmetry. 

\subsubsection{Polar decompositions, Viro patchworking and Morse theory}\label{sec:RelatedWorkPolarPatchworkingMorse}

The first nontrivial example of a polar decomposition arises from Weierstrass' complex-analytic parametrization $\CC/\ZZ+\tau\ZZ\to E_{\tau}$.
Here $\tau=a+bi \in \CC$ and $a,b \in \RR$ with $b>0$ and $E_{\tau}$ is the elliptic curve given by the Weierstrass equation 
\begin{equation*}
    y^2=4x^3-g_{2}(\tau)x-g_{3}(\tau).
\end{equation*}
To see the explicit radial and angular parts, we use Tate's multiplicative perspective.
Let $q=e^{2\pi{i}\tau}$, so that $\ell:=-\text{log}|q|=2\pi{b}$. We then find an isomorphism 
\begin{equation*}
    \CC \big/ (\ZZ+\tau\ZZ) \to \CC^{\times}/\langle q\rangle, \qquad x \mapsto e^{2\pi{i}x}.
\end{equation*}
We have a tropicalization map
\begin{equation*}
    \trop: \CC^{\times}/\langle q \rangle \to \mathbb{R}/\ell\ZZ, \qquad x \mapsto -\log|x|.
\end{equation*}
We view the circle $\RR/\ell\ZZ$ as the tropicalization or radial part of the elliptic curve, see \cite{BPR16,CS13} for a Berkovich-theoretic interpretation. Correspondingly, we view $\ell=2\pi{b}$ as the \emph{radial period} of the elliptic curve.
It plays an important role in the theory of archimedean local heights on elliptic curves, see \cite[Chapter 5, Theorem 3.4]{Silverman94}. The tropicalization map induces an exact sequence 
\begin{equation*}
    1 \to S^{1}\to E_{q}\to \RR/\ell\ZZ\to 1,
\end{equation*}
where $S^{1}$ is embedded through $z\mapsto z\bmod \langle q\rangle $.
We can choose a splitting of this exact sequence by choosing a lift $\theta\in \RR$ such that $q/|q|=e^{i\theta}$.
Explicitly, the splitting is given by 
\[
s_\theta(x\bmod \ell\ZZ)
:= \exp\left({-x+ix\dfrac{\theta}{\ell}}\right)  \mod \langle q \rangle.
\]
The element $\theta=2\pi{a}$ gives a lift, and we view this as an \emph{angular period} of $E_{q}$. The corresponding splitting traces out a logarithmic spiral in $\CC^{\times}$ from $1$ to $q$. In $\CC/\ZZ+\tau\ZZ$, this gives the natural straight path from $0$ to $\tau$.
In principle one obtains an infinite range of lifts by adding $2\pi\ZZ$. For each such lift, we obtain an identification $E_{q}\simeq S^{1}\times \RR/\ell\ZZ$. This is the desired polar decomposition for the elliptic curve $E_{q}$. Note that it is not possible to lift this splitting over $\{q\in \CC:0<|q|<1\}$, since going around the origin $q=0$ gives an action of $\ZZ$ on $H_{1}(E_{q},\ZZ)$ through the monodromy matrix 
$\begin{pmatrix}
    1 & 1 \\ 0 & 1
\end{pmatrix}.$ 
Explicitly, the path from $0$ to $\tau$ is transformed to the path from $0$ to $\tau+1$. 

\begin{remark}
From the perspective of this paper, we can reconstruct the topology of $E_{q}$ through a diagram of topological spaces over  $\trop(E_{q})=\RR/\ell\ZZ$ with a specific CW complex structure. Explicitly, we endow $\trop(E_{q})$ with the CW structure consisting of a single vertex $v$ and a single edge $e$. We can first simply assign the constant diagram where every space is $S^{1}$, and the comparison maps are the identity. The colimit is then $S^{1}$, but the homotopy colimit is $(S^{1})^{2}$, which recovers the structure of the elliptic curve. 

On the other hand, the Tate curve construction tells us that there is a more natural diagram that behaves well with respect to the complex structure. To see this, write $z=e^{-t+iu}\in \CC^{\times}$, so that $(t,u)$ form polar coordinates. In terms of these coordinates, the complex number $q$ with $q/|q|=e^{i\theta}$ acts as $(t,u)\mapsto (t+\ell,u+\theta)$. We then find \[
E_q\simeq
\frac{[0,\ell]\times{S^1}}
{(\ell,u)\sim(0,u-\theta)}.
\] 
This is the homotopy colimit of the twisted diagram
\[
S^1_e
\underset{R_{-\theta}}{\overset{\mathrm{id}}{\rightrightarrows}}
S^1_v, \qquad \text{where } R_{-\theta}:S^{1}_{e}\to S^{1}_{v}, \quad  u\mapsto u-\theta.
\]
This is the natural diagram of spaces arising from the Tate uniformization. The pair $(\theta,\ell)$ of the angular period $\theta$ together with the radial part $\ell$ encompass the full complex structure on $E_{q}$. 
Explicitly, we have $\tau=\dfrac{\theta+i\ell}{2\pi}$. Note that $\theta$ is only well defined up to a multiple of $2\pi$. If we take $\theta+2\pi{n}$ instead, then this changes $\tau$ to $\tau+n$. These are equivalent in $\mathbb{H}/\ZZ$, where $\mathbb{H}$ is the complex upper half plane. In order to remove this ambiguity, we would have to view the twisting map as coming from a map of universal covering spaces $\RR_{e}\to \RR_{v}$. 

Overall, we view $\ell$ as giving rise to metric data on the diagram, and $\theta$ giving rise to twisting data for the comparison maps through $R_{-\theta}$. Even though the constant diagram is equivalent to the Tate diagram from the perspective of homotopy theory, the Tate diagram allows us to see the complex structure and the connection to the upper half plane more clearly.
\end{remark}

The theory of Teichm\"{u}ller spaces generalizes this decomposition. Consider a compact orientable real surface $\overline{X}$ with $D\subset \overline{X}$ a set of marked points. We assume that $2-2g-|D|<0$, so that $X=\overline{X}\backslash D$ is hyperbolic.
There then exists a pair-of-pants decomposition of ${X}$, i.e., an embedding of a disjoint union of copies of $\mathbb{P}^{1}\backslash\{0,1,\infty\}$ into $X$ such that its complement is a disjoint union of circles. 
This is one of the main ingredients in the Fenchel-Nielsen parametrization of the corresponding Teichm\"{u}ller space of all complex structures on $\overline{X}$. 
More precisely, endowing $X$ with the structure of a Riemann surface is equivalent to specifying a point in $\mathbb{R}^{3g-3+|D|}\times \mathbb{R}^{3g-3+|D|}_{>0}$, up to an equivalence under the mapping class group. 
Here $3g-3+|D|$ is exactly the number of copies of $S^{1}$, and a point in $\mathbb{R}^{3g-3+|D|}\times \mathbb{R}^{3g-3+|D|}_{>0}$ corresponds to an assignment of a length $\ell \in \RR_{>0}$ and a twisting factor $\phi \in \mathbb{R}$ for each $S^1$ in the decomposition. 
This data is sufficient to obtain a full complex structure. 
When $g=1$ and $|D|=1$, these factors can be interpreted as the real and imaginary parts of $\tau=a+bi$ in the elliptic curve discussion above. 
Note that the copies of $\mathbb{P}^{1}\backslash\{0,1,\infty\}$ are purely angular objects, see Theorem~\ref{thm:MainTheoremKummerCoveringAngleHE}. 
We view these pair-of-pants decompositions as giving a polar decomposition of hyperbolic Riemann surfaces.

These pair-of-pants decompositions also arise naturally when tropicalizing plane curves, see Section~\ref{sec:ClassicalExamplesD01} for an explicit example. An analogue of this decomposition in higher dimensions 
was introduced for smooth complex projective hypersurfaces by Mikhalkin  \cite{MIKHALKIN2004} who showed that any smooth projective hypersurface $X = V(f)$ in $\CC\PP^{n+1}$ admits a decomposition into complex $n$-dimensional \emph{pairs-of-pants}, i.e. spaces of the form $\CC\PP^{n} \setminus |\mathcal{A}|$ where $\mathcal{A}$ is a generic arrangement of $n+2$ hyperplanes in $\CC\PP^n$. 

After Mikhalkin's paper, Kerr and Zharkov \cite{KerrZharkov} showed that, for a generic smooth complex hypersurface $X \subset (\CC^\times)^n$ with a fixed Newton polytope, one can construct a so-called \emph{phase-tropical hypersurface} $X_{\rm phase} \subset \RR^n \times (S^1)^n$ such that $X$ is homeomorphic to $X_{\rm phase}$. The topological space $X_{\rm phase}$ is built from the tropical hypersurfaces associated to $X$ together with phase (or angle) data, which results in a polyhedral complex with torus fibers. 
This shows that, for generic hypersurfaces in the torus, one can recover the topology of the complex variety from a combinatorial model endowed with angle data. Note that the homotopy colimit also appears through the product construction in \cite{KerrZharkov}. Moreover, one finds a description of the gluing problem (see Section~\ref{sec:TriangulationsRCIM}) in the case where all the tropical charts are compatible. See Examples~\ref{ex:refinedCoxeterArrangement1} and  \ref{exa:DworkFamilyK3}, where more complicated charts are needed for constructing CW complex structures. A similar approach using phase structures and Viro's patchworking method can be found in \cite{BLR24,RRS22,RRS25}.

Our paper can be seen as a natural continuation of these types of patchworking constructions, but from a purely homotopic point of view. This allows us to construct explicit CW complex structures on algebraic varieties. Our approach moreover emphasizes the angular perspective, so that these constructions can be seen as suitable polar decompositions. We note here that our varieties are not assumed to be real. Moreover, we can deal with varieties that are not hypersurfaces, and we do not restrict ourselves to the compact or open case.

Our method carries many resemblances to Morse theory. There, one fixes a suitable function $f: X\to \mathbb{R}$, and one then reconstructs the topology of $X$ from the critical points. Our analogue of the choice of a Morse function corresponds to a choice of a strictly semistable degeneration 
\begin{equation*}
    \mathcal{X}\to U\subset \Spec(\CC[t])
\end{equation*}
such that $X$ occurs as a suitable fiber. These strictly semistable degenerations are explicitly provided by the theory of rational tropical manifolds. Rather than study all of the fibers as in Morse theory, we only consider the special fiber given by $t=0$. The main difference between the two approaches is that the local objects we obtain are explicit combinatorial objects. In particular, the formulas for the $f$-vectors we obtain for complete intersections are stated in terms of the polyhedral structure of the tropicalization, and the structure of the standard angular building block. This has the advantage of being very explicit, and of being very closely related to the corresponding complex structures. On the other hand, constructing these tropical degenerations can be somewhat harder than finding Morse functions in certain circumstances.   

\subsubsection{Angle sets of very affine linear spaces and hyperplane complements}\label{sec:RelatedWorkAngleSets}

Our main result, Theorem~\ref{thm:MainTheoremKummerCoveringAngleHE}, on angle sets of very affine linear spaces can be seen as a natural continuation of \cite{Deligne72,Hattori75,Salvetti87,BjornerZiegler92} on constructing combinatorial models of hyperplane complements. 
Similar matroidal techniques (in the real case) were also used in \cite{RRS22,SY25,MY25}, with the latter focusing on central hyperplane arrangements. In a sense, our model in terms of angle sets can be seen as a \emph{continuous} version of the Salvetti complex, the Bj\"{o}rner-Ziegler complex, and their matroidal relatives. 

We discuss the connections between the discrete approach, using Salvetti complexes and Bj\"{o}rner-Ziegler signs, and our continuous approach using angle sets here.

\begin{enumerate}
    \item The notion of an admissible sign function from \cite[Definition 2.2]{BjornerZiegler92} can be seen as subdividing the tori $(S^{1})^{k}=\RR^{k}/2\ZZ^{k}$ in terms of boxes. That is, we record a certain sign pattern if and only if the corresponding box intersects the angle set nontrivially. As such, we can view these sign patterns as slicing or triangulating the angle set. See Figure~\ref{fig:BZ_Signs2} for an instructive example. 
    \item From a tropical point of view, we can view the hyperplanes in an arrangement as corresponding to coordinates in a torus, see Remark~\ref{rem:EquivalenceVALSandEPHC}. However, the corresponding complex sign patterns are usually incompatible with the sign patterns coming from a different basis of the ambient character lattice. That is, they lead to
    different subdivisions, see Figure~\ref{fig:BZ_Signs2} for instance.
    Angle sets on the other hand are completely independent of a choice of basis, see our definition in Section~\ref{sec:AngleSets}. This is what allows us to easily glue angles from different cells of the tropicalization, where
    different bases of character lattices are used. In principle, there is a way to bypass this problem by taking finer and finer sign patterns on initial degenerations\footnote{These finer sign patterns are considered in \cite[Example 2.3]{BjornerZiegler92} for instance. }, but we shall not take this approach here.  
    \item Although angle sets allow for intrinsic gluing maps as we will see in Section~\ref{sec:AngleSets}, in order to obtain explicit CW complex structures on general algebraic varieties, one has to make sure that the cellulations of the angle sets are chosen properly. This will be discussed in more detail in Section~\ref{sec:CWStructure}. In order to do this, it suffices to take certain toric subdivisions of angle sets. See Examples~\ref{ex:refinedCoxeterArrangement1} and  \ref{exa:DworkFamilyK3}. One can interpret these subdivisions as generalizations of the sign decompositions provided by the Salvetti and Bj\"{o}rner-Ziegler complex. 
\end{enumerate}

    \begin{figure}[ht]
    \begin{center}
    \scalebox{0.5}{
    \definecolor{red_color}{rgb}{1,0,0}
    \definecolor{gray_color}{rgb}{0.3,0.3,0.3}
    \definecolor{blue_color}{rgb}{0.3,0.3,1}
    \definecolor{green_color}{rgb}{0.1,0.7,0.3}
    
    \begin{tikzpicture}[line cap=round,line join=round,x=1cm,y=1cm]
    \clip(4,7) rectangle (15,15);
    \fill[line width=2.8pt,color=red_color,fill=red_color,fill opacity=0.5] (10,8) -- (13,11) -- (13,8) -- cycle;
    \fill[line width=2.8pt,color=red_color,fill=red_color,fill opacity = 0.5] (7,14) -- (10,14) -- (7,11) -- cycle;
    
    \draw [line width=2.8pt] (7,11)-- (10,14);
    \draw [line width=2.8pt] (10,8)-- (13,11);

    \draw [line width=2pt, loosely dashed] (7,11)-- (13,11);
    \draw [line width=2pt, loosely dashed] (10,8)--(10,14);
    \draw [line width=2pt, loosely dashed] (7,8)--(13,14);
    \draw [line width=2.8pt] (10,14)-- (13,14);
    \draw [line width=2.8pt] (7,8)-- (10,8);
    \draw [line width=2.8pt] (10,8)-- (13,8);
    \draw [line width=2.8pt] (13,8)-- (13,14);
    \draw [line width=2.8pt] (7,8)-- (7,14);
    \draw [line width=2.8pt] (7,14)-- (10,14);
    
    \begin{scriptsize}

    \draw [fill=blue_color] (7,8) circle    (6pt);
    \draw [fill=blue_color] (13,14) circle  (6pt);
    \draw [fill=blue_color] (13,8) circle   (6pt);
    \draw [fill=blue_color] (7,14) circle   (6pt);
    \draw [fill=green_color] (10,8)  circle (6pt);
    \draw [fill=green_color] (10,14) circle (6pt);
    \draw [fill=yellow] (7,11) circle (6pt);
    \draw [fill=yellow] (13,11) circle (6pt);
    
    \draw (7.65, 13.35) node [anchor=north west][inner sep=0.75pt] [xscale = 2, yscale = 2] {$ij$};
    \draw (11.65, 9.2) node [anchor=north west][inner sep=0.75pt] [xscale = 2, yscale = 2] {$ji$};

    \draw (13.25,8) node [anchor=north west][inner sep=0.75pt] [xscale = 2, yscale = 2] {$++$};
    \draw (6,14.75) node [anchor=north west][inner sep=0.75pt] [xscale = 2, yscale = 2] {$++$};
    \draw (13.25,14.75) node [anchor=north west][inner sep=0.75pt] [xscale = 2, yscale = 2] {$++$};
    \draw (5.75,8) node [anchor=north west][inner sep=0.75pt] [xscale = 2, yscale = 2] {$++$};
    
    \draw (13.25,11.25) node [anchor=north west][inner sep=0.75pt] [xscale = 2, yscale = 2] {$+-$};
    \draw (5.5,11.25) node [anchor=north west][inner sep=0.75pt] [xscale = 2, yscale = 2] {$+-$};

    \draw (9.5,14.75) node [anchor=north west][inner sep=0.75pt] [xscale = 2, yscale = 2] {$-+$};
    \draw (9.5,7.75) node [anchor=north west][inner sep=0.75pt] [xscale = 2, yscale = 2] {$-+$};

    \end{scriptsize}
    
    \end{tikzpicture}
     }
    \end{center}
    \caption{
    The angle set of $V(1+x+y)\subset (\CC^{\times})^{2}$ gives an interpretation of the Bj\"{o}rner-Ziegler complex. Namely, the complex signs arise by subdividing the torus into four squares, and by keeping track when a particular angle occurs. For instance, the angles in the top-left correspond to the sign pattern $ij$, since $x$ lies in the upper half-plane and $y$ consequently lies in the lower half-plane.
    The diagonal lines correspond to the same type of decomposition, but coming from a different choice of dehomogenizing hyperplane. For instance, if we choose $x_{1}=1/x$, $y_{1}=y/x$, then the natural decomposition of the angle set of the torus given by the complex matroid approach corresponds to the diagonal lines together with the vertical lines.
    }
    \label{fig:BZ_Signs2} 
    \end{figure} 
We note here that the results by Hattori in \cite{Hattori75} on generic hyperplane complements can be directly obtained from our results, since both results are phrased in terms of subspaces of tori $\mathbb{R}^{n}/\ZZ^{n}$. For instance, for $V(1+x+y+z)\subset (\CC^{\times})^{3}$ (corresponding to the hyperplane complement $\mathbb{P}^{2}$ minus four generic hyperplanes), Hattori obtains the $2$-skeleton of $\mathbb{R}^{3}/\ZZ^{3}$. Our angle set deformation retracts onto this set, as we obtain a complement of $\mathbb{R}^{3}$ minus a set of open zonotopes. We do not fully work out the details here to obtain Hattori's results just yet, but it should be clear from this example that they can be directly obtained, once we know the angle sets of generic hypersurfaces.

The angle sets of these  generic hypersurfaces were for instance studied in \cite{NS13NonArchimedean}, where it was shown that they are essentially equivalent to $\mathbb{R}^{n}/\ZZ^{n}$ minus an open zonotope. Additionally, various structural results on the angle sets of affine linear spaces were obtained in \cite{NS13NonArchimedean,NS13PhaseLimit,NS24}, where it was shown that the closure of the angle set of an affine linear space $Z = V(I) \subset (\CC^\times)^n$ decomposes into the union of angle sets of its initial degenerations $V({\rm in}_{w}(I))$, thereby relating the angle behavior at infinity to the tropical fan of the affine linear space in question. For us, a similar decomposition directly follows from the Kato-Nakayama set-up, see Section~\ref{sec:GluingAngleSets}. Moreover, they prove that the phase-limit set of the complement of a hyperplane arrangement can be described as a finite union of products of angle sets of hyperplane complements indexed by flags of flats, and in fact is closely related to the associated Bergman fan of the arrangement. In our paper, we will explicitly use their determination of the closed angle set of a generic hypersurface to make our pair-of-pants decompositions for complete intersections explicit.

\subsubsection{Hodge structures and mirror symmetry}\label{sec:RelatedWorkHyperGeometricHodgeMirror}

Every rational tropical manifold $(X,\Sigma,D)$ comes with a strictly semistable model $\mathcal{X}\to U\subset \Spec(\CC[t])$. The special fiber $\mathcal{X}_{s}$ admits a natural stratification in terms of
Kummer coverings of essential projective hyperplane complements. One can consider the de Rham cohomology groups $H^{m}_{{\rm dR}}(\mathcal{X}_{t},\CC)$ of the fibers. If we want to include $t=0$, then we can take the induced family of Kato-Nakayama spaces 
\begin{equation*}
    \mathcal{X}^{\log}\to U^{\log}.
\end{equation*}
Here $\mathcal{X}$ carries the log structure from the special fiber $\mathcal{X}_{s}$, and $U$ carries the log structure from $t=0$. This gives rise to a vector bundle $\mathcal{H}^{m}_{{\rm dR}}$ with a logarithmic Gauss-Manin connection $\nabla$. We can  consider flat sections of this vector bundle. These satisfy a differential equation known as the Picard-Fuchs equation.
\begin{example}
For instance, consider the Dwork cubic $\overline{X}=V(F)\subset \mathbb{P}^{2}_{\CC(t)}$ defined by 
\begin{equation*}
    F=3XYZ-t^3(X^3+Y^3+Z^3).
\end{equation*}
The scheme $\overline{X}/\CC(t)$ gives rise to a semistable family $\mathcal{X}\to U$ over some Zariski open subset $U\subset \Spec(\CC[t])$ of $t=0$. The corresponding Picard-Fuchs differential equation is 
\[
\frac{d^2w}{dt^2} + \frac{4t^3 - 1}{t(t^3 - 1)} \frac{dw}{dt} + \frac{2t}{t^3 - 1}w = 0.
\]
This has the holomorphic solution
\[
w(t) = \sum_{n=0}^{\infty} \frac{(3n)!}{3^{3n}(n!)^3} t^{3n}={}_{2}F_{1}\left(\frac{1}{3},\frac{2}{3};1;t^3\right),
\]
where the latter is the classical Gaussian hypergeometric function. This can be seen as describing the periods\footnote{One can integrate over many cycles in $H_{1}(\mathcal{X}_{t},\ZZ)$, but they are all related to each other by a linear transformation. } $\int_{\gamma}\varpi$ corresponding to the holomorphic $1$-form $$\varpi(t) =  -\frac{tdy}{3tx^2 - y}.$$ 

From a tropical point of view, the dehomogenized Laurent polynomial $f=3xy-t^3(x^3+y^3+1)$ does not define a smooth tropicalization. However, the initial degenerations are Kummer coverings of very affine linear spaces, see Example~\ref{exa:DworkFamily1}. We thus see that $\overline{X}$ naturally corresponds to a rational tropical manifold $(X,\Sigma,D)$. 
\end{example}

We now note that these periods have various connections to point-counting over finite fields, see \cite{Dwork69}. Indeed, we have the formula  $\lambda_{1}=w(t)/w(t^{p})$ for one of the eigenvalues of Frobenius of the elliptic curve $\mathcal{X}_{t}$\footnote{One has to assume that $\mathcal{X}_{t}$ is an elliptic curve with ordinary reduction here. Here we are considering a natural extension $\mathcal{X}\to U\subset \Spec(\ZZ[t])$ of the Dworkian family. }. 

One can in fact use the semistable reduction type of this Dwork pencil at $t=0$ to compute the periods, and one then uses the corresponding $p$-adic differential equations to find the number of points of a distant elliptic curve modulo a prime. This is the basis of Lauder's deformation technique, see \cite{Lauder04} and \cite[Chapter 20]{Kedlaya10}.

We view our rational tropical manifolds as giving combinatorial prototypes that generalize these Dworkian examples. We note here that the approach using logarithmic geometry in Sections~\ref{sec:RTM} and \ref{sec:TropicalDiagrams} naturally generalizes to non-archimedean fields such as $\QQ_{p}$.  For instance, an analogue of our homotopy colimit formula in Theorem~\ref{thm:MainTheoremHCKNFunctor} over arbitrary non-archimedean fields can be found in \cite[Section 3.6]{AHLS26}.

From a modern perspective, the Dworkian methods rely on a proper understanding of the corresponding Hodge-theoretic invariants of these families. The connection between these and tropical geometry is studied in many papers, see \cite{HK12,BK24,IKMZ19} for instance. Here we note that the first focuses on general sch\"{o}n subvarieties, the second focuses on smooth tropicalizations and the third on smooth tropicalizations with a certain compactification condition. Our rational tropical manifolds are naturally sch\"{o}n, so they fall under the material developed in \cite{HK12}. See \cite{GRW17} for a similar treatment in terms of Berkovich spaces. In contrast to these methods, our approach is purely topological. In particular, we will not need  limiting mixed Hodge structures or analogues of the Steenbrink spectral sequence to define tropical homology groups. This method can thus also compute the integral homology groups of a nearby fiber. In a sense, the tropical homology groups in \cite{IKMZ19} seem closer to duals to de Rham cohomology groups, whereas the homology groups we introduce here arise from the integral homology of a CW complex.

We note that we also obtain a spectral sequence of cohomology groups converging to the cohomology of the general fiber, as in \cite{IKMZ19}. Moreover, the local cohomology groups are also given by Orlik-Solomon algebras in the smooth case. It would be interesting to see if the tropical Hodge spectral sequence can somehow be obtained by integrating differentials over the homology cycles coming from our approach.

Finally, we discuss the connection to the Strominger–Yau–Zaslow (SYZ) viewpoint in mirror symmetry, see \cite{Gross01}, \cite{Gross09}, \cite{GS03}, \cite{Yamamoto22} and \cite{Arguz21}. Our approach in Section~\ref{sec:CYMonodromy} is closest  to the approach taken in the last two papers. Our computations for the Ap\'{e}ry-Fermi pencil give an elementary approach to the monodromy of integral affine structures, as defined formally in \cite{Yamamoto22}. In particular, we are able to reproduce this monodromy from a purely categorical and topological perspective. Namely, we have a diagram of Kato-Nakayama spaces, and walking around the $2$-sphere corresponds to applying the arrows in our angle functor. In particular, this approach shows that the singular points in fact should be considered as infinitesimal points on the cube (as in the theory of Zariski-Riemann spaces) in order to keep them invariant under ramified base extensions $\CC[t]\subset \CC[t^{1/n}]$. The approach also shows that there is a natural $p$-adic analogue: we simply replace the homology groups of the local Kato-Nakayama spaces by the logarithmic prime-to-$p$ \'{e}tale cohomology of the various spaces. The monodromy then arises from exactly the same types of calculations, but with Kummer coverings replacing our angular considerations.

\subsection{Outline}

We start in Section~\ref{sec:RTM} by defining the notion of a rational tropical manifold, denoted by $(X, \Sigma, D)$, through Kummer coverings of very affine linear spaces. 
We emphasize here that we allow both open and closed tropicalizations, so that we are not necessarily dealing with the compact case. 
We define an associated tropical poset $\Pscr(\Sigma_{D})$ for each rational tropical manifold $(X,\Sigma,D)$.
This gives rise to the base category for our diagram of spaces. 
It can be thought of as containing the cells of the rational tropical manifold. 
We also review how these manifolds give rise to strictly semistable models $\mathcal{X}\to U\subset \Spec(\CC[t])$ for some Zariski open subset $U$ containing $t=0$. 

In Section~\ref{sec:TropicalDiagrams}, we enhance these rational tropical manifolds with the structure of a diagram of spaces, denoted by $\Pscr(\Sigma_{D})\to {\rm Top}$. 
The idea is to associate a Kato-Nakayama space to each cell. For cells in the interior, the comparison maps are given by inclusion maps. 
For cells near the boundary, we obtain quotients by real tori. These are not  cofibrations, so it is not immediate that the homotopy colimit of this diagram is homotopy equivalent to its colimit. 
We define a natural cofibrant replacement of the diagram and show that the induced space is homotopy equivalent to the original colimit, so that the diagram gives the desired homotopy type. 
A special case of this result for smooth toric varieties can be found in \cite[Section 5.3]{WZZ99}.   

In Section~\ref{sec:AngleSets}, we prove that the angle map associated to a Kummer covering of a very affine linear space is a homotopy equivalence. 
Moreover, we show how to associate angle spaces to Kato-Nakayama spaces in a functorial way, so that we obtain a diagram of spaces that is equivalent to the original diagram. We call this diagram of spaces the angle functor.

In Section~\ref{sec:CompleteIntersections}, we make the angle functor explicit for complete intersections. 
We show that the angle map continues to give a homotopy equivalence on the Kato-Nakayama space of an snc-compactification associated to a Kummer covering. 
We then describe a cellulation of these angle sets in terms of Coxeter arrangements, and we determine their $f$-vectors in terms of Stirling numbers of the second kind, see Section~\ref{sec:FVectorCI}. This cellulation can be seen as a first step towards obtaining CW complex structures on complete intersections. 

In Section~\ref{sec:CWStructure}, we define the general notion of a cellulated tropical diagram. 
The homotopy colimit of a tropical diagram then has an associated CW complex. 
This gives an explicit CW complex structure on the homotopy colimit, which in turn is homotopy equivalent to a nearby fiber $\mathcal{X}_{t}\backslash \mathcal{D}_{t}$ of an associated strictly semistable model. 
We prove a formula for the $f$-vector of this CW complex in terms of the $f$-vectors of the spaces in the cellulation, see Theorem~\ref{thm:CWComplexStructureTriangulatedDiagram}. In Section~\ref{sec:TriangulationsRCIM}, we show how to explicitly triangulate rational complete intersection manifolds. In Section~\ref{sec:FVectorRefinedArrangements}, we extend our earlier formula for $f$-vectors of Coxeter arrangements to  $f$-vectors of refined Coxeter arrangements. This allows us to give a concrete $f$-vector formula for the CW complex structure associated to complete intersections in toric varieties.  

Finally, we show how our techniques work in practice through a number of explicit examples, see Section~\ref{sec:Examples}. 
We start with rational tropical manifolds of dimensions $0$ and $1$ in Section~\ref{sec:ClassicalExamplesD01}, where they reduce to mixed volumes and pair-of-pants decompositions. 
We then consider the case of toric varieties in Section~\ref{sec:ToricVarieties}, where we obtain an explicit formula for a cellulation of a smooth toric variety in terms of our cellulation of the torus. In Section~\ref{sec:FermatAndDwork}, we give a formula for the tropical $f$-vector of Fermat hypersurfaces in terms of their degree. We also show how our techniques work in the case of the Dwork K3 surface, where an additional desingularization is needed on the tropicalization. 
Finally, in Section~\ref{sec:CYMonodromy}, we study the Ap\'{e}ry-Fermi pencil, and we show how to obtain the monodromy of integral affine structures in a purely topological way from our tropical diagram of spaces. This in particular shows that this monodromy also occurs over $p$-adic fields.

\subsection{Notation and conventions}
We give a short list containing the general notation we shall use throughout this paper: 
\begin{enumerate}
    \item For a category $\mathcal{C}$, we write $\mathcal{C}^{\rm opp}$ for the opposite category.
    
    \item Standard categories such as the category of topological spaces will be denoted by ${\rm Top}$.
    
    \item Functors will be denoted by Gothic letters such as $\mathcal{F}$ and $\mathcal{G}$.
    
    \item  A variety over a field $K$ is a reduced separated scheme of finite type over $K$. It will be denoted by $X$, $Y$ or $Z$. If $K$ is algebraically closed, we will sometimes identify these with the classical varieties obtained by taking their $K$-valued points.
    
    \item Projective space over a field $K$ will be denoted by $\mathbb{P}^{n}={\rm Proj}(K[X_{0},...,X_{n}])$.
    
    \item Affine space over a field $K$ will be denoted by $\mathbb{A}^{n}={\rm Spec}(K[x_{1},...,x_{n}])$.     \item Let $M$ be a finitely generated abelian group. Elements of the group ring $K[M]$ will be written as $\sum_{i=1}^{r}c_{i}\chi^{m_{i}}$, where $c_{i}\in{K}$ and  $m_{i}\in{M}$. We view these as functions $f: M\to K$ such that there are finitely many elements $m\in{M}$ with $f(m)\neq{0}$. In the notation above, we have $f(m_{i})=c_{i}$ for the given $m_{i}\in M$, and $f(m)=0$ elsewhere.
    
    \item Let $M$ be a lattice. We write $T$ or $\mathbb{T}$ for the corresponding torus $\Spec(K[M])$ over a field $K$. We will often pick a basis for $M$, so that we can identify $T$ with $\Spec(K[x_{1}^{\pm},..., x_{n}^{\pm}])$. The corresponding variety is $(K^{\times})^{n}$ for $K^{\times}=K\backslash\{0\}$. 

    \item Models of varieties over a base scheme will be denoted by Gothic letters. We will often write $\phi: \mathcal{X}\to U$ for such a model. Here $U$ is thought of as a parameter space.
    
    \item Combinatorial objects such as simplicial or polyhedral complexes will usually be denoted by Greek letters such as $\Sigma$, $\Gamma$, and $\Delta$. We will generally use the conventions adopted in \cite{Gubler2013, MS15}.
    
    \item The circle $S^{1}$ will always be viewed as a subspace of $\CC$. That is, $S^{1} := \{z\in\CC:|z|=1\}$.
    \item For any vector space $V$ with a fixed basis $v_{i}$, we write $\mathds{1}_n$ for the vector $\sum_{i=1}^{n}v_{i}$, the all-ones vector. For instance, if $V=\RR^{n}$ and $v_{i} = e_{i}$ is the standard basis, then this yields $\mathds{1}_{n}=(1,...,1)$. 
    \item Let $X$ be a variety over $\CC$. We say a certain property holds for general $t\in X(\CC)$ if it holds for all $t\in U(\CC)$, where $U$ is a Zariski dense subset.  
    \item We write $\mathbb{D}:= \{z \in \CC \colon |z|\leq 1\}$ for the closed unit disk in $\CC$.

    \item We write $\sh^{\ell}$ for the shift operator acting on functions $\mathfrak{g}: \ZZ \to \ZZ$ as follows:
    \[  
        \sh^{\ell}(\mathfrak{g})(m) := \mathfrak{g}(m - \ell), \quad \text{for all } \ell, m \in \ZZ.
    \]
\end{enumerate}

\subsection*{Acknowledgments}
    We thank Piotr Achinger, Tadashi Ashikaga, Jeffrey Giansiracusa, Keita Goto, Andrew Harder, Joachim Jelisiejew, Masanori Kobayashi, Arne Kuhrs, Marcin Lara, Sukjoo Lee, Gaku Liu, Ryota Mikami, Franciszek Nowak, Yuji Odaka, Yue Ren, Frank Sottile, Bernd Sturmfels, Yuki Tsutsui, Kazuhiko Yamaki, Yuto Yamamoto and Chi Ho Yuen for the many inspiring conversations, questions, suggestions and comments. We also thank the attendees of the various conferences where preliminary versions of the paper were presented for their helpful comments and questions.

    YE was partially supported by Deutsche Forschungsgemeinschaft (DFG, German Research Foundation) SFB-TRR 195 “Symbolic Tools in Mathematics and their Application”. 
    PH was supported by the JSPS Postdoctoral Fellowship 23769 and KAKENHI 23KF0187 as a Postdoctoral Fellow at the University of Tsukuba and Tohoku University.

\subsection*{Declaration of the use of AI}

This paper was written over the course of the past two years (Jul.~2024 - Oct.~2026) during most of which we did not use any AI. Towards the end of this project, we ended up using the basic versions of ChatGPT and Google Gemini in the following ways: 

\begin{enumerate}
    \item We realized that we needed a more general version of Theorem~\ref{thm:FVectorCI} for refined Coxeter arrangements. After we provided our statement and proof of Theorem~\ref{thm:FVectorCI} to ChatGPT, it was able to produce the statement and most of the steps in the proof of Theorem~\ref{thm:refined_Coxter_F_vect}. 
    \item ChatGPT was similarly used to find the formula in Theorem~\ref{thm:OrderComplexTropicalizationFVector} and for computations in Section~\ref{sec:FermatAndDwork}. 
    \item We also used ChatGPT and Google Gemini to initially help with suggestions for the material in Sections~\ref{sec:CofibrantReplacements}, \ref{sec:TriangulatedTropicalDiagrams} and \ref{sec:CYMonodromy}. 
\end{enumerate}
The exact statements and proofs in all cases were fully written and verified by the authors. We take full responsibility for the correctness and integrity of the present work.

\section{Rational tropical manifolds}\label{sec:RTM}

In this section, we introduce the notion of a \emph{rational tropical manifold}. This will be based on the notion of \emph{Kummer coverings} of very affine linear spaces. These for instance arise naturally from triangulations of Newton polytopes. We show that the stable intersection of two rational tropical manifolds is again a rational tropical manifold in Proposition~\ref{pro:RationalManifoldStableIntersection}. In the last section, we recall how these rational tropical manifolds naturally come equipped with a set of strictly semistable models. These will be used to recover the topology of fibers close to the special fiber.

\subsection{Toric and tropical geometry}\label{sec:ToricTropicalGeometry}

In this section, we review some concepts from toric and tropical geometry, and we establish the notation we will use throughout the rest of this paper. For more background information, we refer the reader to \cite{MS15}, \cite{Gubler2013} and \cite{Payne09}. Throughout this section, $K$ is a field with a non-archimedean valuation $\val: K^{\times}\to \overline{\RR} := \RR \cup \{+ \infty\}$. We allow $\val(\cdot)$ to be the trivial valuation. We call the pair $(K,\val(\cdot))$ a valued field, and we will denote it by $K$ if the valuation is clear, and write $R$ for the valuation ring of $K$.

\begin{remark}\label{rem:ValuedFields}
The valued fields we use throughout this paper will be limited to the following:
\begin{enumerate}
    \item The field of complex numbers $\CC$ with the trivial valuation. 
    \item The field of Laurent power series $\CC(\!(t)\!)$ with the $t$-adic valuation. This sends $t$ to $1$. 
    \item The field of Puiseux series $\CC\{\!\{t\}\!\}$ with $\val(t)=1$. This is an algebraic closure of $\CC(\!(t)\!)$. Note that it is not complete.  
    \item The rational function field $\CC(t)\subset \CC(\!(t)\!)$ with the valuation induced from $\CC(\!(t)\!)$. 
    \end{enumerate}
\end{remark}   
\begin{definition}[Toric varieties]\label{def:ToricVarieties}
Let $M\simeq \ZZ^{n}$ be a lattice of rank $n$ with group ring $A=K[M]$ and torus $T=\Spec(A)$. We write $N :={\rm Hom}(M,\ZZ)$ for the cocharacter lattice. Let $\Delta\subset N_{\RR} := N\otimes_{\ZZ}\RR$ be a rational polyhedral fan. Then to each cone $\sigma\in\Delta$, we can associate a variety $Y_{\sigma} := \Spec(K[\sigma^{\vee}\cap{M}])$ where
\begin{equation*}
        \sigma^{\vee} := \big\{v\in{M}_\RR : \langle v,w \rangle \geq 0\text{ for all }w\in\sigma \big\}
\end{equation*}
is the dual cone to $\sigma$. For any inclusion of cones $\sigma\subset \tau$ in $\Delta$, we obtain open immersions $Y_{\sigma}\to Y_{\tau}$. Gluing all the varieties $Y_{\sigma}$ along these immersions  yields the toric variety $Y(\Delta)$. We write $Y = \Spec(K[M])\subset Y(\Delta)$ for its dense torus.         
\end{definition} 

We now define extended tropicalizations. Our treatment here includes tropicalizations of closed subvarieties of $Y(\Delta)$, rather than just the dense torus $Y$, see \cite[Section 3]{Payne09}. 

\begin{definition}
      [Tropicalization]
      Let $K$ and $\Delta$ be as in Definition~\ref{def:ToricVarieties}. For any $\sigma\in\Delta$, we define the local tropicalization map
      \begin{equation*}
          Y_{\sigma}(K)={\rm Hom}_{{\rm Mon}}(\sigma^{\vee}\cap{M},K)\to {\rm Hom}_{{\rm Mon}}(\sigma^{\vee}\cap{M},\RRi)
      \end{equation*}
      by composing $\psi:\sigma^{\vee}\cap{M}\to K$ with the valuation map $\val:K\to \RRi$. Here $\RRi$ is considered as a monoid with respect to tropical multiplication. By gluing these topological spaces ${\rm Hom}_{{\rm Mon}}(\sigma^{\vee}\cap{M},\RRi)$ together (here we use the order topology on $\RRi$), we obtain the space  $\trop(Y(\Delta))$. We have a natural global tropicalization map 
      \begin{equation*}
          Y(\Delta)(K)\to \trop(Y(\Delta)).
      \end{equation*}
      Note that this is in general not surjective. We write $\trop(Y)$ for the tropical torus, which can be identified with $\RR^{n}$ after choosing a basis of $M$.  

      Suppose that $X\subset Y(\Delta)$ is a closed subvariety. Let $K \subset K'$ be a valued field extension, where $K'$ is algebraically closed and the valuation on $K'$ is not trivial.  We define the tropicalization of $X$ to be the closure of $\trop(X(K'))$ inside $\trop(Y(\Delta))$. This is independent of the field extension $K'\supset K$. Similarly, suppose that $X$ is a closed subvariety of the torus $Y$. We define the tropicalization of $X$ to be the closure of $\trop(X(K'))$ inside $\trop(Y)$.
\end{definition}

\subsection{Initial degenerations and compactifications}\label{sec:InitialDegenerations}
We now recall initial degenerations of ideals and review their relation to compactifications.

\begin{definition}
    Let $K := \CC\{\!\{ t\}\!\}$ be the field of Puiseux series over $\CC$. Let $M$ and $N$ be as in Definition~\ref{def:ToricVarieties}. 
    For any weight vector $w \in N_{\RR}$ the initial form of a Laurent polynomial $f = \sum\limits_{m \in{M}} c_{m} \chi^{m} \in K[M]$ with respect to $w$ is defined by
    \begin{equation*}
        \mathrm{in}_{w}(f) := \sum_{ \substack{m \in M \text{ such that} \\ \val(c_m) + \langle m, w \rangle \text{ is minimal} } }  \overline{t^{-\val(c_m)} c_m} \ \chi^m  \in \CC[M].
    \end{equation*}
    For an ideal $I \subset K[M]$, we define the initial ideal with respect to $w\in N_{\RR}$ to be the ideal generated by all the initial forms $\mathrm{in}_{w}(f)$ for $f\in{I}$. 
    We denote it by ${\rm in}_{w}(I)\subset \CC[M]$. 
    If $X=V(I)$, then the vanishing set of the initial ideal is denoted by $\mathrm{in}_{w}(X)$. If $K$ is one of the other valued fields in Remark~\ref{rem:ValuedFields}, then we define the initial form of a polynomial $f\in K[M]$ at $w\in N_{\RR}$ to be the initial form in $K'[M]$, where $K'$ is the field of Puiseux series. We similarly define initial ideals. 
\end{definition}

    \begin{remark}
        One can also define the initial ideal in terms of models as follows. For simplicity, suppose $w=0$, we can then consider the ideal in $R[M]$ given by
        \begin{equation*}
            \mathcal{I}:=I\cap R[M].
        \end{equation*}
        The initial ideal ${\rm in}_{w}(I)$ is then $\mathcal{I}\otimes_{R[M]}\CC[M]$, which is naturally an ideal in $\CC[M]$. 
        Here we used the natural residue map $R[M]\to \CC[M]$. This is the point of view taken in \cite[Section 5]{Gubler2013}. 
    \end{remark}

    \begin{definition}\label{def:GroebnerComplex}
        Let $K$ be one of the valued fields in Theorem~\ref{rem:ValuedFields} and let $I\subset K[M]$ with $V(I)\subset \Spec(K[M]) =: Y$. For $w\in N_{\RR}$ we define the following subset of $N_{\RR}$
        \begin{equation*}
            C_{I}[w] = \Big\{ w' \in N_{\RR} : \mathrm{in}_{w}(I)=\mathrm{in}_{w'}(I)\Big\}.
        \end{equation*}
        The $\overline{C_{I}[w]}$ for varying $w\in N_{\RR}$ form a polyhedral subdivision of $N_{\RR}$.
        By considering $w\in \trop(V(I))$, this endows
        $\trop(V(I))$ with the structure of a polyhedral complex. This is the Gr\"{o}bner complex structure. If $K$ is trivially valued, then the various $\overline{C_{I}[w]}$ define a complete fan in $N_{\RR}$. We call this the Gr\"{o}bner fan. A polyhedral structure $\Sigma$ on $\trop(X)$ is \emph{admissible} if it refines the Gr\"{o}bner complex. We will say that the polyhedral structure is unimodular if for every polyhedron $\sigma\in\Sigma$, the primitive ray generators of the homogenized cone 
        \begin{equation*}
            \overline{\RR_{>0}(\sigma\times\{1\})}\subset N_{\RR}\times \RR
        \end{equation*}
        form a basis of the lattice $N\times \ZZ$. 
    \end{definition}
    
    \begin{remark}
        Our reason for introducing the Gr\"{o}bner complex here is that the induced tropical compactifications satisfy various pleasant properties, see \cite{Gubler2013,GRW16,GRW17}. In particular, they are necessary for connecting our rational tropical manifolds back to algebraic geometry, see Proposition~\ref{pro:StrictlySemistableModels} and Theorem~\ref{thm:KNFunctorHomotopyColimit3}. 
    \end{remark}

    \subsubsection{The trivially valued case}\label{sec:TriviallyValuedCaseInitials}
    When $K=\CC$, one can also interpret the initial degenerations without the parameter $t$. 
    For more details, see \cite[Section 3.1]{Fulton93} and  \cite{NS13PhaseLimit}.
   
    Let $I\subset \CC[M]$ be an ideal with variety $X=V(I)$. We choose a fan structure $\Delta$ on $N_{\RR}$ that refines the Gr\"{o}bner fan of $I$, as in Definition~\ref{def:GroebnerComplex}. We write $Y(\Delta)$ for the induced toric variety and $\overline{X}$ for the closure of $X$ inside $Y(\Delta)$. For any cone $\tau\in \Delta$, we write
    \begin{align*}
        \tau^{\perp}& :=\{v\in{M}_\RR : \langle v,w \rangle = 0 \text{ for all }w\in\tau\}.
    \end{align*}
    We have a morphism of semigroups 
    \begin{equation*}
        \tau^{\vee}\to \tau^{\perp}, \quad v \mapsto \begin{cases}
            v, \quad \text{if } v \in \tau^\perp \\
            0, \quad \text{otherwise }
        \end{cases}.
    \end{equation*}
    We define $I(\tau)$ to be the image of the ideal $I \cap \CC[\tau^{\vee} \cap M]$ under the map 
    \[
        \CC[\tau^{\vee}\cap M]\to \CC[\tau^{\perp}\cap M]
    \]
    induced by the map of semigroups discussed above. If we write $\mathcal{O}_{\tau}=\Spec(\CC[\tau^{\perp}\cap M])$ for the open torus orbit corresponding to $\tau$, then $I(\tau)$ generates the ideal of $\overline{X} \cap \mathcal{O}_{\tau}$ inside $\mathcal{O}_{\tau}$, so that $\overline{X}\cap \mathcal{O}_{\tau}=V(I(\tau))$.    
    \begin{lemma}
        The initial ideal $\mathrm{in}_{w}(I)$ for $w \in \relint(\tau)$ coincides with the ideal generated by the image of $I(\tau)$ under the inclusion map $\CC[\tau^{\perp}\cap M]\to \CC[\tau^{\vee}\cap M]$. 
    \end{lemma}
    \begin{proof}
        See \cite[Lemma 12]{NS13PhaseLimit}.
    \end{proof}
    
    We now describe the connection between the initial ideal and the intersection $\overline{X}\cap \mathcal{O}_{\tau}$. The torus $\mathbb{T}=\Spec(\CC[M])$ acts on  $\mathcal{O}_{\tau}$. We define $\mathbb{T}_{\tau}$ to be the sub-torus of $\T$ that acts trivially on $\mathcal{O}_{\tau}$. 
    Explicitly ${\mathbb{T}}_{\tau}=\Spec(\CC[M_{\tau}])$, where $M_{\tau} := M / (\tau^{\perp}\cap M)$. It is clear that $M_{\tau}$ is indeed a lattice, because $\tau^{\perp}\cap M$ is saturated. Note that this torus fits in an exact sequence
      \begin{equation*}
            1\to \mathbb{T}_{\tau}\to \mathbb{T}\to \mathcal{O}_{\tau}\to 1
      \end{equation*}
    of algebraic groups. Moreover, ${\mathbb{T}}_{\tau}$ acts on $\mathrm{in}_{w}(X)$ and the quotient by this action is exactly $V(I(\tau))\subset \mathcal{O}_{\tau}$ by the following:
            
    \begin{lemma}\label{lem:InitialDegenerationsTorusQuotient}
        Let $\mathbb{T}_{\tau} = \Spec(\CC[M_{\tau}])\subset \mathbb{T}$ be the torus that acts trivially on the torus orbit $\mathcal{O}_{\tau}$, so that $\mathbb{T}/\mathbb{T}_{\tau}\simeq \mathcal{O}_{\tau}$. 
        The action of the torus $\mathbb{T}$ on itself induces an action of $\mathbb{T}_{\tau}$ on $\mathrm{in}_{w}(X)\subset \mathbb{T}$, and this action is free. Moreover, the quotient $\mathrm{in}_{w}(X)/\mathbb{T}_{\tau}$  can be identified with $\overline{X}\cap \mathcal{O}_{\tau}=V(I(\tau))$.
    \end{lemma}
    \begin{proof}
        See \cite[Corollary 13]{NS13PhaseLimit}. 
    \end{proof}
    
    We can thus view the initial ideal as a twisted version of the local ideal $I(\tau)$ of $\overline{X}$ in~$\mathcal{O}_{\tau}$. 
    
    \begin{remark}
        By choosing a splitting of 
        \begin{equation*}
            1\to \mathbb{T}_{\tau}\to \mathbb{T}\to \mathcal{O}_{\tau}\to 1,
        \end{equation*}
        we can write 
        \begin{equation*}
            \mathrm{in}_{w}(X)\simeq (\mathbb{C}^{\times})^{\mathrm{dim}(\tau)}\times (\overline{X}\cap \mathcal{O}_{\tau}).
        \end{equation*}
    \end{remark}

\begin{definition}
    Let $\tau \in\Delta$. Let $\mathcal{O}_{\tau}$ be the corresponding torus orbit. We write
    \[
    \overline{\mathcal{O}}_{\tau} = \bigsqcup_{\sigma \supset \tau} \mathcal{O}_{\sigma},
    \]
    for the closure of $\mathcal{O}_{\tau}$ in $Y(\Delta)$. 
    Suppose that $\tau$ intersects $\trop(X)$. Then we write $Z_{\tau} =  \overline{\mathcal{O}}_{\tau} \cap \overline{X}$ for the corresponding closed subset of $\overline{X}$.
\end{definition}
 
\subsection{Rational tropical manifolds}

In this section, we define the notion of an embedded rational tropical manifold. We start by recalling various existing notions of smoothness in tropical geometry.
See \cite{IKMZ19}, \cite{HK12}, \cite{GRW17} and \cite{AP20} for more background information.   

\begin{definition}\label{def:Schon}
    Let $K$ be the field of complex Puiseux series. 
    Let $X = V(I)$ be a closed subscheme of the torus $Y := \Spec(K[M])$ with ideal $I$ and let $w\in N_{\RR}$ be a weight vector with initial ideal $\mathrm{in}_{w}(I)$. We say that $X$ is sch\"{o}n at $w$ if the variety $V(\mathrm{in}_{w}(I))$ is smooth over $\CC$. We say that $X$ is sch\"{o}n if it is sch\"{o}n at every $w\in   N_{\RR}\simeq \RR^{n}$. 
\end{definition}

\begin{definition}[Tropical smoothness]\label{def:FirstDefinition}\label{def:TropicalSmoothness}
Let $X=V(I)\subset \Spec(\CC[M])$ be a very affine variety. We say that $X$ is a \emph{very affine linear space} if the ideal $I$ 
is \emph{affine linear} with respect to a basis of $M$. That is, there is a basis $m_{1},...,m_{n}\in {M}$ and a set of generators $f_{1},...,f_{r}$ of $I$ such that 
    \begin{equation*}
        f_{i}=c_{i,0}+\sum_{j=1}^{n} c_{i,j} \chi^{m_{j}} \in \CC[M], \quad \text{for some } c_{i,j} \in \CC, \text{ where }1\leq{i}\leq{r} \text{ and }0\leq{j}\leq{n}.
    \end{equation*}
    Let $X=V(I)\subset \Spec(K[M])$ be a very affine variety over the field of complex Puiseux series. We say that $X$ is tropically smooth at $w\in N_{\RR}$ if $V(\mathrm{in}_{w}(I))$ is a very affine linear space. 
     We say that $X$ is tropically smooth if it is tropically smooth at every $w\in N_{\RR}$. In particular, the base change of a very affine linear space $X=V(I)\subset \CC[M]$ over $\Spec(K)\to \Spec(\CC)$ is tropically smooth everywhere. We will also call these base changes very affine linear spaces.      
\end{definition}
\begin{remark}
    In particular, the toric variety $X=\Spec(K[M])=V(\langle 0\rangle)$ and the empty variety $\varnothing = V(\langle 1\rangle)$ are very affine linear.
\end{remark}

\begin{remark}\label{rem:EquivalenceVALSandEPHC}
    By \cite[Section 4.1]{MS15}, every $m$ dimensional affine linear space in the torus $(\CC^\times)^n$ is the complement of an essential projective hyperplane arrangement $D$ in $\mathbb{P}^{m}$, where $D=\bigcup_{i=0}^{n} D_{i}$ and $D_{i}=V(H_{i})$ for some homogeneous $H_{i}$. To illustrate this, consider $ \CC[\ZZ^{4}]$, and write $e_{i}$ for $i=0,1,2,3$ for the generators of $\ZZ^{4}$. We now consider the homogeneous ideal $J = \langle \chi^{e_{0}}+\chi^{e_{1}}+\chi^{e_{2}}+\chi^{e_{3}}\rangle \subset \CC[\ZZ^4]$. We now consider $I=J\cap \CC[M]$. To find an explicit presentation, we dehomogenize with respect to $e_0$ and obtain $I = \langle 1 + x + y +z \rangle \subset \CC[M]$, where 
    \[ 
     M := \{(a,b,c,d)\in \ZZ^{4}:a+b+c+d=0\},
    \quad x = \chi^{e_{1}-e_{0}}, \quad y = \chi^{e_{2}-e_{0}},
    \quad \text{and} \quad z=\chi^{e_{3}-e_{0}}.
    \] 
    The map
    \begin{equation*}
        \mathbb{P}^{2}\to V(J) \subset \mathbb{P}^{3}, \quad [X:Y:Z]\mapsto [X:Y:Z:-X-Y-Z].
    \end{equation*}
    gives an isomorphism $\mathbb{P}^{2}\backslash D\to V(I)$ for $D = V\left(XYZ(-X-Y-Z)\right)$.
\end{remark}

\begin{corollary}\label{cor:SchoenImpliesSmooth}
    A sch\"{o}n variety $X$ is smooth over $K$. In particular, tropically smooth varieties $X$ are smooth over $K$.  
\end{corollary}
\begin{proof}
    This follows from \cite[Proposition 2.7 and Lemma 2.9]{KatzStapledon12}. Note that affine linear varieties are automatically smooth over $K$, so that tropically smooth varieties are sch\"{o}n.   
\end{proof}

We now give a different, perhaps more familiar definition of a smooth tropicalization in terms of Bergman fans of matroids. For any rational polyhedron $\sigma \subset N_\RR$ we denote by $N(\sigma)$ the linear span of $\sigma$ i.e. the vector space
    \[
        N_\RR(\sigma) := {\rm span}_\RR( \{y - x \colon y \in \sigma\}), \quad \text{for any } x \in \sigma.
    \]
    Note that $N_\RR(\sigma)$ does not depend on the choice of $x \in \sigma$. We also write $N(\sigma) := N \cap N_{\RR}(\sigma)$.
    
\begin{definition}\label{def:SecondDefinition}
    A tropical variety in $N_{\RR}$ of dimension $d$ consists of a pair $(\Sigma, w)$  of a pure $d$-dimensional rational polyhedral complex $\Sigma$ and a weight function
    \[
        w \colon \{ \sigma \in \Sigma \colon \dim(\sigma) = d \} \longrightarrow \ZZ_{>0}, \quad \sigma \longmapsto w(\sigma),
    \]
    such that $\bm{\Sigma}$ is \emph{balanced}. That is, 
    for any $\tau \in \Sigma$ such that $\dim(\tau) = d-1$ we have
    \[
        \sum_{\substack{\sigma \supset \tau \\ \dim(\sigma) = d}} w(\sigma) \  u_{\sigma / \tau} \in N_\RR(\tau),
    \]
    where $u_{\sigma / \tau} \in N$ are any vectors such that $N(\sigma) = N(\tau) \oplus \ZZ u_{\sigma / \tau}$ for all $\sigma \supset \tau$ in $\Sigma$.
    
    We say that a tropical variety $(\Sigma, w)$ is smooth if its facets have weights equal to $1$ and $|\Sigma|$ is locally isomorphic to the support of the Bergman fan of a loopless matroid. A tropical variety $(\Sigma, w)$ is realizable if it is the tropicalization of a closed subscheme $X$ of the torus $Y$.
\end{definition}
\begin{proposition}\label{prop:EquivalentDefinitions}
    Suppose that $\Sigma$ is realizable by an algebraic variety $X$. Then the definitions of tropical smoothness given in Definition~\ref{def:FirstDefinition} and Definition~\ref{def:SecondDefinition} coincide.   
\end{proposition}
\begin{proof}
    If $\Sigma$ is realizable by a tropically smooth variety in the sense of Definition~\ref{def:FirstDefinition}, then the local structure of the tropical variety is of the desired type by \cite[Chapter 4]{MS15}. Conversely, suppose that $\Sigma$ is locally the Bergman fan of a loopless matroid. By \cite[Proposition 4.2]{KatzPayne2011}, we know that the closure of $V(\mathrm{in}_{w}(I))$ in projective space is a linear space. In particular, the initial ideal is affine linear.
\end{proof}

\subsubsection{Kummer coverings and rational tropical manifolds}

We now define Kummer coverings of very affine linear spaces. By gluing these together, we obtain the notion of a \emph{rational tropical manifold}, which can be seen as a tropical manifold with a certain nilpotent structure.

\begin{definition}[Kummer coverings] \label{def:KummerSmooth}
    Let $M$ be a lattice and let $\Lambda\subset M$ be a subgroup of finite index. We have a map of tori 
    \[
    \kappa_{M,\Lambda}:\Spec(\CC[M])\to \Spec(\CC[\Lambda]).
    \]
    Let $X\subset \Spec(\CC[M])$. We say that $X$ is a Kummer covering of a very affine linear space if there exists a lattice inclusion $\Lambda\subset M$ and a very affine linear space $Z\subset \Spec(\CC[\Lambda])$ such that $X=\kappa^{-1}_{M,\Lambda}(Z)$.
\end{definition}

\begin{remark}
    There is an action of the finite group scheme $\Spec(\CC[M/\Lambda])$ on $\Spec(\CC[M])$ and $X$. Their respective quotients are $\Spec(\CC[\Lambda])$ and $Z$. Since we are in characteristic zero, these maps are finite \'{e}tale, so they are Galois coverings with Galois group $M/\Lambda$. The degree of this covering is $[M:\Lambda]$.
\end{remark}
    
    \begin{example}\label{exa:KummerCoverings}\label{exa:KummerCoverings2}
    Consider the lattice $\Lambda=\ZZ^{2}$ with polynomial 
$f=1+\chi^{e_{1}}+\chi^{e_{2}}$
    and very affine linear space 
    \[ 
        Z=V(f)\subset \Spec(\CC[\Lambda]). 
    \]
    If we write $\chi^{e_{1}}=x$ and $\chi^{e_{2}}=y$, then $Z=V(1+x+y)\subset \Spec(\CC[\Lambda])$. 

    Now let $m\in \mathbb{N}_{\geq{1}}$ and consider the lattice $M=\dfrac{1}{m}\ZZ^{2}$. The inverse image $Z_{m}$ is now $V(1+z^{m}+w^{m})$, where $z=\chi^{e_{1}/m}$ and $w=\chi^{e_{2}/m}$. This gives the standard Fermat hypersurface. The map $Z_{m}\to Z$ is induced from the maps 
    \begin{equation*}
        \CC[\Lambda]\to \CC[M].
    \end{equation*}
    In local coordinates, this sends $x$ to $z^m$ and $y$ to $w^{m}$. 
    
        For a second example, consider the lattice $M=\langle(1/2,1/2),(0,1),(1,0)\rangle$. Let $v_{1}=(1/2,1/2)$ and $v_{2}=(0,1)$. Note that these form a basis of $M$. Write 
        \[
        s_{1}=\chi^{v_{1}},\quad s_{2}=\chi^{v_{2}}.
        \]
        We then have that $f=1+\chi^{e_{1}}+\chi^{e_{2}}=1+s_{1}^{2}/s_{2}+s_{2}$ and  $X=V(f)\subset \Spec(\CC[M])\simeq\Spec(\CC[s^{\pm}_{1},s^{\pm}_{2}])$. Note that there is a difference between considering $V(f)$ in $\Spec(\CC[\Lambda])$ and $\Spec(\CC[M])$. 
    \end{example}
    
    \begin{remark}\label{rem:TropicalizationKummerCoveringInvariant}
        Note that the underlying set of the tropicalization of a Kummer covering is the same. In terms of groups, we have restriction maps 
        \begin{equation*}
            {\rm Hom}(M,\RR)\to {\rm Hom}(\Lambda,\RR)
        \end{equation*}
        The subsets $\trop(X)$ and  $\trop(Z)$ are  identified under these maps, since they are fans. The multiplicities however necessarily change by \cite[Proposition 4.2]{KatzPayne2011}. It is thus natural to view Kummer coverings as tropical analogues of non-reduced schemes.  
    \end{remark}

\begin{remark}\label{rem:InitialDegenerationsKummer}
    Note that initial degenerations of Kummer coverings are again Kummer. Indeed, we can write down a tropical basis of the corresponding very affine linear space. The injection of these polynomials into $\Lambda$ then gives the corresponding tropical basis for $X$. Here we used Remark~\ref{rem:TropicalizationKummerCoveringInvariant}.     
\end{remark}

    We now use this to define the notion of a rational tropical manifold.
    \begin{definition}[Rational tropical manifolds]
        Let $I\subset K[M]$ be an ideal with vanishing set $X=V(I)$. We say that $X$ is rationally smooth at $w\in N_{\RR}$ if $\mathrm{in}_{w}(X)\subset \Spec(\CC[M])$ is a Kummer covering of a very affine linear space, as in Definition~\ref{def:KummerSmooth}. We say $X$ is rationally smooth if it is rationally smooth at every $w\in N_{\RR}$. 

        A rational tropical manifold $(X, \Sigma,D)$ consists of the following:
        \begin{enumerate}
            \item A rationally smooth $X=V(I)\subset Y=\Spec(K[M])$. 
            \item A choice of a unimodular admissible polyhedral structure $\Sigma$ on $\trop(X)$ with smooth toric variety $Y(\Delta)$. See Definition~\ref{def:GroebnerComplex}. Here $\Delta$ is the recession fan of $\Sigma$.
            
            \item A subset $D = \{D_1, \dots, D_s\}$ of irreducible components of $Y(\Delta) \setminus Y$.
        \end{enumerate}
       If $X$ is tropically smooth, then we say that $(X, \Sigma, D)$ is a tropical manifold.
    \end{definition}
    \begin{remark}
        We emphasize here that rational tropical manifolds in this paper are \emph{embedded} and \emph{realizable}.
    \end{remark}

    \begin{definition}
    \label{def:OpenVsClosed}
    Let $(X,\Sigma,D)$ be a rational tropical manifold with $D=\{D_{1},...,D_{s}\}$ and corresponding rays $\rho_{i}\in \Delta$. The choice of a polyhedral structure $\Sigma$ gives rise to a compactification $\overline{X}\subset Y(\Delta)$ with tropicalization $\overline{\Sigma}$. We define 
    \[
        X_{D} := \overline{X} \Big\backslash \bigcup_{i=1}^{s}D_{i}.
    \]
    This is the partial compactification of $X$ induced by $D$. There are two special cases: 
\begin{enumerate}
   \item $D=D_{{\rm full}}$, the set of all irreducible components of $Y(\Delta)$. In this case $X_{D}=X$. We denote the corresponding rational tropical manifold by $(X,\Sigma,D_{{\rm full}})$.   
   \item $D = \varnothing$. In this case, $X_{D}=\overline{X}$. We will write the corresponding rational tropical manifolds as $(X,\Sigma,\varnothing)$.
\end{enumerate}
Note that $D$ encodes the components at infinity that are removed from $\overline{X}$.
    
    \end{definition}

\begin{remark}
    Note that rationally smooth $X$ are automatically sch\"{o}n, and thus smooth over $K$ in the ordinary sense, see \cite[Proposition 2.7 and Lemma 2.9]{KatzStapledon12}. Here we used that the residue characteristic of $K$ is zero. In residue characteristic $p$, it is better to demand the degrees of the Kummer coverings to be coprime to $p$.     
\end{remark}

\begin{example}
    Suppose that we take $I = (0)$. Then $X=Y$ is the dense torus, and $\Sigma$ corresponds to a choice of a smooth toric compactification $\overline{X}\supset X$.
    In particular, the tropicalization of any smooth toric variety is a tropical manifold in our sense.      
\end{example}

\begin{example}\label{exa:FermatFamily1}
    Consider the ideal $I = \left \langle 1+\sum_{i=1}^{n}x^{m}_{i} \right \rangle$ for the Fermat hypersurface of degree $m\in\mathbb{N}$. Then $X := V(I)$ is rationally smooth and the tropicalization is dual to an $n$-simplex. Note that $\Delta$ can be taken to be the
    fan of $\mathbb{P}^{n}$. Note that this variety is not tropically smooth.
\end{example}

\begin{example}\label{exa:DworkFamily1}
    Consider the Dwork family $I=\langle f\rangle $ where 
    \begin{equation*}
        f=\prod_{i=1}^{n}x_{i}+t+t\sum_{i=1}^{n}x^{m_i}_{i}
    \end{equation*}
    for some $m_i>1$. To see that $X=V(I)$ is rationally smooth,
    first note that the tropicalization is dual to a cellulation of the Newton polytope with central vertex $(1,...,1)$. The statement then follows from our upcoming Lemma~\ref{lem:KummerLemma}. 

    For instance, consider the Dwork elliptic pencil $f=xy+t(1+x^3+y^3)$. 
    Then $\trop(V(f))$ consists of a triangle with vertices $w_{1} =(0,1)$, $w_{2} =(1,0)$ and $w_{3} = (-1,-1)$. 
    For instance, we have 
    \[
        \mathrm{in}_{w_{1}}(f) = 1 + xy + x^3.
    \] 
    Consider the lattice $\Lambda=\langle (1,1),(3,0)\rangle\subset M$. We can view $\mathrm{in}_{w_{1}}(f)$ as an element of $\CC[\Lambda]$. Note that $V(\mathrm{in}_{w_{1}}(f))\subset \Spec(\CC[\Lambda])$ is a very affine linear space. Its inverse image gives $V(\mathrm{in}_{w_{1}}(f))\subset \Spec(\CC[M])$, which is then a Kummer covering of a very affine linear space. Note that the degree is $3$.
\end{example}

\begin{example}\label{exa:ZeroDimensionalRTM}
    Suppose that $X$ is a zero-dimensional rational tropical manifold. Then locally at $w\in N_{\RR}$, we have that $X$ is a Kummer covering of a linear space consisting of a point. The degree of the Kummer covering is then exactly the index $[M:\Lambda]$ of the corresponding lattice map $M\supset \Lambda$. This in turn is exactly the number of geometric points in the initial degeneration (note that there are no multiplicities since we have an \'{e}tale covering of a point).
\end{example}

Regular unimodular triangulations of lattice polytopes give rise to smooth tropical hypersurfaces by duality. For rational tropical manifolds, this turns into a statement about arbitrary lattice polytopes.

\begin{proposition}
    Let $S\subset M$ be a finite set and consider the space $K[S]\subset K[M]$ of all Laurent polynomials with fixed monomial supports $S\subset M$. Let $\mathcal{C}$ be the convex hull of $S$.  Then there is an $f$ in $K[S]$ with Newton polytope $\mathcal{C}$ whose vanishing set $X=V(f)$ is a rational tropical manifold with respect to some polyhedral structure.
\end{proposition}
\begin{proof}
    We can find a regular triangulation of the convex hull $\mathcal{C} = {\rm conv}(S)$ by \cite[Proposition 2.2.4]{LRS10}. Consider a simplex in the triangulation. We can place one of the vertices at $0$, and the other vertices then give rise to a basis of the rational span of $\mathcal{C}$ inside $M\otimes_{\ZZ}\QQ$. We conclude using Lemma~\ref{lem:KummerLemma}.
\end{proof}
\begin{lemma}\label{lem:KummerLemma}
    Suppose $m_{0},...,m_{n}\in M$ are such that the $m_{i}-m_{0}$ are linearly independent over $\QQ$.
    Then $V(f)\subset \Spec(\CC[M])$ for $f=\sum_{i=0}^{n}\chi^{m_{i}}$ is a rational tropical manifold.   
\end{lemma}
\begin{proof}
We can assume without loss of generality that the $m_{i}-m_{0}$ span $M\otimes_{\ZZ}\QQ$. Consider the homomorphism
\begin{equation*}
    \psi:\ZZ^{n+1}\to M
\end{equation*}
given by $\psi(e_{i})=m_{i}$.\footnote{This will be formalized in terms of affine simplices in Definition~\ref{def:AffineSimplex}.} Note that the induced map $\Xi:={\rm ker}(\mathds{1}_{n+1})\to \Lambda=\psi(\Xi)$ is an isomorphism. Moreover, $\Lambda$ is of finite index in $M$. In particular, $Z=V(f)\subset \Spec(\CC[\Lambda])$ is a very affine linear space and $X=V(f)\subset \Spec(\CC[M])$ is a Kummer covering of $Z$.
\end{proof}
\begin{remark}
    These results explain why we call these objects rational tropical manifolds. Namely, rather than being related to bases of the lattice $M$, they are related to bases of the rational vector space $M\otimes_{\ZZ}\QQ$.  
\end{remark}

\begin{example}\label{rem:DifferentKummerPresentations}
    We remark here that there might be multiple Kummer presentations of a given ideal $I$ with different Kummer degrees. For instance, consider the ideal 
    \begin{equation*}
        I=\langle x^2+y^5+z^7\rangle\subset \CC[\ZZ^{3}]
    \end{equation*}
    Note there is the obvious Kummer covering from $V(I)$ to $V(z_{0}+z_{1}+z_{2})$ of degree $70$. On the other hand, consider the elements $z_{0}=y^{5}/x^2$, $z_{1}=z^7/x^2$, $z_{2}=x/yz^2$ of $\CC[M]$. The corresponding vectors in $M$ span $M$. Moreover, we can express $I$ as 
    \begin{equation*}
        I=\langle 1+z_{0}+z_{1}\rangle. 
    \end{equation*}
    We thus see that $I$ is in fact a tropical manifold, although this is not quite clear from the given presentation.
\end{example}

In general, the stable intersection of two tropical manifolds is not a tropical manifold. For rational tropical manifolds, the stable intersection however is again a rational tropical manifold by the following.

\begin{proposition}\label{pro:RationalManifoldStableIntersection}
Suppose that $X_{1},..,X_{r}$ are rational tropical manifolds that intersect transversely. Then $\bigcap_{i=1}^{r}X_{i}$ is a rational tropical manifold. 
\end{proposition}
\begin{proof}
We use the definition of transversality from \cite[Section 3]{BJSST07} here. That is, if $\sigma_{1}$ and $\sigma_{2}$ are the cones of $X_{1}$ and $X_{2}$ which contain $w$ in their relative interior, then $\RR{\sigma_{1}}+\RR{\sigma_{2}}=N_{\RR}$. We now recall that the initial ideal of a transversal intersection is given by the sum of the initial ideals, see \cite[Lemma 3.4.10]{MS15}. In other words:
\begin{equation*}
    \mathrm{in}_{w}(I_{1})+\mathrm{in}_{w}(I_{2})=\mathrm{in}_{w}(I_{1}+I_{2}).
\end{equation*}
It now suffices to show that the stable intersection of two tropical manifolds is a rational tropical manifold. We will show that the initial ideals at vertices of the stable intersection are rationally smooth. Note that these vertices arise from pairs $(\sigma_{1},\sigma_{2})$ of polyhedral cells from $\Sigma_{1}$ and $\Sigma_{2}$ whose dimensions add up to $n$, and whose spans intersect transversely in the point. This will be sufficient for showing that the entire tropicalization is rationally smooth, since the initial degeneration of a Kummer covering is again a Kummer covering.

We write $N_{i,\RR}$ for their tangent spaces, so that 
\begin{equation*}
    N_{\RR}=N_{1,\RR}\oplus N_{2,\RR}.
\end{equation*}
Dually, we then find 
\begin{equation*}
    M_{\RR}=M_{1,\RR}\oplus M_{2,\RR}
\end{equation*}
for their dual character vector spaces. We now claim that we can find affine linear presentations for $I_{1}$ and $I_{2}$ whose characters land in $M_{1,\RR}$ and $M_{2,\RR}$ respectively. It suffices to do this for $I_{1}$. Choose a basis $\{v_{i}\}$ of $M$ such that $I_{1}$ is affine linear with respect to this basis. Let $f_{1},...,f_{k}$ be the corresponding generators. We now let $B$ be the $k\times (n+1)$-matrix that describes the homogenizations of the $f_{j}$ with respect to the $\{v_{i}\}$, as in \cite[Section 4.1]{MS15}. That is, the rows describe the different $f_{j}$ as follows: the 
zeroth column records the constant coefficient of the $f_{j}$, and the $i$-th column records the coefficient of $f_{j}$ with respect to $\chi^{v_{i}}$.      
We can assume without loss of generality that the first entry $a_{01}$ is nonzero. Indeed, if there is a nonzero constant coefficient, then we can permute the $f_{j}$. If there is no nonzero coefficient, then we can divide by some $\chi^{v_{i}}$ to obtain the desired nonzero coefficient. Note that this does not change the ideal. We can now apply elementary row operations to make all of the entries in the first column nonzero. This gives generators $g_{1},...,g_{k}$ of the same ideal. 

Now note that for any monomial $\chi^{m}$ appearing in $g_{i}$, we have that $\chi^{m}$ is invariant under $N_{1,\RR}$. Indeed, otherwise we would obtain a contribution of some power of $t$, which we cannot divide by since every generator has a nonzero constant term. Since $\chi^{m}$ is invariant, we conclude that $\langle m,n\rangle=0$ for every $n\in N_{1,\RR}$. In other words, $m\in M_{1,\RR}$.

We now have equations for $I_{1}$ and $I_{2}$ in terms of monomials that live in distinct vector spaces. Taking the basis for $I_{1}$ and $I_{2}$, we then find a Kummer covering as in Lemma~\ref{lem:KummerLemma}. Explicitly, if $g_{i,k}=c_{i,0,k}+\sum_{j=1}^{m} c_{i,j,k}\chi^{m_{j}}$, 
then we consider the linear space 
\begin{equation*}
Z=V(\tilde{g}_{i,1},\tilde{g}_{i,2}),
\end{equation*}
where $\tilde{g}_{i,k}=c_{i,0,k}+\sum c_{i,j,k}z_{j}$. Note that the $c_{i,0,k}$ are nonzero by construction. Since the tropicalizations intersect transversely, we find that the monomial supports of the given generators of $I_{1}$ and $I_{2}$ are disjoint (apart from the constant term). As in Lemma~\ref{lem:KummerLemma}, we then find a Kummer covering that maps to both.

\end{proof}

\begin{remark}
The stable intersection of two tropical varieties can be seen as a tropical variant of a generic intersection. The statement of Proposition~\ref{pro:RationalManifoldStableIntersection} can then be seen as an analogue of 
   Bertini's theorem in algebraic geometry for rational tropical manifolds. 
\end{remark}

\begin{example}
Suppose that the intersection of the rational tropical manifolds is zero-dimensional. By Proposition~\ref{pro:RationalManifoldStableIntersection}, we then have a zero-dimensional rational tropical manifold. The local Kummer degree is then exactly the tropical intersection multiplicity by Example~\ref{exa:ZeroDimensionalRTM}. One can use this to give a quick proof of Bernstein's theorem, which says that the generic number of intersection points of $n$ hypersurfaces with fixed monomial supports is the mixed volume of the Newton polytopes. See \cite{MS15} and \cite{HR2025} for similar tropical approaches. We will also review Bernstein's theorem in Section~\ref{sec:ClassicalExamplesD01} to see that it gives an example of a zero-dimensional polar decomposition. 
\end{example}

\subsubsection{Regular CW complex structure on rational tropical manifolds}

We now endow $\overline{\Sigma}$ with the structure of a regular CW complex for rational tropical manifolds. This in particular allows us to interpret the cells in $\overline{\Sigma}$ in terms of a poset for any rational tropical manifold. This poset will be used later on as the base category for our diagram of spaces. We start by defining this poset. 
   
\begin{definition}[Tropical Poset] \label{def:TropicalPoset}
Let $(X,\Sigma,D)$ be a rational tropical manifold with recession fan $\Delta={\rm rec}(\Sigma)$.
    The poset $\Pscr(\overline{\Sigma})$ consists of pairs $(\tau,\sigma)$ of cells $\tau\in \Sigma$ and $\sigma\in \Delta$ such that the recession cone of $\tau$ contains $\sigma$.
    The partial order on $\Pscr(\overline{\Sigma})$ is given by the following:
    \begin{enumerate}
        \item For any $\tau \in \Sigma$ we have $(\tau,\sigma_{1})\leq (\tau,\sigma_{2})$ whenever $\sigma_{2} \subset \sigma_{1}$.
        \item For any $\sigma \in \Delta$ we have $(\tau_{1},\sigma)\geq (\tau_{2},\sigma)$ whenever $\tau_{1}\supset \tau_{2}$.
    \end{enumerate} 
\end{definition}
\begin{remark}
    Note that the poset $\Sigma$ can be identified with the subposet generated by $(\sigma,0)$ in $\Pscr(\overline{\Sigma})$. 
\end{remark}

\begin{proposition}\label{pro:CWComplexTropicalization}
    Let $\Sigma$ be a rational tropical manifold. Then $\Pscr(\overline{\Sigma})$ is the poset of a regular CW complex structure on $\overline{\Sigma}$.
\end{proposition}
\begin{proof}
    Locally, we have that the tropicalization is modeled after the tropicalization of a toric substratum of $\mathbb{A}^{n}$. 
    The strata for these are exactly indexed by the elements of the poset $\Pscr(\overline{\Sigma})$. Since $\overline{\Sigma}$ is constructed by gluing these, we find the desired CW complex structure.    
\end{proof}

\begin{example}
    Consider the fan $\Delta$ associated to $\mathbb{A}^{2}$. Then the cells of $\overline{\Delta}$ are indexed by elements of $\scrP(\overline{\Delta})$ as follows:

    \medskip
    
    \begin{center}
    \begin{tikzpicture}

        \begin{scope}[shift={(-6,0)}]
        \coordinate[nvert] (a)  at (0,0);
        \coordinate        (a1) at (3,0);
        \coordinate        (a2) at (0,3);
        \coordinate        (a3) at (3,3);

        \coordinate[label=below left: $0$] (v) at (0,0);
        \coordinate[label=below: $\rho_1$] (rho1) at (1.5,0);
        \coordinate[label=left:  $\rho_2$] (rho2) at (0,1.5);
        \coordinate[label=below left: $\sigma$] (sigma) at (1.5,1.5);

        \draw[] (a) -- (a1);
        \draw[] (a) -- (a2);

        \path[pattern color = blue, line width = 2, pattern=north west lines, opacity = 0.5] (a) -- (a1) -- (a2) -- plot[domain=1:3.5, smooth] (a) -- cycle;

        \node at (1.5,-1) {$\Delta$};
        
        \end{scope}

        \begin{scope}{shift = {(6,0)}}
            \coordinate[nvert]   (a)  at (0,0);
            \coordinate[nvert]   (a1) at (3,0);
            \coordinate[nvert]   (a2) at (0,3);
            \coordinate[nvert]   (a3) at (3,3);

            \coordinate[label=below left: {$(0,0)$}] (v) at (0,0);
            \coordinate[label=below: {$(\rho_1,0)$}] (rho1) at (1.5,0);
            \coordinate[label=left:  {$(\rho_2,0)$}] (rho2) at (0,1.5);
            \coordinate[label=below left: {$(\sigma, 0)$}] (sigma) at (2,2);
            \coordinate[label=above left:  {$(\rho_2,\rho_2)$}] (rho2) at (0,3);
            \coordinate[label=above:  {$(\sigma,\rho_2)$}] (rho2) at (1.5,3);
            \coordinate[label=above right:  {$(\sigma,\sigma)$}] (rho2) at (3,3);
            \coordinate[label= right:  {$(\sigma,\rho_1)$}] (rho2) at (3,1.5);
            \coordinate[label= below right:  {$(\rho_1,\rho_1)$}] (rho2) at (3,0);

            \draw[] (a) -- (a1);
            \draw[] (a) -- (a2);
            \draw[] (a2) -- (a3);
            \draw[] (a1) -- (a3);

             \path[pattern color = blue, line width = 2, pattern=north west lines, opacity = 0.5] (0,0) -- (3,0) -- (3,3) -- (0,3) -- plot[domain=1:3.5, smooth] (0,0) -- cycle;

             \node at (1.5,-1) {$\overline{\Delta}$};
        \end{scope}

    \end{tikzpicture}
        
    \end{center}
\end{example}

\begin{definition}

    We define $\Pscr(\Sigma_{D})\subset \Pscr(\overline{\Sigma})$ to be the subset of all $(\tau,\sigma)$ such that $\sigma$ does not contain any of the $\rho_{i}$. We call this a partial compactification of $X$ and $\Sigma$.
\end{definition}

\begin{remark}
    In this paper, we will mostly consider the open case $\Pscr(\Sigma_{D}) = \Pscr(\Sigma)$ and the closed case $\Pscr(\Sigma_{D}) = \Pscr(\overline{\Sigma})$. In principle, one could also work with any subset 
    \begin{equation*}
        \Pscr(\Sigma)\subset \Pscr\subset \Pscr(\overline{\Sigma}),
    \end{equation*}
    but we will not elaborate on this.
\end{remark}

\begin{remark}
    Let $\Sigma=\trop(X)$ with compactification $\overline{\Sigma}$. The category associated to the poset  $\Pscr(\overline{\Sigma})$ will serve as the base category for our polar decompositions. In Section~\ref{sec:TropicalDiagrams} we will assign a set of spaces to these cells. On a topological level, these spaces for $(\tau,\sigma)\in \Pscr(\overline{\Sigma})$ are given by the initial degenerations of the ideal $I$ at $w \in {\rm relint}(\tau)$ with respect to the ambient torus $\Spec(\CC[\sigma^{\perp}\cap M])$.

    These initial degenerations unfortunately do not directly give rise to  comparison maps. Moreover, it is not clear exactly how these initial degenerations relate to nearby smooth fibers. To that end, we will use the theory of Kato-Nakayama spaces. This will allow us to find a diagram whose homotopy colimit captures the topology of a nearby smooth fiber.    
    
    To construct explicit CW complex structures on this homotopy colimit, we will replace the Kato-Nakayama spaces by their angle sets in Section~\ref{sec:AngleSets}. To show that no information is lost in this process, we prove that the angle map is a homotopy equivalence for the local Kummer spaces, see Theorem~\ref{thm:MainTheoremAText}. In the remaining sections, we explore how to obtain explicit CW complex structures on these angle sets for complete intersections.      
\end{remark}

\subsection{Semistable models for rational tropical manifolds}\label{sec:ModelsForRTMs}

In this section, we recall how every rational tropical manifold  comes with a natural semistable model after taking a cellulation.

\begin{definition}
    Consider a rational tropical manifold $(X, \Sigma, D)$. By the construction in \cite{Gubler2013}, we obtain a scheme $\mathcal{Y}_{\Sigma}$ over an open subset $U\subset \Spec(\CC[t])$ that contains $t=0$.  
    Its generic fiber is the toric variety $Y(\Delta)$. Here $\Delta$ is the recession fan of $\Sigma$. We take the closure $\mathcal{X}$ of $X$ in $\mathcal{Y}_{\Sigma}$. This is an integral scheme that is flat, proper and of finite presentation over $U$.  We call this the model associated to the rational tropical manifold. 

    Write $\mathcal{D}$ for the closure of the union $\bigcup_{i=1}^{s}D_{i}$ in $\mathcal{X}$. Since $\mathcal{X}$ is regular, we can consider this as an effective Cartier divisor on $\mathcal{X}$. We call the pair $(\mathcal{X}\to U,\mathcal{D})$ a model of the rational tropical manifold $(X,\Sigma,D)$. If $D=D^{\rm{full}}$, then we write $\mathcal{D}^{\rm{full}}$ for the horizontal divisor, and if $D=D^{\varnothing}$, then we write $\mathcal{D}^{\varnothing}$ for the horizontal divisor. 
\end{definition}

\begin{definition}
 Let $\mathcal{X}$ be a regular scheme that is flat, proper and of finite presentation over $R$. We say that $\mathcal{X}$ is strictly semistable if 
 \begin{enumerate}
     \item The generic fiber $\mathcal{X}_{\eta}$ is smooth. 
     \item The special fiber $\mathcal{X}_{s}$ is a strictly normal crossings divisor in $\mathcal{X}$.
 \end{enumerate}
 Let $\mathcal{D}$ be a horizontal effective Cartier divisor on $\mathcal{X}$.
 We say that $(\mathcal{X},\mathcal{D})$ is strictly semistable if $\mathcal{X}_{s}+\mathcal{D}$ is a strictly normal crossings divisor in $\mathcal{X}$. 
\end{definition}

The following useful proposition allows us to construct semistable models of varieties. 
\begin{proposition}\label{pro:StrictlySemistableModels}
Let $(X,\Sigma,D)$ be a rational tropical manifold and let $(\mathcal{X}\to U,\mathcal{D})$ be its associated model. Then $(\mathcal{X} \to U,\mathcal{D})$ is a strictly semistable model. Moreover, we have $\mathcal{X}_{\eta}=\overline{X}$ and $X$ is the complement of the full toric boundary divisor in $\mathcal{X}_{\eta}$.
\end{proposition}
\begin{proof}
Note that since $X$ is rationally smooth, it is automatically sch\"{o}n. The unimodularity then implies that we obtain a regular scheme whose special fiber is a strictly normal crossings divisor, see \cite[Proposition 7.1.1]{BK24}. The proof relies on \cite[Propositions 2.3 and 3.9]{HK12}. We also refer the reader to \cite[Proposition 9.9]{GRW17}. Note that $\mathcal{X}_{s}$ intersects $\mathcal{D}$ transversely since $\Sigma$ refines the Gr\"{o}bner complex, see \cite[Theorem 12.3 and Proposition 12.6]{Gubler2013}. This implies that $\mathcal{X}_{s}+\mathcal{D}$ is a strictly normal crossings divisor.    
\end{proof}

\begin{remark}
    Our assumption that the polyhedral complex $\Sigma$ be unimodular ensures that $\mathcal{X}\to U$ is regular. Without the unimodularity, the model in general will not be regular, and thus also not strictly semistable. Note that we can always find a unimodular triangulation after perhaps extending the base using $\CC[t]\subset \CC[t^{1/M}]$ for some $M$. 
\end{remark}

\begin{remark}
    Note that $\mathcal{X}$ by itself is strictly semistable as well. We view the strictly semistable model $(\mathcal{X} \to U, \varnothing)$ as corresponding to the closed case in Definition~\ref{def:OpenVsClosed}, and $(\mathcal{X} \to U, \mathcal{D}^{\rm full})$ as corresponding to the open case. 
\end{remark}

\begin{remark}\label{rem:RemarkOverIntegers}
Suppose that the tropicalization is smooth. Then most of the constructions given above in fact work over $\Spec(\ZZ[\![t]\!])$. This is reminiscent of Raynaud's construction of the Tate curve, albeit in an embedded algebraic form. For instance, if we take a tropically smooth elliptic curve, then this gives rise to a proper semistable model $\mathcal{X}\to \Spec(\ZZ[\![t]\!])$. The reason for this is that the polyhedral constructions surrounding tropical smoothness do not require any primes to be inverted. Moreover, we can check the semistability prime by prime, where we can use the classical result. Note that the induced map to the moduli stack of elliptic curves does not give a local chart, since the orders of the component groups of the special fibers of the N\'{e}ron are larger than $1$. Instead, tropically smooth curves can be seen as giving local charts for certain modular curves at cusps. We can similarly view tropically smooth varieties as giving charts for more complicated stacks at the cusps.   
\end{remark}
\begin{remark}
In this paper, our strictly semistable models are derived from schemes ${X}/\CC(t)$ embedded into some toric variety. In principle, there is no obstruction to defining an abstract variant where we locally glue together Kummer coverings of very affine linear spaces, or snc-compactifications thereof. See Example~\ref{exa:PairOfPantsg2} for an example of this.

Finally, we note that one can define a notion of quasi-hyperbolicity using tropical techniques. That is, if we can describe the topology of a complex variety $X$ using an angle functor $\mathcal{F}_{\Theta}$ whose maps $\mathcal{F}_{\Theta}(z_{1})\to \mathcal{F}_{\Theta}(z_{0})$ for $z_{1}>z_{0}$ are injections, then $X$ can be considered quasi-hyperbolic. For instance, for compact Riemann surfaces $\overline{X}$ of genus $g$ with $k$ marked points $D\subset \overline{X}$, this is equivalent to $\chi(X\backslash D)\geq{0}$. From this point of view, all rational tropical manifolds $(X,\Sigma,D_{\text{full}})$ lead to quasi-hyperbolic varieties.
\end{remark}

\section{Tropical diagrams of topological spaces}\label{sec:LogGeometry}\label{sec:TropicalDiagrams}

In this section, we enhance the rational tropical manifolds from the previous section with a set of topological spaces. 
The idea will be to use the associated strictly semistable model to construct a Kato-Nakayama space. 
Near the toric boundary, we will use the natural quotient map on the angles. 
We will define a cofibrant replacement of the corresponding diagram of Kato-Nakayama spaces to show that the homotopy colimit of this diagram recovers the homotopy type of a nearby smooth fiber.

We start by setting our notation for diagrams of spaces and homotopy colimits in Section~\ref{sec:DiagramsSpaces}. 
We then review Kato-Nakayama spaces in Section~\ref{sec:LogSchemes}, and we show how these spaces can be interpreted in terms of initial degenerations for rational tropical manifolds in Section~\ref{sec:InitialDegenerationsKNSpaces}. This allows us to define the corresponding Kato-Nakayama functor in Section~\ref{sec:TropicalDiagramsViaKN}.
We show in Theorem~\ref{thm:KNFunctorHomotopyColimit3} how the homotopy colimit of this diagram recovers the topology of a nearby smooth fiber.

\subsection{Diagrams of topological spaces}\label{sec:DiagramsSpaces}

We review the concept of a diagram of topological spaces. These diagrams give rise to a glued space through the categorical construction of a homotopy colimit. We recall this construction together with some basic results and examples, including the case where the indexing diagram is a poset. 
We refer the reader to \cite{WZZ99}, \cite{Kozlov2008} and \cite{Segal68} for more background.

\begin{definition}
Let $\mathcal{I}$ be a category.
A \emph{diagram of topological spaces} over $\mathcal{I}$ is a covariant functor $\mathcal{F}:\mathcal{I}\to \mathrm{Top}$. We will also refer to this as a \emph{diagram of spaces}. For an object $\sigma \in\mathrm{Ob}(\mathcal{I})$, we write $\mathcal{F}(\sigma)$ for the corresponding topological space, and for a morphism $f:\sigma \to \tau$, we write $\mathcal{F}(f):\mathcal{F}(\sigma)\to \mathcal{F}(\tau)$ for the corresponding morphism of topological~spaces.
\end{definition}

\begin{definition}(\v{C}ech diagram)\label{def:GoodCovering}
Let $\mathcal{U}=\{U_{i}\}$ be a collection of subsets of a topological space $X$, where $i$ ranges over an indexing set $I$. Note that we do not require the $U_{i}$ to be open. We endow each $U_{i}$ with the induced subspace topology. We say that $\mathcal{U}$ forms a covering of $X$ if $X = \bigcup_{i \in I} U_i$. We say that $\mathcal{U}$ is locally finite if for every $x\in{X}$, there are finitely many $i \in I$ such that $x\in U_{i}$. Suppose that $\mathcal{U}=\{U_{i}\}$ is a covering of a topological space $X$, and let $\mathcal{I}$ be the category whose objects are subsets $S$ of $I$ such that $\bigcap_{i\in{S}}U_{i}$ is non-empty, and whose morphisms are inclusions. The \emph{\v{C}ech functor} is the following diagram of topological spaces 
    \begin{equation*}
        \mathcal{F}_{\check{C}}:\mathcal{I}\to \mathrm{Top}, \quad \quad S \mapsto \ \bigcap_{i\in{S}}U_{i}.
    \end{equation*}
\end{definition}

\begin{definition}(Poset diagrams)\label{def:PosetDiagram}
Let $\Pscr$ be a poset. We also write $\Pscr$ for the natural category associated to $\Pscr$, which is defined as follows. Its objects are the elements of $\Pscr$. For every $\sigma\geq{\tau}$, we have exactly one morphism $\sigma\to \tau$. A poset diagram is a covariant functor 
    \begin{equation*}
        \Pscr\to {\rm Top}.
    \end{equation*}
\end{definition}
\begin{remark}
    Note that a \v{C}ech diagram is in fact an example of a poset diagram. Indeed, we can define a poset consisting of all $S\subset I$ such that $\bigcap_{i \in S} U_{i}$ is non-empty. 
    Due to the way we assigned a category to a poset, we then have to take the following partial order: if $S_{1}\subset S_{2}$, then $S_{2}\geq S_{1}$. The two notions of diagrams then coincide.  
\end{remark}

\begin{remark}
    Our posets will mostly come from polyhedral complexes $\Sigma$. That is, we endow $\Sigma$ with the following poset structure: 
    for $\sigma,\tau\in \Sigma$, $\sigma\geq{\tau}$ if $\sigma \supset \tau$. This is the reason why we use polyhedral notation for elements of posets throughout this paper. See Definition~\ref{def:TropicalPoset} for a more general poset structure arising from tropical geometry.  
\end{remark}

\begin{definition}[Order complex]
Let $\Pscr$ be a finite poset. The order complex $\Delta(\Pscr)$ of $\Pscr$ is the simplicial complex whose $k$-dimensional cells are given by chains $z=(z_{0}<...<z_{k})$ of length $k$. We denote the set of chains of length $k$ by $\Delta(\Pscr)_{k}$. The vertices of a simplex $(z_{0}<...<z_{k})$ are the chains of length zero given by $z_{i}$ for  $0\leq{i}\leq{k}$. 
\end{definition}
We now discuss the notion of a homotopy colimit, as in \cite[Section 15.2]{Kozlov2008} and \cite{WZZ99}. We will state the main results regarding the homotopy colimit for general diagrams, but only give the definition when the indexing category is a poset.

We first associate an explicit simplex in $\RR^{k+1}$ to each chain $z=(z_{0}<...<z_{k})\in \Delta(\Pscr)$. Let $[k]=\{0,...,k\}$. Then ${\rm Hom}_{{\rm Sets}}([k],\RR_{\geq{0}})=\RR^{k+1}_{\geq{0}}$. We will write elements of ${\rm Hom}_{{\rm Sets}}([k],\RR_{\geq{0}})$ as $t=(t_{0},...,t_{k})$.   
Let 
\begin{equation*}
    \Delta(z) = \left\{\phi\in {\rm Hom}_{{\rm Sets}}([k],\RR_{\geq{0}}):\sum_{i=0}^{k}\phi(i)=1\right\}.
\end{equation*}
In coordinates, we have $\Delta(z)= \{(t_{0},...,t_{k})\in\RR^{k+1}_{\geq{0}}:t_{0}+...+t_{k}=1\}$. In particular, for $k=0$ we find that $\Delta(z)$ is the single point $\{1\}\subset \RR$. We identify subchains of $z=(z_{0}<...<z_{k})$ with non-empty subsets of $[k]$. We write these as $z_{S}\in \Delta(\Pscr)$ for $S\subset [k]$.

Let $\mathcal{F}:\Pscr\to {\rm Top}$ be a functor. We consider the space 
\[
B=\bigsqcup_{z=(z_{0}<...<z_{k})\in \Delta(\Pscr)}\Delta(z)\times \mathcal{F}(z_{k}).
\]
We denote an element in $B$ by $w=[z \ ; \ \phi \ ; \ x]$, where 
\[
    z=(z_{0}<...<z_{k})\in \Delta(\Pscr), \quad \phi\in \Delta(z)\subset {\rm Hom}_{{\rm Sets}}([k],\RR_{\geq{0}}) \quad  \text{and} \quad x\in \mathcal{F}(z_{k}).
\] 
We now introduce an equivalence relation on $B$.

\begin{definition}[Equivalence relation for homotopy colimits]
Let $\Pscr$ be a finite poset and let 
    \[
        \mathcal{F}:\Pscr\to \text{Top}
    \] 
be a poset diagram, as in Definition~\ref{def:PosetDiagram}. 
Let $w=[z \ ; \ \phi \ ; \ x]\in B$, where $z=(z_{0}<...<z_{k})$, $\phi\in \Delta(z)\subset {\rm Hom}_{{\rm Sets}}([k],\RR_{\geq{0}})$ and $x\in \mathcal{F}(z_{k})$. Consider a non-empty subset $S\subset [k]$ containing $\{i\in[k]:\phi(i)\neq{0}\}$. Write $s$ for the maximal element in $S$. Then $w\sim w_{S}=[z_{S} \ ; \ \phi|_{S} \ ; \ x_{S}]$, where $z_{S}$ is the chain obtained from $z$ by only keeping $z_{i}$ with $i\in S$, $\phi|_{S}$ is the restriction of $\phi$ to $S$, and $x_{S}$ is the image of $x$ under $\mathcal{F}(z_{k})\to \mathcal{F}(z_{s})$. We write $R\subset B\times B$ for the equivalence relation this generates. For $w_{1}$ and $w_{2}$, we again write $w_{1}\sim w_{2}$ if $(w_{1},w_{2})\in R$.    
\end{definition}
In other words, if one of the coordinates $t_{i}$ is zero, then we can remove it and obtain an equivalent vector with fewer entries. 
We however have to replace $x\in \mathcal{F}(z_{k})$ by its image in $\mathcal{F}(z_{s})$ under the natural map $\mathcal{F}(z_{k})\to \mathcal{F}(z_{s})$, where $s = \max(S)$. We now use this to give an easy description of the equivalence classes in $B$. This will allow us to work with concrete representatives in the homotopy colimit, which eventually will lead to explicit finite CW structures in Section~\ref{sec:CWStructure}.

\begin{lemma}\label{lem:RepresentativesHomotopyColimit2}
Every $y\in B$ is equivalent to exactly one 
$w=[z \ ; \ \phi \ ; \ x]\in B$ with $z=(z_{0}<...<z_{k})$, 
$\phi(i)>0$ for all $0\leq i \leq k$ and $x\in \mathcal{F}(z_{k})$.
\end{lemma}
\begin{proof}
    This is immediate from the definition of the equivalence relation. 
\end{proof}

\begin{definition}[Homotopy colimit]\label{def:HomotopyColimit}
Let $\Pscr$ be a finite poset and let $\mathcal{F}:\Pscr\to \text{Top}$ be a diagram of spaces. Let $B$ be as above. The homotopy colimit of $\mathcal{F}$ is the space ${\rm hocol}(\mathcal{F})=B/\sim$. We endow ${\rm hocol}(\mathcal{F})$ with the quotient topology.
\end{definition}
\begin{remark}
    In terms of homotopy theory, it is better to speak of \emph{a} homotopy colimit of the diagram, as there are many ways to construct it. The final space is however independent of these choices, in the sense that they are all weakly homotopy equivalent. See \cite{Dugger25,Kozlov2008}.  
\end{remark}

If the diagram arises from a suitable covering $\mathcal{U}=\{U_{i}\}$ of a space $X$, then the projection map from the homotopy colimit to the colimit is a homotopy equivalence.

\begin{theorem}[Projection lemma]\label{thm:ProjectionLemma} 
Let $\{U_{i}\}$ be a \v{C}ech covering of a paracompact topological space $X$ with diagram of spaces $\mathcal{F}_{\check{\mathcal{C}}}$ as in Definition~\ref{def:GoodCovering}. Let 
\begin{equation*}
    \mathrm{hocol}(\mathcal{F}_{\check{\mathcal{C}}})\to \mathrm{colimit}(\mathcal{F}_{\check{\mathcal{C}}})
\end{equation*}
be the induced projection map. 
We then have the following:
\begin{enumerate}
    \item If the $U_{i}$ are open and the covering is locally finite, then the projection map is a homotopy equivalence. Moreover, the colimit is $X$.
    \item If the $U_{i}$ are closed and the covering is finite, then the projection map is a homotopy equivalence and the colimit is $X$.
\end{enumerate}

\end{theorem}
\begin{proof} 
   See \cite[Theorem 15.19 and Remark 15.20]{Kozlov2008} for instance. For a more general result in terms of model categories, see \cite[Proposition 14.2]{Dugger25} and \cite[Proposition 10.10]{Hirschhorn25}. For the last part, we refer the reader to \cite[Lemma 4.5]{WZZ99}.
\end{proof}

\begin{remark}
    There is an analogue for this projection lemma in the context of poset diagrams. 
    Consider a poset diagram $\mathcal{F}:\Pscr\to {\rm Top}$.
    We can then consider the latching map
    \[
    L_{\sigma
    }(\mathcal{F})={\rm colimit}_{\sigma<\tau}(\mathcal{F}(\tau))\to \mathcal{F}(\sigma).
    \] 
    If these latching maps are closed cofibrations\footnote{Note that for a CW complex $X$, any closed subspace $Z$ that is also a CW complex automatically defines a closed cofibration $Z\to X$. We will not need the general definition of a cofibration, as all of our examples are of this form. }, then the map from the homotopy colimit to the colimit is again a homotopy equivalence, see \cite[Proposition 3.1]{WZZ99}.

    In general, it is not enough for the maps $\mathcal{F}(\tau)\to \mathcal{F}(\sigma)$ to be cofibrations. For an example that is pertinent to our polar decompositions, see Example~\ref{exa:PairOfPantsg2}. 
\end{remark}

\begin{proposition}[Theorem 15.12 in \cite{Kozlov2008}]\label{prop:HomotopyEquivalentFunctors}
    Let $\mathcal{F}_{1}\to \mathcal{F}_{2}$ be a natural transformation of diagrams of topological spaces such that $\mathcal{F}_{1}(\sigma)\to \mathcal{F}_{2}(\sigma)$ is a homotopy equivalence for each object $\sigma$. Then the induced map
    \begin{equation*}
        \mathrm{hocol}(\mathcal{F}_{1})\to \mathrm{hocol}(\mathcal{F}_{2})
    \end{equation*}
    is a homotopy equivalence. 
\end{proposition}

\begin{remark}\label{rem:SpectralSequences}
   Suppose that we have a diagram of spaces $\mathcal{F}: \Pscr \to {\rm Top}$ such that the projection map $\mathrm{hocol}(\mathcal{F})\to \mathrm{colimit}(\mathcal{F})=X$ is a homotopy equivalence. One can then express the cohomology groups of any sheaf on $X$ in terms of cohomology groups on the $\mathcal{F}(\sigma)$ using a spectral sequence, see \cite[Proposition 15.9]{Dugger25}. More precisely, suppose that we have a diagram of spaces $\mathcal{F}:\Pscr\to {\rm Top}$ such that the projection map $\mathrm{hocol}(\mathcal{F})\to \mathrm{colimit}(\mathcal{F})=X$ is a homotopy equivalence. 
   For any $n\geq{0}$, we have a functor 
    \begin{equation*}
       \mathcal{E}^{n}: Y\mapsto H^{n}(Y,\ZZ)
    \end{equation*}
    from the opposite of the category of topological spaces to abelian groups.  
    If we have a diagram of spaces $\mathcal{F}:\Pscr\to {\rm Top}$, then we can compose $\mathcal{F}$ with this cohomology theory $\mathcal{E}^{n}$ to obtain a diagram of abelian groups $\mathcal{E}^{n}(\mathcal{F}):\Pscr^{\rm opp}\to {\rm Ab}$. This sends
    \begin{equation*}
        \sigma \mapsto H^{n}(\mathcal{F}(\sigma),\ZZ).
    \end{equation*}
Consider the abelian category ${\rm Ab}^{\Pscr^{\rm opp}}$ of diagrams $\Pscr^{\rm opp}\to {\rm Ab}$. This has enough injectives, so that we can define the right derived functors of any left exact functor. We will take the limit functor 
\[
    {\rm Ab}^{\Pscr^{\rm opp}}\to {\rm Ab}
\]
sending $\mathcal{G}$ to $\mathrm{lim}(\mathcal{G})$.
This is left exact, and we write $H^{q}(\Pscr^{\rm opp},\mathcal{G})$ for the corresponding cohomology groups. We can now precisely state how homotopy colimits allow us to glue together local cohomology groups. 

\begin{proposition}\label{prop:SpectralSequence}
    Let $\mathcal{F}:\Pscr\to {\rm Top}$ be a diagram of spaces such that the projection map 
    \begin{equation*}
        {\rm hocol}(\mathcal{F})\to {\rm colimit}(\mathcal{F})=X
    \end{equation*}
    is a homotopy equivalence. 
    Let $\mathcal{E}^{n}:Y\mapsto H^{n}(Y,\ZZ)$ be the $n$-th cohomology functor with $\ZZ$-coefficients. There is then a spectral sequence 
    \begin{equation*}
        E^{p,q}_{2}=H^{p}(\Pscr^{\rm opp},\mathcal{E}^{q}(\mathcal{F}))\Rightarrow H^{p+q}(X,\ZZ). 
    \end{equation*}
\end{proposition}
 \begin{proof}
     See \cite[Proposition 15.9]{Dugger25}. 
 \end{proof}   
    We view this as a generalization of the \v{C}ech-to-derived spectral sequence. In the case of a smooth tropicalization, these local cohomology groups are given by Orlik-Solomon algebras. 
   \end{remark}

\subsection{Log schemes and Kato-Nakayama spaces}\label{sec:LogSchemes}

Throughout this paper, we will use various concepts from logarithmic geometry to pass from singular to smooth fibers. For more background, we refer the reader to \cite{OgusBook18,Kato89,Usui01}. We start by defining log schemes and their Betti realizations or \emph{Kato-Nakayama spaces}. 
These are also referred to as \emph{real oriented blow-ups} in the literature.

\begin{definition}
    Let $X$ be a scheme over $\CC$. A pre-log structure on $X$ is a sheaf of monoids $\mathcal{M}$ on $X$, together with a morphism of sheaves $\rho:\mathcal{M}\to \mathcal{O}_{X}$. Here $\mathcal{O}_{X}$ is viewed as a sheaf of monoids with respect to multiplication, so that $1\in\mathcal{O}_{X}(U)$ is the unit. We say that $\rho$ is a log structure if the induced map $\rho^{-1}(\mathcal{O}^{\times}_{X})\to \mathcal{O}^{\times}_{X}$ is a bijection. Here $\mathcal{O}^{\times
    }_{X}(U)=(\mathcal{O}_{X}(U))^{\times}$ is the sheaf of invertible elements on $X$. A log scheme $(X,\mathcal{M})$ consists of a scheme $X$ over $\CC$ together with a log structure $\mathcal{M}$.       
\end{definition}

\begin{definition}
    Let $(X,\mathcal{M})$ be a log-scheme over $\CC$.
    The Kato-Nakayama space $X^{\log}$ associated to $(X,\mathcal{M})$ is the set of all pairs $(x,h)$ of points  $x \in X(\mathbb{C})$ and group homomorphisms $h \colon \mathcal{M}^{\rm grp}_{X,x} \to S^1$
    satisfying
    \[
         h(f) = f(x)/ |f(x)| \quad \text{whenever } f \in \mathcal{O}_{X,x}^\times\subset \mathcal{M}^{\rm grp}_{X,x}.
    \]
   We write 
   \begin{equation*}
       \pi:X^{\log}\to X(\mathbb{C})
   \end{equation*}
   for the forgetful map sending $(x,h)\mapsto x$. If $Z\subset X$ is a subset, we write $Z^{\log}$ for the inverse image of $Z$ under the forgetful map. 
   We endow $X(\mathbb{C})$ with the natural complex topology. The topology on $X^{\log}$ is the weakest one such that $\pi:X^{\log}\to X(\CC)$ is continuous, and for any open $U\subset X$ and any $m\in \mathcal{M}_{X}(U)$ the angle map 
   \begin{equation*}
       \mathrm{ang}(m):U^{\log}\to S^{1} 
   \end{equation*}
   sending $(x,h)$ to $h(m_{x})$ 
   is continuous. Here $m_{x}$ is the image of $m$ in $\mathcal{M}^{{\rm grp}}_{X,x}$. For a morphism of log-schemes $f: (X,\mathcal{M})\to (Y,\mathcal{N})$, we obtain a continuous map $f^{\log}: X^{\log}\to Y^{\log}$ as follows. On the level of points, this sends $x\mapsto f(x)$. For the group homomorphism, we compose $\mathcal{M}^{\rm grp}_{X,x}\to S^{1}$ with $\mathcal{N}^{\rm grp}_{Y,f(x)}\to \mathcal{M}^{\rm grp}_{X,x}$.  
\end{definition}

We give several examples of log schemes with their Betti realizations here. We note that it suffices to work with a \emph{chart} for the log structure to determine the points in the Betti realization, see \cite[Section 2.2]{OgusBook18}. This will be used often in the upcoming examples.

\begin{notation}
    When $P$ is a monoid and $X$ is a variety over $\CC$, we will also write $P$ for the constant sheaf on $X$ associated to $P$.
\end{notation}

\begin{example} \label{exa:DivisorialLogStructure}
    Consider a complex variety $X$ with a closed subset $D\subset X$ and corresponding open $V=X\backslash{D}$.
    We endow $X$ with the log structure associated to $D$. The corresponding sheaf of monoids on $X$ is defined by 
    \[
        \mathcal{M}(U) := \{ f \in \mathcal{O}_X(U) \colon f_{| U \cap V} \in \mathcal{O}_{X}(U\cap V)^\times \}.
    \] 
    The morphism $\mathcal{M}(U)\to \mathcal{O}_{X}(U)$ is the inclusion map. 
    We denote the associated log scheme by $(X,\mathcal{M})$. 
    We will only consider this log structure for $D$ an effective Cartier divisor on $X$. In this case, we call this the divisorial log structure associated to $D$. 
    
    Suppose for instance that $X=\Spec(\CC[x])$, $D=V(x)$. We then have the pre-log structure 
    \begin{equation*}
        {\mathbb{N}}(U)\to \mathcal{O}_{X}(U)
    \end{equation*}
    sending $1$ to $x$ for every open $U\subset X$. The associated log structure is $\mathcal{M}$. In other words, this pre-log structure defines a chart for $(X,\mathcal{M})$.
    
    We now determine the Betti realization. We have $X(\CC)=\CC$. For $0\in X(\CC)$, we again obtain all group homomorphisms $\ZZ\to S^{1}$, so we find a circle. For every other point in $X(\CC)$  however, there is only one option.  
    We thus see that $X^{\log}=S^{1}\times \mathbb{R}_{\geq{0}}$. The forgetful map $X^{\log}\to \mathbb{C}$ can be seen as a real oriented blow-up. In a sense, we have extended the usual polar decomposition to $x=0$ by adding a circle, see Figure~\ref{fig:C_KN_space}. 

    \begin{figure}[ht]
    \centering
    \scalebox{1}{
     \begin{tikzpicture}       
       \filldraw[fill=gray!80, draw=black,opacity=0.5] (0,0) -- (3,-1) .. controls (3.5,-0.5) and (3.5,0.5) .. (3,1) -- cycle;
       \draw[black,dashed] (1.5,-0.5) .. controls (1.8,-0.3) and (1.8,0.3) .. (1.5,0.5) ;
       
       \filldraw[fill=blue!20, draw=black, opacity=0.7] (0,0) -- (3,-1) .. controls (2.7,-0.5) and (2.7,0.5) .. (3,1) -- cycle;
       \draw[black] (1.5,-0.5) .. controls (1.2,-0.3) and (1.2,0.3) .. (1.5,0.5) ;
       \draw [fill=black] (0,0) circle (1.5pt);
       
       \draw[->] (0,-1.5) -- (3.5,-1.5) node[right] [xscale = 0.8, yscale = 0.8] {$r$};     
       \node at (1.5, 1) [xscale = 0.8, yscale = 0.8] {$X$};

       \fill[fill=gray!80, opacity=0.5] (5,-0.4) .. controls (5.2,-0.1) and (5.2,0.1) .. (5,0.4) -- (8,1) .. controls (8.5,0.5) and (8.5,-0.5) .. (8,-1) -- cycle  ;
       \draw[dashed, line width = 1.5] (5,-0.4) .. controls (5.2,-0.1) and (5.2,0.1) .. (5,0.4);
       \draw[black,dashed] (6.5,-0.7) .. controls (6.8,-0.3) and (6.8,0.3) .. (6.5,0.7) ;
       \draw (8,1) .. controls (8.5,0.5) and (8.5,-0.5) .. (8,-1);

       \filldraw[fill=blue!20, draw=black, opacity=0.7] (5,-0.4) .. controls (4.8,-0.1) and (4.8,0.1) .. (5,0.4) -- (8,1) .. controls (7.5,0.5) and (7.5,-0.5) .. (8,-1) -- cycle  ;
       \draw[draw=black, line width = 1.5] (5,-0.4) .. controls (4.8,-0.1) and (4.8,0.1) .. (5,0.4);

       \draw[black] (6.5,-0.7) .. controls (6.2,-0.3) and (6.2,0.3) .. (6.5,0.7) ;
       \draw[->] (5,-1.5) -- (8.5,-1.5) node[right] [xscale = 0.8, yscale = 0.8] {$r$};
       \node at (6.5, 1) [xscale = 0.8, yscale = 0.8] {$X^{\log}$};
    \end{tikzpicture}
    }
    \caption{
    The complex plane $X = \Spec(\CC[x])$ is presented on the left in polar coordinates. Here, the point $0$ is the divisor.
    On the right, we depict the corresponding Kato-Nakayama space $X^{\log}$. Note that the point $0$ in $X$ on the left is replaced by an entire circle $S^1$ in $X^{\log}$ on the right. }
    \label{fig:C_KN_space}
  \end{figure}
    
    For a more complicated example, suppose $X = \PP^1$ and $D = (0) + (1) + (\infty)$. The associated Kato-Nakayama space $X^{\log}$ is obtained from $\mathbb{P}^{1}(\mathbb{C})$ by replacing the three points in $D$ with three circles. This is  topologically a sphere with three holes, i.e. a \emph{pair-of-pants}. We stress that the boundaries of these holes are also included in $X^{\log}$. In particular, $X^{\log}$ is \emph{compact}.
\end{example}

Recall that, by Proposition~\ref{pro:StrictlySemistableModels}, a rational tropical manifold $(X, \Sigma, D)$ comes with a strictly semistable model $(\mathcal{X}\to U,\mathcal{D})$.
For the remainder of the paper, we will assume that the strictly semistable models arising from Proposition~\ref{pro:StrictlySemistableModels} are defined over a small enough Zariski open neighborhood $U\subset \Spec(\CC[t])$ containing $t=0$. In other words, we have a proper flat morphism of finite presentation $\mathcal{X}\to U$, where $\mathcal{X}$ is integral and regular. These arise naturally from $X$ defined over the rational function field $\CC(t)$, which is equipped with the $t$-adic valuation.

\begin{definition}\label{def:KNStructureRTM}
Let $(\mathcal{X} \to U,\mathcal{D})$ be a strictly semistable model over $U\subset \Spec(\CC[t])$ associated to a rational tropical manifold $(X, \Sigma, D)$ by Proposition~\ref{pro:StrictlySemistableModels}. We write $\mathcal{X}^{\log}_{\mathcal{D}}$ for the Kato-Nakayama space associated to the strictly normal crossings divisor $\mathcal{X}_{s} + \mathcal{D}$. We endow $U$ with the log structure coming from $t = 0$. 

Note that the fiber of $U^{\log}\to U$ over $0$ can be identified with $S^{1}$. For any such $\zeta\in S^{1}$, we write $\mathcal{X}^{\log,\zeta}_{\mathcal{D}}$ for the fiber under the natural map $\mathcal{X}^{\log}_{\mathcal{D}}\to U^{\log}$.
This is the Kato-Nakayama fiber associated to the rational tropical manifold.
\end{definition}

We now recall a result by Usui, which says that $\mathcal{X}^{\log,\zeta}_{\mathcal{D}}$ captures the topology of nearby smooth fibers. This will be the main tool that allows us to connect objects in tropical and logarithmic geometry to objects in algebraic geometry.  

\begin{theorem}\label{thm:UsuisTheorem2}
    Let $(\mathcal{X}\to U,\mathcal{D})$ be a strictly semistable model over an open $U\subset \Spec(\CC[t])$ that contains $0$. We endow $\mathcal{X}$ and $U$ with the log structures induced by the divisors $\mathcal{X}_{s}+\mathcal{D} $ and $\{0\}$ respectively. Then $\mathcal{X}_{\mathcal{D}}^{\log}\to U^{\log}$ is topologically trivial over $U^{\log}$. In particular, $\mathcal{X}^{\log,\zeta}_{\mathcal{D} ,s}$ is homeomorphic to $\mathcal{X}_{t}\backslash\mathcal{D}_{t}$ for general $t\in U(\CC)$.
\end{theorem}
\begin{proof}
    The first incarnation of this theorem can be found in 
    \cite[Theorem 5.4]{Usui01} (although one takes a slightly different log structure there, but it yields the same result). This result was later generalized by Ogus and Nakayama in \cite[Proposition 5.1]{NO10}. We refer the reader to \cite[Theorem 6.2.3]{BK24} for the current version.
\end{proof}

\subsection{Initial degenerations and Kato-Nakayama spaces}\label{sec:InitialDegenerationsKNSpaces}

In this section we study the initial degenerations associated to rational tropical manifolds, and we interpret these in terms of Kato-Nakayama spaces.

\begin{definition}\label{def:StrataAtInfinity}
    Let $\sigma \in \Delta$ be a cone in the recession fan of $\Sigma$ and $V(\sigma)$ the corresponding toric variety $V(\sigma) := \overline{\Ocal}_\sigma \subset Y(\Delta)$. 
    We denote by $Z_\sigma \subset \mathcal{X}_s$ the closed set obtained by intersecting $\mathcal{X}_s$ with the closure of $\overline{X} \cap V(\sigma)$ in $\mathcal{X}$. Equivalently we have
    \[
        Z_{\sigma} := \bigcup_{\substack{\tau \in \Sigma \\  \sigma \ \subset \ \rec(\tau)}} Z_{\tau}.
    \]
\end{definition}

\begin{definition}[Logarithmic initial degenerations]
\label{def:LogarithmicInitialDegenerations}
    Consider a rational tropical manifold $(X,\Sigma, D)$ with a cell $\tau \in \Sigma$ corresponding to a closed subspace $Z_{\tau}\subset \mathcal{X}_{s}$ by the torus-orbit correspondence, see \cite[Section 7.9]{Gubler2013}. We write $Z^{\circ}_{\tau}\subset Z_{\tau}$ for the corresponding open subspace that lies in the torus $\Spec(\CC[\tau^{\perp}\cap M])$.
    
    We write $Z^{\log}_{\tau}$ and $Z^{\log,\circ}_{\tau}$ for the inverse images of $Z_{\tau}$ and $Z^{\circ}_{\tau}$ respectively under the forgetful map
    \begin{equation*}
        \mathcal{X}^{\log,\zeta}_{\mathcal{D}, s}\to \mathcal{X}_{s}.
    \end{equation*}
    We call $Z^{\log}_{\tau}$ the \emph{logarithmic initial degeneration} of $Z_{\tau}$, and $Z^{\log,\circ}_{\tau}$ the \emph{open logarithmic initial degeneration} associated to $\tau$. 

    Similarly, suppose we have a pair $(\tau,\sigma)\in \Pscr(\Sigma_{D})$ with inclusion $Z_{\sigma}\supset Z_{\tau}$ of closed subsets, as in Definition~\ref{def:StrataAtInfinity}. Note that $Z_{\sigma}$ similarly has a log structure coming from its toric boundary. We write $Z^{\sigma,\log}_{\sigma}$ for the corresponding Kato-Nakayama space. We have a map 
    \begin{equation*}
        Z^{\sigma,\log}_{\sigma}\to Z_{\sigma},
    \end{equation*}
    and for any $\tau$ as above, we write $Z^{\sigma,\log}_{\tau}$ for the inverse image of $Z_{\tau}$ under this map. We write $Z^{\sigma,\log,\circ}_{\tau}$ for the points that map to $Z^{\circ}_{\tau}$. We call $Z^{\sigma,\log}_{\tau}$ the \emph{relative logarithmic initial degeneration} of $Z_{\tau}$ with respect to $\sigma$, and $Z^{\sigma,\log,\circ}_{\tau}$ the \emph{relative open logarithmic initial degeneration} of $Z_{\tau}$ with respect to $\sigma$.
\end{definition}

To make these spaces more concrete, we will now explain how they can be interpreted in terms of initial degenerations. 
Let $X=V(I)\subset \Spec(K[M])$. For any cone $\sigma \subset N_{\RR}$, we denote the lattice orthogonal to $\sigma$
\[
    M(\sigma) := \sigma^{\perp}\cap M.
\]
Its dual lattice is 
\[
    N / N_{\sigma} \simeq {\rm Hom}(M(\sigma), \ZZ).
\]

Recall that the torus  $\mathbb{T}_{\tau}=\Spec(\CC[M/M(\tau)])$ acts on $V(\mathrm{in}_{w}(I))$ for $w\in \mathrm{relint}(\tau)$. Similarly, for any $\sigma\subset \tau$, 
consider the exact sequence
\begin{equation*}
    0 \to M(\sigma)/M(\tau)\to M/M(\tau)\to M/M(\sigma)\to 0.
\end{equation*}
Then by writing
\[
    \mathbb{T}_{\sigma} := \Spec(\CC[M/M(\sigma)]),\quad 
    \mathbb{T}_{\tau}   := \Spec(\CC[M/M(\tau)]),\quad \text{and} \quad \mathbb{T}_{\sigma\subset \tau} := \Spec(\CC[M(\sigma)/M(\tau)]),
\]
we obtain an exact sequence of tori:
\begin{equation*}
    0 \to \mathbb{T}_{\sigma}\to \mathbb{T}_{\tau}\to \mathbb{T}_{\sigma\subset \tau}\to 0.
\end{equation*}
The tori $\mathbb{T}_{\sigma}$ and $\mathbb{T}_{\tau}$ then similarly act on $\mathrm{in}_{w}(I)$ for $w\in {\rm relint}(\tau)$. 
We write $\tau_\sigma$ for the cone in $N/N_{\sigma}$ defined by
\[
    \tau_\sigma := (\tau + N_{\sigma}) / N_{\sigma}.
\]
Its dual cone in $M(\sigma)$ is
\begin{equation*}
    \tau^{\vee}_{\sigma} := \{m\in M(\sigma):\langle m,n\rangle\geq{0}\,\,\forall n \in \tau_{\sigma} \},
\end{equation*}
so that 
\begin{equation*}
    M(\tau) \subset \tau^{\vee}_{\sigma} \subset M(\sigma).
\end{equation*}
We write
\[
    \overline{\tau_{\sigma}^{\vee}} := \tau_{\sigma}^{\vee} / M(\tau)_{\RR} \subset M(\sigma)_{\RR} / M(\tau)_{\RR}.
\] 
This has a corresponding positive real torus 
\begin{equation*}
    \mathbb{T}_{\sigma\subset \tau}(\RR_{>0}) 
    = 
    \mathrm{Hom} \Big(\overline{\tau_{\sigma}^{\vee}},\RR_{>0}\Big).
\end{equation*}
We view this as a subset of $\mathbb{T}_{\sigma\subset \tau}(\CC)$. 
If $\sigma=0$, then we simply write $\mathbb{T}_{\tau}(\RR_{>0})$. We now interpret the corresponding Kato-Nakayama spaces as a quotient under the action of this real positive torus. 
\begin{lemma}\label{lem:InitialDegenerationsVSKNSpaces}
Consider a pair $(\tau,\sigma)\in \Pscr(\Sigma_{D})$. We then have homeomorphisms
\[
    {\rm in}_{w}(X)_{\sigma}/\mathbb{T}_{\sigma\subset \tau}(\RR_{>0}) \to Z^{\sigma,\log, \circ}_{\tau}, \quad 
    {\rm in}_{w}(X)/\mathbb{T}_{\tau} \to {\rm in}_{w}(X)_{\tau}\to Z^{\circ}_{\tau}, \quad \text{and} \quad 
        {\rm in}_{w}(X)/\mathbb{T}_{\sigma} \to {\rm in}_{w}(X)_{\sigma}.
\]
    In particular, for $\sigma=0$, the first diagram becomes 
    \begin{equation*}
        {\rm in}_{w}(X)/\mathbb{T}_{\sigma\subset \tau}(\RR_{>0})=Z^{\log,\circ}_{\tau}. 
    \end{equation*}
    The quotient map ${\rm in}_{w}(X)_{\sigma}\to {\rm in}_{w}(X)_{\sigma}/\mathbb{T}_{\sigma\subset \tau}(\RR_{>0})$ is a homotopy equivalence. 
\end{lemma}
\begin{proof}
        This immediately follows from the fact that the toric boundary of $Z_{\sigma}$ forms an snc-divisor, and the fact that $\mathrm{in}_{w}(X)/\mathbb{T}_{\tau} = Z^{\circ}_{\tau}$, see Lemma~\ref{lem:InitialDegenerationsTorusQuotient}.
\end{proof}

\begin{remark}\label{rem:InitialDegenerationsDirectQuotientKNSpaces}
    We can also directly express the Kato-Nakayama space $Z^{\sigma,\log,\circ}_{\tau}$ as a quotient of $\mathrm{in}_{w}(X)$. To that end, we need to define a mixed torus $W$ using the exact sequence 
    \begin{equation*}
        0 \to \mathbb{T}_{\sigma}(\CC)\to W\to \mathbb{T}_{\sigma\subset \tau}(\mathbb{R}_{>0}) \to 0
    \end{equation*}
    arising from 
    \begin{equation*}
    0 \to \mathbb{T}_{\sigma}\to \mathbb{T}_{\tau}\to \mathbb{T}_{\sigma\subset \tau}\to 0.
\end{equation*}
    That is, we define $W=\{P\in \mathbb{T}_{\tau}(\CC):\pi(P)\in \mathbb{T}_{\sigma\subset \tau}(\RR_{>0})\}$. Note that in local coordinates, we can always find an isomorphism
    \begin{equation*}
        W\simeq \mathbb{R}^{r}_{>0}\times (\CC^{\times})^{m}.
    \end{equation*}
    We then find $\mathrm{in}_{w}(X)/W=Z^{\sigma,\log,\circ}_{\tau}$. This allows us to bypass the double-quotient construction. 

    Similarly, we can also associate groups to the morphism of Kato-Nakayama spaces
    \begin{equation*}
Z^{\log,\circ}_{\tau}\to Z^{\sigma,\log,\circ}_{\tau}.
    \end{equation*}
 Namely, the associated group is $\mathbb{T}_{\sigma}(S^{1})=\mathrm{Hom}(M/M(\sigma),S^{1})\to \mathrm{Hom}(M/M(\sigma),\CC)=\mathbb{T}_{\sigma}(\CC)$. This acts on the first space $Z^{\log,\circ}_{\tau}$, and the quotient is exactly $Z^{\sigma,\log,\circ}_{\tau}$.
\end{remark}
\begin{remark}
    We remark here that $Z^{\sigma,\log}_{\tau}$ is a manifold with boundary, and $Z^{\sigma,\log,\circ}_{\tau}$ is its interior. As such, the inclusion map is a homotopy equivalence $Z^{\sigma,\log,\circ}_{\tau}\to Z^{\sigma,\log}_{\tau}$. The proofs of these facts follow quickly from the fact that  $\mathcal{X}_{s}+\mathcal{D}$ is a Cartier divisor, so that all the strata are smooth, and they intersect transversally.    
\end{remark}

\subsubsection{Cofibrant replacements}\label{sec:CofibrantReplacements}
Let $w\in {\rm relint}(\tau)$ for some $(\tau,\sigma)\in \Pscr(\Sigma_{D})$. 
Consider the quotient map 
\begin{equation*}
    {\rm in}_{w}(X) 
    \to {\rm in}_{w}(X)_{\sigma} 
    := {\rm in}_{w}(X)/\mathbb{T}_{\sigma}.
\end{equation*}
This is not injective for $\sigma\neq (0)$, and thus not a cofibration. We now define a cofibrant replacement of this map. 
For this, we will use a toric version of the Borel construction for group actions. Consider a cone $\sigma\subset N_{\RR}$ with orthogonal lattice $M(\sigma)=\sigma^{\perp}\cap M$. We define 
\[
  \overline{\sigma^{\vee}} := (\sigma^{\vee} + M(\sigma)_{\RR}) / M(\sigma)_{\RR} \subset M_\RR / M(\sigma)_{\RR} 
  \quad \text{and} \quad 
  \overline{\sigma^{\vee}}_{\ZZ} := \overline{\sigma^{\vee}} \cap M/M(\sigma). 
\]
We write $V_{\sigma}$ for the toric variety 
$V_{\sigma} := \Spec(\CC[\overline{\sigma^{\vee}}_{\ZZ}])$.
Note that we have a dense torus inclusion 
\begin{equation*}
    \Spec(\CC[M/M(\sigma)])\to \Spec(\CC[\overline{\sigma^{\vee}}_{\ZZ}]),
\end{equation*} 
and the $\CC$-points of this toric variety can be identified with morphisms of monoids
\[
V_{\sigma}(\CC) \simeq
{\rm Hom}_{\rm Mon}
\big(\overline{\sigma^\vee}_{\ZZ},\CC\big).
\]
Here $\CC$ is thought of as a multiplicative monoid. We will also work with the space $V_{\sigma}(\mathbb{D})$, where $\mathbb{D}\subset \CC$ is the closed unit disk. This is also considered as a multiplicative monoid with unit $1$. Our main reason for using the unit disk here is that it is compact. This will allow us to use the Smale-Vietoris theorem.

\begin{lemma}
    The spaces $V_{\sigma}(\CC)$ and $V_{\sigma}(\mathbb{D})$ are contractible.
\end{lemma}
\begin{proof}
    This follows since the toric part is removed. 
    See \cite[Page 29]{Fulton93}. 
\end{proof}

We now suppose we have two cones $\sigma_{1}\subset \sigma_{2}\subset N_{\RR}$. 
The inclusion $M(\sigma_{2})\subset M(\sigma_{1})$ induces a surjective homomorphism
\[
    \pi:M/M(\sigma_2)\to M/M(\sigma_1).
\]
This restricts to a homomorphism of monoids
\[
\pi:
\overline{\sigma^\vee}_{2, \ZZ}
\to
\overline{\sigma^\vee}_{1, \ZZ}.
\]
Using the identification 
$
V_{\sigma_i}(\mathbb{D}) =
{\rm Hom}_{\rm Mon}
\big(\overline{\sigma^\vee}_{i, \ZZ},\mathbb{D}\big)$, 
we obtain a natural map
\[
V_{\sigma_1}(\mathbb{D})
\to
V_{\sigma_2}(\mathbb{D}),
\qquad
\rho\mapsto \rho\circ\pi.
\]
These maps are functorial with respect to inclusions of cones, i.e., if $\sigma_1\subset\sigma_2\subset\sigma_3$, 
then the map $V_{\sigma_1}(\mathbb{D})\to V_{\sigma_3}(\mathbb{D})$ is the composite
\[
V_{\sigma_1}(\mathbb{D})
\to
V_{\sigma_2}(\mathbb{D})
\to
V_{\sigma_3}(\mathbb{D}).
\]
\begin{example}
Let  $M=\ZZ^{2}$. Let $\sigma_{1}=\mathbb{R}_{\geq{0}}e_{1}$ and $\sigma_{2}=\mathbb{R}_{\geq{0}}e_{1}+\mathbb{R}_{\geq{0}}e_{2}$. As monoids, we have
\[
\sigma^{\perp}_{1} = \langle \pm{}e_{2}\rangle, \quad
\sigma^{\perp}_{2} = \langle 0\rangle, 
\quad \text{and} \quad
\sigma^{\vee}_{1} = \langle e_{1},\pm{}e_{2}\rangle, \quad
\sigma^{\vee}_{2} = \langle e_{1},e_{2}\rangle.
\]
The quotient $M/M(\sigma_{1})$ is generated by $\pm e_{1}$. The image of $\sigma^{\vee}_{1}$ is the same as the image of $\sigma^{\vee}_{2}$. The map 
\begin{equation*}
V_{\sigma_{1}}(\mathbb{D})=\mathbb{D}\to V_{\sigma_{2}}(\mathbb{D})=\mathbb{D}^{2}
\end{equation*} 
now sends $x\mapsto (x,1)$. 
\end{example}

\begin{definition}\label{def:CofibrantReplacements1}
Let $\sigma\subset \tau$ and $w\in {\rm relint}(\tau)$. 
We define 
\begin{equation*}
Z^{\sigma,\log,\circ}_{\tau,c}:=(Z^{\log,\circ}_{\tau}\times V_{\sigma}(\mathbb{D}))/\mathbb{T}_{\sigma}(S^{1}).
\end{equation*}

Here $\mathbb{T}_{\sigma}(S^{1})$, with $\mathbb{T}_{\sigma} :=\Spec(\CC[M/M(\sigma)])$, acts diagonally. 
We will use this space as a cofibrant replacement of the quotient map $Z^{\log,\circ}_{\tau}\to Z^{\sigma,\log,\circ}_{\tau}=Z^{\log,\circ}_{\tau}/\mathbb{T}_{\sigma}(S^{1})$. 
We write an equivalence class in this set as $[\phi,\psi]$ for 
$\phi\in Z^{\log,\circ}_{\tau}$ and $\psi\in V_{\sigma}(\mathbb{D})$. 
We have an injective map 
\[
    Z^{\log,\circ}_{\tau} \hookrightarrow Z^{\log,\circ}_{\tau,c}, 
    \quad x\mapsto [x,1_{\sigma}].
\] 
Here $1_{\sigma}$ is the point $V_{\sigma}(\mathbb{D}) = {\rm Hom}_{\rm Mon}(\overline{\sigma^\vee}_\ZZ, \mathbb{D})$ given by
\[
    1_{\sigma} \colon \overline{\sigma^\vee}_\ZZ \to \mathbb{D}, \quad m \mapsto 1.
\]
If we have an inclusion of cones $\sigma_{1}\subset \sigma_{2}\subset \tau$, the map $V_{\sigma_{1}}(\mathbb{D})\to V_{\sigma_{2}}(\mathbb{D})$ 
is equivariant with respect to the natural inclusion
\[
\mathbb{T}_{\sigma_{1}}(S^1)
\hookrightarrow
\mathbb{T}_{\sigma_{2}}(S^1).
\]
Indeed, both maps are induced by composition with the projection map
\[
\pi:M/M(\sigma_{2})\to M/M(\sigma_{1}).
\]
Thus, if $\lambda\in\mathbb{T}_{\sigma_{1}}(\mathbb{D})$ and
$\rho\in V_{\sigma_{1}}(\mathbb{D})$, then
\[
(\lambda\rho)\circ\pi
=
(\lambda\circ\pi)(\rho\circ\pi).
\]
Hence the map $V_{\sigma_{1}}(\mathbb{D})\to V_{\sigma_{2}}(\mathbb{D})$ induces a map 
\[
    Z^{\log,\circ}_{\sigma_{1},c} 
    \to 
    Z^{\log,\circ}_{\sigma_{2},c}.
\]
\end{definition}

\begin{lemma}\label{lem:RelativeInitialDegenerationsCofibrationsHE}
    Suppose that $\sigma_{1}\subset \sigma_{2}\subset \tau$. Then the maps 
  \[
  Z^{\log,\circ}_{\tau}\to Z^{\sigma_{1},\log,\circ}_{\tau,c}\to Z^{\sigma_{2},\log,\circ}_{\tau,c}
  \]
  are all closed cofibrations. The natural maps $Z^{\sigma_{i},\log,\circ}_{\tau,c}\to Z^{\sigma_{i},\log,\circ}_{\tau}$ are homotopy equivalences.  

\end{lemma}
\begin{proof}
One easily sees that all the maps are semi-algebraic and injective. Since these are closed subsets, and all sets in question are CW complexes (by semi-algebraic triangulation), we find that the maps are closed cofibrations.

We now show that the forgetful maps are homotopy equivalences. Fix an integral point $v \in {\rm relint}(\sigma)\cap N$. Since every element of $M(\sigma)$ vanishes on $\sigma$, the pairing with $v$
descends to a homomorphism
\[
\ell_v:M/M(\sigma)\to \ZZ, \quad m \mapsto \langle v, m\rangle.
\]
Its restriction to $\overline{\sigma^\vee}_{\ZZ}$ satisfies
\[
\ell_v(m)\geq 0, \quad \text{for all } m \in \overline{\sigma^\vee}_{\ZZ},
\]
with equality if and only if $m=0$. Using the identification $V_{\sigma}(\mathbb{D}) = {\rm Hom}_{\rm Mon} \bigl(\overline{\sigma^\vee}_{\ZZ},\mathbb{D}\bigr)$, we define a homotopy map
\[
    H:V_{\sigma}(\mathbb{D})\times[0,1]\to V_{\sigma}(\mathbb{D})
\]
by setting
\[
H(\psi,t)(m)
:=
t^{\ell_v(m)}\psi(m), \quad \text{for each } \psi \in V_\sigma(\mathbb{D}), t \in [0,1] \text{ and } m \in \overline{\sigma^{\vee}}_{\ZZ}.
\]
This is well defined because $\ell_v$ is additive. Indeed, we have the following:
\[
H(\psi,t)(m_1+m_2)
=
t^{\ell_v(m_1+m_2)}\psi(m_1+m_2)
=
H(\psi,t)(m_1)H(\psi,t)(m_2), \quad \text{for all } m_1, m_2 \in \overline{\sigma^\vee}_{\ZZ}.
\]
Moreover,
\[
H(\psi,1)=\psi,
\quad \text{and} \quad
H(\psi,0)=0_{\sigma},
\]
where $0_{\sigma}\in V_{\sigma}(\mathbb{D})$ is the torus-fixed point defined by
\[
0_{\sigma}(0)=1,
\qquad
0_{\sigma}(m)=0
\quad\text{for }0\neq m \in \overline{\sigma^\vee}_{\ZZ}.
\]
This homotopy is $\mathbb{T}_{\sigma}(S^{1})$-equivariant. Indeed, if
$\lambda\in\mathbb{T}_{\sigma}(S^{1})$, then
\[
H(\lambda\psi,t)(m)
=
t^{\ell_v(m)}\lambda(m)\psi(m)
=
\lambda(m)H(\psi,t)(m).
\]
It therefore induces a homotopy
\[
Z^{\sigma,\log,\circ}_{\tau,c}\times[0,1]
\to
Z^{\sigma,\log,\circ}_{\tau,c} \quad \text{and} \quad 
([\phi,\psi],t)
\longmapsto
[\phi, H(\psi,t)].
\]
At $t=0$, this has image
\[
\bigl(Z^{\log,\circ}_{\tau}\times\{0_{\sigma}\}\bigr)/
\mathbb{T}_{\sigma}(S^{1})
\simeq
Z^{\sigma,\log,\circ}_{\tau}.
\]
Hence the natural map
\[
Z^{\sigma,\log,\circ}_{\tau,c}
\to
Z^{\sigma,\log,\circ}_{\tau}
\]
is a homotopy equivalence.
\end{proof}

\begin{definition}\label{def:CofibrantReplacements2}
    The cofibrant replacement of $Z^{\sigma,\log}_{\tau}$ is 
    \begin{equation*}
        Z^{\sigma,\log}_{\tau,c} := 
        \Big( 
        Z^{\log}_{\tau}\times V_{\sigma}(\mathbb{D})
        \Big)
        /
        \mathbb{T}_{\sigma}(S^{1}). 
    \end{equation*}
\end{definition}

\begin{remark}
      As before, we have that the map $Z^{\sigma,\log}_{\tau,c}\to Z^{\sigma,\log}_{\tau}$ is a homotopy equivalence, and all the maps 
    \[
    Z^{\sigma_{1},\log}_{\tau,c}\to Z^{\sigma_{2},\log}_{\tau,c}
    \]
    are cofibrations. This will be used later on to prove that the natural diagram of spaces associated to a rational tropical manifold faithfully represents the topology of a nearby fiber, see Theorem ~\ref{thm:KNFunctorHomotopyColimit3}. 
    
\end{remark}

\begin{example}\label{exa:MainSealedExample3}
Let $X=Y=\CC^{\times}$, and let $\Sigma=\Delta$ be the fan of $\PP^1$ as a toric variety. 
The compactified tropicalization $\mathbb{T}\mathbb{P}^{1}$ as a set is then $\{-\infty\}\cup \mathbb{R}\cup \{\infty\}$. 
The poset $\Pscr(\mathbb{T}\mathbb{P}^{1})$ consists of $5$ elements and the corresponding Kato-Nakayama functor is as follows. 
On each affine chart $U_{i}=\CC$, we have $U^{\log}_{i} = \mathbb{R}_{\geq{0}}\times S^{1}$. Consider $Z^{\log}_{i}=S^{1}\subset U^{\log}_{i}=\mathbb{R}_{\geq{0}}\times S^{1}$. 
We then obtain the diagram of spaces
\begin{equation*}
    Z_{1}\leftarrow Z_{1}^{\log}\to X^{\log}\leftarrow Z_{2}^{\log}\to Z_{2}.
\end{equation*}
Here the outer two maps collapse the circles onto two points. 

We compute a natural cofibrant replacement of this diagram here. We replace $Z_{i}$ by $(Z^{\log}_{i}\times \mathbb{D})/S^{1}$. 
Note that the map $Z^{\log}_{i}\to (Z^{\log}_{i}\times \mathbb{D})/S^{1}$ is an inclusion, as it can be identified with the standard inclusion $S^{1}\to \mathbb{D}$. 
We thus obtain a cofibrant diagram  
\begin{equation*}
    \mathbb{D} \leftarrow S^{1}\to X^{\log}\rightarrow S^{1}\to \mathbb{D}.
\end{equation*}
It is easy to see that the colimit of this diagram is homotopy equivalent to $\mathbb{P}^{1}(\CC)$. In other words, the homotopy colimit of the original diagram is homotopy equivalent to the colimit.  
This is a general fact, as we will see in Theorem~\ref{thm:KNFunctorHomotopyColimit3}. To prove this, we will use the cofibrant replacements introduced in this section. 
\end{example}

\subsection{Constructing tropical diagrams using Kato-Nakayama spaces}\label{sec:TropicalDiagramsViaKN}

In this section we use the material from the previous sections to associate a diagram of spaces to a rational tropical manifold. We call this a tropical diagram of spaces. This diagram consists of a topological space for every cell $(\tau,\sigma)\in\Pscr(\Sigma_{D})$ in the rational tropical manifold $(X, \Sigma, D)$. We will show in Theorem~\ref{thm:KNFunctorHomotopyColimit3}
that the homotopy colimit of this diagram recovers the homotopy type of a nearby fiber. 

\begin{remark}
    Recall that a rational tropical manifold $(X, \Sigma, D)$ comes with a strictly semistable model $(\mathcal{X}\to U,\mathcal{D})$.  
    For the remainder of the paper, we will assume that the strictly semistable models arising from Proposition~\ref{pro:StrictlySemistableModels} are defined over a small enough Zariski open neighborhood $U\subset \Spec(\CC[t])$ containing $t=0$. In other words, we have a proper flat morphism of finite presentation $\mathcal{X}\to U$, where $\mathcal{X}$ is integral and regular. These arise naturally from $X$ defined over the rational function field $\CC(t)$, which is equipped with the $t$-adic valuation. 
\end{remark}

We now define the Kato-Nakayama functor $\mathcal{F}_{\mathrm{KN}}$ of the rational tropical manifold $(X, \Sigma, D)$. The base category of this functor is the tropical poset $\Pscr(\Sigma_{D})$ introduced in Definition~\ref{def:TropicalPoset}.

\begin{definition}[Kato-Nakayama functor]\label{def:KNFunctor}
Let $\Pscr(\Sigma_{D})$ be the poset associated to the rational tropical manifold $(X,\Sigma,D)$. Let $(\tau,\sigma)\in \Pscr(\Sigma_{D})$ and let $Z^{\sigma,\log}_{\tau}$ be the relative logarithmic initial degeneration from Definition~\ref{def:LogarithmicInitialDegenerations}. We set 
\begin{equation*}
    \mathcal{F}_{\rm KN}(\tau,\sigma)=Z^{\sigma,\log}_{\tau}.
\end{equation*}
For $\sigma_{1}\subset \sigma_{2}\subset \tau$, we have a quotient map 
\begin{equation*}
    Z^{\sigma_{1},\log}_{\tau}\to Z^{\sigma_{2},\log}_{\tau}.
\end{equation*}
For a fixed $\sigma\in \Delta$, and $\tau_{1}\supset \tau_{2}$, we obtain an inclusion map
\begin{equation*}
    Z^{\sigma,\log}_{\tau_{1}}\subset Z^{\sigma,\log}_{\tau_{2}}.
\end{equation*}
This defines the corresponding Kato-Nakayama functor $\mathcal{F}_{\rm KN}:\Pscr(\Sigma_{D})\to {\rm Top}$.
\end{definition}

\begin{example}\label{exa:P2ExampleHC}
    Consider the fan in $\RR^{2}$ associated to $\mathbb{P}^{2}$ with tropical manifold $\overline{\Sigma}=\mathbb{T}\mathbb{P}^{2}$. Note that the poset $\Pscr(\overline{\Sigma})$ consists of $19$ elements. In this case, every space is homotopic to a torus $(S^{1})^{k}/S^{1}$ for some $k$. We will describe the diagram up to homotopy.
    
    For every cell $\sigma$ in the interior, we assign $(S^{1})^{3}/S^{1}$ with the identity map everywhere. For the vertex at the boundary of a ray, we assign $(S^{1})^{2}/S^{1}$. We similarly assign  $(S^{1})^{2}/S^{1}$ to the adjacent rays.  The boundary maps from the interior to the rays are
    \begin{align*}
        [\theta_{X}:\theta_{Y}:\theta_{Z}]\mapsto [\theta_{X}:\theta_{Y}],\qquad 
        [\theta_{X}:\theta_{Y}:\theta_{Z}]\mapsto [\theta_{X}:\theta_{Z}],\qquad 
        [\theta_{X}:\theta_{Y}:\theta_{Z}]\mapsto [\theta_{Y}:\theta_{Z}].
    \end{align*}
    Note that this torus gives the homotopy type of $Z^{\sigma,\log}_{\sigma}=\mathbb{P}^{1,\log}$ which is homotopy equivalent to $\mathbb{P}^{1}\backslash\{0,\infty\}$. The angle map arises by composition from the group homomorphism $M(\sigma)\to M$, where $M$ is the character lattice of $\mathbb{P}^{2}$. Similarly, for every boundary vertex in the triangle $\mathbb{T}\mathbb{P}^{2}$, we obtain the trivial space. Note that this diagram is the same as the one in \cite[Section 5.3]{WZZ99}.
\end{example}

\begin{remark}
    Note that $Z^{\sigma,\log}_{\tau}$ is a manifold with boundary, and its interior is $Z^{\sigma,\log,\circ}_{\tau}$. This is homotopy equivalent to $\mathrm{in}_{w}(X)_{\sigma}$ for $w\in {\rm relint}(\tau)$. That is, if we are only interested in the homotopy types of the spaces $\mathcal{F}_{\rm KN}(\tau,\sigma)$, we can simply work with the corresponding initial degenerations and change the ambient torus.

    Note however that these initial degenerations do not give rise to natural horizontal comparison maps for pairs $(\tau_{1},\sigma)$ and $(\tau_{2},\sigma)$ with $\tau_{1}\subset \tau_{2}$. For Kato-Nakayama spaces, we however do obtain natural inclusions of subspaces $Z^{\sigma,\log}_{\tau_{2}}\subset Z^{\sigma,\log}_{\tau_{1}}$. This, together with Usui's Theorem~\ref{thm:UsuisTheorem2}, is our main reason for working with Kato-Nakayama spaces.

\end{remark}

We will now prove that the homotopy colimit of this diagram preserves the topology of $\mathcal{X}_{s}$. To that end, we will need the cofibrant replacements from Definitions~\ref{def:CofibrantReplacements1} and \ref{def:CofibrantReplacements2}.

\begin{definition}
    The sealed Kato-Nakayama functor is defined by $\mathcal{F}_{\rm SKN}(\tau,\sigma)=Z^{\sigma,\log}_{\tau,c}$, see Definition~\ref{def:CofibrantReplacements2}. The projection maps 
    \begin{equation*}
        {Z}^{\sigma,\log}_{\tau,c}\to {Z}^{\sigma,\log}_{\tau}
    \end{equation*} give rise to a natural  transformation 
    \begin{equation*}
        \mathcal{F}_{\rm SKN}\to \mathcal{F}_{\rm KN}. 
    \end{equation*}
    
\end{definition}
\begin{lemma}\label{lem:SealedFunctorHE2}
    The natural transformation $\mathcal{F}_{\rm SKN}\to \mathcal{F}_{\rm KN}$ is a homotopy equivalence. 
\end{lemma}
\begin{proof}
    This follows from Lemma~\ref{lem:RelativeInitialDegenerationsCofibrationsHE} and Proposition~\ref{prop:HomotopyEquivalentFunctors}.
\end{proof}
\begin{lemma}\label{lem:SealedLatchingCofibration2}
    The latching maps for $\mathcal{F}_{\rm SKN}$ are cofibrations. 
\end{lemma}
\begin{proof}
Note that the latching maps describe the spaces coming from the toric boundary. By construction of the cofibrant replacements, these spaces on the boundary inject into the given space, so we conclude. 
\end{proof}

\begin{theorem}\label{thm:KNFunctorHomotopyColimit3}
    Consider the Kato-Nakayama functor $\mathcal{F}_{\rm KN}:\Pscr(\Sigma_{D}) \to {\rm Top}$. There is a homotopy equivalence
    \begin{equation*}
        \mathrm{hocol}(\mathcal{F}_{\mathrm{KN}})\to \mathcal{X}^{\log,\zeta}_{\mathcal{D},s}.
    \end{equation*}
\end{theorem}
\begin{proof}
By Lemma~\ref{lem:SealedFunctorHE2}, we can work with the sealed functor. Using Lemma~\ref{lem:SealedLatchingCofibration2} and \cite[Proposition 3.1]{WZZ99}, we find that the map from the homotopy colimit to the colimit is a homotopy equivalence. 

We now have to construct a deformation retract from the colimit of $\mathcal{F}_{\rm SKN}$ to $\mathcal{X}^{\log,\zeta}_{\mathcal{D}, s}$. We have a map from this colimit to the Kato-Nakayama space by projecting every point in $V_{\sigma}(\mathbb{D})$ to $0_{\sigma}$. We now immediately find that the fibers are contractible. Since the spaces in question are compact,
we can apply Theorem~\ref{thm:SmaleVietoris}.

\end{proof}
\begin{corollary}\label{cor:KNHomotopColimitUsuisTheorem}
    Let $\mathcal{F}_{\rm KN}$  be as in Theorem~\ref{thm:KNFunctorHomotopyColimit3}. Then there is a homotopy equivalence 
    \begin{equation*}
        {\rm hocol}(\mathcal{F}_{\rm KN})\to \mathcal{X}_{t}\backslash \mathcal{D}_{t}
    \end{equation*}
    for general $t\in U(\CC)$.
\end{corollary}
\begin{proof}
This follows from Theorem~\ref{thm:KNFunctorHomotopyColimit3} and Theorem~\ref{thm:UsuisTheorem2}.
\end{proof}

\begin{example}
    Note that if we apply Theorem~\ref{thm:KNFunctorHomotopyColimit3} to $\mathbb{P}^{1}$ with $D=\{0,\infty\}$ and the standard fan $\Sigma$ of $\mathbb{P}^{1}$ with $\overline{\Sigma}=\mathbb{T}\mathbb{P}^{1}$, then the induced Kato-Nakayama diagram is the one considered in Example~\ref{exa:MainSealedExample3}. 

    More generally, if we use this on the fan $\Sigma=\Delta$ of a smooth toric variety, then we obtain the homotopy colimit from \cite[Section 5.3]{WZZ99}. To see this, we need to replace the Kato-Nakayama spaces by real tori $(S^{1})^{n}$.  This idea of replacing Kato-Nakayama spaces by angle sets will work more generally for arbitrary rational tropical manifolds, as we will see in Section~\ref{sec:AngleSets}. In Section~\ref{sec:ToricVarieties}, we will use this proposition to find CW complex structures on toric varieties. We note as in \cite{WZZ99} that this proposition automatically recovers Danilov's spectral sequence for the cohomology of a toric variety, see \cite[Chapter 3, Paragraph 12]{Danilov78}. 
    
   Theorem~\ref{thm:KNFunctorHomotopyColimit3} can also be used to recover the cohomology of any toric snc-compactification of an essential projective hyperplane complement. In this case the local components are given by the cohomology of local hyperplane complements, which are given by affine Orlik-Solomon algebras. It is not clear to us whether this recovers the de Rham exact sequence, as in the case of a toric variety. More generally, it might be possible to show locally that all the topological cohomology groups can be interpreted in terms of differential-theoretic cohomology groups. This would allow one to generalize the theory of tropical homology to the case of open rational tropical manifolds.
\end{example}

\begin{example}\label{exa:SurfaceTropicalDiagram}
    Consider the space $X=V(1+x+y+z)$ with tropicalization $\Sigma$ dual to the $3$-simplex. Note that the $3$-simplex also gives the tropicalization of the toric variety $Y(\Delta)=\mathbb{P}^{3}$. We determine the diagram of spaces here. 

    \begin{figure}[ht]
        \centering
        \includegraphics[scale=0.25]{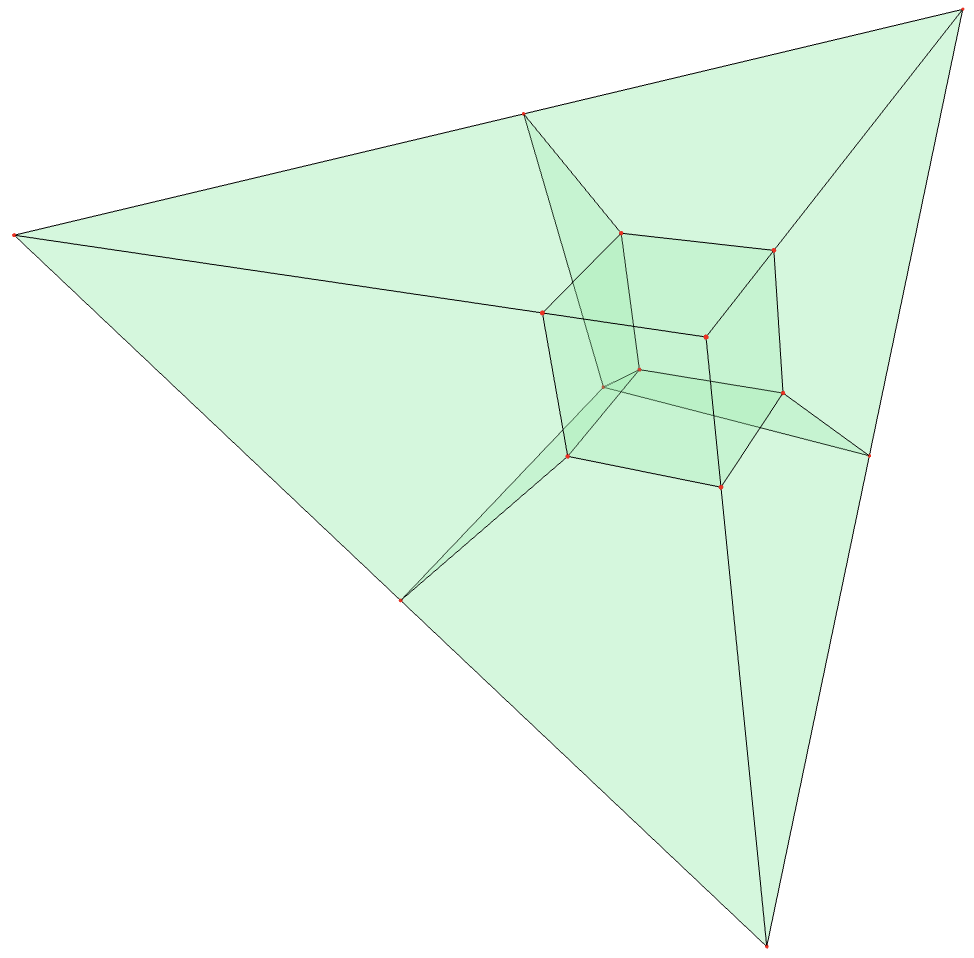}
        \caption{The tropical plane $\Sigma$ compactified in the tropical projective space $\mathbb{TP}^3$.}
        \label{fig:placeholder}
    \end{figure}

    We obtain three different types of initial ideals, up to an automorphism of the torus: 
    \begin{align*}
        I_{1} = \langle 1+x+y+z\rangle, \quad
        I_{2} = \langle 1+x+y\rangle,   \quad
        I_{3} = \langle 1+x\rangle.          
    \end{align*}
    The first ideal only occurs in the torus $\Spec(\CC[\ZZ^3])$. 
    It describes the space at the central vertex of the tropicalization. 
    Its vanishing set is isomorphic to a hyperplane complement $\mathbb{P}^2\backslash D$, where $D$ is the union of $4$ generic hyperplanes. 
    By \cite{Hattori75}, the topology of this space can be described by the $2$-skeleton of $\RR^{3}/\ZZ^3$. 
    We will investigate this space in more detail using angle sets in Section~\ref{sec:AngleSets}, see Example~\ref{exa:HattoriAngleSet}.
    
    For $I_{2}$, we can consider this ideal either in the torus $\Spec(\CC[\ZZ^3])$ or in the torus $\Spec(\CC[\ZZ^2])$. This gives two different topological spaces. The second is isomorphic to $\mathbb{P}^{1}\backslash\{0,1,\infty\}$, a standard pair of pants. The first is the product $\mathbb{P}^{1}\backslash\{0,1,\infty\}\times \CC^{\times}$. Note that the first space is associated to a ray $(\tau,0)$. The second space is associated to $(\tau,\tau)$, which is the vertex at infinity.  
    
    For $I_{3}$, we have three choices for the ambient torus: $\Spec(\CC[\ZZ^{i}])$ for $1\leq{i}\leq{3}$. This gives the spaces $(S^{1})^{i}$. The space $(S^{1})^{2}$ is associated to a $2$-cell $(\tau,0)$. For any of the adjacent rays at infinity $(\tau,\sigma)$, we then obtain a copy of $S^{1}$. Finally, the space for  $(\tau,\tau)$ is a point.  
\end{example}

\begin{example}\label{exa:K3SurfaceTropicalDiagram}
    Consider the surface given by the polynomial
    \begin{equation*}
        f=1+t(x+y+z+x^{-1}+y^{-1}+z^{-1}).
    \end{equation*}
    Its tropicalization can be found in Figure~\ref{fig:AperyK3}. If we take the compactification of $X$ in the corresponding toric variety, we obtain a K3-surface. 

    The spaces attached to the cells are the same as in Example~\ref{exa:SurfaceTropicalDiagram}. In fact, for any smooth tropical surface, the corresponding diagram can be described purely in terms of those spaces. We have added the corresponding angle space in Figure~\ref{fig:AperyK3}. We will see later on that these faithfully capture the topology of the different $\mathcal{F}_{\rm KN}(\tau,\sigma)$. See Figure~\ref{fig:AperyK32} for the corresponding diagram of angle spaces.   

    \begin{figure}[ht]
        \centering
        \scalebox{0.8}{\includegraphics[width=0.5\linewidth]{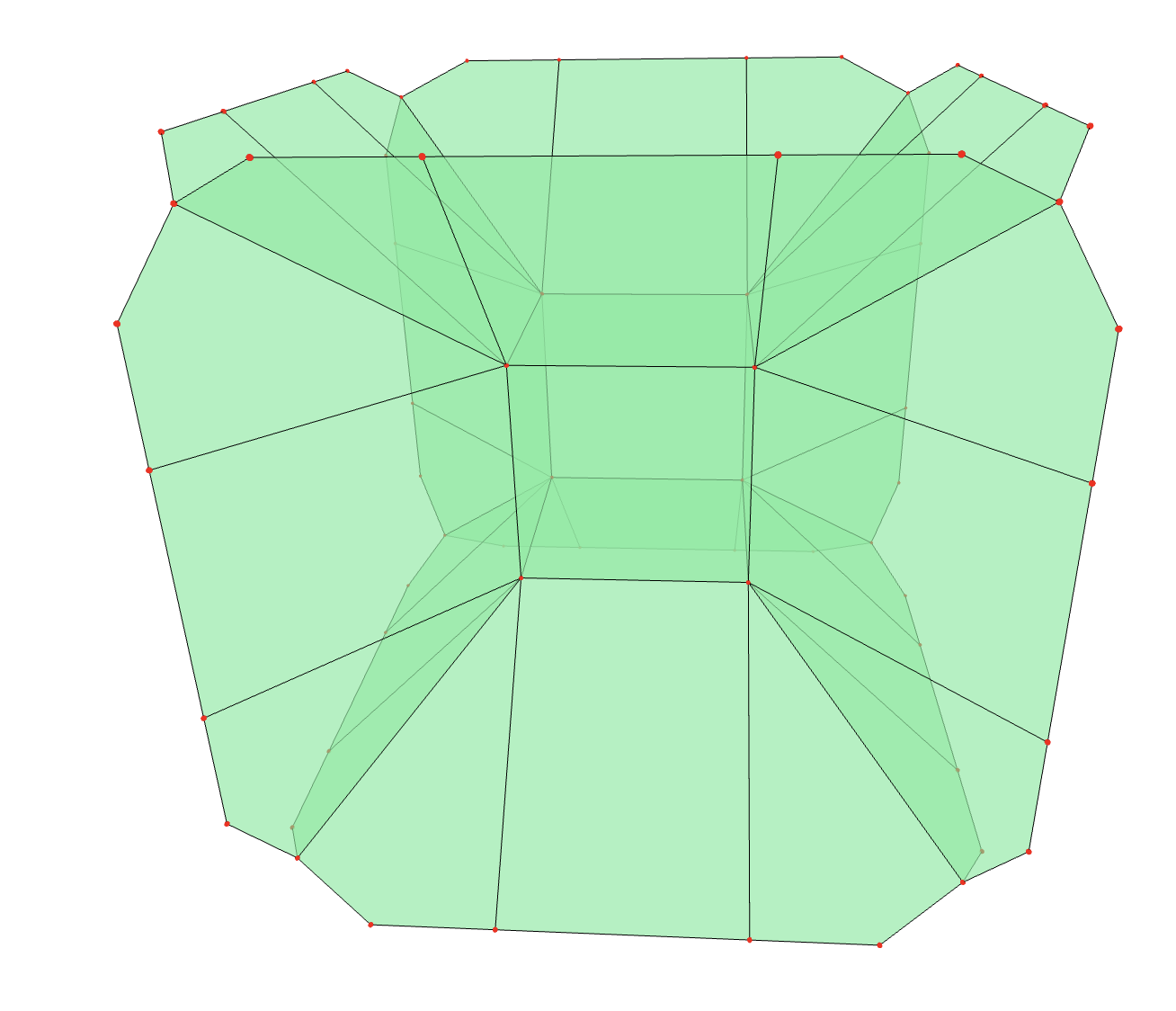}}
        \caption{The tropicalization of the Ap\'{e}ry family of K3-surfaces. The diagram is made up of $6$ different topological spaces. }
        \label{fig:AperyK3}
    \end{figure}
\end{example}
\begin{example}
    Consider the Dwork family 
    \begin{equation*}
        f=xyz+t(1+x^4+y^4+z^4).
    \end{equation*}
    The induced subdivision of the Newton polytope of $f$ is a triangulation with four simplices. 
    We conclude that the tropicalization is a rational tropical manifold by Lemma~\ref{lem:KummerLemma}. Note that the compactified variety $\overline{X}$ is a K3-surface. Moreover, the topological spaces at the boundary are not always connected. For instance, when the initial ideal is $1+x^4$, we have four points. 

    Similarly, we can consider the Fermat family 
    \begin{equation*}
        f=1+x^4+y^4+z^4.
    \end{equation*}
    This also defines a rational tropical manifold, although it is not as degenerate as the Dwork manifold. In a sense, the various topological spaces have increased in complexity for the Fermat family, but the tropicalization has decreased in complexity. 

    Similarly, if we consider a full unimodular triangulation of the simplex $4\cdot \Delta_{1}$, then the tropicalization is quite complex, but the topological spaces are simply essential projective hyperplane complements.   
\end{example}

For future reference, we now define the notion of a tropical diagram. In these, we have topological spaces that are homotopy equivalent to the $Z^{\sigma,\log}_{\tau}$. This is enough to capture the topology of a nearby fiber by Theorem~\ref{thm:KNFunctorHomotopyColimit3}. 

\begin{definition}\label{def:TropicalDiagram}
    Let $(X,\Sigma,D)$ be a rational tropical manifold. A tropical diagram is a functor
    \begin{equation*}
        \mathcal{F}:\Pscr{(\Sigma_{D})}\to {\rm Top}
    \end{equation*}
    together with a natural transformation $\mathcal{F}_{\rm KN}\to \mathcal{F}$ that is a homotopy equivalence on objects. 
\end{definition}
\begin{corollary}
    Let $(X,\Sigma,D)$ be a rational tropical manifold and let $\mathcal{F}$ be a tropical diagram. Then for any strictly semistable model $(\mathcal{X}\to U,\mathcal{D})$ associated to $(X,\Sigma,D)$, we have a homotopy equivalence 
    \begin{equation*}
        {\rm hocol}(\mathcal{F})\to \mathcal{X}_{t}\backslash \mathcal{D}_{t}
    \end{equation*}
    for general $t\in U(\CC)$.
\end{corollary}
\begin{proof}
    See Theorem~\ref{thm:KNFunctorHomotopyColimit3}, Proposition~\ref{prop:HomotopyEquivalentFunctors}, and Corollary~\ref{cor:KNHomotopColimitUsuisTheorem}.  
\end{proof}
\begin{remark}
In the upcoming sections, we will for instance define certain angle functors $\mathcal{F}_{\Theta}$, see Definitions~\ref{def:AngleFunctorGeneral} and \ref{def:AngleFunctorCI}. The angle map in these cases automatically induces a natural transformation $\mathcal{F}_{\rm{KN}}\to \mathcal{F}_{\Theta}$.
\end{remark}

\section{Angle sets}\label{sec:AngleSets}

    In this section, we study angle sets of linear spaces in the complex torus $(\CC^{\times})^n$ and their Kummer coverings. 
    Namely we show that any Kummer covering of a linear space in $(\CC^{\times})^n$ is homotopy equivalent to its angle set in $(S^1)^{n}$.
    
    \subsection{The angle set}
        We start by defining the \emph{angle set} or \emph{coamoeba}. We refer the reader to \cite{NS13NonArchimedean} and \cite{NS24} for more background on these angle sets. In this article, we shall however define these angle sets in the coordinate-free approach commonly used in toric geometry. 
        Note that this is similar to the approach used in Section~\ref{sec:ToricTropicalGeometry} to define the compactified tropicalization. 
        
        \begin{definition} 
            Let $M$ be a lattice and let $A=\CC[M]$ be the group ring over $\CC$ with torus $T=\Spec(A)$. We have 
            \begin{equation*}
                T(\CC)=\mathrm{Hom}(M,\mathbb{C}^{\times}),
            \end{equation*}
            where the latter is the set of group homomorphisms from $M$ to $\mathbb{C}^{\times}$. Consider the group homomorphism
            \[
            \begin{matrix}
            \ang  : & \mathbb{C}^{\times} & \longrightarrow & S^{1} \\
                              &         z      &  \longmapsto    &  z/|z|
            \end{matrix} \quad ,
            \]
            where $|z|=\sqrt{x^2+y^2}$ is the modulus of the complex number $z=x+iy$. By composition, we then obtain a morphism 
            \[
            T(\CC)=\mathrm{Hom}(M,\mathbb{C}^{\times})\to \mathrm{Hom}(M,S^{1}).
            \]
            This is the angle map for $T$. It is surjective and continuous. Here we endow $S^{1}$ and $\CC$ with the Euclidean topology, and $\mathrm{Hom}(M,S^{1})$ with the product topology.
        \end{definition}

    \begin{remark}
    After fixing a basis $e_{1}, \dots, e_n$ of $M$, we can identify $T(\CC)$ with $(\CC^\times)^n$ and $\mathrm{Hom}(M,S^{1})$ with $(S^{1})^{n}$.
    With respect to this basis, the angle map is given by 
    \[
        x \in T(\CC) \cong (\CC^\times)^{n} \to {\rm Hom}(M,S^1) \cong (S^1)^n , \quad  (z_j)_{j = 1}^{n} \mapsto (z_j/|z_{j}|)_{j=1}^{n}.
    \]
    Note in particular that the angle map is \emph{semi-algebraic}.
    \end{remark}
    \begin{example}\label{exa:ProjectiveAngles}
       Let $M := {\rm ker}(\mathds{1}_{n+1})$ be the character lattice of $\mathbb{P}^{n}$, where $\mathds{1}_{n+1}$ is the vector in $\ZZ^{n+1}$ whose coordinates are all equal to $1$. Explicitly, 
        \[
            M := \left\{ m =(m_j)_{j=0}^{n} \in \ZZ^{n+1} \colon m_0 + \dots + m_n = 0 \right\}.
        \]
        We then have an exact sequence 
        \begin{equation*}
            (0)\to M\to \ZZ^{n+1}\to \ZZ \to (0).
        \end{equation*}
        By applying $\mathrm{Hom}(-,\CC^{\times})$ to this, we obtain the exact sequence 
        \begin{equation*}
             (0)\to \CC^{\times}\to (\CC^{\times})^{n+1}\to \mathrm{Hom}(M,\CC^{\times}) \to (0).
        \end{equation*}
        This allows us to view the dense torus in $\mathbb{P}^{n}(\CC)$ as the quotient $(\CC^{\times})^{n+1}/\CC^{\times}$. Similarly, we can 
        identify $\mathrm{Hom}(M,S^{1})$ with $(S^{1})^{n+1}/S^{1}$ where $S^1$ acts diagonally, or equivalently sets of $n+1$ angles, up to simultaneous scaling. We denote these equivalence classes by $[\theta_{0}:...:\theta_{n}]$, in keeping with the usual notation for elements of projective space. 
        Similarly, if we are given an element $x=[x_{0}:...:x_{n}]\in\mathbb{P}^{n}(\CC)$ in the open dense torus, we then write 
        \begin{equation*}
            \mathrm{ang}(x)=[\ang(x_{0}):...:\ang(x_{n})].
        \end{equation*}
        for the induced angle map. Note that if $I^\ast$ is a homogeneous ideal in $\CC[\ZZ^{n+1}]$, then the angle set $\hat{\Theta} \subset (S^1)^{n+1}$ of the affine cone $\widehat{X} = V(I^\ast) \cap (\CC^{\times})^{n+1}$ is invariant under the action of $S^{1}$, and its quotient is the angle set of $X = \widehat{X} / \CC^{\times}$. If we explicitly want to find a presentation of $X$, then we have to find generators of the ideal $I = I^\ast \cap \CC[M]$. These can be found by dehomogenizing the ideal~$I^\ast$.      
    \end{example}
    
    \begin{remark}
    It will also be useful to think of the following transcendental approach for the angle map. Consider the group homomorphism $\mathbb{R}\to S^{1}$ defined by 
    \begin{equation*}
        \theta \mapsto e^{ i \theta}.
    \end{equation*}
    This induces an isomorphism $\mathbb{R}/2\pi\mathbb{Z}\simeq S^{1}\subset \CC^{\times}$. For a nonzero complex number, we then write $\theta(z)\in \RR$ for any angle such that $e^{i\theta(z)}=z/|z|$. Similarly, we write $\overline{\theta(z)}$ for its image in the quotient $\RR/2\pi\ZZ$.       %
    In terms of these coordinates, the angle map can be defined as follows: 
    \begin{equation*}
        (\mathbb{C}^{\times})^{n} \to (\mathbb{R}/2\pi\mathbb{Z})^{n}, \quad (z_{j})_{j=1}^{n}\mapsto (\overline{\theta(z_{j}}))_{j=1}^{n}.
    \end{equation*}
    This can be viewed as the argument function from complex analysis. 
    \end{remark}
    
    \begin{definition}\label{def:AngleMap}
    Let $A=\CC[M]$, and 
        let $X\to T=\Spec(A)$ be a closed embedding. We define the angle set of $X$ to be the image of $X(\CC)$ in $\mathrm{Hom}(M,S^{1})$ under 
        \begin{equation*}
            X(\CC)\to T(\CC)=\mathrm{Hom}(M,\mathbb{C}^{\times})\to \mathrm{Hom}(M,S^{1}). 
        \end{equation*}
        We denote it by $\Theta$.
    \end{definition}
    \begin{remark}
        Note that $\Theta$ is a semi-algebraic set.  Indeed, the angle map is semi-algebraic and $X$ is semi-algebraic, so this follows immediately from the Tarski-Seidenberg theorem.
    \end{remark}

\begin{figure}[ht]
\scalebox{0.9}{
    \begin{minipage}[c]{0.45\textwidth}
        \centering
        \includegraphics[scale=0.24]{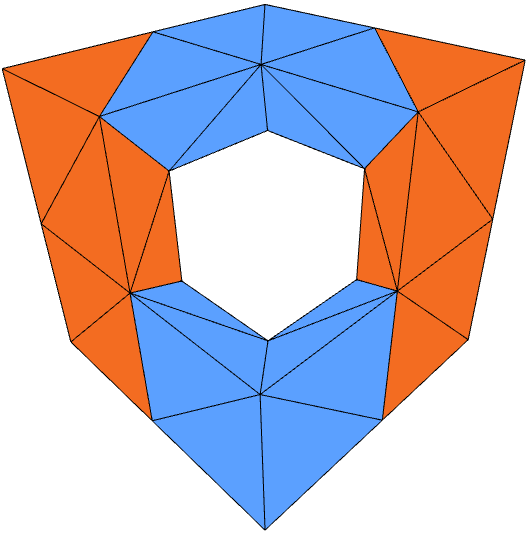}
        \end{minipage}
    \begin{minipage}[c]{0.45\textwidth}
        \centering
        \includegraphics[scale=0.45]{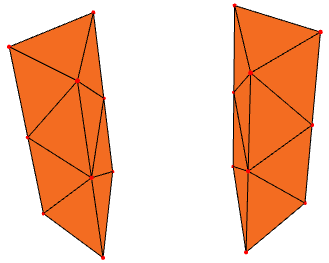}
        \end{minipage}
        }
        \caption{The completed angle set of $V(1+x+y+z)\subset (\CC^{\times})^{3}$ and the corresponding angle set of the initial degeneration $V(1+x+y)\subset(\CC^{\times})^{3}$. Note that the second is a subcomplex of the first. We will see in Section~\ref{sec:CompleteIntersections} that this corresponds to an inclusion of Kato-Nakayama spaces $Z^{\log}_{\sigma}\subset Z^{\log}$.}
        \label{fig:AngleSetInitialDegeneration}
\end{figure}

    \begin{example}\label{exa:HattoriAngleSet}
\label{ex:cube3}
    Consider the lattice $M=\ker{\mathds{1}_{4}}\subset \ZZ^{4}$. Let $e_{0}, \dots, e_{3}$ be the standard basis of $\ZZ^{4}$. Consider the homogeneous linear ideal $I^\ast=\langle \chi^{e_{0}} + \chi^{e_{1}}+\chi^{e_{2}}+\chi^{e_{3}} \rangle \subset \CC[\ZZ^{4}]$ with $I=\langle 1+x+y+z\rangle=I^\ast \cap \CC[M]$. Here $x = \chi^{e_1-e_0}$, $y = \chi^{e_2-e_0}$ and $z = \chi^{e_3-e_0}$. The space $X=V(I)$ is isomorphic to the complement of four generic hyperplanes in $\mathbb{P}^{2}$. Indeed, the map 
    \begin{equation*}
        [X:Y:Z]\mapsto [X:Y:Z:-X-Y-Z] 
    \end{equation*}
    realizes the hyperplane complement $Y=\mathbb{P}^{2}\backslash V\big(XYZ(-X-Y-Z)\big)$ as the intersection of the dense torus of $\mathbb{P}^{3}$ with the vanishing set of $I^\ast$. By \cite[Theorem 1]{Hattori75}, the complement $Y$ is homotopy equivalent to the $2$-skeleton of $\mathbb{R}^{3}/\mathbb{Z}^{3}$. Moreover, we have $\pi_{1}(Y)\simeq \mathbb{Z}^{3}$. The same is true for the angle set of $Y$, as guaranteed by Theorems~\ref{thm:MainTheoremAText} and \ref{thm:HEforKNSpaces}.
    
    We will study angle sets in more detail in Section~\ref{sec:AngleSetsCI}. For now, we note that the closure of the angle set is the complement of the following open zonotope
    \begin{equation*}
        Z := \Big\{\theta \in (-\pi, \pi)^3 : |\theta_{i}-\theta_{j}|<\pi \text{ for all }1 \leq i<j\leq{3}\Big\},
    \end{equation*}
    in the cube $C = [-\pi,\pi)^3$, by the argument in \cite[Lemma 3.1]{NS13NonArchimedean}. We can also view $Z$ as the Minkowski sum of the four line segments $(0,\pi w_j)$ with $1 \leq j \leq 3$ and the segment $(0, -\pi(w_{1}+w_{2}+w_{3}))$ where $w_1,w_2,w_3$ is the standard basis of $\RR^3$. Note that these vectors define the rays of the tropicalization of $X$. This open polytope $Z$ is the rhombic dodecahedron, which has $12$ faces, $24$ edges, and $14$ vertices. We will give a general formula for the $f$-vector of the angle sets of these types of very affine linear spaces in Theorem~\ref{thm:FVectorCI}.   
    
    Denote the closure of the complement of $Z$ in $[-\pi,\pi)^{3}$ by $\Theta^{\log}$. A pictorial representation of this set in $\mathbb{R}^{3}$ can be found in Figure~\ref{fig:AngleSetInitialDegeneration}. It is polyhedral, and it comes with a natural simplicial subdivision. This subdivision keeps track of all the different subspaces arising from initial degenerations.
    \end{example}

\subsection{Homotopy equivalences}\label{sec:HomotopyEquivalences}

In this section we recall two theorems that allow us to deduce global homotopy equivalences from local ones. We start by recalling a theorem due to McCord, which states that a continuous map $f:X\to Y$ that is a homotopy equivalence on a suitable open cover of $Y$, is also globally a homotopy equivalence. This will be one of the main ingredients in the proof of Theorem~\ref{thm:MainTheoremAText}. 
\begin{definition}\label{def:BaseLike}
    Let $\mathcal{U}$ be an open cover of a topological space $X$. We say that $\mathcal{U}$ is basis-like if the following holds. Suppose $x\in U\cap V$ for $U,V\in\mathcal{U}$. Then there exists a $W\in\mathcal{U}$ such that $x\in{W}\subset U\cap V$.  
\end{definition}

\begin{theorem}\label{thm:McCord}
    Suppose that $f:X\to Y$ is a continuous map of topological spaces, and suppose that $\mathcal{U}$ is a basis-like open cover such that the induced restrictions $f^{-1}(U)\to U$ are weak homotopy equivalences for all $U\in \mathcal{U}$. Then $f$ is a weak homotopy equivalence. 
\end{theorem}
\begin{proof}
  See \cite[Theorem 6]{McCord66}.
\end{proof}

We now quickly discuss the Smale-Vietoris theorem. This similarly allows us to deduce a global homotopy equivalence by imposing certain local conditions. In this case however, we assume that the fibers are contractible.

\begin{theorem}[Smale-Vietoris Theorem] \label{thm:SmaleVietoris}
Let $X$ and $Y$ be compact topological spaces that can be endowed with the structure of a simplicial complex. 
Let $f:X\to Y$ be a continuous map such that $f^{-1}(y)$ is contractible for every $y$ (in particular, $f$ is surjective). Then $f$ is a homotopy equivalence.
\end{theorem}
\begin{proof}
Note that $X$ and $Y$ are locally contractible by \cite[Appendix 4 p. 523]{Hatcher2001}. We then conclude using \cite{SmaleVietoris57}.
\end{proof}

\begin{remark}
    There are more general versions of this theorem that only require the morphism to be proper, and the complexes to be locally finite CW complex. We refer to  \cite[Page 17]{Lacher77} for a version of this.
\end{remark}

As a consequence, we deduce the following for semi-algebraic sets.
\begin{corollary}\label{cor:SmaleVietorisSemialgebraic}
    Let $f:X\to Y$ be a continuous morphism of compact semi-algebraic sets such that $f^{-1}(y)$ is contractible for every $y$. Then $f$ is a homotopy equivalence.
\end{corollary}

\begin{proof}
    This follows from Theorem~\ref{thm:SmaleVietoris} and the fact that semi-algebraic sets are simplicial complexes, see \cite[Theorems 2 and 3]{Lojasiewicz64}. 
\end{proof}

\begin{example}
    We give a simple example to show that the analogous statement for non-compact spaces is not true. Let $X=V(u(uv-1))\subset \RR^{2}$ and let $Y=\RR$. We consider the projection map $X\to Y$ given by $(u,v)\mapsto u$. Then the fibers are contractible, since they are all polyhedra. However, the map $X\to Y$ is not a homotopy equivalence, since $X$ has $3$ connected components, and $Y$ has only one. The naive unbounded version of the Smale-Vietoris theorem thus does not hold. 

    Note that the fiber over $u = 0$ is contractible, but the fiber over any $(-\epsilon,\epsilon)$ for $\epsilon>0$ small is not contractible. In our main theorem, this will not be the case.
\end{example}

\subsection{Proof of \texorpdfstring{Theorem~\ref{thm:MainTheoremAText}}{Theorem 4.32}}

    To prove that the angle map is a homotopy equivalence, we will find an open basis-like covering of the angle set such that the induced angle maps are homotopy equivalences. In fact, the opens in our covering will be contractible by a construction using the local conical structure theorem, see Lemma~\ref{lem:OpenContractibleNeighborhood}. Their inverse images will be intersections of cones with an affine linear space, so these will be convex and thus contractible, see Lemma~\ref{lem:ConvexityLemma}. We then conclude using McCord's theorem, see \cite[Theorem 6]{McCord66}.
    
    Throughout this section, we will use the following notation. 
    Let $\norm{\cdot}$ be a semi-algebraic norm on $\RR^n$ i.e. a norm whose unit ball is a semi-algebraic subset of $\RR^n$. We denote by $B(a,\epsilon)$ and $S(a,\epsilon)$ the open ball and the sphere of center $a$ and radius $\epsilon$ respectively.

    \begin{remark}
        The norm $\norm{\cdot}$ could for example be the usual Euclidean norm $\norm{\cdot}_2$, or the supremum norm $\left\lVert x \right\rVert_{\infty} := \max\limits_{1 \leq i \leq n} |x_{i}|$. Note that the latter is again semi-algebraic and continuous, since $\max(\cdot, \cdot)$ and $|\cdot|$ are both semi-algebraic and continuous. 
     \end{remark}
    
    \subsubsection{Local conical structure}\label{sec:LocalConicalStructure}
    We start by collecting some useful results from real semi-algebraic geometry.
    \begin{definition}
        Let $Y\subset \RR^{n}$ be any semi-algebraic set and let $a\in\RR^{n}$. The \emph{cone $C(a,Y)$ of $Y$ over $a$} is defined as follows:
        \[ 
            C(a, Y) := \Big\{ \lambda{a}+(1-\lambda)y  \colon \text{for } \lambda \in [0,1]  \text{ and } y \in Y \Big\}.
        \]
    \end{definition}
    
    \begin{lemma}\label{lem:InjectivityCone}
    Suppose that $Y=S(a, \epsilon)\cap{A}$ for some semi-algebraic set $A \subset \RR^n$. Then the map $[0,1)\times Y\to C(a,Y)$ given by $(\lambda,y)\mapsto \lambda{a}+(1-\lambda)y$ is injective. 
    \end{lemma}
    \begin{proof}
        Suppose that $z=\lambda{a}+(1-\lambda){y}$. Then $z-a=(1-\lambda)(y-a)$. Taking absolute values, we  find that $\left\lVert z-a\right\rVert=|1-\lambda|\cdot \epsilon$. In particular, $\lambda$ is unique for $z$, provided that $\lambda\neq{1}$. If $y_{1}\neq{y_{2}}$ with $z=\lambda{a}+(1-\lambda){y}_{i},$ then subtracting yields $(1-\lambda)(y_{1}-y_{2})=0$, a contradiction. We conclude that the map is injective.       
    \end{proof}
    
    \begin{remark}
    Note that the map above is almost surjective. Namely, the only missing point is the cone point $a$. 
    \end{remark}
    For the following theorem, we recall that a non-isolated point of a semi-algebraic set $A$ is a point $a$ such that for every $\epsilon>0$, we have that $B(a, \epsilon)\cap A$ contains points other than $a$. In other words, $\{a\}$ is not open.   
    
    \begin{theorem}[\textbf{Local conical structure}, Theorem 9.3.6 in \cite{BochnakCosteRoy}] \label{thm:LocalConicalStructure}
       Suppose that $a$ is a non-isolated point of $A$  and let $\left\lVert\cdot\right\rVert$ be any semi-algebraic norm on $\RR^{n}$. Then there exists an $\epsilon>0$ such that 
        \begin{equation*}
            \overline{B}(a, \epsilon)\cap A\simeq C(a,S(a, \epsilon)\cap{A}).
        \end{equation*}
        That is, the closed $\epsilon$-neighborhood of $A$ around $a$ is a cone over $a$ with base $S(a, \epsilon)\cap{A}$, the closed disk of radius $\epsilon$ around $a$.  
    \end{theorem}
    
    \begin{corollary}\label{cor:LocalContractibility}
        The set $\overline{B}(a, \epsilon)\cap A$ is contractible for small enough $\epsilon$. The same is true for $B(a,\epsilon)\cap A$. 
    \end{corollary}
    \begin{proof}
    If $a$ is an isolated point, then this is trivial. We assume $a$ is not an isolated point. 
        By Theorem \ref{thm:LocalConicalStructure}, we only have to show that $C(a,Y)$ is contractible for $Y = S(a, \epsilon)\cap{A}$. It suffices to show that it is a \emph{star domain}. That is, there is a point $P$ such that for any other point $Q$, the line segment between $P$ and $Q$ lies in $C(a,Y)$. We choose $P = a \in C(a,Y)$. 
        Then the line segment between $a$ and any other point in $C(a,Y)$ is automatically in $C(a,Y)$. Indeed, we calculate 
        \begin{align*}
            \mu{a}+(1-\mu)(\lambda{a}+(1-\lambda)y)=a(\mu+\lambda-\mu\lambda)+(1-\mu)(1-\lambda)y.
        \end{align*}
        But $1-(\mu+\lambda-\mu\lambda)=(1-\mu)(1-\lambda)$. Moreover, one easily checks that these values lie in between $0$ and $1$. We conclude that the line segment is in $C(a,Y)$, as desired. The open case follows in exactly the same way.
    \end{proof}

    \begin{definition}
    Let $a,b$ be real numbers with $a < b$. We define four maps
    \begin{align*}
        &\phi_{(1,\pm)}: (a,b) \times \RR_{>0}\to \RR^{2}, \qquad (s,r) \mapsto \pm\left(\dfrac{r}{\sqrt{1+s^2}},\dfrac{sr}{\sqrt{1+s^2}}\right), \\
        &\phi_{(2,\pm)}: (a,b) \times \RR_{>0}\to \RR^{2}, \qquad  (s,r) \mapsto \pm\left(\dfrac{sr}{\sqrt{1+s^2}},\dfrac{r}{\sqrt{1+s^2}}\right).
    \end{align*}
    We denote by $C_{(j,\pm)}(a,b)$ for $j=1,2$ the four open convex polyhedral cones in $\RR^2 \cong \CC$ defined as follows:
    \begin{align*}
        C_{(1,\pm)}(a,b) &:= \pm {\rm relint} \Big({\rm cone}\big((1,a), (1,b)\big) \Big), \\
       C_{(2,\pm)}(a,b)  &:= \pm {\rm relint}\Big( {\rm cone}\big((a,1), (b,1)\big) \Big). 
    \end{align*}
\end{definition}

\begin{lemma}
    The maps $\phi_{(j,\pm)} \colon (a,b) \times \RR_{>0} \longrightarrow C_{(j,\pm)}(a,b)$ are semialgebraic homeomorphisms.
\end{lemma}
\begin{proof}
    We define an inverse for each of the maps $\phi_{(j,\pm)}$ as follows. 
    \[
    \psi_{(1,\pm)}(x,y)= \left( \frac{y}{x}, \pm\sqrt{x^2+y^2}\right) \quad \text{and}  \quad \psi_{(2,\pm)}(x,y) = \left(\frac{x}{y}, \pm\sqrt{x^2+y^2}\right).
    \]
    Note that these maps are semi-algebraic and well defined on the respective $C_{(j,\pm)}(a,b)$. Moreover one can easily check that the maps $\psi_{(j,\pm)}$ are indeed inverses of the maps $\phi_{(j,\pm)}$.
\end{proof}

\begin{definition}\label{def:Accounting}
    Let $a_{\ell}< b_{\ell}\in \RR$ for $\ell=1,...,n$, with vectors $a=(a_{\ell})_{\ell=1}^{n}$ and $b=(b_{\ell})_{\ell=1}^{n}$. 
    Let $p_{\ell}\in \{1,2\}\times\{+, -\}$, with $p=(p_{\ell})_{\ell=1}^{n}$. We then write 
    \begin{equation*}
        C_{p}(a,b) := \prod_{\ell=1}^{n} C_{p_{\ell}}(a_{\ell},b_{\ell}).
    \end{equation*}
    As $a$, $b$, and $p$ vary, the sets $C_p(a,b)$ form a basis for the topology of $(\CC^\times)^n$. 
    Note that the semi-algebraic isomorphisms $\phi_{p_{\ell}}$ induce a semi-algebraic isomorphism
    \begin{equation*}
        \phi_{p}:\prod_{\ell=1}^{n} \big( (a_{\ell},b_{\ell})\times \RR_{>0} \big) \longrightarrow C_{p}(a,b).
    \end{equation*}

\end{definition}

\begin{definition}
    Consider the space $\prod_{i=1}^{n}\Big((a_{i},b_{i})\times \RR_{>0}\Big)\subset \RR^{2n}$. We endow $\RR^{2n}$ with the supremum norm
    \begin{equation*}
        \left\lVert x \right\rVert_{\infty} := \max_{1 \leq i \leq 2n} |x_{i}|, \quad \text{for }  x = (x_{i})_{i=1}^{2n}  \in  \RR^{2n}.
    \end{equation*}
    Note that the map $\RR^{2n}\to \RR$ given by $x\mapsto \left\lVert x\right\rVert_{\infty}$ is semi-algebraic. This is because $\max(\cdot)$ and $|\cdot|$ are semi-algebraic. 
\end{definition}

\begin{definition}

   For $i=1,...,n$, we have the following slope and radius projection maps 
   \[
   \pi_{i,s}:\prod_{\ell = 1}^{n} \Big((a_{\ell},b_{\ell})\times \RR_{>0}\Big)\to (a_{i},b_{i}) \quad \text{and} \quad \pi_{i,r}:\prod_{i=1}^{n}\Big((a_{i},b_{i})\times \RR_{>0}\Big)\to \RR_{>0}.
   \]
   We call these the \emph{slopes} and \emph{radii} projections. 
   Consider the semi-algebraic homeomorphism 
     \[
     \phi_{p}:\prod_{\ell = 1}^{n}\Big((a_{\ell},b_{\ell})\times \RR_{>0}\Big)\to C_{p}(a,b).
     \]
    The slopes (resp. radii) of a point $z \in C_{p}(a,b)$ are the slopes (resp. radii) of $\phi^{-1}_{p}(z)$.  
\end{definition}

\begin{definition}
    Let $p$, $a$ and $b$ be as in Definition \ref{def:Accounting}. We define the following
    \[
    \Theta_{p}(a,b) := \Theta\cap C_{p}(a,b) \quad \text{and} \quad \mathcal{S}_{p}(a,b):=\phi_{p}^{-1}(\Theta_{p}(a,b)).
    \]
    Note that we have a semi-algebraic homeomorphism $\phi_{p} : \mathcal{S}_{p}(a,b)\to \Theta_{p}(a,b)$.
\end{definition}

\begin{lemma}\label{lem:OpenContractibleNeighborhood}
    For every $z\in \Theta$, there are vectors $a$, $b$ and $p$ such that $\Theta_{p}(a,b)$ is a contractible open neighborhood of $z$.
\end{lemma}
\begin{proof}
First, consider any open neighborhood $C_{p}(a,b)$ of $z$ and write
\[
z'=\phi_{p}^{-1}(z).
\]
Note that the radii of $z'$ are all equal to $1$. We choose $\epsilon>0$ sufficiently small such that
\[
0<\epsilon<1
\quad
\text{and}
\quad
a_{i}<\pi_{i,s}(z')-\epsilon
<
\pi_{i,s}(z')+\epsilon<b_{i}
\quad \text{for all } i=1,\ldots,n.
\]
Consider the open ball $B(z',\epsilon)$ with respect to the supremum norm. Explicitly, we have
\[
B(z',\epsilon)
=
\prod_{\ell=1}^{n}
\Big(
\big(\pi_{\ell,s}(z')-\epsilon,
\pi_{\ell,s}(z')+\epsilon\big)
\times
(1-\epsilon,1+\epsilon)
\Big).
\]
By Corollary \ref{cor:LocalContractibility}, after making $\epsilon$ smaller if necessary, the set $B(z',\epsilon)\cap\mathcal{S}_{p}(a,b)$ is contractible. 

Let $s(z') := \big(\pi_{i,s}(z')\big)_{i=1}^{n}$. Since every point of $\mathcal{S}_{p}(a,b)$ maps under $\phi_{p}$ to a point of
\[
\Theta\subset (S^{1})^{n},
\]
all of its radius coordinates are equal to $1$. Hence the radius conditions defining $B(z',\epsilon)$ are automatically satisfied on $\mathcal{S}_{p}(a,b)$, and intersecting with $B(z',\epsilon)$ only restricts the slope coordinates. We therefore have
\[
\phi_{p}
\Big(
B(z',\epsilon)\cap\mathcal{S}_{p}(a,b)
\Big)
=
\Theta_{p}
\Big(
s(z')-\epsilon,
s(z')+\epsilon
\Big),
\]
where, as before, we write $\epsilon$ for
$\epsilon\cdot(1,\ldots,1)$. Since $\phi_{p}$ is a homeomorphism, it follows that
\[
\Theta_{p}
\Big(
s(z')-\epsilon,
s(z')+\epsilon
\Big)
\]
is contractible. Moreover, it is an open neighborhood of $z$ in $\Theta$. This proves the result.
\end{proof}

\begin{lemma}\label{lem:SlopeCoveringBaseLike}
Let $\mathcal{U}$ be the collection of all contractible subsets of $\Theta$ of the form $\Theta_p(a,b)$. 
Then $\mathcal{U}$ is a basis-like open covering
of $\Theta$.
\end{lemma}

\begin{proof}
By Lemma \ref{lem:OpenContractibleNeighborhood}, every point
$z\in\Theta$ is contained in some contractible set $\Theta_p(a,b)$. 
Therefore $\mathcal{U}$ is an open covering of $\Theta$. 
Now suppose that
\[
    z \in \Theta_p(a,b) \cap \Theta_q(c,d)
\]
for two elements of $\mathcal{U}$. 
Since the sets $C_p(a,b)$, as $a$, $b$, and $p$ vary, form a basis for the topology of $(\CC^\times)^n$, there exists $C_r(a',b')$ such that
\[
    z\in C_r(a',b') \subset C_p(a,b)\cap C_q(c,d).
\]
By Lemma \ref{lem:OpenContractibleNeighborhood}, we may choose
$r,a',b'$ sufficiently small so that $\Theta_r(a',b')$ is
contractible. 
Intersecting with $\Theta$, we obtain
\[
    z \in\Theta_r(a',b') \subset \Theta_p(a,b)\cap\Theta_q(c,d).
\]
Thus $\mathcal{U}$ is basis-like.
\end{proof}

\begin{definition}
    We define $X_{p}(a,b) := \ang^{-1}\Big(\Theta_{p}(a,b)\Big)$ to be the inverse image of $\Theta_{p}(a,b)$ under the angle map $\ang:X\to \Theta$. 
\end{definition}

\begin{lemma}\label{lem:CompatibilityLemmaLocalAngles}
    We have $X_{p}(a,b)=X \cap C_{p}(a,b)$. 
\end{lemma}
\begin{proof}
Recall that $C_p(a,b)$ is invariant under coordinate-wise multiplication by positive real numbers. Let $z = (z_i)_{i=1}^{n} \in X$. Since
$ \ang(z) = \left(z_1/|z_1|, \dots, z_n/|z_n| \right)$,
we have
\[
z\in C_p(a,b)
\quad \text{if and only if} \quad
\ang(z)\in C_p(a,b).
\]
Suppose first that $z\in X_p(a,b)$. 
Then $\ang(z)\in\Theta_p(a,b)\subset C_p(a,b)$, and hence $z\in C_p(a,b)$. 
Therefore
\[
    X_p(a,b)\subset X\cap C_p(a,b).
\]

Conversely, suppose that $z\in X\cap C_p(a,b)$. Then $\ang(z)\in C_p(a,b)$.
Since $z\in X$, we also have $\ang(z)\in\Theta$. 
Therefore we obtain
\[
\ang(z)\in\Theta\cap C_p(a,b)=\Theta_p(a,b),
\]
so $z\in X_p(a,b)$. We conclude that $X_p(a,b)=X\cap C_p(a,b)$.
\end{proof}

\begin{lemma}\label{lem:ConvexityLemma}
    The set  $X_{p}(a,b)$ is convex. 
\end{lemma}
\begin{proof}
    Let $F_{1},...,F_{k}$ be the affine linear equations that cut out $X$. We write $\overline{X}$ for the affine linear space cut out by the $F_{i}$ in $\CC^{n}$. Note that $\overline{X}$ is convex. Indeed, it is the vanishing set of the affine linear forms ${\rm Re}(F_{i})$ and ${\rm Im}(F_{i})$. Since $C_{p}(a,b)$ is convex, we deduce that $X_{p}(a,b) = C_{p}(a,b) \cap \overline{X} \subset (\CC^\times)^n$ is convex by Lemma \ref{lem:CompatibilityLemmaLocalAngles}, as desired.
\end{proof}

\begin{theorem}\label{thm:MainTheoremAText}
   Let $X \subset (\CC^{\times})^n$ be a very affine linear space. Then the \emph{angle map}
\[
    \ang: X \longrightarrow (S^1)^{n}, \qquad x = (x_1, \dots, x_n) \longmapsto \Big(x_1/|x_1|, \dots, x_n/|x_n|\Big)
\]
is a homotopy equivalence from $X$ onto its image $\Theta := \ang(X)$.
\end{theorem}
\begin{proof}
By Lemma \ref{lem:OpenContractibleNeighborhood}, for every $z\in\Theta$, there are $a$, $b$ and $p$ such that $\Theta_{p}(a,b)$ is an open contractible neighborhood of $z$. These form a basis-like open cover by contractible neighborhoods by Lemma \ref{lem:SlopeCoveringBaseLike}. Moreover, note that the inverse images $X_{p}(a,b)$ are convex by Lemma \ref{lem:ConvexityLemma} and thus contractible (since it is non-empty). We now conclude using Theorem \ref{thm:McCord} that $\ang \colon X \to \Theta$ is a weak equivalence. Since $X$ and $\Theta$ are semi-algebraic sets, they are CW complexes, and thus a weak equivalence implies a homotopy equivalence by Whitehead's theorem.
\end{proof}

\subsubsection{The homotopy equivalence for Kummer coverings}

We now show how the same technique also gives a homotopy equivalence for Kummer coverings. The main idea is simply to use the fact that angles work well with power maps. 

\begin{definition}
Let $\kappa_{m}:\Spec(\CC[M])\to \Spec(\CC[M])$ be the $m$-th Kummer map associated to the multiplication-by-$m$ map on $M$. 
We say that $C_{p}(a,b)$ is $m$-admissible for $m \in \ZZ_{\geq 1}$ if its inverse image under the map $\kappa_{m}$ is a disjoint union of open sets, each of which
maps homeomorphically onto $C_{p}(a,b)$ under $\kappa_{m}$.
\end{definition}

\begin{lemma}\label{lem:SmallMAdmissible}
Let $m\in\mathbb{N}$ be fixed. 
    Let $z=(z_{i})_{i=1}^{n}\in (S^{1})^{n}$. Then any sufficiently small open neighborhood $C_{p}(a,b)$ of $z$ is $m$-admissible.    
\end{lemma}
\begin{proof}
    Choose $C_p(a,b)$ sufficiently small so that, in each coordinate, the angle width is strictly less than $2\pi/m$. Then the map $z\mapsto z^m$ has $m$ disjoint inverse branches on each coordinate sector, each mapping homeomorphically onto that sector. 
    Taking products over the $n$ coordinates, we find that $\kappa_m^{-1}(C_p(a,b))$ is a disjoint union of $m^n$ open sets, each of which maps homeomorphically onto the open set $C_p(a,b)$.
\end{proof}

\begin{theorem}[Theorem~\ref{thm:MainTheoremKummerCoveringAngleHE}]\label{thm:KummerHomotopyEquivalence}\label{rem:KummerCoveringGeneralHomotopyEquivalence}
    Let $X \to Z$ be a Kummer covering of a very affine linear space $Z$. Then the angle map $\ang\colon X \to \Theta$ is a homotopy equivalence. 
\end{theorem}
\begin{proof}
    We start by showing the result for the standard Kummer covering $X = Z_m$. 
    We take as our covering of $\Theta_{m}$ the inverse image under $\Theta_{m}\to \Theta$ of all $U_{\alpha}\subset \Theta$ from the proof of Theorem \ref{thm:MainTheoremAText} that are $m$-admissible. It follows from Lemma \ref{lem:SmallMAdmissible} that these cover $\Theta_{m}$.
    Note that, for each $U_{\alpha}$, its inverse image is a disjoint union of $m^{n}$ copies of $U_{\alpha}$.
    
    By our earlier observations, we find that the map $\ang^{-1}(U_{\alpha})\to U_{\alpha}$ is a homotopy equivalence. 
    Moreover, the map 
    \[
        \kappa^{-1}_{m}(\ang^{-1}(U_{\alpha}))\to \ang^{-1}(U_{\alpha})
    \]
    is a covering map of degree $m^{n}$, in the ordinary sense. 
    Since the set $\ang^{-1}(U_{\alpha})$ is contractible, we find that
    $\kappa^{-1}_{m}(\ang^{-1}(U_{\alpha}))$ is a disjoint union of
    $m^{n}$ copies of $\ang^{-1}(U_{\alpha})$. 
    Moreover, by $m$-admissibility, the inverse image of $U_{\alpha}$ in $\Theta_{m}$ is a disjoint union of $m^{n}$ open sets, and the angle map sends each component of $\kappa^{-1}_{m}\big(\ang^{-1}(U_{\alpha})\big)$ onto the corresponding component of this inverse image. In particular, the angle map on each component automatically yields
    a homotopy equivalence onto its image.
    
    The resulting open covering of $\Theta_m$ is basis-like. 
    Indeed, if two lifted components intersect at a point $z\in\Theta_m$, then their images in $\Theta$ contain the image of $z$. 
    Since the covering of $\Theta$ is basis-like, we may choose a smaller $U_{\beta}$ contained in their intersection, and after shrinking $U_{\beta}$ further if necessary, we may assume that it is $m$-admissible. 
    The component of its inverse image containing $z$ is then contained in the intersection of the two given lifted components.

    Finally, for a general Kummer covering there exists $m \geq 1$ and a subgroup $H$ of $M/mM$ such that $X = Z_m/H$. So the angle set of $X$ is the quotient $\Theta_m / H$.
    The above $m$-admissible covering is $H$-equivariant, so after taking quotients the same argument gives a basis-like open covering of $\Theta_m/H$ whose inverse images in $X$ are contractible. 
    Hence McCord's theorem applies in exactly the same way.
    We then conclude $X\to \Theta_{m}$ is a homotopy equivalence.
\end{proof}

   \begin{example}
        In order to calculate with the space $\Theta_{m}$, it suffices to determine the angle set of $\Theta$. Indeed, let $\widetilde{\Theta}$ be the pullback of $\Theta$ under the universal covering map 
        \begin{equation*}
            \mathbb{R}^{n}\to (S^{1})^{n} \cong \mathbb{R}^{n}/\mathbb{Z}^{n}.
        \end{equation*}
        Note that this is obtained from a single copy of $\Theta$ by gluing infinitely many copies along $\ZZ^{n}$. Any $\Theta_{m}$ is then obtained from $\tilde{\Theta}$ by quotienting by the subgroup $m\ZZ^{n}$. In terms of fundamental domains, if $[0,1]^{n}$ is the fundamental domain for $\Theta$, then $[0,m]^{n}$ is the fundamental domain for $\Theta_{m}$. 
        We call $\tilde{\Theta}$ the infinite Kummer covering of $X=V(I)$. See Figure \ref{fig:Examplen=2} for the infinite Kummer covering for $f=1+x+y$. We stopped at $m=5$ here, so we can view this as determining the topology of $V(1+x^5+y^5)\subset (\CC^{\times})^{2}$. Note that its smooth closure defines the Fermat curve of genus $6$.   
        \begin{figure}[ht]
        \centering    
            \begin{tikzpicture}[scale=0.5, line width = 0.5, clip=true]
                \def\mincoord{-5}
                \def\maxcoord{5}

                \foreach \i in {-2,...,2} {
                    \foreach \j in {-2,...,2} {
                        \pgfmathsetmacro\xshift{2*\i}
                        \pgfmathsetmacro\yshift{2*\j}
            
                        \filldraw[fill=orange!70, draw=black] 
                            ({-1 + \xshift}, {0 + \yshift}) --
                            ({0 + \xshift}, {1 + \yshift}) --
                            ({-1 + \xshift}, {1 + \yshift}) -- cycle;
            
                        \filldraw[fill=orange!70, draw=black] 
                            ({0 + \xshift}, {-1 + \yshift}) --
                            ({1 + \xshift}, {-1 + \yshift}) --
                            ({1 + \xshift}, {0 + \yshift}) -- cycle;
                    }
                }

                \foreach \m in {-2,...,2} {
                    \draw[black, line width = 0.2] (2*\m, \mincoord) -- (2*\m, \maxcoord);
                }

                \foreach \m in {-3,...,2} {
                    \draw[blue, line width = 2] (2*\m + 1, \mincoord) -- (2*\m + 1, \maxcoord);
                }

                \foreach \m in {-2,...,2} {
                    \draw[black, line width = 0.2] (\mincoord, 2*\m) -- (\maxcoord, 2*\m);
                }
                \foreach \m in {-3,...,2} {
                    \draw[green, line width = 2] (\mincoord, 2*\m + 1) -- (\maxcoord, 2*\m + 1);
                }

                \foreach \m in {0,...,5} {
                    \draw[black, line width = 0.2] 
                        (\maxcoord - 2*\m, \mincoord) -- 
                        (\maxcoord, \mincoord + 2*\m);
                }
                \foreach \m in {0,...,5} {
                    \draw[black, line width = 0.2] 
                        (\mincoord, \mincoord + 2*\m) -- 
                        (\maxcoord - 2*\m, \maxcoord);
                }

                \foreach \m in {0,...,4} {
                    \draw[red, line width = 2] 
                        (\mincoord, \mincoord + 2*\m + 1) -- 
                        (\maxcoord - 2*\m - 1, \maxcoord);
                }
                \foreach \m in {0,...,4} {
                    \draw[red, line width = 2] 
                        (\maxcoord - 2*\m-1, \mincoord) -- 
                        (\maxcoord, \mincoord + 2*\m+1);
                }

                \draw[->, black, line width = 1] (-5, 0) -- (5.5, 0) node[right] {$x$};
                \draw[->, black, line width = 1] (0, -5) -- (0, 5.5) node[above] {$y$};

                \draw[black, line width = 3] (1,1) -- (-1,1) -- (-1,-1) -- (1,-1) -- cycle;

                \draw (0.5,-6) node {$\Sigma_{[2]} \subset \mathcal{K}_{[2]}$};

                 \foreach \m in {-5,...,5} {
                        \foreach \n in {-5,...,5} {
                                \draw [fill=black] (\m, \n) circle (2pt);
                    }
                }
                
            \end{tikzpicture}
            \caption{The lifting of the angle set of $V(1+x+y)\subset (\CC^{\times})^{2}$ to $\RR^{2}$. The current fundamental domain determines the topology of the open Fermat curve $V(1+x^5+y^5)\subset (\CC^{\times})^{2}$.}
        \label{fig:Examplen=2}
    \end{figure}
    \end{example}

\subsection{Gluing angle sets}\label{sec:GluingAngleSets}

    Recall from Section \ref{sec:InitialDegenerationsKNSpaces} that $Z^{\sigma,\log,\circ}_{\tau}$ is the relative logarithmic initial degeneration. We view this as a subset of the $\CC$-valued points of $\Spec(\CC[\sigma^{\perp}\cap M])$. We also have the initial degeneration ${\rm in}_{w}(X)_{\sigma}\subset \Spec(\CC[\sigma^{\perp}\cap M])$, and there is a natural quotient map 
    \begin{equation*}
    \mathrm{in}_{w}(X)_{\sigma}\to Z^{\sigma,\log,\circ}_{\tau},
    \end{equation*} 
    see Lemma \ref{lem:InitialDegenerationsVSKNSpaces}. 
    \begin{lemma}
    The angle sets of $\mathrm{in}_{w}(X)_{\sigma}$ and $Z^{\sigma,\log,\circ}_{\tau}$ inside $\Spec(\CC[M(\sigma)])$ coincide. 
    There is moreover an action of $\mathbb{T}_{\sigma}(S^{1})$ on the angle set of ${\rm in}_{w}(X)$, and its quotient is the angle set of ${\rm in}_{w}(X)_{\sigma}$. 
    \end{lemma}
    \begin{proof}
    This follows from Lemma \ref{lem:InitialDegenerationsVSKNSpaces}. 
    \end{proof}
    
    \begin{definition}
 The relative open logarithmic angle set is the angle set of $\mathrm{in}_{w}(X)_{\sigma}\subset \Spec(\CC[M(\sigma)])$ for $w\in {\rm relint}(\tau)$. We denote it by $\Theta^{\sigma,\log,\circ}_{\tau}\subset {\rm Hom}(M(\sigma),S^{1})$.     
    \end{definition}
    
    We now construct an analogue of $Z^{\sigma,\log}_{\tau}$ on the level of angles. Recall that $Z_{\sigma}$ is a smooth variety inside a toric variety $Y(\star(\sigma))$ with dense torus $\Spec(\CC[\sigma^{\perp}\cap M])$. Moreover, the toric boundary is an snc divisor, and we have a natural log structure. 
    \begin{definition}
    Let $\mathcal{M}_{Z_{\sigma},x}$ be the monoid at $x$ associated to the log structure on $Z_{\sigma}$ coming from its toric boundary. Note that $\sigma^{\perp}\cap M=:M(\sigma)\subset \mathcal{M}^{{\rm gp}}_{Z_{\sigma},x}$. We define its angle set as follows. Let $(x,h)\in Z^{\sigma,\log}_{\tau}$, where $h:\mathcal{M}^{{\rm gp}}_{Z_{\sigma},x}\to S^{1}$ is a group homomorphism. By restricting $h$ to $M(\sigma)$, we obtain an element of ${\rm Hom}(M(\sigma),S^{1})$. This defines the angle map 
    \begin{equation*}
    Z^{\sigma,\log}_{\tau}\to \mathrm{Hom}(M(\sigma),S^{1})
    \end{equation*}    
    and we define its angle set $\Theta^{\sigma,\log}_{\tau}$ to be the image under this map. 
    \end{definition}  
\begin{remark}
As before, there is an action of $\mathbb{T}_{\sigma}(S^{1})$ on $\Theta^{\log}_{\tau}$, and its quotient is $\Theta^{\sigma,\log}_{\tau}$.  
\end{remark}

\begin{proposition}
The angle map defines a homotopy equivalence
\begin{equation*}
Z^{\sigma,\log,\circ}_{\tau}\to \Theta^{\sigma,\log,\circ}_{\tau}.
\end{equation*} 
\end{proposition}
\begin{proof}
This follows from Theorem \ref{thm:KummerHomotopyEquivalence}.
\end{proof}
   
   \begin{remark}
   We do not know at the moment whether the angle map gives a homotopy equivalence
   \begin{equation*}
 Z^{\sigma,\log}_{\tau}\to   \Theta^{\sigma,\log}_{\tau}.
   \end{equation*}
   We will only prove this when the compactification is of $\mathbb{P}^{n}$-type, see Theorem \ref{thm:RelativeAngleFunctor}. This is our reason for using the open angle sets for the upcoming angle functor. 
   \end{remark}

We now define the angle functor. The definition will use the comparison maps on the level of Kato-Nakayama spaces. 
   
   \begin{definition}[Angle functor]\label{def:AngleFunctorGeneral}
   For $(\tau,\sigma)\in \Pscr(\Sigma_{D})$, we set
   \begin{equation*}
   \mathcal{F}_{\Theta} \Big((\tau,\sigma)\Big) 
   := \Theta^{\sigma,\log,\circ}_{\tau}.
   \end{equation*}
   The comparison maps are as follows. Recall that $Z^{\sigma,\log,\circ}_{\tau}\to Z^{\sigma,\log}_{\tau}$ is a homotopy equivalence. For each such pair, we fix a homotopy inverse $Z^{\sigma,\log}_{\tau}\to Z^{\sigma,\log,\circ}_{\tau}$. 
   Similarly, we fix a homotopy inverse $\Theta^{\sigma,\log,\circ}_{\tau}\to Z^{\sigma,\log,\circ}_{\tau}$ of the angle map.   
   
    For any $(\tau_{1},\sigma_{1})\to (\tau_{2},\sigma_{2})$ with map $Z^{\sigma_{1},\log}_{\tau_{1}} \to Z^{\sigma_{2},\log}_{\tau_{2}}$, we now define the map $\Theta^{\sigma_{1},\log,\circ}_{\tau_{1}} \to \Theta^{\sigma_{2},\log,\circ}_{\tau_{2}}$ to be the composite of the following maps:
    
            \begin{equation*}
\begin{tikzcd}
	Z^{\sigma_{1},\log}_{\tau_{1}} & Z^{\sigma_{2},\log}_{\tau_{2}} \\
	Z^{\sigma_{1},\log,\circ}_{\tau_{1}} & Z^{\sigma_{2},\log,\circ}_{\tau_{2}} \\
	\Theta^{\sigma_{1},\log,\circ}_{\tau_{1}} & \Theta^{\sigma_{2},\log,\circ}_{\tau_{2}} 
	\arrow[from=1-1, to=1-2]
	\arrow[from=2-1, to=1-1]
	\arrow[from=1-2, to=2-2]
	\arrow[from=3-1, to=2-1]
    \arrow[from=2-2, to=3-2]
\end{tikzcd}
\end{equation*}
For each pair $(\tau,\sigma)$, we have fixed a set of homotopy inverses, so that, up to homotopy, we can identify $Z^{\sigma,\log}_{\tau}$ and $\Theta^{\sigma,\log,\circ}_{\tau}$. Moreover, the maps between the different $\Theta^{\sigma,\log,\circ}_{\tau}$ are exactly induced from the maps of the $Z^{\sigma,\log}_{\tau}$. By construction, we thus have an induced angle functor $\mathcal{F}_{\Theta}$. This is the angle functor associated to the rational tropical manifold. Strictly speaking, it relies on our  choices of homotopy inverses, but we will omit this dependence. 

      We have a natural transformation 
     \begin{equation*}
     \mathcal{F}_{\rm KN}\to \mathcal{F}_{\Theta}
     \end{equation*} 
     sending $Z^{\sigma,\log}_{\tau}$ to $\Theta^{\sigma,\log,\circ}_{\tau}$. Note again that the commutativity on the level of morphisms follows exactly because we are changing the objects using a fixed homotopy equivalence for every object.        
   \end{definition}
    
    \begin{theorem}\label{thm:MainThmBInText}
    Let $(X, \Sigma, D)$ be a rational tropical manifold with strictly semistable model $(\mathcal{X} \to U,\mathcal{D})$. 
    There is a homotopy equivalence 
    \begin{equation*}
        {\rm hocol}(\mathcal{F}_{\Theta})\to \mathcal{X}_{t}\backslash \mathcal{D}_{t}
    \end{equation*}
    for general $t\in U(\CC)$. 
    \end{theorem}
    \begin{proof}
    Consider the map of diagrams $\mathcal{F}_{\rm KN}\to \mathcal{F}_{\Theta}$. 
    Since all of the maps on objects are homotopy equivalences by Theorem \ref{thm:KummerHomotopyEquivalence}, we find a homotopy equivalence of homotopy colimits using Proposition \ref{prop:HomotopyEquivalentFunctors}. But the homotopy colimit of $\mathcal{F}_{\rm KN}$ is homotopy equivalent to $\mathcal{X}^{\log,\zeta}_{\mathcal{D}, s}$ by Theorem \ref{thm:KNFunctorHomotopyColimit3}, so we find the desired statement using Theorem \ref{thm:UsuisTheorem2}.  
    \end{proof}
    
    \begin{remark}
    This theorem gives us our first polar decomposition for general complex varieties. Here we view the tropicalizations as containing the radial parts. In the homotopy colimit, we then take products of the radial and angular parts, so that we locally obtain a product decomposition. This is globally not true however. We will see many examples of this.  
    \end{remark}

\begin{remark}
Since the homotopy inverses used in the definition of the angle functor are not directly explicit, this theorem does not directly allow us to find $f$-vectors for instance. To that end, we will study more explicit angle functors in the next section for complete intersections. For those, the angle map $Z^{\sigma,\log}_{\tau}\to \Theta^{\sigma,\log}_{\tau}$ will be a homotopy equivalence, so that we can bypass the homotopy inverses.

In practice, the maps of angle sets can simply be obtained by perturbing the angle sets so that they become subsets of the larger ambient angle set. This makes the map of angle sets a cofibration, but note that the latching maps are in general not cofibrations.   
\end{remark}

\begin{remark}
Suppose that $(X,\Sigma,D)$ is a curve. When the associated generic fiber is hyperbolic, this decomposition we give here is known as a pair-of-pants decomposition. See Section \ref{sec:ClassicalExamplesD01} for explicit examples.  When the curve is elliptic, then it is the torus decomposition, and when we have a projective or affine line, it is a decomposition of a toric variety. Our polar decompositions can thus be seen as interpolating between these different notions.
\end{remark}

\section{Complete intersections}
\label{sec:CompleteIntersections}

In this section, we show how to make the angle functor explicit for complete intersections in tori, as well as their associated smooth toric compactifications. This can be seen as generalizing Bernstein's theorem to higher-dimensional varieties. To describe this angle functor, we show that the angle map for the Kato-Nakayama spaces associated to complete intersections is a homotopy equivalence. The angle functor can then be constructed without relying on tubular neighborhoods, as in Definition \ref{def:AngleFunctorGeneral}. 

We then review the description of the angle set of a complete intersection given in \cite{NS13NonArchimedean} from the point of view of Coxeter arrangements. This gives a triangulation of the angle set, and we calculate its $f$-vector. This will be used in the upcoming sections to describe the $f$-vector of a general complete intersection.      
 
\subsection{Rational complete intersection manifolds}

We first define our notion of a complete intersection. 
\begin{definition}\label{def:RationalCIManifold}
Let $X = V(I) \subset \Spec(K[M])$ be a subvariety of the torus over $K = \CC(t)$ and let $r = {\rm codim}(X)$ be its codimension.
We say that $X$ (or $\trop(X)$) is a \emph{rationally smooth complete intersection} if there exist $r$ generators $f_{i} \in I$ such that the tropical hypersurfaces $\trop(V(f_{i}))$ are rationally smooth, and intersect transversally everywhere. 
We say that $X$ or $\trop(X)$ is a \emph{smooth complete intersection} if $I$ can be generated by $r$ elements $f_{i}$ whose tropicalizations $\trop(V(f_{i}))$ are smooth, intersect transversally and properly, and their stable intersection is tropically smooth.

Let $(X,\Sigma,D)$ be a rational tropical manifold. We say that $(X,\Sigma,D)$ is a rational complete intersection manifold if $X$ is a rationally smooth complete intersection.
Suppose that $(X,\Sigma,D)$ is a rational tropical manifold such that $I$ is generated by a single polynomial $f$. We then call  
$X=V(f)$ or $(X,\Sigma,D)$ a rationally smooth hypersurface. We similarly define the notion of a smooth complete intersection $(X,\Sigma,D)$.  
\end{definition}

\begin{remark}
By Proposition \ref{pro:RationalManifoldStableIntersection}, we have that the stable intersection of rational tropical manifolds is again a rational tropical manifold. The same is of course not true for tropical manifolds, hence the added condition in that case.
\end{remark}

\begin{remark}\label{rem:CILocalComponents}
Let $X$ be a rational complete intersection manifold with generators $f_{i}$. By the argument in Lemma \ref{lem:KummerLemma}, the local linear space is given by a stable intersection of hypersurfaces of the form 
\begin{equation*}
g_{j}=\sum_{i=0}^{m}c_{i,j}x_{i}.
\end{equation*}  
As in the proof of that lemma, we can assume that the $g_{j}$ have a nonzero constant coefficient. 
The nonconstant monomial supports are then disjoint by transversality. 
We call these  \emph{standard rational complete intersection spaces}. 

Note that the corresponding angle set is simply the product of the angle sets. 
In particular, we see that it is sufficient to determine the angle set of a single linear form. This will be done in Section \ref{sec:AngleSetsCI}. 
\end{remark}

\begin{definition}\label{def:UniversalBernsteinFamily}
Suppose that $S_{i}\subset \ZZ^n$ are finite subsets with convex hulls $\mathcal{C}_{i}$, for $i=1,...,r$. Consider the family of all tuples of Laurent polynomials $f_{i}\in K[M]$ whose Newton polytopes are $\mathcal{C}_{i}$. We call this the universal Bernstein family for $r$ polynomials with fixed polynomial supports $S_{i}$. If we write $A=\CC[\mathbb{N}^{m}]$ for $m=\sum |S_{i}|$, then the universal Bernstein family can be viewed as a morphism 
\begin{equation*}
    \mathcal{X}\to U,
\end{equation*}
where $U$ is an open subset of $\Spec(A)$. Explicitly, we have $\mathcal{X}=V(f_{1},...,f_{r})$, where $f_{i}\in A[M]$. Moreover, any point in $U(K)$ corresponds to a set of coefficients for the $f_{i}$  such that the induced Newton polytopes are all $\mathcal{C}_{i}$.
\end{definition}

\begin{corollary}
The universal Bernstein family $\mathcal{X}\to U$ is generically smooth. The homotopy type of a general fiber $\mathcal{X}_{t}\backslash \mathcal{D}_{t}$ is constant and equal to   ${\rm hocol}(\mathcal{F}_{\Theta})$ for a choice of a rational tropical manifold in the family.   
\end{corollary}
\begin{proof}
Generic smoothness follows from the generic smoothness theorem, after
restricting $U$ if necessary. The statement about the homotopy type of
a general fiber then follows from Theorem \ref{thm:MainThmBInText}.
\end{proof}

We now specialize to a very specific family. Let $\Delta_n$ denote the standard simplex of dimension $n$ in $\RR^n$ defined as follows 
\[
    \Delta_n := \conv(0, e_1, ..., e_n),
\] 
where $e_1,\dots, e_n$ are the standard basis vectors. We set
\[
    S_{i} := d_i \Delta_n \cap \ZZ^n, \quad \text{for any } 1 \leq i \leq r.
\] 
In other words, these are the integer points of the simplices $d_{i} \cdot \Delta_{n}\subset \RR^{n}$. 
Note that any collection of generic polynomials $f_i$ of degree $d_i$ have monomial supports $S_i$. 

We will show that for some choice of coefficients in $U$, the family $\mathcal{X} = V(f_{1},...,f_{r})\to U$ from Definition \ref{def:UniversalBernsteinFamily} admits a smooth tropicalization, so that we describe the topology of these spaces in terms of tropical manifolds. 
In this case, we do not need the more general notion of a rational tropical manifold, although this does lead to a more efficient CW complex structure, see Section \ref{sec:Examples}. 

For any $1 \leq r \leq n$ and positive integers $d_1 , \dots, d_r \geq 1$, the \emph{Cayley polytope} of the lattice polytopes $d_i \Delta_n$ is the lattice polytope in $\RR^{n + r}$ defined as follows:
\[
    \mathscr{C}(d_1 \Delta_n, \dots, d_r \Delta_n) := \conv\left( \bigcup_{i = 1}^{r} (d_i \Delta_n) \times e_{n+i} \right) \subset \RR^{n + r}.
\]

We recall the following from \cite[Section 4.6]{MS15}. Let $f_1, \dots, f_r \in K[x_1^{\pm}, \dots, x_n^{\pm}]$ be polynomials such that ${\rm Newt}(f_i) = d_i \Delta_n$, in particular we have $\deg(f_i) = d_i$. Each polynomial $f_i$ induces a height function on the lattice points of $d_i \Delta_n$. These heights induce a height function on the lattice points of $\mathscr{C}(d_1 \Delta_n, \dots, d_r \Delta_n)$ and hence a regular polyhedral subdivision of $\mathscr{C}(d_1 \Delta_n, \dots, d_r \Delta_n)$. When this subdivision is a unimodular triangulation then the tropical complete intersection $\trop(V(f_1)) \cap_{\rm st} \dots \cap_{\rm st} \trop(V(f_r))$ is smooth. See \cite[Page 212]{MS15} in particular.

On the other hand, we note that for each $1 \leq i \leq r$ the polytope $d_i \Delta_n$ can be identified with the slice
\begin{equation}\label{eq:SliceCayleyPolytope1}
    \mathscr{C}(d_1 \Delta_n, \dots, d_r \Delta_n)  \cap \{ x_{n+i} = 1 \}.
\end{equation}
So, any regular unimodular triangulation of $\mathscr{C}(d_1 \Delta_n, \dots, d_r \Delta_n)$ induces regular unimodular triangulations of the polytopes $d_i \Delta_n$ by slicing $\mathscr{C}(d_1 \Delta_n, \dots, d_r \Delta_n)$ as in \eqref{eq:SliceCayleyPolytope1}. Moreover the corresponding tropical complete intersection is smooth.

\begin{lemma}\label{lem:multiplicityOfStableIntersection}
    Let $f_1, \dots, f_r \in K[x_1^{\pm}, \dots, x_n^{\pm}]$ be polynomials of degrees $d_1, \dots, d_r$ and suppose that the induced regular subdivision of the Cayley polytope is a unimodular triangulation. Then the following hold:
    \begin{enumerate}
        \item the hypersurfaces $Z_{i}:=\trop(V(f_{i}))$ are smooth,
        \item the hypersurfaces $Z_{i}$ intersect transversely,
        \item the local multiplicities of $Z=\bigcap\limits_{i=1}^{r} Z_{i}$ are all equal to $1$.
    \end{enumerate}
\end{lemma}    
\begin{proof}
    See \cite[Proposition 7.4]{BB13}.
\end{proof}

\begin{corollary}\label{cor:TropicallySmoothIntersection}
    Let $f_{i}$ be as in Lemma \ref{lem:multiplicityOfStableIntersection} with vanishing set $X=\bigcap_{i=1}^{r}V(f_{i})$. Then there is a $\Sigma$ such that $(X,\Sigma,D)$ is a tropical manifold for any $D$.
\end{corollary}
\begin{proof}
As in the proof of Proposition \ref{pro:RationalManifoldStableIntersection}, we can find local coordinates so that the monomial supports of the hypersurfaces have at most $1$ in common. Since the multiplicity is $1$, we find that the monomials of the hypersurfaces generate a basis of $M$, so that we indeed obtain a tropical manifold. 
\end{proof}

Therefore, to find a smooth tropical complete intersection of hypersurfaces with degrees $d_1, \dots, d_r$ it is enough to find a regular unimodular triangulation of $\mathscr{C}(d_1 \Delta_n, \dots, d_r \Delta_n)$. However, lattice polytopes are not always guaranteed to have regular unimodular triangulations, see \cite{OhsugiHibi}. Fortunately, in this case we do have a positive result.

\begin{proposition}\label{prop:SmoothStableIntersection}
    The Cayley polytope $\mathscr{C}(d_1 \Delta_n, \dots, d_r \Delta_n)$ admits a regular unimodular triangulation. In particular, for any integers $d_1,\dots, d_r \geq 1$, there exist polynomials $f_1, \dots, f_r \in K[x_1^{\pm}, \dots, x_n^{\pm}]$ of respective degrees $d_1, \dots, d_r$ such that the stable intersection of the tropical hypersurfaces $\trop(V(f_i)) \subset \RR^n$ is a tropical manifold. 
\end{proposition}
\begin{proof}
    The existence of a unimodular triangulation of $\mathscr{C}(d_1 \Delta_n, \dots, d_r \Delta_n)$ is a particular case of \cite[Lemma 2.10]{AdiprasitoLiuTemkin} and \cite[Remark 3.5]{BackmanLiu}. The existence of the polynomials $f_i$ such that the stable tropical intersection
    \[
        \trop(V(f_1)) \cap_{\rm st} \dots  \cap_{\rm st}  \trop(V(f_r))
    \]
    is smooth follows from Corollary \ref{cor:TropicallySmoothIntersection}. 
\end{proof}

\subsection{Special compactifying varieties}

In this section, we specialize to Kummer coverings of very affine linear spaces that admit $\mathbb{P}^{n}$ as an snc-compactifying variety. For instance, the local spaces arising from rational complete intersection manifolds are of this type. We will show that the associated angle map on the level of Kato-Nakayama spaces is a homotopy equivalence.  
We note here that for this result, we do not need the affine result Theorem \ref{thm:MainTheoremAText}, as we can use the classical Smale-Vietoris theorem. We will however need Theorem \ref{thm:MainTheoremAText} in Proposition \ref{pro:ChangingToricCompactification} to conclude that the extended angle map is again a homotopy equivalence for any other compactifying variety.

We start by defining the corresponding very affine linear spaces. 
\begin{definition}\label{def:PnSNCCompactification}
Let $X = V(I)\subset \Spec(\CC[M])=Y$ be a very affine linear space. We say that $X$ is of $\mathbb{P}^{n}$ type if there is an embedding $Y\to \mathbb{P}^{n}$ such that the closure $\overline{X}$ of $X$ in $\mathbb{P}^{n}$ is an snc-compactification. That is, $\overline{X}$ is smooth and the boundary $\partial{X}=\overline{X}\backslash X$ is an snc divisor. We write $X^{\log}$ for the corresponding Kato-Nakayama space. 

Suppose that $X$ is a Kummer covering of a very affine linear space. We say that $X$ is of $\mathbb{P}^{n}$ type if there is a finite Kummer covering $X\to Z$ to a very affine linear space $Z$ of $\mathbb{P}^{n}$-type. We write $\overline{X}$ for the corresponding induced compactification in this case. We write $X^{\log}$ for the corresponding Kato-Nakayama space.        
\end{definition}

\begin{remark}
    We will focus on the very affine linear case in this section.  That is, we will prove that the angle map $X^{\log}\to \Theta^{\log}$ is a homotopy equivalence. The Kummer case easily follows from this.
\end{remark}

\begin{example}
    Suppose that $\overline{X} = V(x_{0} + \dots + x_n)$.
    Then $\mathbb{P}^{n}$ is an snc-compactification for $\overline{X} \cap (\CC^\times)^{n+1} / \CC^\times$. 
    Indeed, all the boundary loci of $\mathbb{P}^{n}$  intersect $\overline{X}$ transversely. For instance, if we take $x_{0}=0$, then we obtain the homogeneous vanishing set of $x_{1}+...+x_{n}$, which is smooth of the expected dimension.
\end{example}

\begin{example}\label{exa:CounterExample}
    Suppose that $I=(w_{1}-(x-y),w_{2}-(x-z),w_{3}-(y-z))$. Then $\mathbb{P}^{5}$ is not an snc-compactification for $X=V(I)$. Indeed, the intersection $\overline{X}\cap V(w_{1})\cap V(w_{2})\cap V(w_{3})$ is of codimension $2$ in $\overline{X}$, rather than $3$.
\end{example}

We now show that the angle map is a homotopy equivalence on the level of Kato-Nakayama~spaces. 

\begin{theorem}\label{thm:HEforKNSpaces}
    Let $X \subset (\CC^{\times})^n$ be a linear space of $\PP^n$-type. 
    Then, the induced map $X^{\log}\to \Theta^{\log}$ is a homotopy equivalence. 
\end{theorem}
\begin{proof}
We begin by defining convenient coordinate charts. 
Let $\CC[X_{0},...,X_{n}]$ be the homogeneous coordinate ring of $\mathbb{P}^{n}$, and let $U_{i} := D_{+}(X_{i})\simeq \CC^{n}$ be the standard open cover of $\PP^n$. We use the following notation
\begin{equation*}
    B_{j}:= \Big\{[x_0: \dots : x_n] \in U_{j}:|x_i/x_j| \leq 1 \text{ for all } 0 \leq i \leq n \Big\}.
\end{equation*}
We note that the sets $B_0, \dotsm B_n$ cover all of $\PP^n$. Similarly, for any $\varnothing \neq I \subset [n]$, we define 
\[
    B_{I} :=\bigcap_{j\in{I}}B_{j}\subset \bigcap_{j\in{I}}U_{j} 
    \quad \text{and} \quad 
    X^{\log}_{I} := \{(P,\psi)\in{X^{\log}}:P\in{B_{I}}\}.
\]
This is a compact subset of $X^{\log}$. We write $\Theta^{\log}_{I}$ for the induced angle set, which is a compact subset of $\Theta^{\log}$. We claim that the restricted angle map 
\[
    \ang_{I} \colon X^{\log}_{I}\to \Theta^{\log}_{I}
\]
is a homotopy equivalence for any nonempty $I \subset [n]$. Indeed, since $\overline{X}$ is the compactification of a linear space in $(\CC^{\times})^n$, we can write $\overline{X}=V(F_{1},...,F_{\ell})$ for homogeneous linear polynomials $F_{j}$. We fix $i\in {I}$. On the open $D_{+}(X_{i})$, we then have 
\begin{equation*}
    f_{j,i}:=F_{j}/X_{i}=c_{0,j}+\sum_{k} c_{k,j}x_{k},
\end{equation*}
where $x_{j}=X_{j}/X_{i}$. We now write $x_{j}=r_{j}s_{j}$, where $r_{j}=|x_{j}|$ and $s_{j}\in S^{1}$. Note that the $r_{j}$ and $s_{j}$ give coordinates on the corresponding Kato-Nakayama space.
We fix the angles $s_j$. For each $1\leq j\leq \ell$, the equation
$f_{j,i}=0$ gives two real linear conditions on the radii, obtained by
taking its real and imaginary parts.

Note that the conditions defining $B_I$ are given by
\[  
    0\leq r_k\leq 1 \quad \text{ for all } 1 \leq k \leq n 
    \quad  \text{ and } \quad 
    r_j=1 \quad  \text{for } j\in I\backslash\{i\}.
\]
Hence the fibers of $\ang_I$ are polyhedra. In particular, the fibers are contractible. Since $X^{\log}_{I}\to \Theta^{\log}_{I}$ is a proper map of compact semi-algebraic sets (note that $X^{\log}_{I}$ is a closed subset of a compact space and $\Theta^{\log}_{I}$ is the image of a compact set), we find that the angle map $\ang_{I} \colon X^{\log}_{I}\to \Theta^{\log}_{I}$ is a homotopy equivalence by Theorem \ref{thm:SmaleVietoris}.

The sets $X^{\log}_{I}$ and $\Theta^{\log}_{I}$ give a finite compact covering of $X^{\log}$ and $\Theta^{\log}$ respectively. The angle map then gives a natural transformation of diagrams of spaces. Note that these diagrams of spaces are closed \v{C}ech covers. Since the natural transformation is a homotopy equivalence on objects, by Proposition \ref{prop:HomotopyEquivalentFunctors}, we find a homotopy equivalence of homotopy colimits. By Theorem~\ref{thm:ProjectionLemma} (here we use that the two diagrams are closed \v{C}ech covers), these are equal to their colimits, which are $X^{\log}$ and $\Theta^{\log}$ respectively. This gives the desired result.
\end{proof}

\begin{example}
    Consider the cone $\sigma\subset \RR^{3}$ given by $v(z)=0$, $v(x)\geq{0}$, $v(y)\geq{0}$. This belongs to the tropicalization of $V(1+x+y+z)$. We now subdivide this cone to obtain a different compactification. For instance, we can subdivide this into four cones, obtained by adding $w_{0}=(1,2)$, $w_{1}=(1,1)$ and $w_{2}=(2,1)$. Note that the induced cones are all smooth, since the determinants are $\pm{1}$. We calculate the induced equation over the cone spanned by $(1,0)$ and $(1,2)$. We find that $v(y)\geq{v(x^{2})}$, so the coordinate ring is given by 
    \begin{equation*}
        A=\CC[x,y_{1},z^{\pm}]=\CC[x,y/x^2,z^{\pm}].
    \end{equation*}
    We have $y=y_{1}\cdot x^2$. In particular, the new equation becomes 
\begin{align*}
    f=1+x+y_{1}x^2+z.
\end{align*}
For Kato-Nakayama spaces, we thus find 
\begin{equation*}
f=1+r_{x}s_{x}+s_{y_{1}}s^{2}_{x}r_{y_{1}}r_{x}^{2}+r_{z}s_{z}.
\end{equation*}
For fixed angles, we find that the equations for the radii are non-linear. In particular, we cannot use the argument in the proof of Theorem \ref{thm:HEforKNSpaces}. We will see in Proposition \ref{pro:ChangingToricCompactification} that we can still deduce that the angle map is a homotopy equivalence for these toric compactifications. 
\end{example}

We now extend 
Theorem \ref{thm:HEforKNSpaces} to show that the angle map is a homotopy equivalence for arbitrary toric compactifications. Note that it already follows from Theorem \ref{thm:HEforKNSpaces} that the two spaces involved are homotopy equivalent, but it does not directly follow that the angle map is a homotopy equivalence.  
\begin{proposition}\label{pro:ChangingToricCompactification}
    Let $X$ be as in Definition \ref{def:PnSNCCompactification}, with $\Delta$ the fan of $\mathbb{P}^{n}$. Suppose that $\Delta'$ is a refinement of $\Delta$, and let $\overline{X}_{\Delta'}$ be the corresponding compactification with Kato-Nakayama space $X^{\log}_{\Delta'}$. 
    Then the compactified angle map $X^{\log}_{\Delta'}\to \Theta^{\log}$ is a homotopy equivalence.  
\end{proposition}
\begin{proof}
    We first consider the compactification $\overline{X}$ coming from $\mathbb{P}^{n}$. Consider the diagram 
    \begin{equation*}
        \begin{tikzcd}
        	X & {X}^{\log}\\
        	\Theta& \Theta^{\log}
        	\arrow[from=1-1, to=1-2]
        	\arrow[from=1-1, to=2-1]
        	\arrow[from=1-2, to=2-2]
        	\arrow[from=2-1, to=2-2].
        \end{tikzcd}
    \end{equation*}
    Note that $\Theta\to \Theta^{\log}$ is a homotopy equivalence. Indeed, $3$ of the $4$ maps in the diagram above are already homotopy equivalences by Theorem \ref{thm:HEforKNSpaces}, \ref{thm:MainTheoremAText} and the fact that $X^{\log}$ is a manifold with boundary, with interior $X$. 

    We now change our compactification to the one in the main theorem. Note that the angle set $\Theta^{\log}$ does not change. We now again obtain a diagram as above. Note that $X\to X^{\log}_{\Delta'}$ is a homotopy equivalence, since $X^{\log}_{\Delta'}$ is a manifold with boundary and $X$ is the interior. 
    By Theorem \ref{thm:MainTheoremAText}, we have that $X\to \Theta$ is a homotopy equivalence, so that we conclude using Whitehead's theorem for instance that $X^{\log}_{\Delta'}\to \Theta^{\log}$ is a homotopy equivalence.
    \end{proof}

We now show how the angle sets of the initial degenerations fit into $\Theta^{\log}$ in a compatible way.

\begin{theorem}\label{thm:RelativeAngleFunctor}
    Let $X\subset (\CC^{\times})^{n}$ be a very affine linear space and suppose that it admits $\mathbb{P}^{n}$ as an snc-compactifying variety. Then all arrows in the following diagram are homotopy equivalences:
              \begin{equation*}
    \begin{tikzcd}
    	X & {X}^{\log}\\
    	\Theta& \Theta^{\log}
    	\arrow[from=1-1, to=1-2]
    	\arrow[from=1-1, to=2-1]
    	\arrow[from=1-2, to=2-2]
    	\arrow[from=2-1, to=2-2].
    \end{tikzcd}
    \end{equation*}
    Moreover, the same is true for any $Z^{\log}_{\tau}\subset X^{\log}$, and we have a commutative diagram
                  \begin{equation*}
    \begin{tikzcd}
    	Z^{\log}_{\tau_{1}} & {Z}^{\log}_{\tau_{2}} & X^{\log}\\
    	\Theta^{\log}_{\tau_{1}}& \Theta^{\log}_{\tau_{2}} & \Theta^{\log}
    	\arrow[from=1-1, to=1-2]
        \arrow[from=1-2, to=1-3]
    	\arrow[from=1-1, to=2-1]
    	\arrow[from=1-2, to=2-2]
        \arrow[from=1-3, to=2-3]
        \arrow[from=2-2, to=2-3]
    	\arrow[from=2-1, to=2-2]
    \end{tikzcd}
    \end{equation*}
    for $\tau_{1}\supset \tau_{2}$.
\end{theorem}

\begin{proof}
    The first part is Theorem \ref{thm:HEforKNSpaces}. For the second part, note that all the initial ideals $\mathrm{in}_{w}(I)$ of the ideal $I$ of $X$ satisfy the same snc-condition. 
    Indeed, on the homogeneous side, the relevant initial degenerations are obtained by passing to suitable coordinate strata determined by the weight vector. Since the boundary divisors meet $\overline{X}$ with the expected codimensions, the same remains true after passing to these strata and to further intersections.
    But this yields the correct codimension by assumption, so we find that it also yields the correct codimension for further intersections, as desired. In particular, we find that the $Z^{\log}_{\tau} \to  \Theta^{\log}_{\tau}$ are all homotopy equivalences. 
    Moreover, it is clear from the construction of the angle map that this yields a commutative diagram of angle sets. This finishes the proof.        
\end{proof}

\begin{remark}
     It is unclear to the authors whether the same holds for snc-compactifications of arbitrary hyperplane complements. Intuitively, the angles of infinite directions can collide as in  Figure \ref{fig:AngleSet3PointsinP1}, as the angles for the three points at infinity all intersect, and this could potentially cause the topology of $\Theta^{\log}$ to be different from that of $\Theta$. In our examples, the topology doesn't change however, as we have suitably nice compactifying varieties.    
\end{remark}

\begin{remark} 
For a Kummer covering of a very affine linear space here, we can apply the technique from Theorem \ref{thm:KummerHomotopyEquivalence} and obtain
the appropriate results in a very similar way. We leave the details to the reader.   
\end{remark}

\subsection{The angle functor for complete intersections}\label{sec:ChangingCompactifications}\label{sec:GluingForCompleteIntersections}

We now combine these results and define an angle functor that captures the topology of a nearby fiber for a complete intersection. 

\begin{definition}[Angle functor for complete intersection manifolds]\label{def:AngleFunctorCI}
    Let $X=V(I)$ be a rational complete intersection manifold. We define the associated angle functor $\mathcal{F}_{\Theta, \mathrm{CI}}$ as follows:
        \begin{equation*}
        \mathcal{F}_{\Theta,\mathrm{CI}}(\tau,\sigma) := \Theta^{\sigma,\log}_{\tau}.
        \end{equation*}
    The induced maps are either inclusions of angle sets, or quotients of angle sets as in Lemma \ref{lem:InitialDegenerationsVSKNSpaces} and Remark \ref{rem:InitialDegenerationsDirectQuotientKNSpaces}. We have an induced natural transformation 
        \begin{equation*}
        \mathcal{F}_{\rm KN}\to \mathcal{F}_{\Theta,{\rm CI}},
        \end{equation*}
    sending $Z^{\sigma,\log}_{\tau}$ to $\Theta^{\sigma,\log}_{\tau}$.
    
    We note that, unlike the general angle functor $\mathcal{F}_{\Theta}$ defined in Definition \ref{def:AngleFunctorGeneral}, whose comparison maps are defined using chosen homotopy inverses and whose objects are the open angle spaces $\Theta^{\sigma,\log,\circ}_{\tau}$, the complete-intersection angle functor uses the closed angle spaces $\Theta^{\sigma,\log}_{\tau}$ and has canonical comparison maps given by inclusions and quotients.
\end{definition}

\begin{corollary}\label{cor:AngleFunctorCI}
Let $(X,\Sigma,D)$ be a rational complete intersection manifold with strictly semistable model $(\mathcal{X} \to U, \mathcal{D})$. 
There is a homotopy equivalence 
\begin{equation*}
{\rm hocol}(\mathcal{F}_{\Theta,{\rm CI}})\to \mathcal{X}_{t}\backslash \mathcal{D}_{t}
\end{equation*}
   for general $t\in U(\CC)$. 
\end{corollary}
\begin{proof}
This follows from Theorem \ref{thm:RelativeAngleFunctor}, together with its Kummer analogue described above, applied objectwise to the Kato--Nakayama functor.
\end{proof}

\begin{remark}
The main difference between Theorem \ref{thm:MainTheoremAText} and Corollary \ref{cor:AngleFunctorCI} is that we have an explicit triangulation of each $\Theta^{\sigma,\log}_{\tau}$, and we have canonical gluing maps. This will allow us to convert the homotopy colimit to an explicit CW complex structure for a general nearby fiber. Moreover, the $f$-vector formula for the local parts will allow us to find a global formula for the $f$-vector of the homotopy colimit. These formulas  will be given in Section \ref{sec:CWStructure}.
\end{remark}

\begin{remark}\label{rem:KummerSmoothRemark2}
 From a computational point of view, it seems one can use Corollary \ref{cor:AngleFunctorCI} to define a topological version of  
 homotopy continuation, see \cite{HuberSturmfels1995} for the zero-dimensional case, and \cite{Jensen16,DR24} for tropical versions. Namely, one can introduce a point to track for each corner point of the angle cellulation. By tracking these corner points, one can then induce CW complex structures on general CI algebraic varieties. 
\end{remark}

\subsection{Angle sets of complete intersections}\label{sec:AngleSetsCI}

In this section, we recall the description of the angle sets of linear hypersurfaces from \cite{NS13NonArchimedean}. 
These angle sets are obtained from affine Coxeter arrangements by deleting certain open zonotopes. 

We use the notation $[n] := \{0,\ldots,n\}$. For every nonempty subset $I \subset [n]$, let
\[
    V_I:=\RR^I/\RR\mathds{1}_I,
    \qquad
    L_I:=\ZZ^I/\ZZ\mathds{1}_I,
\]
where $\mathds{1}_I$ is the all-ones vector. We have a covering map
\[
    V_I\longrightarrow (S^1)^I/S^1,
    \qquad
    [u_i]_{i\in I}\longmapsto [e^{i\pi u_i}]_{i\in I},
\]
whose kernel is $2L_I$. When $I=[n]$, we simply write $V=V_{[n]}$ and
$L=L_{[n]}$.

\begin{definition}\label{def:CoxeterArrangement}
For every nonempty $I \subset [n]$, let $\mathscr{A}_I$ be the affine
Coxeter arrangement in $V_I$ consisting of the hyperplanes
\[
    H_{ij}^{(\ell)}
    =\{x\in V_I:x_i-x_j=\ell\},
    \qquad i,j\in I,\quad i\ne j,\quad \ell\in\ZZ.
\]
We denote by $\mathcal{A}_I$ the induced simplicial complex.
For $\varnothing\ne S \subset  I$, define
\[
    Z_{I,S}
    :=\{x\in V_I:|x_i-x_j|<1\text{ for all }i,j\in S\},
    \qquad
    U_{I,S}:=\bigcup_{\gamma\in2L_I}(\gamma+Z_{I,S}).
\]
We then define the deleted Coxeter complex
\[
    \mathcal{B}_{I,S}
    :=\{\tau\in\mathcal{A}_I:
         \relint(\tau)\not \subset  U_{I,S}\}.
\]
We also use the conventions $Z_{I,\varnothing}=U_{I,\varnothing}=\varnothing$ and $\mathcal{B}_{I,\varnothing}=\mathcal{A}_I$.

Since the boundaries of the sets $Z_{I,S}$ lie in the hyperplanes of
$\mathscr{A}_I$, each cell of $\mathcal{A}_I$ has relative interior either
contained in or disjoint from $U_{I,S}$. In particular,
$\mathcal{B}_{I,S}$ is a subcomplex of $\mathcal{A}_I$.
When $S=I$, we write $Z_I:=Z_{I,I}$ and
$\Sigma_I:=\mathcal{B}_{I,I}$. We also use the abbreviations
\[
\mathcal{A}_n:=\mathcal{A}_{[n]}, \quad 
\mathcal{B}_{n,S}:=\mathcal{B}_{[n],S}, \quad
Z_S:=Z_{[n],S} \quad \text{and} \quad \Sigma_n := \Sigma_{[n]}.
\]
\end{definition}

For any pair of sets $\varnothing\ne I \subset  J \subset [n]$, restriction of coordinates induces a map
\[
    p_{J,I}:V_J\longrightarrow V_I,
    \qquad
    [x_j]_{j\in J}\longmapsto[x_i]_{i\in I}.
\]
and we have $p_{J,I}(2L_J) = 2L_I$ and therefore descend to the
corresponding angle tori. We shall see in Proposition~\ref{prop:CompatibilityQuotientTriangulation} that, for
$S \subset  I$, they also induce simplicial maps
$\mathcal{B}_{J,S}\to\mathcal{B}_{I,S}$.

Now consider the hyperplane
\[
    X=V(x_0+\cdots+x_n)\subset (\CC^{\times})^{n+1}/ \CC^\times.
\]
Its tropicalization $\Delta$ is the fan with cones
\[
    \sigma_I = \cone(\overline e_i : i\in I),
    \qquad I  \subset [n],\quad |I| \le n-1,
\]
where $\overline e_i$ is the image of the $i$-th standard basis vector in
$V$. For $w\in\relint(\sigma_I)$, we have
\[
    {\rm in}_w(X) = V\left(\sum_{i\in[n]\setminus I} x_i\right) \subset (\CC^{\times})^{n+1}/ \CC^\times.
\]
Thus the complexes $\mathcal{B}_{[n],S}$ naturally describe the angle set of $X$ as well as the angle sets of its degenerations, as the following proposition explains.

\begin{figure}[ht]
        \centering    
        
        \begin{minipage}{.35\textwidth}
        \scalebox{0.8}{
            \begin{tikzpicture}

                \coordinate (a) at (0, 0);
                \coordinate (a1) at (3, 0);
                \coordinate (a2) at (0, 3);
                \coordinate (a3) at (-2,-2);
                
                \coordinate (e1) at (1, 0);
                \coordinate (e2) at (0, 1);
                \coordinate (e3) at (-1,-1);

                \draw (1,-0.5) node   {$\bar{e}_1$};
                \draw (0.5,1)  node   {$\bar{e}_2$};
                \draw (-1.5,-1) node  {$\bar{e}_3$};

                \draw[] (a) -- (a1);
                \draw[] (a) -- (a2);
                \draw[] (a) -- (a3);
                \draw[->, color = green, line width = 1.5] (a) -- (e1);
                \draw[->, color = blue, line width = 1.5] (a) -- (e2);
                \draw[->, color = red, line width = 1.5] (a) -- (e3);
                \draw [fill=black] (a) circle (2pt);

                \draw (0.5,-2) node {$\trop(X)$};
                
            \end{tikzpicture}}
        \end{minipage}
       \begin{minipage}{.55\textwidth}
            \begin{tikzpicture}[scale=0.5, line width = 0.5, clip=true]
                \def\mincoord{-5}
                \def\maxcoord{5}

                \foreach \i in {-2,...,2} {
                    \foreach \j in {-2,...,2} {
                        \pgfmathsetmacro\xshift{2*\i}
                        \pgfmathsetmacro\yshift{2*\j}
            
                        \filldraw[fill=orange!70, draw=black] 
                            ({-1 + \xshift}, {0 + \yshift}) --
                            ({0 + \xshift}, {1 + \yshift}) --
                            ({-1 + \xshift}, {1 + \yshift}) -- cycle;
            
                        \filldraw[fill=orange!70, draw=black] 
                            ({0 + \xshift}, {-1 + \yshift}) --
                            ({1 + \xshift}, {-1 + \yshift}) --
                            ({1 + \xshift}, {0 + \yshift}) -- cycle;
                    }
                }

                \foreach \m in {-2,...,2} {
                    \draw[black, line width = 0.2] (2*\m, \mincoord) -- (2*\m, \maxcoord);
                }

                \foreach \m in {-3,...,2} {
                    \draw[blue, line width = 2] (2*\m + 1, \mincoord) -- (2*\m + 1, \maxcoord);
                }

                \foreach \m in {-2,...,2} {
                    \draw[black, line width = 0.2] (\mincoord, 2*\m) -- (\maxcoord, 2*\m);
                }
                \foreach \m in {-3,...,2} {
                    \draw[green, line width = 2] (\mincoord, 2*\m + 1) -- (\maxcoord, 2*\m + 1);
                }

                \foreach \m in {0,...,5} {
                    \draw[black, line width = 0.2] 
                        (\maxcoord - 2*\m, \mincoord) -- 
                        (\maxcoord, \mincoord + 2*\m);
                }
                \foreach \m in {0,...,5} {
                    \draw[black, line width = 0.2] 
                        (\mincoord, \mincoord + 2*\m) -- 
                        (\maxcoord - 2*\m, \maxcoord);
                }

                \foreach \m in {0,...,4} {
                    \draw[red, line width = 2] 
                        (\mincoord, \mincoord + 2*\m + 1) -- 
                        (\maxcoord - 2*\m - 1, \maxcoord);
                }
                \foreach \m in {0,...,4} {
                    \draw[red, line width = 2] 
                        (\maxcoord - 2*\m-1, \mincoord) -- 
                        (\maxcoord, \mincoord + 2*\m+1);
                }

                \draw[->, black, line width = 1] (-5, 0) -- (5.5, 0) node[right] {$x$};
                \draw[->, black, line width = 1] (0, -5) -- (0, 5.5) node[above] {$y$};

                \draw[black, line width = 3] (1,1) -- (-1,1) -- (-1,-1) -- (1,-1) -- cycle;

                \draw (0.5,-6) node {$\Sigma_{2} := \mathcal{B}_{[2],[2]} \subset \mathcal{A}_{[2]}$};

                 \foreach \m in {-5,...,5} {
                        \foreach \n in {-5,...,5} {
                                \draw [fill=black] (\m, \n) circle (2pt);
                    }
                }
            \end{tikzpicture}
            \end{minipage}
        \caption{The affine Coxeter arrangement $A_2$ in the case $n = 2$. The zonotope $Z_{[2]}$ and its $2L$-translates are the hexagons in white. The complex $\Sigma_{[2]}$ is the complex consisting of the triangles colored in orange. The square delimited in black is a fundamental domain for the torus 
        $V / 2 L  \cong  (S^1)^{[2]} / S^1 \cong (S^1)^2$ . The black dots in the right side of the figure represent the lattice $L$ in $V$. Here, to draw the figure, the vector space $V$ is identified with $\RR^2$ via the map $[x_0:x_1:x_2] \mapsto (x_1-x_0, x_2-x_0)$.}
        \label{fig:KummerExamplen=2}
    \end{figure}

We also have compatibility with the quotient maps obtained by forgetting
coordinates.

\begin{proposition}\label{prop:CompatibilityQuotientTriangulation}
Let $\varnothing\ne I \subset  J \subset [n]$, and let
$p_{J,I} : V_J\to V_I$. 
For every cell $\tau\in\mathcal{A}_J$, the image $p_{J,I}(\tau)$ is a cell of $\mathcal{A}_I$, and
\[
    p_{J,I}\bigl(\relint(\tau)\bigr) = \relint\bigl(p_{J,I}(\tau)\bigr).
\]
Moreover, for every $S \subset  I$,
\[
    \tau\in\mathcal{B}_{J,S} 
    \quad \text{if and only if} \quad 
    p_{J,I}(\tau)\in\mathcal{B}_{I,S}.
\]
Consequently, coordinate restriction induces simplicial maps
\[
    p_{J,I}:\mathcal{A}_J\longrightarrow\mathcal{A}_I,
    \qquad
    p_{J,I}:\mathcal{B}_{J,S}\longrightarrow\mathcal{B}_{I,S}.
\]
These maps also descend to the quotients by any finite-index lattices
$H_J \subset 2L_J$ and $H_I \subset 2L_I$ satisfying
$p_{J,I}(H_J) \subset  H_I$.
\end{proposition}

\begin{proof}
Fix $i_0\in I$, and use the coordinates $z_i=x_i-x_{i_0}$ on $V_J$.
The maximal cells of $\mathcal{A}_J$ are integral translates of the
simplices
\[
    0\le u_{j_1}\le\cdots\le u_{j_d}\le1,
    \qquad d=|J|-1,
\]
where $(j_1,\ldots,j_d)$ is an ordering of $J\setminus\{i_0\}$.
The map $p_{J,I}$ simply forgets the coordinates outside $I$. Since the
vertices of these simplices are the monotone $0$--$1$ vectors,
coordinate restriction sends them to the vertices of the corresponding
simplex for $I$, possibly with repetitions. It follows that every face
of $\mathcal{A}_J$ maps onto a face of $\mathcal{A}_I$.
Since the relative interior of a simplex is given by positive barycentric
coordinates at its vertices, the same argument also gives
$p_{J,I}(\relint(\tau)) = \relint(p_{J,I}(\tau))$.

Now fix $S \subset  I$. By the defining inequalities, we have
\[
    Z_{J,S} = p_{J,I}^{-1}(Z_{I,S})
\]
when $S$ is nonempty; for $S=\varnothing$, the same holds by our
convention. Since $p_{J,I}(2L_J)=2L_I$, we deduce
\[
    U_{J,S}=p_{J,I}^{-1}(U_{I,S}).
\]
Therefore,
\[
    \relint(\tau) \subset  U_{J,S}
    \quad \text{if and only if} \quad
    \relint\bigl(p_{J,I}(\tau)\bigr) \subset  U_{I,S},
\]
which gives the desired equivalence for the deleted complexes.
The assertions about the induced maps follow immediately.
\end{proof}

\begin{proposition}\label{pro:TriangulationCI}
Let $\varnothing\ne S \subset  I \subset [n]$, and let $H\subset 2L_I$
be a lattice of finite index. The logarithmic angle set of the Kummer
covering of
\[
    X_{I,S}:=V\Big(\sum_{i\in S}x_i\Big)
\]
corresponding to $H$ is
\[
    |\mathcal{B}_{I,S}|/H.
\]
The arrangement $\mathcal{A}_I$ induces a cellulation of this angle set. In particular, for $I=S=[n]$, we recover the angle set $|\Sigma_n|/H$.

For $S_1 \subset S_2 \subset I$, we have an inclusion of subcomplexes $\mathcal{B}_{I,S_1} \subset \mathcal{B}_{I,S_2}$. 
Thus the cellulation is compatible with initial degenerations. More precisely,
for every cone $\sigma_I$ of $\Delta$, the logarithmic angle set of the
corresponding initial degeneration is
\[
     |\mathcal{B}_{[n],[n]\setminus I}|/H \subset |\Sigma_n|/H.
\]
\end{proposition}

\begin{proof}
By \cite[Lemma 3.1]{NS13NonArchimedean}, the angle set of the
linear hyperplane with coordinates in $I$ is obtained by removing the
$2L_I$-translates of $Z_{I,I}$. 
The coordinates in $I\setminus S$ are free, so the same description in $V_I$ is precisely $|\mathcal{B}_{I,S}|$. Passing to the covering
associated to $H$ gives the first assertion. Since the boundaries of
the deleted regions lie in the Coxeter arrangement, the remaining cells
give a cellulation.

For $w\in\relint(\sigma_I)$, the initial equation is
$\sum_{i\in[n]\setminus I}x_i=0$, so the description of its angle set follows by taking $S=[n]\setminus I$. More generally, $S_1\subset S_2$ implies
$Z_{I,S_2}\subset Z_{I,S_1}$ and therefore
$\mathcal{B}_{I,S_1}\subset\mathcal{B}_{I,S_2}$.
This gives the claimed inclusion of subcomplexes.
\end{proof}

\subsection{The \texorpdfstring{$f$}{f}-vector for the building blocks of complete intersections}
\label{sec:FVectorCI}

We now count the cells in $\mathcal{B}_{n,S}$. These counts give
the local $f$-vectors used for complete intersections. We begin with the following lemma.

\begin{lemma}\label{lem:liftedSpoints}
    Let $0 \leq t_1 < \dots < t_k < 1$ represent points on the circle $\RR/\ZZ$ of length $1$ and let $\epsilon = (\epsilon_1, \dots, \epsilon_k) \in \{0,1\}^k$. Then the lifted points $\tilde{t}_i := t_{i} + \epsilon_i \in \RR/2\ZZ$ are all within a segment of length $1$ if and only if $\epsilon$ is of the following form
    \[
     1^m 0^{k-m}:= (\underbrace{1,\dots,1}_{m},
     \underbrace{0,\dots,0}_{k-m})
    \qquad\text{or}\qquad
    0^m 1^{k-m}:=  (\underbrace{0,\dots,0}_{m},
     \underbrace{1,\dots,1}_{k-m})
    \]
\end{lemma}
\begin{proof}
Suppose first that
\[
    \epsilon = 1^m 0^{k-m}
\]
with $0< m <k$. Then the lifted points are
\[
    t_{m+1} < \dots < t_k < 1+t_1 < \dots <1+t_m \in [0,2).
\]
On the circle $\RR/2\ZZ$ they are all contained in the arc
\[
    [t_{m+1},1+t_m],
\]
whose length is
\[
    1+t_m-t_{m+1}
    =
    1-(t_{m+1}-t_m)
    <1.
\]
The cases $m=0$ and $m=k$ are immediate, since all the points then lie in
an interval of length $t_k-t_1<1$. Note that the vectors $\epsilon$ of the form
\[
    (\underbrace{0,\dots,0}_{m}, \underbrace{1,\dots,1}_{k-m})
\]
are obtained from the preceding ones by adding $1$ simultaneously to all
the points in $\RR/2\ZZ$, and hence the lifted points with these vectors also lie in an interval of size less than $1$.

Conversely, suppose that the vector  $\epsilon = (\epsilon_1, \dots, \epsilon_k)$ changes value at least twice when reading the vector from left to right. Then
there exist indices $i<j<\ell$ such that
\[
    (\epsilon_i,\epsilon_j,\epsilon_\ell)
    =
    (1,0,1)
\quad 
\text{or}
\quad
    (\epsilon_i,\epsilon_j,\epsilon_\ell)
    =
    (0,1,0).
\]
Consider first the pattern $(1,0,1)$. Recall that $0 \leq t_i < t_j < t_\ell < 1$. The corresponding three lifted
points are
\[
    t_j < 1+t_i < 1+t_\ell \in [0,2).
\]
In cyclic order on $\RR/2\ZZ$, the three gaps between these points have lengths
\[
    1-(t_j-t_i) < 1,\qquad
    t_\ell-t_i < 1,\qquad
    1-(t_\ell-t_j) < 1.
\]
Note that if three points in $\RR/2\ZZ$ lie in an interval of length less than $1$, then the distance between the outer two points is greater than $1$. Since the distances between the three points $t_{j},1+t_{i},1+t_{\ell}$ are all less than $1$ by the above, we conclude that they cannot lie in an interval of length less than $1$.
The case $(0,1,0)$ is identical.
\end{proof}

\begin{theorem}\label{thm:FVectorCI}
For any lattice $H \subset 2L$ let $\Delta_{n,S,H}$ be the complex in $V / H$ induced by $\mathcal{B}_{n,S}$. Then the $f$-vector of the CW complex $\Delta_{n,S,H}$ is given by 

\[
    f_{d} =  [2L : H] \ d! \ g_{n,s}(d) \quad \text{for all } 0 \leq  d \leq n
\]
where
\[
 g_{n,s}(d) 
 =  2^n
     \left\{ \begin{matrix} n+1 \\ d+1 \end{matrix} \right\} - 2^{n+1-s} \sum_{\ell = 1}^{n + 1 - d} \left[ \binom{n+1}{\ell} - \binom{n+1-s}{\ell} \right] \left\{ \begin{matrix} n + 1 - \ell \\ d \end{matrix} \right\}.
 \]

Here $\left\{\begin{matrix} m \\ r \end{matrix}\right\}$ denotes the Stirling number of the second kind.
\end{theorem}

\begin{proof}
First we assume that $H = 2L$.
Fix $0 \leq d \leq n$ and let us first count the $d$-dimensional cells    
in the complex $\mathcal{A}_{n}/L \subset (S^1)^{n+1} / S^1$. 
A cell of dimension $d$ in the arrangement $\mathcal{A}_n/L$ is determined by the following data:
     \begin{itemize}
         \item a partition of $\{0,1,\dots,n\}$ into $d+1$ non-empty blocks recording the coordinates that have the same fractional part, and 
         \item a cyclic order on the fractional parts or equivalently on the $d+1$ blocks of the above partition.
    \end{itemize}
    Since the number of partitions of $\{0,1,\dots,n\}$ into $d+1$ cyclically ordered non-empty blocks is
    \[
        \frac{(d+1)!}{d+1} \ \left\{ \begin{matrix} n + 1 \\ d + 1 \end{matrix} \right\} = d! \ \left\{ \begin{matrix} n + 1 \\ d + 1 \end{matrix} \right\},
    \]
    we deduce that the number of $d$-cells in $\Delta_{n,\varnothing, 2L}$ (or equivalently in $\mathcal{A}_{n} \cap [-1,1)^n$) is exactly
    \[
        2^{n} \ d! \ \left\{ \begin{matrix} n + 1 \\ d + 1 \end{matrix} \right\}.
    \]

 The next step is to compute the deleted $d$-dimensional cells. As we saw, each cell in $V/L$ has $2^n$ lifts in $V/2L$. So a $d$-dimensional cell in $V/2L$ is given by the data of 
    \begin{itemize}
        \item A cyclically ordered partition $P = (P_0, \dots, P_d)$ of $\{0,\dots,n\}$ into non-empty blocks 
        \item a choice of a lift for each coordinate, i.e. lifting $\overline{y} \in \RR/\ZZ$ to either $\overline{y}$ or to $\overline{y+1}$ in $\RR/ 2\ZZ$. 
    \end{itemize}
    We now count the number of cells of dimension $d$ in $\Delta_{n,\varnothing, 2L}$ whose relative interior is contained in $Z_S$. To do so, let us first fix a cell of dimension $d$ in $\mathcal{A}_{n} \cap [0,1)^n$ corresponding to a cyclically ordered partition of $[n]$ into $d+1$ non-empty blocks
    \[
        P = (P_0, \dots, P_d).
    \]
    Let
    \[
        C_1,\dots,C_k
    \]
    be the blocks of $P$ that meet $S$, listed according to the increasing
    order of their fractional parts after choosing representatives in
    $[0,1)$. Write these fractional parts as
    \[
        0\leq t_1<t_2<\cdots<t_k<1.
    \]
    Thus
    \[
        k=k(P):=
        \#\{B\text{ a block of }P:B\cap S\neq\varnothing\}.
    \]
    
A lift of the cell from $V/L$ to $V/2L$ is obtained by choosing, for each
coordinate, whether its representative is lifted by $0$ or by $1$.
We first analyze the choices on the coordinates indexed by $S$.

If two elements of $S$ belong to the same block $C_j$, then in an
admissible lift they must receive the same $0/1$-shift. Indeed, if one
were shifted by $0$ and the other by $1$, the two lifted coordinates
would differ by exactly $1$, contrary to the defining strict inequality
for $Z_S$. Consequently, the lift on the $S$-coordinates is encoded by a
vector
\[
    \epsilon=(\epsilon_1,\dots,\epsilon_k)\in\{0,1\}^k,
\]
where $\epsilon_j$ is the common shift assigned to the elements of
$C_j\cap S$.

By Lemma~\ref{lem:liftedSpoints}, such a lift falls in the interior of $Z_S$ precisely when $\epsilon$ changes value at most once. The vectors $\epsilon$ that achieve that 
are exactly:
\[
    1^m0^{k-m}
    \quad\text{and}\quad
    0^m1^{k-m},
    \qquad 0\leq m\leq k.
\]
There are $2k$ distinct such vectors: each of the two families contains
$k+1$ vectors, and the two constant vectors $0^k$ and $1^k$ occur in both
families.

The coordinates indexed by $\{0,\dots,n\}\setminus S$ are unrestricted,
so they may be lifted independently. This contributes
\[
    2^{n+1-s}
\]
choices.

At this point we have therefore obtained
\[
    2^{n+1-s}\cdot 2k(P)
\]
binary lift vectors. However, two binary vectors that differ by
simultaneously changing every $0$ to $1$ and every $1$ to $0$ determine
the same lift in $V/2L$: indeed, this operation amounts to adding the
diagonal vector $\bm 1_{n+1}$ to all coordinates, which is trivial in
$V=\RR^{n+1}/\RR\bm 1_{n+1}$. Thus every lift is represented by exactly
two complementary binary vectors. It follows that the number of
admissible lifts of the fixed cell $P$ is
\[
    \frac{2^{n+1-s}\cdot 2k(P)}{2}
    =
    2^{n+1-s}k(P).
\]

We now sum over all cyclically ordered partitions $P$. The total number
of deleted $d$-cells is therefore equal to
\[
    2^{n+1-s}\sum_P k(P).
\]
The sum $\sum_P k(P)$ has the following useful interpretation. Since $k(P)$ is the number of blocks of $P$ that meet $S$, the sum counts pairs $(P,B)$ where $P$ is a cyclically ordered partition of $\{0,\dots,n\}$ into $d+1$ non-empty blocks and $B$ is a distinguished
block of $P$ satisfying
\[
    B\cap S\neq\varnothing.
\]
Once such a block $B$ is fixed, there is a unique way to rotate
the cyclic order and write
\[
    P=(P_0,P_1,\dots,P_d)
\]
so that $P_0=B$. Thus the pairs $(P,B)$ above are in bijection with ordered partitions
\[
    (P_0,\dots,P_d) \quad \text{such that } P_0 \cap S \neq \varnothing.
\]
Conversely, forgetting the distinguished first block of such a rooted
cyclic partition recovers the underlying cyclic partition $P$, while
$P_0$ records the distinguished block meeting $S$. Hence
\[
    \sum_P k(P)
    =
    \#\left\{
        (P_0,\dots,P_d):
        P_0\cap S\neq\varnothing
    \right\}.
\]
Let us count these partitions. The number of choices for $P_0$ of size $\ell$ that meet $S$ is exactly
\[
    \binom{n+1}{\ell} - \binom{n+1-s}{\ell},
\]
since the second term counts precisely the $\ell$-subsets contained
entirely in the complement of $S$.

After choosing $P_0$, the remaining $n+1-\ell$ elements must be
partitioned into $d$ non-empty blocks, which can be done in
\[
    \left\{
    \begin{matrix}
        n+1-\ell\\
        d
    \end{matrix}
    \right\}
\]
ways, and these $d$ blocks can then be ordered after $P_0$ in $d!$ ways.
Consequently,
\[
    \sum_P k(P)
    =
    d!
    \sum_{\ell=1}^{n+1-d}
    \left[
        \binom{n+1}{\ell}
        -
        \binom{n+1-s}{\ell}
    \right]
    \left\{
    \begin{matrix}
        n+1-\ell\\
        d
    \end{matrix}
    \right\}.
\]
Hence the number of deleted $d$-cells in $V/2L$ is
\[
    2^{n+1-s}d!
    \sum_{\ell=1}^{n+1-d}
    \left[
        \binom{n+1}{\ell}
        -
        \binom{n+1-s}{\ell}
    \right]
    \left\{
    \begin{matrix}
        n+1-\ell\\
        d
    \end{matrix}
    \right\}.
\]

Subtracting this from the total number of $d$-cells gives
\[
    f_d
    =
    d!
    \left(
        2^n
        \left\{
        \begin{matrix}
            n+1\\
            d+1
        \end{matrix}
        \right\}
        -
        2^{n+1-s}
        \sum_{\ell=1}^{n+1-d}
        \left[
            \binom{n+1}{\ell}
            -
            \binom{n+1-s}{\ell}
        \right]
        \left\{
        \begin{matrix}
            n+1-\ell\\
            d
        \end{matrix}
        \right\}
    \right) = d! \ g_{n,s}(d).
\]

    Finally, for a general lattice $H\subset 2L$, the quotient
    $V/H$ consists of $[2L:H]$ translates of the complex $\Delta_{n,S,2L}$. Therefore
    \[
        f_d = [2L:H] \ d! \ g_{n,s}(d),
    \]
    which gives the desired formula.
 \end{proof}

 We refer the reader to Theorem \ref{thm:refined_Coxter_F_vect} for a more general result.

\section{CW complex structures on algebraic varieties}
\label{sec:CWStructure}

In this section, we discuss a general scenario in which a tropical diagram can be used to construct CW complex structures. The main idea is to associate a suitable cellulation to every space. We show that such cellulated tropical diagrams naturally lead to CW complex structures on the homotopy colimit, which in turn give rise to a CW complex structure on algebraic varieties by Theorem \ref{thm:CWComplexStructureTriangulatedDiagram}. Our description of the corresponding CW complex in particular gives a formula for the $f$-vector of any open rational complete intersection manifold, see Section \ref{sec:FVectorRefinedArrangements}.

\subsection{Cellulated tropical diagrams}\label{sec:TriangulatedTropicalDiagrams}

In this section we discuss a general notion of a cellulated tropical diagram. Throughout this section, $D^{m}=\{x\in \RR^{m}:||x||\leq{1}\}$ will denote the closed $m$-dimensional ball. Here $||x||=\sqrt{x_{1}^{2}+...+x_{m}^{2}}$ is the standard Euclidean norm. 
\begin{definition} 
    Let $X$ and $Y$ be CW complexes, with $n$-skeletons $X_{n}$ and $Y_{n}$. A cellular map of CW complexes is a continuous map 
    \[
    f:X\to Y
    \]
    such that $f(X_{n})\subset Y_{n}$ for all $n$. We write $\text{CW}$ for the category of CW complexes, and $\text{Finite-CW}$ for the full subcategory of finite CW complexes.   
\end{definition}

\begin{definition}[Cellulations of tropical diagrams]\label{def:TriangulationTropicalDiagram}
    Let $\Pscr(\Sigma_{D})$ be the poset attached to a rational tropical manifold $(X,\Sigma,D)$, and let
    \begin{equation*}
        \mathcal{F}:\Pscr(\Sigma_{D})\to {\rm Top}
    \end{equation*}
    be a tropical diagram, see Definition \ref{def:TropicalDiagram}. That is, we have assigned a space $\mathcal{F}(\tau,\sigma)$ to every $(\tau,\sigma)$ that is homotopy equivalent to the Kato-Nakayama space $Z^{\sigma,\log}_{\tau}$, and this assignment preserves morphisms. We assume here that all of the spaces $\mathcal{F}(\tau,\sigma)$ are compact Hausdorff spaces. 
    
    A cellulation of the tropical diagram consists of a functor 
    \begin{equation*}
    \mathcal{F}^{\triangle}:\Pscr(\Sigma_{D})\to \text{Finite-CW}
    \end{equation*}
    such that the composite with the forgetful functor $\text{Finite-CW}\to {\rm Top}$ is $\mathcal{F}$. In other words, we have endowed the spaces $\mathcal{F}(\tau,\sigma)$ with a finite CW complex structure in such a way that all the comparison maps are cellular. 
    
\end{definition}

\begin{remark}
    We will see in Section \ref{sec:TriangulationsRCIM} that the cellulations from Section \ref{sec:CompleteIntersections} give rise to a cellulated tropical diagram. 
\end{remark}

\begin{definition}\label{def:CellsPerLength}\label{def:ShiftedFVector}
    Consider a cellulated tropical diagram $\mathcal{F}^{\triangle}$. We define the function   
    \begin{equation*}
        \mathfrak{g}_{k}:\ZZ\to\ZZ, \quad   m \mapsto \sum_{(z_{0}<...<z_{k})}|\mathcal{F}^{\triangle}(z_{k})_{m}|. 
    \end{equation*}
    Here the sum runs over all chains of length $k$ in $\Pscr(\Sigma_{D})$, and $\mathcal{F}^{\triangle}(z_{k})_{m}$ denotes the set of all $m$-simplices in $\mathcal{F}^{\triangle}(z_{k})$.
    Note that $\mathfrak{g}_{k}(m)=0$ for $m<0$ and for large $m$. 
    We define the shifted $f$-vector for chains of length $k$, which we denote by $\mathfrak{f}_{k}$, as follows:
    \[
        \mathfrak{f}_k : \ZZ \to \ZZ, \quad m \mapsto \mathfrak{g}_{k}(m - k).
    \] 
    These shifted $f$-vectors will appear in Theorem~\ref{thm:CWComplexStructureTriangulatedDiagram} to work with homotopy colimits. 
\end{definition}

We now recall the definition of the homotopy colimit, as in Definition \ref{def:HomotopyColimit}. We will work with the poset $\Pscr=\Pscr(\Sigma_{D})$ associated to a rational tropical manifold $(X,\Sigma,D)$ here. For a tropical diagram 
\[
\mathcal{F}:\Pscr(\Sigma_{D})\to \text{Top}
\]
as in Definition \ref{def:TriangulationTropicalDiagram},
we define 
\[
B=\bigsqcup_{z=(z_{0}<...<z_{k})\in \Delta(\Pscr(\Sigma_{D}))}\Delta(z)\times \mathcal{F}(z_{k}).
\]
Here $\Delta(z)$ is the $k$-dimensional simplex associated to the chain $(z_{0}<...<z_{k})$. We write elements of $B$ as $[z \  ; \ \phi \ ; \ x]$, where $z=(z_{0}<...<z_{k})\in \Delta(\Pscr(\Sigma_{D}))$, $\phi\in \Delta(z)\subset {\rm Hom}_{\text{Sets}}([k],\RR_{\geq{0}})=\RR^{k+1}_{\geq{0}}$ and $x\in \mathcal{F}(z_{k})$.
Recall from Lemma \ref{lem:RepresentativesHomotopyColimit2} that we have unique representatives for elements $y\in {\rm hocol}(\mathcal{F})$. That is, for every $y\in {\rm hocol}(\mathcal{F})$,  there is a unique $w=[z \ ; \ \phi \ ; x]$ with $z=(z_{0}<...<z_{k})$, $\phi(i)>0$ for all $0\leq{i}\leq{k}$ and $x\in \mathcal{F}(z_{k})$, such that $q(w)=y$. Here $q:B\to {\rm hocol}(\mathcal{F})$ is the quotient map.

\begin{lemma}
    Let $\mathcal{F}$ be as in Definition \ref{def:TriangulationTropicalDiagram}. Then the space ${\rm hocol}(\mathcal{F})$ is a compact Hausdorff space. 
\end{lemma}
\begin{proof}
The space $B$ is automatically a compact Hausdorff space by our assumption on the $\mathcal{F}(z)$.  We then have that ${\rm hocol}(\mathcal{F})$ is the image of a compact space by a continuous map, so it is compact. To see that it is Hausdorff, write $R\subset B\times B$ for the equivalence relation. We need to prove that this is a closed subset. But this easily follows from the fact that the comparison maps are continuous.
\end{proof}

We now use this to endow the homotopy colimit with the structure of a finite CW complex. 

\begin{definition}[Characteristic maps] \label{def:CharacteristicMaps}    
Let $\mathcal{F}^{\triangle}:\Pscr(\Sigma_{D})\to {\rm Finite-CW}$ be a cellulated tropical diagram. We define a set of characteristic maps for ${\rm hocol}(\mathcal{F})$, as in Definition \ref{def:TriangulationTropicalDiagram}. 
For $w\in \Pscr(\Sigma_{D})$, let  $\Phi(w)$ be the set of characteristic maps  $\rho:D^{m}\to \mathcal{F}(w)$ defining the CW structure on $\mathcal{F}(w)$. For every pair $(z,\rho)$ of a chain $z=(z_{0}<...<z_{k})\in \Delta(\Pscr(\Sigma_{D}))$ and a characteristic map $\rho\in \Phi(z_{k})$, we define 

\begin{equation*} 
    \chi_{z,\rho}:\Delta(z)\times D^{m}\to {\rm hocol}(\mathcal{F})   
\end{equation*}
by $\chi_{z,\rho}(t,x)=q([z \ ; \ t \ ; \ \rho(x)])$, where $q$ is the quotient map $$B=\bigsqcup_{z=(z_{0}<...<z_{k})\in \Delta(\Pscr(\Sigma_{D}))}\Delta(z)\times \mathcal{F}(z_{k})\to {\rm hocol}(\mathcal{F}).$$  
Note that the spaces $\Delta(z)\times D^{m}$ are homeomorphic to closed disks of dimension $k+m$.
\end{definition}

\begin{proposition}\label{prop:CWComplexHomotopyColimit}
    Let $\mathcal{F}^{\triangle}$ be a cellulation of a tropical diagram $\mathcal{F}:\Pscr(\Sigma_{D})\to \text{Top}$. The characteristic maps $\chi_{z,\rho}$ define a finite CW complex structure on ${\rm hocol}(\mathcal{F})$. 
\end{proposition}
\begin{proof}
    By Lemma \ref{lem:RepresentativesHomotopyColimit2}, we find that the interiors are disjoint, and they cover ${\rm hocol}(\mathcal{F})$. Moreover, they are homeomorphisms on their interiors. Since the comparison maps are cellular, we easily find that the dimensions of cells in the boundary are of lower dimension. We now conclude using \cite[Proposition A.2]{Hatcher2001}.
\end{proof}

\begin{remark}
    If all of the comparison maps are injective, then it is easy to see that this defines a regular CW complex. This is for instance the case for rational complete intersections $(X,\Sigma,D)$ with $D=D_{\text{full}}$. It is not clear to the authors if this is always the case for arbitrary $D$.
\end{remark}

\begin{definition}[Tropical CW complex] \label{def:TropicalDeltaComplex}
    Let $\mathcal{F}^{\triangle}$ be a cellulated tropical diagram with characteristic maps 
    \[
        \mathscr{S} = \Big\{\chi_{z,\rho} \ \colon \
        z = (z_{0},...,z_{k}) \in \Delta(\Pscr(\Sigma_{D})), 
        \rho\in \Phi(z_{k})\Big\}
    \] 
    as in Definition \ref{def:CharacteristicMaps}. The pair 
    $({\rm hocol}(\mathcal{F}),\mathscr{S})$ 
    is the tropical CW complex associated to 
    $\mathcal{F}^{\triangle}$. 
    We denote it by $\Delta(\mathcal{F}^{\triangle})$. 
    We will also write 
    $|\Delta(\mathcal{F}^{\triangle})|={\rm hocol}(\mathcal{F})$. 
    The pair $\Delta(\mathcal{F}^{\triangle})$ defines a CW complex structure on $|\Delta(\mathcal{F}^{\triangle})|={\rm hocol}(\mathcal{F})$ by Proposition \ref{prop:CWComplexHomotopyColimit}. We will write $|\Delta(\mathcal{F}^{\triangle})_{m}|$ for the number of cells of dimension $m$. The full $f$-vector of this CW complex is defined by 
    \[
        \mathfrak{f}(\Delta(\mathcal{F}^{\triangle}))(m)=|\Delta(\mathcal{F}^{\triangle})_{m}|.
    \]
    We view this as a function from $\ZZ\to \ZZ$, where 
    $\mathfrak{f}(\Delta(\mathcal{F}^{\triangle}))(m) = 0$ if $m < 0$ or if there are no cells of dimension $m$.
\end{definition}

By combining our description of the tropical CW complex together with Theorem \ref{thm:MainTheoremPolarDecomposition}, we now directly find the following.

\begin{theorem}\label{thm:CWComplexStructureTriangulatedDiagram}
Let $(X,\Sigma,D)$ be a rational tropical manifold, and let $(\mathcal{X}\to U,\mathcal{D})$ be an associated strictly semistable model.
Consider a cellulated tropical diagram 
\begin{equation*}
    \mathcal{F}^{\triangle}:\Pscr(\Sigma_{D}) \to {\rm Finite-CW}
\end{equation*}
with associated tropical CW complex $\Delta(\mathcal{F}^{\triangle})$, see Definition \ref{def:TropicalDeltaComplex}. This finite CW complex is homotopy equivalent to $\mathcal{X}_{t}\backslash \mathcal{D}_{t}$ for general $t\in U(\CC)$. Its $f$-vector is given by
    \begin{equation*}
        \mathfrak{f}(\Delta(\mathcal{F}^{\triangle}))
        = \sum_{k \geq 0} \mathfrak{f}_{k}.
    \end{equation*}
    Here $\mathfrak{f}_{k}$ is the shifted $f$-vector associated to $\mathcal{F}^{\triangle}$, see Definition \ref{def:CellsPerLength}.
\end{theorem}

\begin{proof}
    This follows from Theorem \ref{thm:MainTheoremPolarDecomposition} and Proposition \ref{prop:CWComplexHomotopyColimit}, since 
    \begin{equation} |\Delta(\mathcal{F}^{\triangle})_{j}|= \sum_{i=0}^{j}\sum_{(z_{0}<...<z_{i})\in \Delta(\Pscr(\Sigma_{D}))_{i}}|\mathcal{F}^{\triangle}(z_{i})_{j-i}|.
    \end{equation}
\end{proof}

\begin{corollary}
    Let $\mathcal{F}^{\triangle}$ be a cellulated tropical diagram for a connected rational tropical manifold $(X, \Sigma, D)$. 
    Then the fundamental group of $\mathcal{X}_{t}\backslash \mathcal{D}_{t}$ can be generated by 
    \[
        \sum_{z_{0}}|\mathcal{F}^{\triangle}(z_{0})_{1}|+\sum_{z_{0}<z_{1}}|\mathcal{F}^{\triangle}(z_{1})_{0}|-\sum_{z_{0}}|\mathcal{F}^{\triangle}(z_{0})_{0}|+1  
    \] 
    many elements. 
    Moreover, the kernel of this presentation can be normally generated by 
    \[
    |\Delta(\mathcal{F}^{\triangle})_{2}|=\sum_{z_{0}}|\mathcal{F}^{\triangle}(z_{0})_{2}|+\sum_{z_{0}<z_{1}}|\mathcal{F}^{\triangle}(z_{1})_{1}|+\sum_{z_{0}<z_{1}<z_{2}}|\mathcal{F}^{\triangle}(z_{2})_{0}|
    \]
    many elements.  
\end{corollary}

\begin{proof}
    This is immediate from the fact that $\Delta(\mathcal{F}^{\triangle})$ induces a CW complex structure on $\mathcal{X}_{t}\backslash \mathcal{D}_{t}$. 
\end{proof}

\begin{remark}
    We can write down a presentation of the fundamental group as follows. Let $\Delta(\mathcal{F}^{\triangle})_{1}$ be the $1$-skeleton, which is a finite graph. Let $T\subset \Delta(\mathcal{F}^{\triangle})_{1}$ be a spanning tree of the $1$-skeleton, and fix a base point $x_{0}$ in the spanning tree. The edges $e_{i}$ outside the spanning tree then give generators $g_{i}$, so that we can identify
    \[
    \pi_{1}(\Delta(\mathcal{F}^{\triangle})_{1})=\langle g_{1},...,g_{k}\rangle,
    \]
    the free group on $k$ generators. 
    
    Using the van Kampen theorem, we can then obtain a full presentation as follows. For every $2$-cell $\alpha\in \Delta(\mathcal{F}^{\triangle})_{2}$, we have an attaching map
    \[
    \phi_{\alpha}:S^{1}\to \Delta(\mathcal{F}^{\triangle})_{1}.
    \]
    This induces a homomorphism
    \[
    \psi_{\alpha}:\ZZ\simeq \pi_{1}(S^{1})\to \pi_{1}(\Delta(\mathcal{F}^{\triangle})_{1})=\langle g_{1},...,g_{k}\rangle.
    \]
    We then have the presentation
    \[
    \pi_{1}(\Delta(\mathcal{F}^{\triangle}))\simeq \langle g_{1},...,g_{k}\rangle /\langle \psi_{\alpha}(1) \ | \ \alpha\in \Delta(\mathcal{F}^{\triangle})_{2}\rangle.
    \]
    Here we note that the latter denotes the normal subgroup generated by $\psi_{\alpha}(1)$ for the finitely many $2$-cells  $\alpha\in\Delta(\mathcal{F}^{\triangle})_{2}$. We call $\pi_{1}(\Delta(\mathcal{F}^{\triangle}))$ a \emph{tropical fundamental group}. We will think of it through a presentation of the type presented above. We note that this relies on the choice of a spanning tree.  
\end{remark}

\subsection{Cellulated tropical diagrams for rational complete intersection manifolds}\label{sec:TriangulationsRCIM}

In this section, we show how to obtain cellulations of the tropical diagrams associated to rational complete intersection manifolds. 
In the easiest cases (such as Fermat hypersurfaces), one can simply use the cellulation from Proposition \ref{pro:TriangulationCI}, and this is compatible with all the maps. 
In more general cases, we shall need to refine this cellulation for the following reasons:
\begin{enumerate}
    \item Nontrivial gluing maps: for different adjacent vertices in the tropicalization, the local cellulations might not be compatible. That is, there might be \emph{twisting maps} that are not induced from initial degenerations.
    \item Projections of cells on the interior: at the boundary, projected cells might induce further subdivisions.
\end{enumerate}

To solve these problems, we can simply further subdivide the Coxeter arrangement and make sure that all the comparison maps are well defined. 
We give an $f$-vector formula later on in terms of refinements of the graphic matroid.

\subsubsection{Local cellulations}

We start by defining the local cellulations. Throughout this section, we shall use the following notation:
\[
    [n] :=\{0,1,\dots,n\}, \quad V := \RR^{n+1} / \RR \mathds{1}_{n+1} \quad \text{and} \quad L := \ZZ^{n+1}/\ZZ \mathds{1}_{n+1}.
\]
We write $\Xi$ for the dual lattice of $L$. So
\[
 \Xi := {\rm ker}(\mathds{1}_{n+1}) = \{ x \in \ZZ^{n+1} \colon x_0 + \dots + x_n = 0 \} \cong {\rm Hom}(L, \ZZ).
\]
Recall from Example \ref{exa:ProjectiveAngles} that we have an isomorphism of groups 
\begin{equation*}
    \mathrm{Hom}(\Xi,S^{1}) \simeq (S^{1})^{n+1}/S^{1}.
\end{equation*}
Indeed, applying ${\rm Hom}(-,S^{1})$ to the exact sequence 
\begin{equation*}
    0 \to \Xi\to \ZZ^{n+1} \to \ZZ\to 0,
\end{equation*}
we obtain
\begin{equation*}
    0 \to {\rm Hom}(\ZZ,S^{1})\to {\rm Hom}(\ZZ^{n+1},S^{1})\to {\rm Hom}(\Xi,S^{1})\to 0.
\end{equation*}

\begin{remark}
We can make this isomorphism explicit as follows. 
If we have a morphism $\phi: \ZZ^{n+1}\to S^{1}$, then we can restrict it to $\Xi$ to obtain an element  of ${\rm Hom}(\Xi,S^{1})$. Since ${\rm Hom}(\ZZ^{n+1},S^{1})\simeq (S^{1})^{n+1}$ by sending $\phi$ to $(\phi(e_{i}))_{i=0}^{n}$, we also obtain an isomorphism
\[
{\rm Hom}(\Xi,S^{1})\simeq (S^{1})^{n+1}/S^{1}.
\]
Explicitly, the isomorphism sends $\phi:\Xi\to S^{1}$ to the class of any $(s_{0},\ldots,s_{n})\in (S^{1})^{n+1}$ satisfying $\phi(e_{i}-e_{j})=s_{i}s_{j}^{-1}$ for all $i,j\in [n]$. In dehomogenized coordinates with respect to $e_{0}$, we have
\[
    \phi\mapsto [1:\phi(e_{1}-e_{0}):...:\phi(e_{n}-e_{0})].
\]
\end{remark}

We now consider the group homomorphism 
\begin{equation*}
    \RR\to \RR/2\ZZ\simeq S^{1}, \quad x \mapsto e^{\pi{i}x}. 
\end{equation*}
By composition, this gives  
    \[
        {\rm Hom}(\Xi,\RR)\to {\rm Hom(\Xi,S^{1})}.
    \]
We identify this with the natural universal covering map
    \[
        V:=\RR^{n+1}/\RR\mathds{1}\to (S^{1})^{n+1}/S^{1}.
    \] 
Note that the kernel of ${\rm Hom}(\Xi,\RR)\to \mathrm{Hom}(\Xi,S^{1})$ is $\mathrm{Hom}(\Xi,2\ZZ)$,  which we identify as before with $2L$. We thus obtain an exact sequence 
    \[
        0 \to 2L\to V\to (S^{1})^{n+1}/S^{1}.
    \]

Note that every element of $\Xi_{\RR}$ can be identified with a linear form on $V$. Indeed, by applying ${\rm Hom}(-,\RR)$, we obtain isomorphisms 
\begin{equation*}
   \Xi_{\RR}\simeq {\rm Hom}({\rm Hom}(\Xi,\RR),\RR)\simeq {\rm Hom}(V,\RR).
\end{equation*}
Similarly, elements of ${\rm Hom}(V,\RR)$ will be viewed as linear forms on $\RR^{n+1}$ that are invariant under the action of $\RR\mathds{1}$.

\begin{definition}[Refined Coxeter arrangements]
\label{def:refinedCoxeterArrangement}
For $u\in \Xi_{\QQ}$ we define 
\[
    d_u := \min\{ d \in \ZZ_{\geq 1} \colon d u \in \Xi \} \quad \text{and} \quad \hat{u} := d_u  u. 
\]
Let $Q$ be a finite set of non-zero vectors in $\Xi_{\QQ}$. The extended set of Coxeter hyperplanes associated to $Q$ is defined as follows:
\[
    E(Q) := E \cup \{\hat{u} \colon u \in Q\} \subset \  \Xi.
\]
For $\ell \in \ZZ$ and $u \in \Xi$, we define the affine hyperplane
\[
    H_u^{(\ell)} := \{x\in V: \langle u, x\rangle=\ell\}.
\]
The refined Coxeter arrangement $\mathscr{A}_{n}(Q)$ is the infinite arrangement 
\[
    \mathscr{A}_{n}(Q) := \{ H_{u}^{(\ell)} \colon u \in E(Q) \text{ and } \ell \in \ZZ\}.
\] 
We note that the refined Coxeter arrangement is $L$-periodic. 
We write $\mathcal A_n(Q)$ for the polyhedral complex induced by
the hyperplanes in $\mathscr{A}_n(Q)$. Note that $|\mathcal{A}_{n}(Q)|=V$.  
\end{definition}

We now define a cellulation of the universal covering of the angle set of the hyperplane 
\[
    V\Big(\langle x_{0}+x_{1}+...+x_{n}\rangle\cap \CC[\Xi]\Big)\subset \Spec(\CC[\Xi]).
\]  
Note that the angle set of this complex variety can be identified with a subset of $(S^{1})^{n+1}/S^{1}$, and thus the universal covering can be identified with $V=\RR^{n+1}/\RR\mathds{1}$.
\begin{definition}[Refined cellulations of angle sets]
\label{def:RefinedTriangulationsOfAngleSets}

Let $S \subset [n]$ be a subset with $s=|S|\geq2$. We set
\[
    E_S=\{e_i-e_j:i<j,\ i,j\in S\} \subset E
    \quad \text{and} \quad
    Z_S := \{x\in V: |x_i-x_j|<1\text{ for all }i,j\in S\}.
\]
Let $Q\subset \Lambda_{\QQ}$ be a finite set of nonzero vectors. 
Define
\[
    \mathcal B_{n,S}(Q)
    :=
    \Big\{
    \tau \in \mathcal A_n(Q): {\rm relint}(\tau)\not \subset 
    \bigcup_{\gamma\in 2L}(\gamma+Z_S) 
    \Big\}.
\]
Then $\mathcal{B}_{n,S}(Q)$ is the  cellulation of the universal covering of the angle set of
\[
    V\Big(\langle \sum_{i\in{S}}x_{i}\rangle\cap \CC[\Xi]\Big)\subset \Spec(\CC[\Xi])
\] 
induced by $\mathcal{A}_{n}(Q)$. 
Its cells are invariant under the action of vectors $e_{i}$ for $i\notin{S}$.  

Similar to Section~\ref{sec:AngleSetsCI}, for any $I\subset [n]$ we write $\Xi_{I} :=\text{ker}(\mathds{1}_{I})\subset \RR^{I}$ for the corresponding lattice, and $V_{I}=\RR^{I}/\RR\mathds{1}_{I}$. For any finite set $Q_{I}\subset \Xi_{I,\QQ}$ of nonzero vectors, we then have analogous definitions of $\mathscr{A}_{I}(Q_{I})$, $\mathcal{A}_{I}(Q_{I})$, $E_{I,S}$, $Z_{I,S}$, and $\mathcal{B}_{I,S}(Q_{I})$. 
\end{definition}

We now use these Coxeter arrangements to define a CW complex structure on angle sets. 

\begin{definition}[Local cellulations]\label{def:LocalTriangulation}
Let $V=\RR^{n+1}/\RR\mathds{1}$ and let $Q\subset \Xi_{\QQ}$ be a finite subset of nonzero vectors. 
 For any finite-index sublattice $H  \subset  2 L$, let $V/H$ be the corresponding $n$-dimensional torus. 
Set
\[
    \Delta_{n,S,H}(Q)
    :=
    \mathcal B_{n,S}(Q)/H.
\]
Note that $\Delta_{n,S,H}(Q)$ gives a CW complex structure on a subspace of the torus $V/H$. We similarly define $\Delta_{I,S,H}$ for any subgroup $H_{I}\subset 2L_{I}$ of finite index. It gives a CW complex structure on a subspace of the torus $V_{I}/H_{I}$.
\end{definition}

\subsubsection{Gluing cellulations}

We define a general framework for gluing these cellulations. The idea is that the material from the previous section gives a local prototype for a cellulation. We then have to make sure that the induced refined Coxeter arrangements are compatible in a suitable sense. 

\begin{definition}[Affine simplices]\label{def:AffineSimplex}
        Let $\psi: \ZZ^{n+1}\to M$ be a group homomorphism whose image $\Lambda={\rm Im}(\psi)$ is a subgroup of finite index in $M$, and the restricted homomorphism $\psi \colon \Xi \to \Lambda$ is an isomorphism. 
        We call this an \emph{affine simplex} in $M$. If $\Lambda = M$, then this is called a \emph{unimodular affine simplex} in $M$.
\end{definition}
\begin{remark} We note the following:
    \begin{enumerate}
        \item for any fixed $j$, the $n$ elements $\{\psi(e_{i}-e_{j}) \in \Lambda\subset M \colon 0 \leq i\neq j \leq n\}$ form a rational basis of $M$. 
        If $\psi$ defines a unimodular affine simplex, then these form a $\ZZ$-basis for $M$. 
        
        \item We have inclusion maps $\Lambda\to M\to \dfrac{1}{m}\Lambda$, where $m = [M:\Lambda]$.
    \end{enumerate}
\end{remark}

\begin{definition}[Affine simplices compatible with polynomials]\label{def:AffineSimplexCompatiblePolynomial}
    Let $f\in \CC[M]$. We say that an affine simplex $\psi:\ZZ^{n+1}\to M$ is compatible with $f$ if ${\rm Supp}(f)\subset \{\psi(e_{i})\}$. In particular, this implies that $V(f)$ is Kummer smooth, see Lemma \ref{lem:KummerLemma}.  
\end{definition}

\begin{example}\label{exa:AffineSimplexExample}
    Suppose that $f=x^2+xy^3+y^5$, with $\Lambda=\ZZ^{2}$. We then obtain a unimodular affine simplex 
    \begin{equation*}
        \psi:\ZZ^{3}\to M, \quad e_0 \mapsto x^2, \ e_1 \mapsto xy^3 , e_2 \mapsto y^5.
    \end{equation*}
    Note that the induced map $\Lambda \to M$ is indeed an isomorphism. For instance, the vectors $e_{0}-e_{2}$ and $e_{1}-e_{2}$ map to $(2,-5)$ and $(1,-2)$, which form a basis of $M$. 
    The corresponding hypersurface is 
    \begin{equation*}
        f/y^5=(x^2/y^5)+(x/y^{2})+1.
    \end{equation*}
\end{example}

    For an affine simplex $\psi:\ZZ^{n+1}\to M$, we write $\psi_{\QQ}:{\rm ker}(\mathds{1}_{n+1})_{\QQ}\to M_{\QQ}$ for the induced isomorphism of rational vector spaces. We then have an inclusion of lattices 
    \[
    \Xi \subset \psi^{-1}_{\QQ}(M)\subset \frac{1}{m} \cdot \Xi
    \]
    for $m=[M:\Lambda]$. 
    In particular, we can identify the various angle sets.
    \begin{definition}\label{def:TriangulationData}[Cellulation data]
        Let $\psi:\ZZ^{n+1}\to M$ be an affine simplex that is compatible with a polynomial $f\in \CC[M]$, as in Definition \ref{def:AffineSimplexCompatiblePolynomial}. Let $\psi_{\QQ}$ be the induced restricted map $\Xi_{\QQ}\to M_{\QQ}$. The lattices $M$ and $\psi^{-1}_{\QQ}(M)$ are then isomorphic.
        The isomorphism of lattices induces an isomorphism of angle sets 
        \[
            \theta_{\psi}: {\rm Hom}(\psi^{-1}_{\QQ}(M),S^{1}) \simeq {\rm Hom}(M,S^{1}).
        \]
        Let $Q\subset \Xi_{\QQ}$ be a finite subset. We call the pair $(\psi,Q)$ a set of cellulation data for $f$. It induces 
        a cellulation on the Kato-Nakayama angle set $\Theta^{\log}_{V(f)}$ of $V(f)\subset \Spec(\CC[M])$ by Definition \ref{def:LocalTriangulation} with 
        $H = {\rm Hom}(\psi^{-1}_{\QQ}(M), 2 \ZZ)$, $n = {\rm rank}(M)$, and $S\subset [n]$ such that $\mathrm{Supp}(f) = \{\psi(e_{i}):i\in{S}\}$. 
         We call this the induced cellulation attached to the pair  $(\psi,Q)$.
    \end{definition}

\begin{remark}
    Suppose we have chosen cellulation data $(\psi,Q)$ for some polynomial $f\in \CC[M]$. The pair $(\psi,Q)$ then also provides a set of cellulation data for every initial form ${\rm in}_{w}(f)$ corresponding to a cell $\tau\in \trop(V(f))$ with $w\in {\rm relint}(\tau)$. Here we take the standard Gr\"{o}bner complex structure on $\trop(V(f))$. 
\end{remark}

\begin{definition}[Compatible cellulations]
     Let $(\psi_{1},Q_{1})$ and $(\psi_{2},Q_{2})$ be two sets of cellulation data for a polynomial $f\in \CC[M]$. We say that these are compatible if they induce the same cellulation on $\Theta^{\log}_{V(f)}$.
\end{definition}

\begin{lemma}
    Suppose that $(\psi_{1},Q_{1})$ and $(\psi_{2},Q_{2})$ are cellulation data for a polynomial $f\in \CC[M]$. Then there is a $(\psi,Q)$ whose induced cellulation of $\Theta^{\log}_{V(f)}$ refines both.
\end{lemma}

\begin{proof}

    We start with $(\psi_{1},Q_{1})$. We then take $Q=Q_{1}\cup \psi^{-1}_{1,\QQ}(\psi_{2,\QQ}(E(Q_{2})))$. That is, we take all of the Coxeter hyperplanes $E(Q_{2})$ that define the cellulation with respect to $\psi_{2}$ and $Q_{2}$, and we add them to $Q_{1}$. The induced cellulation of $\Theta^{\log}_{V(f)}$ is then automatically a refinement of both.
\end{proof}

\begin{example}[$n=2$] \label{ex:refinedCoxeterArrangement1}
    Let $M=\ZZ^2$ and consider the curve cut out by the polynomial
    \[
        f = 1+x+y^2 + ty^{-1} \in \CC(t)[M].
    \]
    Its tropicalization in $\RR^2$ has two vertices $v_1 = (0,1)$ and $v_2 = (0,0)$ and the corresponding initial forms in $\CC[M]$ are
    \[
        {\rm in}_{v_1}(f) = 1+x+y^{-1} \quad \text{and} \quad {\rm in}_{v_2}(f) = 1 + x + y^2.
    \]
    We recall that
    \[
        E := \{e_1 - e_0, e_2 - e_0, e_2 - e_1  \}.
    \]
    We start with the following local cellulation data for the polynomial $f$. At the vertex $v_1$ we have
    \[
        \psi_1 \colon \ZZ^3 \to M, \quad e_0 \mapsto (0,0), \quad e_1 \mapsto (1,0) \text{,} \quad e_2 \mapsto (0,-1) \qquad \text{and} \qquad Q_1 = \varnothing.
    \]
  At the vertex $v_2$ we have
    \[
        \psi_2 \colon \ZZ^3 \to M, \quad e_0 \mapsto (0,0), \quad e_1 \mapsto (1,0) \text{,} \quad e_2 \mapsto (0,2) \qquad \text{and} \qquad Q_2 = \varnothing.
    \]
    To refine this cellulation data we compute
    \[
            \psi_1(E) = \{ (-1,-,1), (1,0), (0,-1) \} 
            \quad \text{and} \quad 
            \psi_2(E) = \{ (-1,2), (1,0), (0,2) \}. 
    \]
    We find  
    \[
    \begin{aligned}
       (\psi_1)_\QQ^{-1}(\psi_2(E))  &= \{(3,-1,-2), (-1,1,0), (2,0,-2)\} \\
       (\psi_2)_\QQ^{-1}(\psi_1(E))  &= \left\{ \left(\frac{3}{2}, -1, -\frac{1}{2}\right), (-1,1,0), \left(\frac{1}{2}, 0, -\frac{1}{2}\right) \right\}.
    \end{aligned}
    \]
    We thus update $Q_1$ and $Q_2$ to obtain:
    \[
    \begin{aligned}
        Q_1' &:= Q_1 \cup (\psi_1)_\QQ^{-1}(\psi_2(E))
        \quad \text{and} \quad
        Q_2' &:= Q_2 \cup (\psi_2)_\QQ^{-1}(\psi_1(E)).
    \end{aligned}
    \]
    Next we calculate the two sublattices $H_1, H_2$ of $\ZZ^3/\ZZ \mathds{1}_3$ and we obtain
    \[
    \begin{aligned}
        H_1 := {\rm Hom}\Big((\psi_1)_\QQ^{-1}(M), 2\ZZ\Big) 
        \cong \Big\langle [0,2,0], [0,0,2] \Big\rangle_{\ZZ} \\
        H_2 := {\rm Hom}\Big((\psi_2)_\QQ^{-1}(M), 2\ZZ\Big)
        \cong \Big\langle [0,2,0], [0,0,4]\Big\rangle_{\ZZ}.
    \end{aligned}
    \]
    The cellulations $(\psi_1, E(Q_1'))$ and $(\psi_2, E(Q_2'))$ of the angle sets of $\mathrm{in}_{v_1}(f)$ and $\mathrm{in}_{v_2}(f)$ are drawn in Figure~\ref{fig:refined-coxeter-arrangements}.

\begin{figure}[htbp]
    \centering

    \begin{tikzpicture}[scale=0.9,line cap=round,line join=round]
        \clip (-3,-3) rectangle (3,3);

        \foreach \i in {-2,...,2} {
            \foreach \j in {-2,...,2} {
                \pgfmathsetmacro\xshift{2*\i}
                \pgfmathsetmacro\yshift{2*\j}

                \filldraw[
                    fill=orange!45,
                    draw=orange!65!black,
                    line width=0.3pt
                ]
                    ({-1+\xshift},{\yshift}) --
                    ({\xshift},{1+\yshift}) --
                    ({-1+\xshift},{1+\yshift}) -- cycle;

                \filldraw[
                    fill=orange!45,
                    draw=orange!65!black,
                    line width=0.3pt
                ]
                    ({\xshift},{-1+\yshift}) --
                    ({1+\xshift},{-1+\yshift}) --
                    ({1+\xshift},{\yshift}) -- cycle;
            }
        }

        \foreach \ell in {-6,...,6} {
            \draw[black!45,line width=0.35pt]
                (\ell,-3) -- (\ell,3);
            \draw[black!45,line width=0.35pt]
                (-3,\ell) -- (3,\ell);
            \draw[black!45,line width=0.35pt]
                (-3,{-3-\ell}) -- (3,{3-\ell});
        }

        \foreach \ell in {-12,...,12} {
            \draw[blue!80!black,line width=0.6pt]
                (-3,{\ell/2}) -- (3,{\ell/2});
        }

        \foreach \ell in {-12,...,12} {
            \draw[red!80!black,line width=0.6pt]
                (-3,{(\ell+3)/2}) --
                (3,{(\ell-3)/2});
        }

        \draw[->,black,line width=0.8pt]
            (-3,0) -- (3.4,0) node[right] {$x$};
        \draw[->,black,line width=0.8pt]
            (0,-3) -- (0,3.4) node[above] {$y$};

        \draw[black,line width=1.8pt]
            (-1,-1) rectangle (1,1);

        \foreach \i in {-3,...,3} {
            \foreach \j in {-3,...,3} {
                \fill[black] (\i,\j) circle (1.15pt);
            }
        }
    \end{tikzpicture}
    \hspace{2cm}
    \begin{tikzpicture}[scale=0.9,line cap=round,line join=round]
        \clip (-3,-3) rectangle (3,3);

        \foreach \i in {-2,...,2} {
            \foreach \j in {-2,...,2} {
                \pgfmathsetmacro\xshift{2*\i}
                \pgfmathsetmacro\yshift{2*\j}

                \filldraw[
                    fill=orange!45,
                    draw=orange!65!black,
                    line width=0.3pt
                ]
                    ({-1+\xshift},{\yshift}) --
                    ({\xshift},{1+\yshift}) --
                    ({-1+\xshift},{1+\yshift}) -- cycle;

                \filldraw[
                    fill=orange!45,
                    draw=orange!65!black,
                    line width=0.3pt
                ]
                    ({\xshift},{-1+\yshift}) --
                    ({1+\xshift},{-1+\yshift}) --
                    ({1+\xshift},{\yshift}) -- cycle;
            }
        }

        \foreach \ell in {-6,...,6} {
            \draw[black!45,line width=0.35pt]
                (\ell,-3) -- (\ell,3);
            \draw[black!45,line width=0.35pt]
                (-3,\ell) -- (3,\ell);
            \draw[black!45,line width=0.35pt]
                (-3,{-3-\ell}) -- (3,{3-\ell});
        }

        \foreach \ell in {-12,...,12} {
            \draw[blue!80!black,line width=0.6pt]
                (-3,{\ell+6}) -- (3,{\ell-6});
        }

        \draw[->,black,line width=0.8pt]
            (-3,0) -- (3.4,0) node[right] {$x$};
        \draw[->,black,line width=0.8pt]
            (0,-3) -- (0,3.4) node[above] {$y$};

        \draw[black,line width=1.8pt]
            (-1,-2) rectangle (1,2);

        \foreach \i in {-3,...,3} {
            \foreach \j in {-3,...,3} {
                \fill[black] (\i,\j) circle (1.15pt);
            }
        }
    \end{tikzpicture}

    \caption{
        Refined Coxeter arrangements associated with the
        cellulation data $(\psi_1,Q_1')$ (left) and
        $(\psi_2,Q_2')$ (right).
        The orange regions indicate the surviving Coxeter cells
        corresponding to the universal covering of the angle sets.
        The blue and red lines represent the additional
        hyperplanes introduced by the refinements.
        The thick black rectangles are fundamental domains
        for the lattices $H_1$ and $H_2$, respectively.
    }
    \label{fig:refined-coxeter-arrangements}
\end{figure}

Note that the initial form of $f$ corresponding to the edge $\tau$ between $v_1$ and $v_{2}$ is 
\[
    {\rm in}_{\tau}(f) = 1 +x,
\]
and its angle set is the circle $S^1$, which is depicted in Figure~\ref{fig:refined-coxeter-arrangements} as the vertical sides of the thick black rectangles. Since we refined the cellulation, this circle is now subdivided into 4 pieces in both angle sets. For the angle set on the left, the circle is attached to two orange triangles, whereas on the right it is attached to four orange triangles. 
\end{example}

\begin{definition}[Cellulations for hypersurfaces]
\label{def:HypersurfaceTriangulation}

    Let $(X,\Sigma,D)$ be a rational tropical hypersurface. 
    We define a cellulation of the angle functor $\mathcal{F}_{\Theta,{{\rm CI}}}:\Pscr(\Sigma_{D})\to {\rm Top}$. Consider the set of all interior vertices $(v,0)\in \Pscr(\Sigma_{D})$. In other words, $v\in \Sigma$ is a vertex of the chosen polyhedral structure. We assign cellulation data $(\psi_{v},Q_{v})$ to $\mathrm{in}_{v}(f)$. We will assume that these are compatible. That is, the induced cellulation data $(\psi_{v},Q_{v})$ for $w\in {\rm relint}(\tau)$ are compatible for all vertices $v\in \tau$. 
    
    We then use induction on $\dim(\sigma)$ to construct more general cellulations. For cells $(\tau,\sigma)$ with $\dim(\sigma)=1$, we consider the cells in $\mathcal{F}_{\Theta,{\rm CI}}(\tau,0)$ as subsets of $\mathrm{Hom}(M,\RR)$. We then take their images under the linear projection map
    \[
    \mathrm{Hom}(M,\RR)\to \mathrm{Hom}(M(\sigma),\RR)
    \]
    and take a cellulation that refines this subdivision. Moreover, we can refine the cellulation to ensure that the cellulation data is compatible, as in the zero-dimensional case. We then inductively apply the same procedure for higher-dimensional $\sigma$ until we have reached all cells. The induced maps of CW complexes are then automatically cellular.
\end{definition}

    \begin{remark}
        In general, this will require us to apply some version of Fourier-Motzkin to explicitly calculate the induced convex polytopes in the quotients. We note that this is possible using larger $Q$ for the infinite strata. When $Q=\varnothing$ everywhere, this is not necessary however by Proposition \ref{prop:CompatibilityQuotientTriangulation}.
    \end{remark}

\begin{remark}
For more general rational complete intersections that are not hypersurfaces, we can cellulate angle sets as follows. We only have to cellulate the angle sets at the vertices, as this also cellulates the angle sets at cells $(\tau,0)$, and we can then again project for cells of the form $(\tau,\sigma)$. 

At a vertex, we have a decomposition
\begin{equation*}
M_{\QQ}=M_{1,\QQ}\oplus ...\oplus M_{k,\QQ},
\end{equation*}
as in Proposition \ref{prop:SmoothStableIntersection}. We make sure that the supports of all the local equations are contained in $M_{i,\QQ}$. We now repeat the construction from Definition \ref{def:TriangulationData}. 
\end{remark}

\begin{remark}
    When the variety $X$ is defined over $\CC$, we do not need any complicated charts, essentially by Proposition \ref{prop:CompatibilityQuotientTriangulation}. This will for instance allow us to find cellulations of toric varieties, Fermat hypersurfaces and Fermat complete intersections.  
\end{remark}

\begin{corollary}\label{cor:CITriangulation}
    Let $\mathcal{F}^{\triangle}$ be the cellulation of the angle functor $\mathcal{F}_{\Theta,{\rm CI}}$ associated to a rational complete intersection manifold $(X,\Sigma,D)$ and a set of charts. Then the geometric realization  of the tropical CW complex $\Delta(\mathcal{F}^{\triangle})$ is homotopy equivalent to $\mathcal{X}_{t}\backslash \mathcal{D}_{t}$ for general $t\in U(\CC)$.
    This endows $\mathcal{X}_{t}\backslash \mathcal{D}_{t}$ with the structure of a finite CW complex, up to homotopy.  
\end{corollary}

\begin{remark}
    Using Theorem~\ref{thm:CWComplexStructureTriangulatedDiagram} and Corollary~\ref{cor:CITriangulation}, we obtain an explicit formula for the $f$-vector of the associated tropical CW complex. Especially when the local charts are simple, this will be easy to calculate. We will go through many explicit examples in Section \ref{sec:Examples}.    
\end{remark}

\subsection{F-vectors for refined Coxeter arrangements}
\label{sec:FVectorRefinedArrangements}

In this section, we generalize the formula from
Theorem~\ref{thm:FVectorCI} to refined Coxeter arrangements defined by
rational vectors.  The main point is that, for a rational vector
$u$, we allow the levels of the corresponding affine hyperplanes to
lie in a suitable discrete subgroup of $\RR$.  With this convention,
the refined arrangement is $L$-periodic, so that we can first count
its cells on $V/L$ and then lift them to $V/2L$, exactly as in the
proof of Theorem~\ref{thm:FVectorCI}.

We use the same notation as in the previous section and we fix a finite set $Q$ of nonzero vectors in $\Xi_{\QQ}$, where $\Xi:={\rm ker}(\mathds{1}_{n+1})$. We recall from Definitions \ref{def:refinedCoxeterArrangement} and \ref{def:RefinedTriangulationsOfAngleSets} that
\[
    E(Q) := \{e_j - e_i:0\leq i<j\leq n\} \cup \{ \hat{u} \colon u \in Q\}.
\]
and for any $S \subset [n]$, with $s:=|S|\geq2$ we have
\[
    E_S := \{e_j - e_i : i<j \in S\}.
\]
Let $H \subset 2L$ be a lattice of finite index. To compute the $f$-vector of the arrangement $\Delta_{n,S,H}(Q)$ on the torus, we shall find it convenient to work with matroids over $\ZZ$ as follows. 
We associate to $E(Q)$ a realizable matroid $\Mscr$ over $\ZZ$ in the sense of \cite{FM16}. We define a matroid $\Mscr$ over $\ZZ$ as follows:

\[
    \Mscr(I) := \Xi \Big / \sum_{u\in I}\ZZ u, \quad \text{for }  I  \subset   E(Q).
\]
For any  $I \subset  E(Q)$ we write
\[
    {\rm rk}_{\Mscr}(I)
    :=
    {\rm rank}_{\ZZ}
    \left(
        \sum_{u\in I}\ZZ u
    \right),
\quad
\text{and}
\quad
    \mu_{\Mscr}(I) := |\Mscr(I)_{\rm tor}|.
\]
The arithmetic Tutte polynomial of $\Mscr$ is
\[
    {\rm T}_{\Mscr}(x,y)
    :=
    \sum_{I  \subset  E(Q)}
    \mu_{\Mscr}(I)
    (x-1)^{n-{\rm rk}_{\Mscr}(I)}
    (y-1)^{|I|-{\rm rk}_{\Mscr}(I)}.
\]
For $I \subset  E(Q)$, the contraction $\Mscr/I$ is the matroid
over $\ZZ$ on the ground set $E(Q) \setminus I$ given by
\[
    (\Mscr/I)(J)
    :=
    \Mscr(I\cup J),
    \qquad
    J \subset  E(Q)\setminus I.
\]
Thus
\[
    {\rm T}_{\Mscr/I}(x,y)
    =
    \sum_{J \subset  E(Q) \setminus I}
    \mu_{\Mscr}(I\cup J) \ 
    (x-1)^{
        n
        -
        {\rm rk}_{\Mscr}(I\cup J)
        }
     \ 
    (y-1)^{
        |J|
        -
        {\rm rk}_{\Mscr}(I\cup J)
        +
        {\rm rk}_{\Mscr}(I)
    }.
\]
We can now state the refined $f$-vector formula.

\begin{theorem}\label{thm:refined_Coxter_F_vect}
Let $Q \subset \Xi_\QQ$ be a finite set of nonzero rational vectors, let $S \subset [n]$ with $s = |S| \geq 2$, and let $H  \subset  2L$ be a finite-index sublattice. For
$0\leq d\leq n$, the number of $d$-dimensional cells of the complex $\Delta_{n,S,H}(Q)$ is
\[
    f_d
    =
    [2L:H]
    \sum_{\substack{
        I  \subset  E(Q)\\
        {\rm rk}_{\Mscr}(I)=n-d}}
    \left(
        2^n
        -
        2^{n+1-s}
        \bigl(
            s-{\rm rk}_{\Mscr}(I\cap E_S)
        \bigr)
    \right)
    {\rm T}_{\Mscr/I}(1,0).
\]
In particular, the formula depends only on the matroid $\Mscr$ over
$\ZZ$ associated to the $\Xi$-integral vectors $E(Q)$, the distinguished
Coxeter subset $E_S$, and the index $[2L:H]$.
\end{theorem}

\begin{proof}
We follow the same proof strategy as in Theorem~\ref{thm:FVectorCI}.
We first count the cells of the refined arrangement modulo $L$.
Then we determine how many of their lifts to $V/2L$ survive the
deletion of the translates of $Z_S$, and finally pass to $V/H$.

\begin{enumerate}
\item[Step 1:]
Since every vector in $E(Q)$ belongs to $\Xi$, the arrangement
$\mathscr A_n(Q)$ is $L$-periodic.
It therefore descends to a toric arrangement on
$T:=V/L\simeq{\rm Hom}(\Xi,\RR/\ZZ)$.
For $u\in E(Q)$, set
\[
  H_u:=\{\chi\in T:\chi(u)=0\},
\]
and, for $I \subset  E(Q)$, let
\[
  H_I:=\bigcap_{u\in I}H_u,
  \qquad H_\varnothing:=T.
\]
By the definition of the matroid over $\ZZ$,
\[
  H_I\simeq
  {\rm Hom}\left(
    \Xi\Big/\sum_{u\in I}\ZZ u,\RR/\ZZ
  \right)
  ={\rm Hom}\bigl(\Mscr(I),\RR/\ZZ\bigr).
\]
In particular, $H_I$ is a disjoint union of
$\mu_{\Mscr}(I)$ real tori, each of dimension
\[
  q=n-{\rm rk}_{\Mscr}(I).
\]
This also holds for $I=\varnothing$, since
$\Mscr(\varnothing)=\Xi$ and
$\mu_{\Mscr}(\varnothing)=1$.

For a cell $\tau$ of the quotient arrangement on $T$,
define its \emph{support} by
\[
  I(\tau):=\{u\in E(Q):\tau \subset  H_u\}.
\]
For a fixed $I \subset  E(Q)$, the union of the cells
whose support is exactly $I$ is the locally closed set
\[
  U_I:=H_I\setminus
  \bigcup_{u\in E(Q)\setminus I}H_u.
\]
Each connected component of $U_I$ is a relatively open
$q$-dimensional polyhedral cell. Indeed, in the universal
cover $V$, the arrangement contains all the affine Coxeter
hyperplanes
\[
  \langle e_j-e_i,x\rangle=\ell,
  \qquad
  0\leq i<j\leq n,\quad \ell\in\ZZ.
\]
Their relatively open faces are convex and bounded, and
the additional hyperplanes subdivide those faces into
further relatively open convex polyhedra.

No such face is identified with itself under a nonzero
translation by $L$. To see this, suppose $x,y$ lie in
one Coxeter face and $y-x\in L$.
For every $i<j$, the integer
$\langle e_j-e_i,y-x\rangle$ must vanish, because the
values of $\langle e_j-e_i,\cdot\rangle$ on a single
Coxeter face are either equal to one integer or lie
strictly between two consecutive integers.
Hence
\[
  \langle e_j-e_i,y-x\rangle=0
  \qquad\text{for every }i<j,
\]
which implies $x=y$ in $V$.
Thus the quotient components are open $q$-balls.

Writing $c_I$ for the number of cells with support exactly
$I$, we obtain
\[
  \chi_c(U_I)=(-1)^q c_I,
\]
where $\chi_c$ denotes compactly supported Euler
characteristic.

By inclusion--exclusion,
\[
  \chi_c(U_I)
  =\sum_{J \subset  E(Q)\setminus I}
    (-1)^{|J|}\chi_c(H_{I\cup J}).
\]
The space $H_{I\cup J}$ is a disjoint union of compact
tori of dimension $n-{\rm rk}_{\Mscr}(I\cup J)$, so
\[
  \chi_c(H_{I\cup J})=
  \begin{cases}
    \mu_{\Mscr}(I\cup J),
      &{\rm rk}_{\Mscr}(I\cup J)=n,\\
    0,
      &{\rm rk}_{\Mscr}(I\cup J)<n.
  \end{cases}
\]
Consequently,
\[
  c_I = (-1)^q\chi_c(U_I) =
  \sum_{\substack{
    J \subset  E(Q)\setminus I\\
    {\rm rk}_{\Mscr}(I\cup J)=n
  }}
  (-1)^{|J|-q}\mu_{\Mscr}(I\cup J)\\
  ={\rm T}_{\Mscr/I}(1,0).
\]
The last equality follows directly from the subset
expansion of the arithmetic Tutte polynomial of the
contraction: at $(x,y)=(1,0)$, only subsets $J$
satisfying ${\rm rk}_{\Mscr}(I\cup J)=n$ contribute,
and the exponent of $-1$ is
\[
  |J|-
  \bigl(
    {\rm rk}_{\Mscr}(I\cup J)
    -{\rm rk}_{\Mscr}(I)
  \bigr)
  =|J|-q.
\]
Hence ${\rm T}_{\Mscr/I}(1,0)$ counts precisely the
cells of the quotient arrangement having support $I$.

\item[Step 2:]
The natural covering $V/2L\to V/L$ has degree $2^n$.
Since the cells of the quotient arrangement are open balls,
every cell has exactly $2^n$ distinct lifts. 

Fix such a cell $\tau$ with support $I$.
The ordinary Coxeter hyperplanes are contained in the
refined arrangement, so equality of the fractional
parts of any two coordinates $x_i,x_j$ is constant on
${\rm relint}(\tau)$.
Partition the indices in $S$ into blocks
$B_1,\ldots,B_k$ according to equality of their
fractional parts.
A Coxeter root $e_i-e_j\in E_S$ belongs to $I$
precisely when $i$ and $j$ lie in the same block.
Therefore
\[
  {\rm rk}_{\Mscr}(I\cap E_S)
  =\sum_{a=1}^k(|B_a|-1)=s-k,
\]
and hence $k = s - {\rm rk}_{\Mscr}(I\cap E_S)$.

We now count the lifts of $\tau$ that are deleted by
the translates of $Z_S$.
By Lemma~\ref{lem:liftedSpoints} and the lifting
calculation in the proof of Theorem~\ref{thm:FVectorCI},
a cell of the ordinary Coxeter arrangement whose
$S$-coordinates have $k$ distinct fractional-part
blocks has exactly
\[
  2^{n+1-s}k
\]
deleted lifts to $V/2L$.

This calculation also applies to the refined cell
$\tau$. Indeed, every cell of the refined arrangement
is contained in a cell of the ordinary Coxeter
arrangement. Moreover, the boundaries of the
translates of $Z_S$ are contained in the ordinary
Coxeter hyperplanes. Consequently, the deletion
status of a lift is unchanged under refinement.

It follows that precisely $2^{n+1-s}k$ of the
$2^n$ lifts of $\tau$ are deleted.
Hence the number of surviving lifts is
\[
\begin{aligned}
  2^n-2^{n+1-s}k
  &=
  2^n-2^{n+1-s}
  \bigl(
    s-{\rm rk}_{\Mscr}(I\cap E_S)
  \bigr).
\end{aligned}
\]
In particular, this number depends only on
$I\cap E_S$, not on the additional normals in $E(Q)$.

\item[Step 3:]
Every cell having support $I$ has dimension
$n-{\rm rk}_{\Mscr}(I)$.
By Step 1, there are ${\rm T}_{\Mscr/I}(1,0)$
such cells in $V/L$; by Step 2, each has
\[
  2^n-2^{n+1-s}
  \bigl(
    s-{\rm rk}_{\Mscr}(I\cap E_S)
  \bigr)
\]
surviving lifts to $V/2L$.
We therefore obtain
\[
\begin{aligned}
 f_d\bigl(\Delta_{n,S,2L}(Q)\bigr)
  =
  \sum_{\substack{
    I \subset  E(Q)\\
    {\rm rk}_{\Mscr}(I)=n-d
  }}
  \left(
    2^n-2^{n+1-s}
    \bigl(
      s-{\rm rk}_{\Mscr}(I\cap E_S)
    \bigr)
  \right)
  {\rm T}_{\Mscr/I}(1,0).
\end{aligned}
\]
\end{enumerate}

Finally, the covering $V/H\to V/2L$ has degree
$[2L:H]$.
Both the refined arrangement and the deleted set are
$2L$-periodic, and every cell has exactly
$[2L:H]$ lifts under this covering.
Multiplying the preceding identity by $[2L:H]$ proves the desired formula
\[
 f_d\bigl(\Delta_{n,S,H}(Q)\bigr)
  =
  [2L:H]
  \sum_{\substack{
    I \subset  E(Q)\\
    {\rm rk}_{\Mscr}(I)=n-d
  }}
  \left(
    2^n-2^{n+1-s}
    \bigl(
      s-{\rm rk}_{\Mscr}(I\cap E_S)
    \bigr)
  \right)
  {\rm T}_{\Mscr/I}(1,0).
  \qedhere
\]
\end{proof}

\begin{remark}\label{rem:RefinedCoxeterEmptyQ}
We show that Theorem~\ref{thm:refined_Coxter_F_vect}
recovers Theorem~\ref{thm:FVectorCI} when $Q=\varnothing$.
In this case,
\[
    E(Q)=E=\{e_j-e_i:0\leq i<j\leq n\},
\]
and $\Mscr$ is the graphic matroid of the complete graph on $n+1$ vertices, with all arithmetic multiplicities $\mu_{\Mscr}(I)$ equal to $1$. 
We note that only flats $I \subset  E$ contribute to the formula, since for a non-flat $I$ the matroid $\Mscr/I$ contains a loop and we get ${\rm T}_{\Mscr/I}(1,0)=0$.
The flats $I_P$ of rank $n-d$ are in one-to-one correspondence with partitions
\[
    P = \{B_0,\ldots,B_d\}
\]
of $[n]$ into $d+1$ nonempty blocks, where $I_P$
consists of the edges joining elements of the same block.
Contracting $I_P$ gives the graphic matroid of the complete
multigraph on $d+1$ vertices. Since parallel elements
do not affect the Tutte evaluation at $(1,0)$, we obtain
\[
    {\rm T}_{\Mscr/I_P}(1,0)=d!.
\]

Let $k(P)$ denote the number of blocks of $P$ meeting $S$.
Then
\[
    {\rm rk}_{\Mscr}(I_P\cap E_S)=s-k(P).
\]
Thus Theorem~\ref{thm:refined_Coxter_F_vect} becomes
\[
    f_d
    =[2L:H]\,d!
    \left(
        2^n\left\{\begin{matrix}n+1\\d+1\end{matrix}\right\}
        -2^{n+1-s}\sum_P k(P)
    \right),
\]
where the sum runs over all partitions of $[n]$
into $d+1$ nonempty blocks.

As in the proof of Theorem~\ref{thm:FVectorCI},
the sum $\sum_P k(P)$ counts partitions with a
distinguished block meeting $S$.
Choosing this block to have size $\ell$ gives
\[
    \sum_P k(P)
    =
    \sum_{\ell=1}^{n+1-d}
    \left[
        \binom{n+1}{\ell}
        -\binom{n+1-s}{\ell}
    \right]
    \left\{\begin{matrix}n+1-\ell\\d\end{matrix}\right\}.
\]
Substituting yields
\[
\begin{aligned}
    f_d &=[2L:H]\,d!\Bigg(
        2^n\left\{\begin{matrix}n+1\\d+1\end{matrix}\right\} 
        -2^{n+1-s}
        \sum_{\ell=1}^{n+1-d}
        \left[
            \binom{n+1}{\ell}
            -\binom{n+1-s}{\ell}
        \right]
        \left\{\begin{matrix}n+1-\ell\\d\end{matrix}\right\}
    \Bigg)\\
    &=[2L:H]\,d!\,g_{n,s}(d),
\end{aligned}
\]
which is precisely the formula of
Theorem~\ref{thm:FVectorCI}.
\end{remark}

\begin{example}[Refined $f$-vectors for $n=2$]
\label{ex:refinedCoxeterFVector2}
We apply Theorem~\ref{thm:refined_Coxter_F_vect} to the two
refinements in Example~\ref{ex:refinedCoxeterArrangement1}.
Here $n=2$, $S=[2]=\{0,1,2\}$, and $E_S=E$. For
$I \subset  E(Q)$, the lifting factor becomes
\[
  w(I):=2^2-2^{3-3}\bigl(3-{\rm rk}_{\Mscr}(I\cap E)\bigr)
       =1+{\rm rk}_{\Mscr}(I\cap E).
\]
Thus the formula reads
\[
 f_d\bigl(\Delta_{2,[2],H}(Q)\bigr)
  =[2L:H]\sum_{\substack{I \subset  E(Q)\\
                 {\rm rk}_{\Mscr}(I)=2-d}}
      w(I)\,{\rm T}_{\Mscr/I}(1,0).
\]
For brevity, write the Coxeter normals, in their original
$\ZZ^3$ coordinates, as
\[
 a := e_1-e_0=(-1,1,0), \quad 
 b := e_2-e_0=(-1,0,1), \quad
 c := e_2-e_1=(0,-1,1).
\]
The two refined ground sets are
\[
 \begin{aligned}
 E(Q_1')&=\{a,b,c,p,q\}, \quad \text{where }
   &u& := (3,-1,-2), &v& := (2,0,-2) = -2b,\\
 E(Q_2')&=\{a,b,c,r,t\}, \quad \text{where }
   &r& := (3,-2,-1), &t& :=(1,0,-1) = -b.
 \end{aligned}
\]
In the second equality we have cleared the denominators
in $Q_2'$ as prescribed by
Definition~\ref{def:refinedCoxeterArrangement}.
In each calculation $\Mscr$ denotes the matroid over $\ZZ$
associated with the displayed ground set.
We evaluate the contractions using
\[
 {\rm T}_{\Mscr/I}(1,0)
 =\sum_{\substack{J \subset  E(Q)\setminus I\\
           {\rm rk}_{\Mscr}(I\cup J)=2}}
   (-1)^{|J|-(2-{\rm rk}_{\Mscr}(I))}
   \mu_{\Mscr}(I\cup J).
\]
All torsion multiplicities can be computed directly from
the displayed triples: for a rank-two subset they are
obtained by taking the greatest common divisor of the
$2\times2$ minors of the matrix of its vectors in $\ZZ^3$.

We start with the cellulation $(\psi_1, Q_1')$. Since $H_1 = 2L$, the index $[2L:H_1]$ is $1$.
The following table lists every support $I \subset  E(Q_1')$
with nonzero contraction evaluation, grouped by cell dimension.
The last column is its contribution to $f_d$.
\[
\begin{array}{|c|c|c|c|c|}
\hline
 d&I&{\rm T}_{\Mscr/I}(1,0)&w(I)&w(I){\rm T}_{\Mscr/I}(1,0)\\ \hline
 0&\{c,p\}&2&2&4\\
 0&\{c,q\}&1&2&2\\
 0&\{a,p,q\}&1&2&2\\
 0&\{a,b,c,p,q\}&1&3&3\\ \hline
 1&\{a\}&2&2&4\\
 1&\{c\}&4&2&8\\
 1&\{p\}&4&1&4\\
 1&\{q\}&2&1&2\\
 1&\{b,q\}&1&2&2\\ \hline
 2&\varnothing&8&1&8\\
 \hline
\end{array}
\]
Summing the contributions in the table, we obtain the $f$-vector
\[ 
    f \bigl(\Delta_{2,[2],H_1}(Q_1')\bigr)=(11,20,8).
\]

Similarly, for the cellulation $(\psi_2, Q_2')$ we have $[2L : H_2] = 2$. We first compute the contributions on $V/2L$:
\[
\begin{array}{|c|c|c|c|c|}
\hline
 d&I&{\rm T}_{\Mscr/I}(1,0)&w(I)&w(I){\rm T}_{\Mscr/I}(1,0)\\ \hline
 0&\{c,r\}&2&2&4\\
 0&\{b,r,t\}&1&2&2\\
 0&\{a,b,c,r,t\}&1&3&3\\ \hline
 1&\{a\}&1&2&2\\
 1&\{c\}&3&2&6\\
 1&\{r\}&4&1&4\\
 1&\{b,t\}&2&2&4\\ \hline
 2&\varnothing&6&1&6\\
 \hline
\end{array}
\]
Multiplying by $[2L:H_2]=2$, we obtain the $f$-vector
\[
    f \bigl(\Delta_{2,[2],H_2}(Q_2')\bigr)=(18,32,12).
\]

Moreover, the two $f$-vectors obtained this way from Theorem~\ref{thm:refined_Coxter_F_vect} are indeed the $f$-vectors of the complexes in Figure~\ref{fig:refined-coxeter-arrangements}.
\end{example}

\section{Examples and explicit formulae}\label{sec:Examples}

In this section, we apply our theorems to explicit examples. We start with several classical examples to show how the theory works in the easiest cases, see Section \ref{sec:ClassicalExamplesD01}. Namely, in dimension $0$, our polar decomposition is related to Bernstein's theorem and in dimension $1$ it is related to pair-of-pants decompositions. Next, we show how our theory works for toric varieties in Section \ref{sec:ToricVarieties}. We then turn to Fermat and Dwork hypersurfaces in Section \ref{sec:FermatAndDwork}. 
Using the angle functor constructed in Section \ref{sec:AngleSetsCI} and our formula for the $f$-vector of the homotopy colimit in Theorem~\ref{thm:CWComplexStructureTriangulatedDiagram}, we are able to explicitly determine their $f$-vectors. Finally, we work out the CW complex structure explicitly for the K3-surface arising from the strictly semistable model associated to the Ap\'{e}ry-Fermi pencil, see Section \ref{sec:CYMonodromy}. This allows us to interpret the monodromy of integral affine structures in the SYZ-program explicitly in terms of angle sets.

\begin{remark}[Notation for $f$-vectors]
    Throughout this section, we will calculate many $f$-vectors, so we fix our notation here. Consider a function $g:\ZZ\to \ZZ$. Let $\ell\geq{0}$ be an integer. Recall that  ${\rm sh}^{\ell}(g)$ is defined by ${\rm sh}^{\ell}(g)(m)=g(m-\ell)$. 

    Consider the complex $\Delta_{n,S,H}$ from Theorem \ref{thm:FVectorCI} with $s = |S|$ and $k=[2L : H]$. We will write $f(\Delta_{n,s,k})$ for the $f$-vector in Theorem \ref{thm:FVectorCI}. That is, 
    \begin{equation*}
        f(\Delta_{n,s,k})(m)=km!g_{n,s}(m).
    \end{equation*}

    We view this as a function $f(\Delta_{n,s,k}):\ZZ\to \ZZ$ that maps all but finitely many $m$ to $0$. 
    Here we define $g_{n,1}$ to be the indicator function at $m=0$, so that $f(\Delta_{n,s,k})$ gives $k$ many points. 
    If $k=1$,  then we will also write $f(\Delta_{n,s})$ for this function. If $k=1$ and $s=0$, then we write $f(\Delta_{n})$ for this function.  If this function is zero outside a certain range, then we will also truncate the vector to the region of interest. For instance, we have $f(\Delta_{2,2})(1)=f(\Delta_{2,2})(2)=2$ and $f(\Delta_{2,2})(m)=0$ for other $m$, so we will also write $f(\Delta_{2,2})=(2,2)$ or $(2,2,0)$ if the other vectors in question have the same size.  
\end{remark}

\subsection{Classical examples: dimension 
\texorpdfstring{$0$}{0} 
and \texorpdfstring{$1$}{1}}\label{sec:ClassicalExamplesD01}

We begin with two classical examples to highlight their polar nature.
\begin{example}[Bernstein's theorem]
We consider the setup in Bernstein's theorem. That is, we have $n$ finite sets $S_{i}\subset \ZZ^{n}$ with convex hulls $\mathcal{C}_{i}$. For generic elements $f_{i}\in \CC[S_{i}]$, Bernstein's theorem states that the cardinality of
\begin{equation*}
    X = V(f_{1},...,f_{n})
\end{equation*}
is finite, and it is equal to the mixed volume of the Newton polytopes $\mathcal{N}_{i}=\mathcal{C}_{i}$ of the $f_{i}$. We view this as a statement about the number of connected components of $X=V(f_{1},...,f_{n})$\footnote{It will follow from our setup that $X$ is generically smooth, see \cite[Theorem 3.5]{HHR24}.}. We choose a set of rationally smooth $f_{i}\in K[S_{i}]$ which intersect transversely\footnote{Here we mean that their tropicalizations are rationally smooth, and that their tropicalizations intersect transversely.} and we calculate their cardinality. Note that this setup can be used to prove Bernstein's theorem, see \cite{HuberSturmfels1995} for a similar approach. 

As in Proposition \ref{pro:RationalManifoldStableIntersection}, we find that the initial ideals at $w\in \trop(X)$ can be generated  by binomials of the form
\begin{equation*}
    g=x^{n_{i}}_{i}-1.
\end{equation*}
The number of zeros in the initial degeneration is then exactly equal to $\prod_{i=1}^{n}n_{i}$. Note that this is equal to the cardinality of the corresponding angle set. This is the polar decomposition in the zero-dimensional case.   
\end{example}

\begin{example}[Pair-of-pants decompositions of Riemann surfaces]
Let $X$ be a smooth algebraic curve over $\CC$ with smooth completion $\overline{X}$. We write $g=g(\overline{X})$ for its genus and $k = |\overline{X}\backslash X|$. If $2 - 2g - k < 0$, then we say that $X$ is hyperbolic. In this case, there exists a pair-of-pants decomposition. That is, we can find a decomposition of $X$ into copies of $\mathbb{P}^{1}\backslash\{0,1,\infty\}$ and $S^{1}$.

Tropically, this can be deduced as follows. We consider a flat family $\mathcal{X}\to U\subset \Spec(\CC[t])$ of Riemann surfaces  such that the fiber over $t=1$ is $X$, and such that the curve for $t=0$ is semistable, with irreducible components birational to $\mathbb{P}^{1}$. Note for instance that a smooth tropical curve leads to a family of this form. Our polar decomposition is then exactly a pair-of-pants decomposition of $X$. See Figure \ref{fig:smooth_cubic2} for the case of an elliptic curve with $9$ punctures. 

   \begin{figure}[ht]
        \centering
        \begin{minipage}{.45\textwidth}
        \scalebox{0.7}{
        \begin{tikzpicture}
            \coordinate (a) at (-2, 2); 
            \node at (-2.7,2.2) [xscale=1, yscale=1]{$(-2,2)$};
            \coordinate (a1) at (-2, 4);
            \coordinate (a2) at (-4, 0);
            
            \coordinate (b) at (-1, 2);
            \node at (-0.3,2.2) [xscale=1, yscale=1]{$(-1,2)$};
            \coordinate (b1) at (-1, 4);
            
            \coordinate (c) at (0, 1);
            \node at (-0.7,1) [xscale=1, yscale=1]{$(0,1)$};
            
            \coordinate (d) at (1, 1);
            \node at (1.7, 1) [xscale=1, yscale=1]{$(1,1)$};
            
            \coordinate (e) at (2, 2);
            \node at (2.5, 2.3) [xscale=1, yscale=1]{$(2,2)$};
            
            \coordinate (e1) at (4, 2);            
            \coordinate (e2) at (2, 4);
            
            \coordinate (f) at (0, 0);
            \node at (-0.7, 0) [xscale=1, yscale=1]{$(0,0)$};
            \coordinate (f1) at (-2,-2);
            
            \coordinate (g) at (1, 0);
            \node at (1.7, 0) [xscale=1, yscale=1]{$(1,0)$};
            
            \coordinate (h) at (2, -1);
            \node at (2.4, -0.7) [xscale=1, yscale=1]{$(2,-1)$};
            \coordinate (h1) at (4, -1);
            
            \coordinate (i) at (2, -2);
            \node at (2.4, -2.4) [xscale=1, yscale=1]{$(2,-2)$};
            \coordinate (i1) at (4, -2);
            \coordinate (i2) at (0, -4);

            \draw [fill=black] (a) circle (2.5pt);
            \draw [fill=black] (b) circle (2.5pt);
            \draw [fill=black] (c) circle (2.5pt);
            
            \draw [fill=black] (d) circle (2.5pt);
            \draw [fill=black] (e) circle (2.5pt);
            \draw [fill=black] (f) circle (2.5pt);
            
            \draw [fill=black] (g) circle (2.5pt);
            \draw [fill=black] (h) circle (2.5pt);
            \draw [fill=black] (i) circle (2.5pt);

             \draw (a) -- (b);
             \draw (a) -- (a1);
             \draw (a) -- (a2);
             
             \draw (b) -- (c);
             \draw (b) -- (b1);
             
             \draw (c) -- (d);
             \draw (d) -- (e);
             \draw (e) -- (e1);
             \draw (e) -- (e2);
             
             \draw (c) -- (f);
             \draw (f) -- (g);
             \draw (f) -- (f1);
             \draw (h) -- (g);
             \draw (d) -- (g);
             \draw (h) -- (h1);
             \draw (i) -- (h);
             \draw (i) -- (i1);
             \draw (i) -- (i2);
        \end{tikzpicture}}
        \end{minipage}
        \begin{minipage}{.45\textwidth}
            \scalebox{0.8}{
            \begin{tikzpicture}[x=0.75pt,y=0.75pt,yscale=-1,xscale=1]
            
            \draw [line width=0.75]    (356.97,121.19) -- (371.57,107.98) ;
            \draw [line width=0.75]    (342.37,107.98) -- (356.97,94.78) ;
            \draw [line width=0.75]    (313.18,107.98) -- (342.37,107.98) ;
            \draw [line width=0.75]    (298.59,94.78) -- (313.18,107.98) ;
            \draw [line width=0.75]    (283.99,107.98) -- (259.66,107.98) ;
            \draw [line width=0.75]    (264.53,90.38) -- (279.12,90.38) ;
            \draw [line width=0.75]    (256.86,23.97) -- (264.53,90.38) ;
            \draw [line width=0.75]    (249.86,23.97) -- (245.07,94.78) ;
            \draw [line width=0.75]    (192.63,147.34) -- (245.07,94.78) ;
            \draw [line width=0.75]    (196.45,154.12) -- (259.66,107.98) ;
            \draw [line width=0.75]    (283.99,107.98) -- (298.59,121.19) ;
            \draw [line width=0.75]    (298.59,121.19) -- (298.59,147.6) ;
            \draw [line width=0.75]    (313.67,122.07) -- (313.67,148.48) ;
            \draw [line width=0.75]    (313.67,122.07) -- (342.86,122.07) ;
            \draw [line width=0.75]    (342.86,122.07) -- (342.86,148.48) ;
            \draw [line width=0.75]    (342.86,148.48) -- (313.67,148.48) ;
            \draw [line width=0.75]    (313.18,160.81) -- (342.37,160.81) ;
            \draw [line width=0.75]    (342.37,160.81) -- (356.97,174.01) ;
            \draw [line width=0.75]    (356.97,147.6) -- (371.57,160.81) ;
            \draw [line width=0.75]    (356.97,121.19) -- (356.97,147.6) ;
            \draw [line width=0.75]    (356.97,174.01) -- (356.97,200.42) ;
            \draw [line width=0.75]    (376.43,178.41) -- (376.43,191.62) ;
            \draw [line width=0.75]    (376.43,191.62) -- (463.83,198.41) ;
            \draw [line width=0.75]    (371.66,213.98) -- (463.83,206.41) ;
            \draw [fill={rgb, 255:red, 175; green, 170; blue, 170 }  ,fill opacity=1 ]   (249.86,23.97) .. controls (250.86,20.97) and (255.86,19.97) .. (256.86,23.97) ;
            \draw    (249.86,23.97) .. controls (248.86,23.97) and (254.86,25.97) .. (256.86,23.97) ;
            \draw [fill={rgb, 255:red, 175; green, 170; blue, 170 }  ,fill opacity=1 ]   (463.83,198.41) .. controls (467.83,199.41) and (466.83,206.41) .. (463.83,206.41) ;
            \draw    (463.83,198.41) .. controls (461.83,200.41) and (461.83,205.41) .. (463.83,206.41) ;
            \draw    (192.63,147.34) .. controls (194.54,146.44) and (197.39,152.78) .. (196.45,154.12) ;
            \draw    (192.63,147.34) .. controls (191.71,151.14) and (193.42,153.59) .. (196.45,154.12) ;
            \draw    (238.86,101.61) .. controls (243.24,101.17) and (251.16,106.4) .. (250.68,114.32) ;
            \draw  [dash pattern={on 0.84pt off 2.51pt}]  (238.86,101.61) .. controls (238.37,109.97) and (245.44,113.59) .. (250.68,114.32) ;
            \draw    (245.68,85.32) .. controls (249.21,81.6) and (259.42,81.03) .. (263.8,85) ;
            \draw  [dash pattern={on 0.84pt off 2.51pt}]  (245.68,85.32) .. controls (249.27,91.37) and (260.09,89.28) .. (263.8,85) ;
            \draw    (271.83,90.38) .. controls (276.47,95.65) and (276.01,102.05) .. (271.83,107.98) ;
            \draw  [dash pattern={on 0.84pt off 2.51pt}]  (271.83,90.38) .. controls (267.23,96.09) and (267.23,102.7) .. (271.83,107.98) ;
            \draw    (305.88,101.38) .. controls (303.72,108.86) and (298.86,113.7) .. (293.99,115.9) ;
            \draw  [dash pattern={on 0.84pt off 2.51pt}]  (305.88,101.38) .. controls (296.91,102.26) and (293.02,107.1) .. (293.99,115.9) ;
            \draw    (328.75,107.98) .. controls (333.4,113.26) and (332.43,115.9) .. (329.24,122.07) ;
            \draw  [dash pattern={on 0.84pt off 2.51pt}]  (328.75,107.98) .. controls (324.16,113.7) and (324.64,116.78) .. (329.24,122.07) ;
            \draw    (298.78,132.06) .. controls (302.31,128.34) and (309.07,127.79) .. (313.45,131.75) ;
            \draw  [dash pattern={on 0.84pt off 2.51pt}]  (298.78,132.06) .. controls (302.37,138.11) and (309.74,136.04) .. (313.45,131.75) ;
            \draw    (328.75,148.48) .. controls (333.4,153.76) and (333.4,156.84) .. (327.78,160.81) ;
            \draw  [dash pattern={on 0.84pt off 2.51pt}]  (328.26,148.48) .. controls (322.21,152.44) and (322.21,157.28) .. (327.29,160.81) ;
            \draw    (363.3,153.32) .. controls (361.13,160.8) and (356.27,165.64) .. (351.4,167.84) ;
            \draw  [dash pattern={on 0.84pt off 2.51pt}]  (363.3,153.32) .. controls (354.32,154.2) and (350.43,159.04) .. (351.4,167.84) ;
            \draw    (343.07,133.06) .. controls (346.6,129.34) and (352.38,129.55) .. (356.75,133.51) ;
            \draw  [dash pattern={on 0.84pt off 2.51pt}]  (343.07,133.06) .. controls (346.66,139.1) and (353.05,137.8) .. (356.75,133.51) ;
            \draw    (357.18,184.12) .. controls (360.71,180.4) and (371.86,180.48) .. (376.24,184.44) ;
            \draw  [dash pattern={on 0.84pt off 2.51pt}]  (357.18,184.12) .. controls (360.77,190.17) and (372.53,188.73) .. (376.24,184.44) ;
            \draw    (351.69,99.62) .. controls (360.51,101.39) and (365.31,106.23) .. (367.26,111.51) ;
            \draw  [dash pattern={on 0.84pt off 2.51pt}]  (351.69,99.62) .. controls (351.2,107.99) and (362.01,110.77) .. (367.26,111.51) ;
            \draw [line width=0.75]    (245.63,200.34) -- (298.59,147.6) ;
            \draw [line width=0.75]    (249.45,207.12) -- (313.18,160.81) ;
            \draw    (245.63,200.34) .. controls (247.54,199.44) and (250.39,205.78) .. (249.45,207.12) ;
            \draw    (245.63,200.34) .. controls (244.66,202.15) and (248,208.49) .. (249.45,207.12) ;
            \draw    (290.86,155.61) .. controls (295.24,155.17) and (303.96,159.37) .. (303.47,167.29) ;
            \draw  [dash pattern={on 0.84pt off 2.51pt}]  (290.86,155.61) .. controls (290.37,163.97) and (298.23,166.55) .. (303.47,167.29) ;
            \draw [line width=0.75]    (304.53,252.99) -- (356.97,200.42) ;
            \draw [line width=0.75]    (308.45,260.12) -- (371.66,213.98) ;
            \draw    (304.63,253.34) .. controls (306.54,252.44) and (309.39,258.78) .. (308.45,260.12) ;
            \draw [fill={rgb, 255:red, 175; green, 170; blue, 170 }  ,fill opacity=1 ]   (304.63,253.34) .. controls (303.66,255.15) and (307,261.49) .. (308.45,260.12) ;
            \draw    (339.86,217.61) .. controls (344.24,217.17) and (352.16,220.68) .. (351.68,228.61) ;
            \draw  [dash pattern={on 0.84pt off 2.51pt}]  (339.86,217.61) .. controls (339.37,225.97) and (346.43,227.87) .. (351.68,228.61) ;
            \draw [line width=0.75]    (371.57,160.81) -- (462.83,165.84) ;
            \draw [line width=0.75]    (376.43,178.41) -- (462.83,173.84) ;
            \draw [fill={rgb, 255:red, 175; green, 170; blue, 170 }  ,fill opacity=1 ]   (462.83,165.84) .. controls (466.28,167.09) and (465.83,173.84) .. (462.83,173.84) ;
            \draw    (462.83,165.84) .. controls (460.83,167.84) and (460.83,172.84) .. (462.83,173.84) ;
            \draw    (386.93,161.75) .. controls (391.58,167.03) and (390.11,171.82) .. (385.93,177.75) ;
            \draw  [dash pattern={on 0.84pt off 2.51pt}]  (386.93,161.75) .. controls (382.33,167.47) and (381.33,172.47) .. (385.93,177.75) ;
            \draw [line width=0.75]    (290.86,23.97) -- (298.59,94.78) ;
            \draw [line width=0.75]    (283.86,23.97) -- (279.12,90.38) ;
            \draw [fill={rgb, 255:red, 175; green, 170; blue, 170 }  ,fill opacity=1 ]   (283.86,23.97) .. controls (284.86,20.97) and (289.86,19.97) .. (290.86,23.97) ;
            \draw    (283.86,23.97) .. controls (282.86,23.97) and (288.86,25.97) .. (290.86,23.97) ;
            \draw    (280.68,82.32) .. controls (284.21,78.6) and (292.42,80.03) .. (296.8,84) ;
            \draw  [dash pattern={on 0.84pt off 2.51pt}]  (280.68,83.32) .. controls (284.27,89.37) and (293.09,88.28) .. (296.8,84) ;
            \draw [line width=0.75]    (376.43,86.62) -- (463.83,93.41) ;
            \draw [line width=0.75]    (371.57,107.98) -- (463.73,100.41) ;
            \draw [fill={rgb, 255:red, 175; green, 170; blue, 170 }  ,fill opacity=1 ]   (463.83,93.41) .. controls (467.83,94.41) and (466.73,100.41) .. (463.73,100.41) ;
            \draw    (463.83,93.41) .. controls (460.35,93.45) and (461.73,99.41) .. (463.73,100.41) ;
            \draw    (391.06,87.94) .. controls (395.7,93.22) and (395.01,100.47) .. (390.83,106.41) ;
            \draw  [dash pattern={on 0.84pt off 2.51pt}]  (391.06,87.94) .. controls (386.46,93.66) and (386.23,101.12) .. (390.83,106.41) ;
            \draw [line width=0.75]    (368.86,22.97) -- (376.43,86.62) ;
            \draw [line width=0.75]    (361.86,22.97) -- (356.97,94.78) ;
            \draw [fill={rgb, 255:red, 175; green, 170; blue, 170 }  ,fill opacity=1 ]   (361.86,22.97) .. controls (362.86,19.97) and (367.86,18.97) .. (368.86,22.97) ;
            \draw    (361.86,22.97) .. controls (360.86,22.97) and (366.86,24.97) .. (368.86,22.97) ;
            \draw    (357.86,79.97) .. controls (361.39,76.25) and (371.09,76.33) .. (375.47,80.29) ;
            \draw  [dash pattern={on 0.84pt off 2.51pt}]  (357.86,79.97) .. controls (361.45,86.02) and (371.76,84.58) .. (375.47,80.29) ;
            \draw  [color={rgb, 255:red, 0; green, 0; blue, 0 }  ,draw opacity=0 ][fill={rgb, 255:red, 175; green, 170; blue, 170 }  ,fill opacity=0.41 ] (283.86,23.97) -- (290.86,23.97) -- (298.59,94.78) -- (313.18,107.98) -- (328.75,107.98) -- (342.37,107.98) -- (356.97,94.78) -- (361.86,22.97) -- (368.86,22.97) -- (376.43,86.62) -- (464.8,93.41) -- (463.73,100.41) -- (371.57,107.98) -- (356.97,121.19) -- (342.86,122.07) -- (313.67,122.07) -- (313.67,148.48) -- (342.86,148.48) -- (343.83,122.07) -- (356.97,121.19) -- (357.94,147.6) -- (372.54,160.81) -- (463.8,165.84) -- (463.8,173.84) -- (377.41,178.41) -- (377.41,191.62) -- (464.8,198.41) -- (464.8,206.41) -- (372.64,213.98) -- (309.42,260.12) -- (305.5,252.99) -- (357.94,200.42) -- (357.94,174.01) -- (342.37,160.81) -- (313.18,160.81) -- (249.45,207.12) -- (245.63,200.34) -- (298.59,147.6) -- (298.59,121.19) -- (283.99,107.98) -- (259.66,107.98) -- (196.45,154.12) -- (192.63,147.34) -- (245.07,94.78) -- (249.86,23.97) -- (256.86,23.97) -- (264.53,90.38) -- (279.12,90.38) -- cycle ;
            \draw  [dash pattern={on 0.84pt off 2.51pt}]  (209.86,129.61) .. controls (209.37,137.97) and (211.44,138.59) .. (216.68,139.32) ;
            \draw    (209.86,129.61) .. controls (214.24,129.17) and (217.16,131.4) .. (216.68,139.32) ;
            \draw    (248.68,42.32) .. controls (252.21,38.6) and (254.3,38.36) .. (258.68,42.32) ;
            \draw  [dash pattern={on 0.84pt off 2.51pt}]  (248.68,42.32) .. controls (252.27,48.37) and (254.97,46.61) .. (258.68,42.32) ;
            \draw  [dash pattern={on 0.84pt off 2.51pt}]  (282.68,43.32) .. controls (286.27,49.37) and (288.97,47.61) .. (292.68,43.32) ;
            \draw    (282.68,43.32) .. controls (286.21,39.6) and (288.3,39.36) .. (292.68,43.32) ;
            \draw    (360.41,44.87) .. controls (363.95,41.15) and (367.09,41.33) .. (371.47,45.29) ;
            \draw  [dash pattern={on 0.84pt off 2.51pt}]  (360.41,44.87) .. controls (364,50.92) and (367.76,49.58) .. (371.47,45.29) ;
            \draw    (258.86,186.61) .. controls (263.24,186.17) and (266.96,187.37) .. (266.47,195.29) ;
            \draw  [dash pattern={on 0.84pt off 2.51pt}]  (258.86,186.61) .. controls (258.37,194.97) and (261.23,194.55) .. (266.47,195.29) ;
            \draw    (315.86,240.61) .. controls (320.24,240.17) and (323.96,241.37) .. (323.47,249.29) ;
            \draw  [dash pattern={on 0.84pt off 2.51pt}]  (315.86,240.61) .. controls (315.37,248.97) and (318.23,248.55) .. (323.47,249.29) ;
            \draw    (392.93,192.75) .. controls (397.58,198.03) and (396.11,206.82) .. (391.93,212.75) ;
            \draw  [dash pattern={on 0.84pt off 2.51pt}]  (392.93,192.75) .. controls (388.33,198.47) and (387.33,207.47) .. (391.93,212.75) ;
            \draw    (433.93,195.75) .. controls (438.58,201.03) and (439.11,202.82) .. (434.93,208.75) ;
            \draw  [dash pattern={on 0.84pt off 2.51pt}]  (433.93,195.75) .. controls (429.33,201.47) and (430.33,203.47) .. (434.93,208.75) ;
            \draw    (432.93,163.75) .. controls (437.58,169.03) and (436.11,169.37) .. (431.93,175.3) ;
            \draw  [dash pattern={on 0.84pt off 2.51pt}]  (432.93,163.75) .. controls (428.33,169.47) and (427.33,170.02) .. (431.93,175.3) ;
            \draw    (435.93,91.03) .. controls (441.93,93.03) and (440.93,99.03) .. (435.93,103.03) ;
            \draw  [dash pattern={on 0.84pt off 2.51pt}]  (435.93,91.03) .. controls (431.33,96.75) and (431.33,97.74) .. (435.93,103.03) ;
            \end{tikzpicture}}
            \end{minipage}    
        \caption{Left: a tropically smooth elliptic curve. Right: a pair of pants decomposition of the intersection of an elliptic curve with the torus in $\PP^2$.}
        \label{fig:smooth_cubic2}
    \end{figure}
     
\end{example}

\begin{example}[Explicit pair-of-pants decomposition for $g=2$]\label{exa:PairOfPantsg2}
We now use our angle machinery to find an explicit pair-of-pants decomposition. That is, we write down a diagram of angle spaces associated to a semistable model, and we find the corresponding CW complex. This will also highlight the necessity of the homotopy colimit in our setup. 

Consider the genus-$2$ curve over $\CC(t)$ given by the equation
\begin{align*}
        y^2&=f(x),\\ 
        f(x)&=x(x-t)(x-1)(x-1-t)(x-2)(x-2-t).
    \end{align*}
    By compactifying the corresponding affine scheme\footnote{Explicitly, take the coordinate rings $R[x,y]/(y^2-f(x))$ and $R[1/x,y/x]/((y/x)^2-g(1/x))$, where $g=f/x^6$ and $R=\CC[t]$. Their spectra glue naturally, and this defines the semistable model. }, we obtain a natural 
    semistable model $\mathcal{X}\to U$. Its special fiber consists of two irreducible components intersecting in $3$ points. The irreducible components of the special fiber are isomorphic to $\mathbb{P}^{1}$. The corresponding open logarithmic initial degenerations $Z^{\log,\circ}_{\tau}$ from the KN-functor can be viewed as $\mathbb{P}^{1}\backslash\{0,1,\infty\}$. The dual intersection graph $\Sigma$ of $\mathcal{X}_{s}$ is the theta graph, which consists of $2$ vertices and $3$ edges.

    \begin{figure}[ht]
    \begin{center}
    \scalebox{0.4}{
    \definecolor{red_color}{rgb}{1,0,0}
    \definecolor{gray_color}{rgb}{0.3,0.3,0.3}
    \definecolor{blue_color}{rgb}{0.3,0.3,1}
    \definecolor{green_color}{rgb}{0.1,0.7,0.3}

    \begin{tikzpicture}[line cap=round,line join=round,x=1cm,y=1cm]
    \clip(4,7) rectangle (15,15);
    \fill[line width=2.8pt,color=red_color,fill=red_color,fill opacity=0.5] (10,8) -- (13,11) -- (13,8) -- cycle;
    \fill[line width=2.8pt,color=red_color,fill=red_color,fill opacity = 0.5] (7,14) -- (10,14) -- (7,11) -- cycle;
    
    \draw [line width=2.8pt] (7,11)-- (10,14);
    \draw [line width=2.8pt] (10,8)-- (13,11);
    
    \draw [line width=2.8pt] (10,14)-- (13,14);
    \draw [line width=2.8pt] (7,8)-- (10,8);
    \draw [line width=2.8pt] (10,8)-- (13,8);
    \draw [line width=2.8pt] (13,8)-- (13,14);
    \draw [line width=2.8pt] (7,8)-- (7,14);
    \draw [line width=2.8pt] (7,14)-- (10,14);
    
    \begin{scriptsize}

    \draw [fill=blue_color] (7,8) circle    (6pt);
    \draw [fill=blue_color] (13,14) circle  (6pt);
    \draw [fill=blue_color] (13,8) circle   (6pt);
    \draw [fill=blue_color] (7,14) circle   (6pt);
    \draw [fill=green_color] (10,8)  circle (6pt);
    \draw [fill=green_color] (10,14) circle (6pt);
    \draw [fill=yellow] (7,11) circle (6pt);
    \draw [fill=yellow] (13,11) circle (6pt);
    
    \draw (5.5,8) node [anchor=north west][inner sep=0.75pt] [xscale = 2, yscale = 2] {$-\pi$};
    \draw (13.5,8) node [anchor=north west][inner sep=0.75pt] [xscale = 2, yscale = 2] {$\pi$};
    \draw (6,14.5) node [anchor=north west][inner sep=0.75pt] [xscale = 2, yscale = 2] {$\pi$};
    
    \end{scriptsize}
    
    \end{tikzpicture}
     }
    \end{center}
    \caption{\label{fig:AngleSet3PointsinP1} The Kato-Nakayama angle set of $V(1+x+y)\subset \Spec(\CC[x^{\pm},y^{\pm}])$. Note that this variety is isomorphic to the hyperplane complement $\mathbb{P}^{1}\backslash\{0,1,\infty\}$.  }
    \end{figure}

    We view $\Sigma$ as an abstract tropical manifold.  We can model the topology of the open components  $\mathbb{P}^{1}\backslash\{0,1,\infty\}$ with the very affine linear space $V(1+x+y)\subset \Spec(\CC[\ZZ^2])$. The angle set of the corresponding Kato-Nakayama space can be found in Figure \ref{fig:AngleSet3PointsinP1}. These give the spaces over the vertices of $\Sigma$.  
    We assign the angle set $S^{1}$ to the three edges. Note that the comparison maps are the natural inclusions. 

    We use Theorem~\ref{thm:CWComplexStructureTriangulatedDiagram} to calculate the $f$-vector of our CW complex. This yields 
    \begin{align*}
        \mathfrak{f}_{0}&=2\cdot f(\Delta_{2,3})+3\cdot f(\Delta_{2,2})=(6,12,4)+(6,6,0)=(12,18,4)\\
        \mathfrak{f}_{1}&=6\cdot {\rm sh}^{1}(f(\Delta_{2,2}))=(0,12,12).
        \end{align*}
        We thus find 
        \begin{equation*}
            \mathfrak{f}(\Delta(\mathcal{F}^{\triangle}))=\mathfrak{f}_{0}+\mathfrak{f}_{1}=(12,30,16).
        \end{equation*}
        We conclude that the Euler characteristic is 
        \[
\chi(\Delta)=12-30+16=-2,
        \]
        which matches the classical $2-2g=2-4=-2$. 

        Note that we cannot take a simple colimit here, even though all of the comparison maps are cofibrations. Indeed, the corresponding $f$-vector would be $(3,6,4)$, with Euler characteristic $1$. The homotopy colimit thus prevents these types of anomalies. The problem here lies in the fact that the three circles overlap, which does not happen in the homotopy colimit by construction.  
\end{example}

\begin{example}
    More generally, let $\mathcal{X}_{s}$ be a  semistable curve of genus $g$ over $\CC$ whose irreducible components are $\mathbb{P}^{1}$. We write $\Sigma$ for the dual intersection graph and $\Pscr(\Sigma)$ for the induced poset.
    We will assume that $\Sigma$ is maximal, in the sense that every vertex has valency $3$. Note that $\beta_{1}(|\Sigma|)=g=1/2\cdot v+1$, so that $2-2g=-v$. As before, we can turn this into a cellulated tropical diagram. Namely, we associate the compactified angle set $\Theta^{\log}$ of $V(1+x+y)\subset \Spec(\CC[\ZZ^{2}])$ to every vertex. On the edges, we have a copy of $S^{1}$, triangulated using two edges and two vertices. The inclusion maps are as in Example \ref{exa:PairOfPantsg2}.  

    We now calculate the Euler characteristic of $|\Delta(\mathcal{F}^{\triangle})|$ using our angle sets. Let $v=|\Sigma_{0}|$ and $e=|\Sigma_{1}|$ be the number of vertices and edges respectively.

    We will need two $f$-vectors here: $f(\Delta_{2,3})=(3,6,2)$ and $f(\Delta_{2,2})=(2,2,0)$. Note that these describe the $f$-vectors of our cellulations of the angle sets of the spaces $V(1+x+y)$ and $V(1+x)$ in $\Spec(\CC[\ZZ^{2}])$ respectively. 
    
    We have
    \begin{align*}
        \mathfrak{f}_{0}&=v\cdot f(\Delta_{2,3})+e\cdot f(\Delta_{2,2})=(3v+2e,6v+2e,2v)\\
        \mathfrak{f}_{1}&=2e\cdot {\rm sh}^{1}(f(\Delta_{2,2}))=(0,4e,4e).
    \end{align*}

    This yields
    \begin{equation*}
        \mathfrak{f}(\Delta(\mathcal{F}^{\triangle})) = (3v+2e,6v+6e,2v+4e).
    \end{equation*}
    Note that $\chi(\Delta(\mathcal{F}^{\triangle})) = 3v+2e-(6v+6e)+2v+4e=-v$, as expected. 

\end{example}

\subsection{Toric varieties}\label{sec:ToricVarieties}

In this section, we combine Theorem \ref{thm:FVectorCI} with Theorem~\ref{thm:CWComplexStructureTriangulatedDiagram} to obtain a formula for the $f$-vector of the cellulated tropical diagram associated to a smooth toric variety. 
We will also apply our machinery to a compactification of a hyperplane complement, and see how it also gives a CW complex structure on its compactification through initial degenerations. 
\begin{remark}
    The homotopy colimit presentation seems to lead to too many auxiliary cells for toric varieties. We will give some examples to show how to modify the diagram.   
\end{remark}

\begin{definition}
    Consider the standard real torus $(S^{1})^{n}$ of dimension $n$. The cellulation of $(S^{1})^{n}$ in Theorem \ref{thm:FVectorCI} has $f$-vector $f(\Delta_{n})$, where 
    \begin{equation*}
        f(\Delta_{n})(d) = 2^{n} \ d! \ \left\{ \begin{matrix} n + 1 \\ d + 1 \end{matrix} \right\}.
    \end{equation*}
    The number of vertices is $2^{n}$, and the number of maximal cells is $2^{n}n!$. 
\end{definition}
\begin{example}
    The first five $f$-vectors are 
    \begin{align*}
        f(\Delta_{0})&=(1),\\
        f(\Delta_{1})&=(2,2),\\
        f(\Delta_{2})&=(4,12,8),\\
        f(\Delta_{3})&=(8,56,96,48),\\
        f(\Delta_{4})&=(16,240,800,960,384).
    \end{align*}
\end{example}
\begin{definition}\label{def:ToricFVectorsDefinition}
    Let $(X,\Sigma,D)$ be a rational tropical manifold with $X=(\CC^{\times})^{n}$. 
    The cosedentarity of an element $z=(\tau,\sigma)\in \Pscr(\Sigma_{D})$ is $n-{\rm dim}(\sigma)$.
    We denote it by $c(z)$. 
\end{definition}

\begin{corollary}\label{thm:FVectorToricVariety}
Consider a rational tropical manifold $((\CC^{\times})^{n},\Sigma,D)$. 
    We have
    \[
    \mathfrak{f}_{k}=\sum_{(z_{0}<...<z_{k})} 
    {\rm sh}^{k} (f(\Delta_{c(z_{k})})) \quad \text{and} \quad \mathfrak{f}(\Delta(\mathcal{F}^{\triangle})) = \sum_{k = 0}^{\dim(\Sigma)} \mathfrak{f}_{k},
    \]
    where the first sum runs over the chains $(z_0 < \dots <z_{k})$ of length $k$ in the poset $\Pscr(\Sigma_D)$.

\end{corollary}
\begin{proof}
    This follows from Theorem \ref{thm:FVectorCI} and Theorem~\ref{thm:CWComplexStructureTriangulatedDiagram}. Indeed, the local spaces are exactly real tori of dimension $c(z_{k})$ for a chain $(z_{0}<...<z_{k})$ of length $k$.  
\end{proof}

\begin{example}\label{exa:ProjectiveLineExample}
    Suppose that $\Sigma=\Delta$ is the fan for $\mathbb{P}^{1}$. This gives the CW complex structure on the right in Figure \ref{fig:SphereTriangulations}. We have five cells in $\Pscr(\mathbb{T}\mathbb{P}^{1})$, so we find using Theorem \ref{thm:FVectorToricVariety} that 
    \begin{equation*}
        \mathfrak{f}_{0}=3f(\Delta_{1,0})+2f(\Delta_{0,0})=(6,6,0)+(2,0,0)=(8,6,0). 
    \end{equation*}

\begin{figure}[ht]
    \centering
    \scalebox{0.6}{
    \begin{tikzpicture}[
      line cap=round,
      line join=round,
      outline/.style={draw=black!75, line width=0.8pt},
      hidden/.style={draw=black!65, dashed, line width=0.7pt},
      shade/.style={fill=gray!55, fill opacity=0.4}]

    \begin{scope}[shift={(-5, 3)}]
      \draw[black, line width=0.5pt] (-5,0) -- (-1,0);
      \draw[black, line width=0.5pt] (4,0) -- (12,0);
      \foreach \x in {-5, -1, 4, 12, 8}
      \fill[blue!75!black] (\x,0) circle (2.2pt);
    \end{scope}
    
    \begin{scope}[shift={(-8,0)}]
    \begin{scope}[rotate=90]
      \fill[shade] (0,2.65) -- (-1.2,0) -- (1.2,0) -- cycle;
      \fill[shade] (-1.2,0) -- (0,-2.65) -- (1.2,0) -- cycle;
    
      \draw[outline] (0,2.65) -- (-1.2,0) -- (0,-2.65)
        -- (1.2,0) -- cycle;
      \draw[hidden] (-1.2,0) arc[start angle=180,end angle=0,
        x radius=1.2, y radius=0.38];
      \draw[outline] (1.2,0) arc[start angle=0,end angle=-180,
        x radius=1.2, y radius=0.38];
    
      \foreach \x/\y in {0/2.65, -1.2/0, 1.2/0, 0/-2.65}
        \fill[black!75] (\x,\y) circle (2pt);
    \end{scope}
    \end{scope}

    \begin{scope}[shift={(8,0)}]
    \begin{scope}[rotate=90]
      \fill[shade] (0,10) -- (-1.25,7.2) -- (1.25,7.2) -- cycle;
      \fill[shade] (-1.25,7.2) rectangle (1.25,4.8);
      \fill[shade] (-1.25,4.8) rectangle (1.25,2.3);
      \fill[shade] (-1.25,2.3) -- (0,-0.5) -- (1.25,2.3) -- cycle;

      \draw[outline] (0,10) -- (-1.25,7.2);
      \draw[outline] (0,10) -- (1.25,7.2);
      \draw[outline] (-1.25,7.2) -- (-1.25,2.3);
      \draw[outline] (1.25,7.2) -- (1.25,2.3);
      \draw[outline] (-1.25,2.3) -- (0,-0.5) -- (1.25,2.3);

      \foreach \y in {7.2,4.8,2.3} {
        \draw[hidden] (-1.25,\y) arc[start angle=180,end angle=0,
          x radius=1.25, y radius=0.42];
        \draw[outline] (1.25,\y) arc[start angle=0,end angle=-180,
          x radius=1.25, y radius=0.42];
      }
    
      \foreach \x/\y in {
        0/10, -1.25/7.2, 1.25/7.2,
        -1.25/4.8, 1.25/4.8, -1.25/2.3, 1.25/2.3, 0/-0.5
      }
        \fill[black!75] (\x,\y) circle (2pt);
    \end{scope}
    \end{scope}
    
    \end{tikzpicture}}
    \caption{\label{fig:SphereTriangulations} Two CW complex structures of a $2$-sphere (bottom row) obtained from two different CW complex structures on $\mathbb{T}\mathbb{P}^{1}$ (top row).}
\end{figure}

    We similarly calculate 
    \begin{equation*}
\mathfrak{f}_{1}=4\cdot{\rm sh}^{1}(f(\Delta_{1,0}))=(0,8,8).
    \end{equation*}
    Together, this yields
    \begin{equation*}
        \mathfrak{f}(\Delta(\mathcal{F}^{\triangle})) = (8,14,8).
    \end{equation*}
    We have $\chi(\mathfrak{f}(\Delta(\mathcal{F}^{\triangle})))=2$, as expected. 
    
    We can also endow $\mathbb{T}\mathbb{P}^{1}$ with the more economic CW complex structure on the left of Figure \ref{fig:SphereTriangulations}. The corresponding tropical diagram then gives  
    \begin{align*}
        \mathfrak{f}_{0}&=f(\Delta_{1,0})+2f(\Delta_{0,0})=(4,2,0)\\
        \mathfrak{f}_{1}&=2\cdot{}{\rm sh}^{1}(f(\Delta_{1,0}))=(0,4,4).
    \end{align*}
    This gives the $f$-vector $(4,6,4)$ of the space on the left in Figure \ref{fig:SphereTriangulations}. %
\end{example}
\begin{example}\label{exa:P2DiagramOfSpaces}
    Consider $X=(\CC^{\times})^{2}$ with compactification $\PP^{2}(\CC)$. This is a toric variety with character lattice $M=\{(a,b,c)\in\ZZ^{3}:a+b+c = 0\}$. See Example \ref{exa:P2ExampleHC}. We have $19$ cells in $\Pscr(\mathbb{T}\mathbb{P}^{2})$. We first calculate 
    \begin{align*}
        \mathfrak{f}_{0}&=7\cdot f(\Delta_{2,0})+9\cdot f(\Delta_{1,0})+3\cdot f(\Delta_{0,0})=(49,102,56,0,0),\\
        \mathfrak{f}_{1}&=30\cdot {\rm sh}^{1}(f(\Delta_{2,0}))+12\cdot{}{\rm sh}^{1}(f(\Delta_{1,0}))=(0,168,504,336,0)\\
        \mathfrak{f}_{2}&=24\cdot {\rm sh}^{2}(f(\Delta_{2,0}))=(0,0,96,288,192).
    \end{align*}
Together, this yields
\begin{equation*}
    \mathfrak{f}=\sum \mathfrak{f}_{k}=(49,270,656,624,192),
\end{equation*}
so that $\chi(\Delta(\mathcal{F}^{\triangle}))=3$. As in the previous example, we can also simply endow $\mathbb{T}\mathbb{P}^{2}$ with the regular CW complex structure with $7$ cells. Here we have $3$ vertices, $3$ edges and one two-dimensional cell. Our formula here gives 
\begin{align*}
    \mathfrak{f}_{0}&=f(\Delta_{2})+3f(\Delta_{1})+3f(\Delta_{0})=(13,18,8,0,0)\\
    \mathfrak{f}_{1}&=(0,36,84,48,0),\\
    \mathfrak{f}_{2}&=(0,0,24,72,48).
\end{align*}
We thus obtain the $f$-vector 
\begin{equation*}
    \mathfrak{f}=(13,54,116,120,48).
\end{equation*}
As before, we can check that $\chi(\mathfrak{f}(\Delta(\mathcal{F}^{\triangle})))=3$. In general, we can use this to construct a more economic CW complex structure on toric varieties.

\end{example}

\begin{example}\label{exa:P2DiagramOfSpaces4H}
    We now consider $X=\mathbb{P}^{2}$ embedded in $\mathbb{P}^{3}$ through 
    \begin{equation*}
        [X:Y:Z]\mapsto [X:Y:Z:-X-Y-Z].
    \end{equation*}
    In other words, it is given as the vanishing set of the homogeneous polynomial 
    \begin{equation*}
        F=X+Y+Z+W.
    \end{equation*}

    Note that this diagram gives another CW complex structure on $\mathbb{P}^{2}$. We will calculate its $f$-vector here. This will require several non-toric $f$-vectors however. We will need the following:  %
    \begin{align*}
        f(\Delta_{3,4})&=(7,42,60,24,0),\\ %
        f(\Delta_{3,3})&=(6,30,36,12,0),\\ %
        f(\Delta_{3,2})&=(4,12,8,0,0),\\
        f(\Delta_{2,3})&=(3,6,2,0,0),\\ %
        f(\Delta_{2,2})&=(2,2,0,0,0),\\ %
        f(\Delta_{1,1})&=(1,0,0,0,0). %
    \end{align*}
    
We now have 
\begin{align*}
    \mathfrak{f}_{0}&=f(\Delta_{3,4})+4f(\Delta_{3,3})+6f(\Delta_{3,2})+4f(\Delta_{2,3})+12f(\Delta_{2,2})+6f(\Delta_{1,1})=(97,282,260,72,0),\\
    \mathfrak{f}_{1}&=48\cdot{\rm sh}^{1}(f(\Delta_{3,2}))+8\cdot{\rm sh}^{1}(f(\Delta_{3,3}))+24\cdot{\rm sh}^{1}(f(\Delta_{2,2}))=(0,288,864,672,96),\\
    \mathfrak{f}_{2}&=48\cdot{\rm sh}^{2}(f(\Delta_{3,2}))=(0,0,192,576,384).
\end{align*}
We then find 
\begin{equation*}
    \mathfrak{f}=\mathfrak{f}_{0}+\mathfrak{f}_{1}+\mathfrak{f}_{2}=(97,570,1316,1320,480).
\end{equation*}
Its Euler characteristic is $3$. Perhaps as expected, choosing a larger initial diagram for $\mathbb{P}^{2}$ here in terms of the embedding $\mathbb{P}^{2}\to \mathbb{P}^{3}$ yields a CW complex with more cells.

\end{example}

\subsection{Fermat and Dwork hypersurfaces}\label{sec:FermatAndDwork}

In this section, we focus on Fermat and Dwork hypersurfaces. We will give a completely explicit formula for the $f$-vector of the tropical complex associated to a Fermat hypersurface. Moreover, we will compute the $f$-vector associated to the Dwork pencil that gives a K3 surface. This in particular makes the construction in Section \ref{sec:TriangulationsRCIM} explicit.   

We start by defining Fermat and Dwork hypersurfaces.
\begin{definition}\label{def:ProjectiveFermat}
    Consider a collection of independent vectors $a_{0},...,a_{n}\in\mathbb{N}^{n+1}\subset \ZZ^{n+1}$ all of total degree $m$. The corresponding projective Fermat hypersurface is $V(F)\subset \mathbb{P}^{n}_{\CC}$, where 
    \begin{equation*}
        F=\sum_{i=0}^{n}\chi^{a_{i}}.
    \end{equation*}
    Let $\langle F \rangle\subset \CC[\ZZ^{n+1}]$ be the ideal generated by $F$. 
    We write $I=\CC[\Xi] \cap \langle F \rangle$ for the corresponding ideal in $\CC[\Xi]$, where $\Xi := {\rm ker}(\mathds{1}_{n+1}) \subset \ZZ^{n+1}$. We define the corresponding toric Fermat hypersurface in $\Spec(\CC[\Xi])$ to be $X=V(I)$. If $a_{i} = me_{i}$ for all $i$, then we call $X$ the standard toric Fermat hypersurface of degree $m$. We denote the corresponding polynomial by $F_{m}$.  
\end{definition}
\begin{remark}
    Note that the standard toric Fermat hypersurface by definition comes with a natural simplex presentation. 
\end{remark}

\begin{definition}
Let $K=\CC(t)$, $a_{0},...,a_{n} \in\mathbb{N}^{n+1}$ and $F$ be as in Definition \ref{def:ProjectiveFermat}. We define the corresponding projective Dwork hypersurface in $\mathbb{P}^{n}_{K}$ to be $V(G)$, where  
    \begin{equation*}
        G= t F + \prod_{i=0}^{n}\chi^{e_{i}}.
    \end{equation*}
    The corresponding toric Dwork hypersurface is $X=V(I)$, where $I=K[\Xi]\cap (G)$.  
\end{definition}

\begin{example}
Suppose $n=3$ and $m=4$. The projective Fermat hypersurface of degree $d$ is cut out in $\PP^3$ by 
\begin{equation*}
    F=X^{4}+Y^{4}+Z^{4}+W^{4}.
\end{equation*}
Here $X=\chi^{e_{0}}$, $Y=\chi^{e_{1}}$, $Z=\chi^{e_{2}}$ and $W=\chi^{e_{3}}$. 
In the torus, we can for instance dehomogenize with respect to $W$ and write 
\begin{equation*}
    f=x^4+y^4+z^4+1.
\end{equation*}
Here $x=\chi^{e_{0}-e_{3}}$, $y=\chi^{e_{1}-e_{3}}$ and $z=\chi^{e_{2}-e_{3}}$. We have $(f)=\CC[\Xi]\cap \langle F \rangle$. Similarly, we have 
\begin{equation*}
    G = t F + XYZW 
\end{equation*}
and 
\begin{equation*}
    g = t(1 + x^4 + y^4 + z^4) + xyz.
\end{equation*}
Note that the zero sets of $F$ and $G$ are both K3-surfaces. We will provide CW complexes for these in the next sections.
\end{example}

\subsubsection{An $f$-vector formula for Fermat hypersurfaces}

We now give an explicit formula for the $f$-vector of a CW complex structure on Fermat hypersurfaces. To that end, we will need a precise description of the polyhedral structure on the tropicalization, and of the order complex of the corresponding tropical poset. We subdivide cells of this order complex according to their sedentarity, depth and chain length.   

We start by defining the corresponding local tropical model. 
\begin{definition}[Linear hypersurfaces]
Let $\Xi={\rm ker}(\mathds{1})\subset \ZZ^{n+1}$. 
Let $F=\sum_{i=0}^{n}\chi^{e_{i}}$ with $X_{\rm proj}=V(F)\subset \mathbb{P}^{n}_{\CC}$. Let $X=V(I)$, where $I = \CC[\Xi]\cap \langle F \rangle$. The tropicalization of $X$ 
is the \emph{standard linear tropical hypersurface}. Let $\Delta := {\rm conv}(\{e_{i}:0\leq{i}\leq{n}\})\subset \RR^{n+1}$. The projection of $\Delta$ in $\RR^{n+1}/\mathbb{R}\mathds{1}_{n+1}$ is a simplex. 
We endow $\trop(X)$ with the natural polyhedral structure whose cells are dual to cells in $\Delta$. We write $\Sigma$ for this polyhedral complex. For a choice $(X,\Sigma,D)$ of a divisor $D$ at infinity, we will write
$\Pscr(\Sigma_{D})$ for the corresponding tropical poset.   
\end{definition}

\begin{definition}[Sedentarity, depth and length]
    Let $z=(\tau,\sigma)\in \Pscr(\overline{\Sigma})$. The sedentarity of $z$ is ${\rm dim}(\sigma)$. The depth of $z$ is ${\rm dim}(\tau)$. Similarly, let $(z_{0}<...<z_{k})$ be a cell of dimension $k$ in the order complex $\Delta(\Pscr(\overline{\Sigma}))$. Then the sedentarity of the chain is the sedentarity of $z_{k}$, and the depth of the chain is the depth of $z_{k}$. The type of a cell $(z_{0}<...<z_{k})$ is the triple $({\rm dim}(\sigma),{\rm dim}(\tau),k)$. Note that ${\rm dim}(\sigma)\leq {\rm dim}(\tau)$ by definition of the tropical poset. 
\end{definition}

\begin{theorem}\label{thm:OrderComplexTropicalizationFVector}
The number of $k$-dimensional cells in the order complex of the standard
tropical linear hypersurface whose type is $(s,d,k)$ is
\[
    N_n(s,d,k)
    =
    \frac{(n+1)!}{(n+1-d)!\,(d-s)!\,s!}
    \sum_{j=0}^{k}(-1)^{k-j}\binom{k}{j}(2j+1)^{d-s}, \quad \text{for} \quad \begin{matrix}
    0 \leq  s \leq d \leq n-1 \\ 
    0 \leq k \leq n - 1.
\end{matrix}
\]
\end{theorem}

\begin{proof}
Write $[n]=\{0,\ldots,n\}$. A cell of the standard tropical
hyperplane is determined by the set $A \subset [n]$ of coordinates
where the minimum is attained. Since the minimum must be attained at
least twice, we must have $|A|\geq 2$. 
If the cell has dimension $d$, then $|A| = n+1-d$. Thus the number of $d$-dimensional cells in the standard tropical hyperplane is
\[
    \binom{n+1}{n+1-d}=\binom{n+1}{d}.
\]

Fix such a cell $\tau=\tau_A$, and set
$B=[n]\setminus A$. Then the recession cone of $\tau_A$ is generated by the $d$ rays
indexed by $B$. An $s$-dimensional face $\sigma$ of this cone is determined by choosing exactly $s$ of those rays. Therefore, the number of
pairs $(\tau,\sigma)$ with $\dim(\tau)=d$ and $\dim(\sigma)=s$ is
\[
    \binom{n+1}{d}\binom{d}{s}
    =
    \frac{(n+1)!}{(n+1-d)!\,(d-s)!\,s!}.
\]

We next count, for one fixed pair $z=(\tau_A,\sigma_C)$, the chains
\[
    (z_0<z_1<\cdots<z_k),
\]
where $z_{k}=z$. 
Here $C \subset  B$ indexes the rays generating $\sigma_C$, and $|C| = s$. Now set $E := B \setminus C$, so $|E| = d - s$. Every element below $z$ in the tropical poset has the form
\[
    \bigl(\tau_{A\cup D},\sigma_{C\cup S}\bigr),
\]
where $D$ and $S$ are disjoint subsets of $E$. Indeed, adding indices
to $A$ passes to a face of $\tau_A$, while adding rays to $C$ enlarges
the cone $\sigma_C$. The two sets must be disjoint so that the enlarged
cone remains in the recession cone of the smaller cell.

For each index in $E$, there are three choices: it can belong to $D$, belong to $S$, or belong to neither. 
In the poset, the first two choices are incomparable and both lie below the third choice. This is exactly the face poset of a segment: its two endpoints lie below the whole segment. 
Making these choices independently for all indices in $E$ shows that the interval of elements below $z$ is the face poset of a $(d-s)$-dimensional cube, with $z$ corresponding to the whole cube.

It remains to count chains of length $k$ whose maximal element is the whole cube. Let $r = d-s$, and let $Q$ be an $r$-dimensional cube. For $j\geq 0$, consider
weak chains of faces
\[
    F_0 \subset F_1 \subset \cdots \subset F_j = Q,
\]
where repetitions are allowed. In each coordinate, the face is either
fixed at $0$, fixed at $1$, or free. Along the chain, a coordinate is
either free from the beginning, or is fixed at one of the two endpoints
and becomes free at one of the $j$ steps. Thus there are $1+2j$ choices
per coordinate, and the number of weak chains is $(2j+1)^r$.

To count strict chains of length $k$, use inclusion-exclusion on the
events $F_{i-1}=F_i$, for $1\leq i\leq k$. If equality is required at
$k-j$ specified steps, collapsing those repeated faces leaves a weak
chain with $j$ steps, counted by $(2j+1)^r$. There are $\binom{k}{j}$
ways to choose the remaining $j$ steps. Hence the number of strict
chains ending at $Q$ is
\[
    \sum_{j=0}^{k}(-1)^{k-j}\binom{k}{j}(2j+1)^r
    =
    \sum_{j=0}^{k}(-1)^{k-j}\binom{k}{j}(2j+1)^{d-s}.
\]
This count is the same for every pair $(\tau,\sigma)$ of dimensions
$(d,s)$. Multiplying it by the number of such pairs gives the stated
formula.
\end{proof}

\begin{definition}\label{def:TropicalComplexFermat}
Let $m \geq{1}$. We construct a tropical complex for the angle functor $\mathcal{F}_{\Theta,CI}$ associated to the Fermat hypersurface $V(F_{m})$. There is one vertex, and it comes with a natural simplex presentation $\psi:\ZZ^{n+1}\to \ZZ^{n}$. We take $Q=\varnothing$, which gives the cellulation data on the interior. The projections at the infinite strata are compatible with the projection maps by Proposition \ref{prop:CompatibilityQuotientTriangulation}, so we take the corresponding cells at infinity.
\end{definition}

\begin{theorem}\label{thm:FermatFVector}
Consider the tropical complex for $(X,\Sigma,D)$ from Definition \ref{def:TropicalComplexFermat}. Its geometric realization is homotopy equivalent to $V(F_{m})$. We have \begin{equation*}
\mathfrak{f}_{k}=\mathrm{sh}^{k} \left(\sum_{0 \leq s \leq d \leq n-1} N_{n}(s,d,k)f(\Delta_{n-s,n+1-d,m^{n-s}})\right).
\end{equation*} 
The $f$-vector for the tropical complex is  
\begin{equation*}
    \mathfrak{f}=\sum_{k=0}^{n-1} \mathfrak{f}_{k}.
\end{equation*}
\end{theorem}
\begin{proof}
This follows from Theorem \ref{thm:MainTheoremPolarDecomposition} and Theorem~\ref{thm:CWComplexStructureTriangulatedDiagram}. 
\end{proof}

\begin{example}\label{exa:FermatK3Example}

   Consider the homogeneous polynomial $F=X^4+Y^4+Z^4+W^4$ with hypersurface $X=V(F)$. This defines a K3 surface. The $f$-vectors for the angle spaces we need are 
    \begin{align*}
        f(\Delta_{3,4,64})&=(448, 2688, 3840, 1536,0),\\
        f(\Delta_{3,3,64})&=(384, 1920, 2304, 768,0),\\
        f(\Delta_{3,2,64})&=(256, 768, 512, 0,0),\\
        f(\Delta_{2,3,16})&=(48, 96, 32, 0,0),\\
        f(\Delta_{2,2,16})&=(32, 32, 0, 0,0),\\
        f(\Delta_{1,1,4})&=(4, 0, 0, 0,0).
    \end{align*}

    Note that the first three $f$-vectors are associated to $V(1+x^4+y^4+z^4)$, $V(1+x^4+y^4)$ and $V(1+x^4)$ inside $(\CC^{\times})^{3}$. The fourth and fifth $f$-vectors are  associated to the angle sets of $V(1+x^4+y^4)$ and $V(1+x^4)$ inside $(\CC^{\times})^{2}$. Finally, the last is associated to the angle set of $V(1+x^4)$ in $\CC^{\times}$. 
Note that the coefficients $N_{3}(s,d,k)$ were already found in Example \ref{exa:P2DiagramOfSpaces4H}. 
        We calculate 
        {\footnotesize
            \begin{align*}
        \mathfrak{f}_{0}&=f(\Delta_{3,4,64})+4f(\Delta_{3,3,64})+6f(\Delta_{3,2,64})+4f(\Delta_{2,3,16})+12f(\Delta_{2,2,16})+6f(\Delta_{1,1,4})=(4120, 15744, 16256, 4608, 0),\\
        \mathfrak{f}_{1}&=8\cdot{\rm sh}^{1}(f(\Delta_{3,3,64}))+48\cdot{\rm sh}^{1}(f(\Delta_{3,2,64}))+24\cdot{\rm sh}^{1}(f(\Delta_{2,2,16}))=(0, 16128, 52992, 43008, 6144),\\
                \mathfrak{f}_{2}&=48\cdot{\rm sh}^{2}(f(\Delta_{3,2,64}))=(0, 0, 12288, 36864, 24576).
            \end{align*}
        }
    Summing these yields the $f$-vector
    \begin{equation*}
        \mathfrak{f}=(4120, 31872, 81536, 84480, 30720),
    \end{equation*}
    with Euler characteristic $24$.

\end{example}

In general, for any $n$, one only has to calculate the $N_{n}(s,d,k)$ once, since this leads to the corresponding coefficients for all Fermat hypersurfaces of degree $m$. We record the corresponding result here for the Fermat surface.
\begin{corollary}\label{cor:FermatSurfaceFVector}
    Consider the Fermat surface $X=V(X^{m}+Y^{m}+Z^{m}+W^{m})$ for $m\geq{1}$. Then the individual contributions of the $f$-vector are given by
    \begin{align*} 
        \mathfrak{f}_{0}
        &= \left( 55 m^3 + 36 m^2 + 6 m, \; 234 m^3 + 48 m^2, \; 252 m^3 + 8 m^2, \; 72 m^3, \; 0 \right),\\
        \mathfrak{f}_{1}
        &= \left( 0, \; 240 m^3 + 48 m^2, \; 816 m^3 + 48 m^2, \; 672 m^3, \; 96 m^3 \right),\\
        \mathfrak{f}_{2}
        &=\left( 0, \; 0, \; 192 m^3, \; 576 m^3, \; 384 m^3 \right).
    \end{align*}
    The resulting $f$-vector is 
    \begin{equation*}
\mathfrak{f} = \left( 55 m^3 + 36 m^2 + 6 m, \; 474 m^3 + 96 m^2, \; 1260 m^3 + 56 m^2, \; 1320 m^3, \; 480 m^3 \right).
    \end{equation*}
    In particular, it is on the order of $O(m^3)$. 
\end{corollary}

\begin{corollary}
    The $f$-vector of the tropical Fermat CW complex is on the order of $O(m^{n})$. 
\end{corollary}
\begin{proof}
    This follows in general since the topological data of the tropicalization is fixed, but the number of cells per angle space grows as $O(m^{n})$.
\end{proof}
\subsubsection{Dwork hypersurfaces}

For the Dwork family, we will focus on one example for explicit formulas. The general case requires a special polyhedral structure to remove the singularities. Our example will again be a K3 surface. In this case, the cellulations of the different adjacent cells will not be compatible, so we will have to refine the Coxeter triangulation as we outlined in Section \ref{sec:TriangulationsRCIM}.

\begin{figure}[ht]
    \centering
    \includegraphics[width=0.35\linewidth]{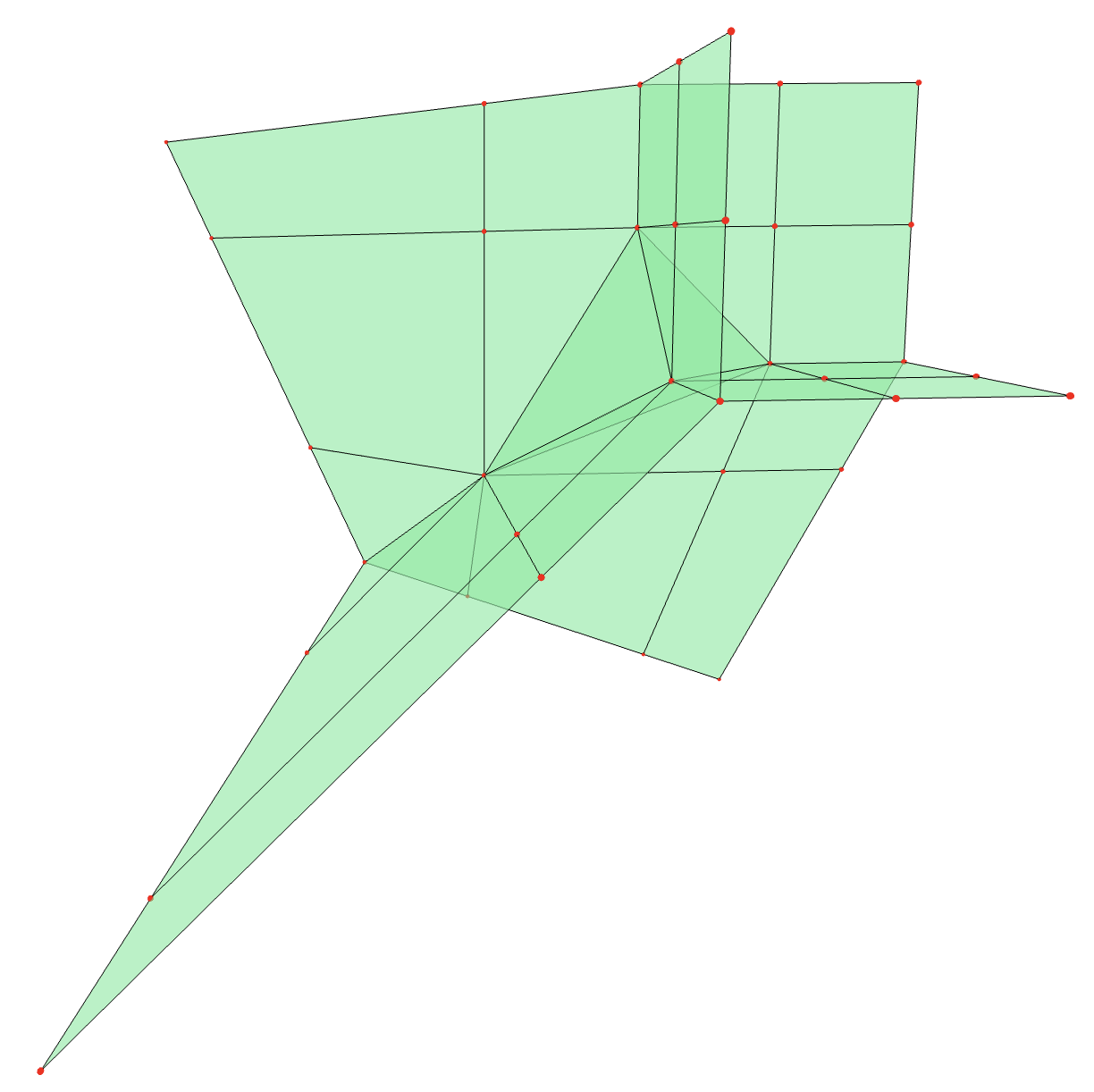}
    \caption{The tropicalization of the Dwork family in Example~\ref{exa:DworkFamilyK3}. The subdivision ensures that we obtain a regular scheme. }
    \label{fig:DworkK3_subdivided}
\end{figure}

\begin{example}\label{exa:DworkFamilyK3}
    Consider the degree-four polynomial 
    \begin{equation*}
        F=XYZW+t(X^4+Y^4+Z^4+W^4).
    \end{equation*}
    Dehomogenizing with respect to $W$, we obtain the polynomial 
    \begin{equation*}
        f=xyz+t(x^4+y^4+z^4+1).
    \end{equation*}
    The induced subdivision of the Newton polytope $\mathcal{N}(f)$ consists of four simplices. This gives rise to four vertices in the tropicalization, see Figure~\ref{fig:DworkK3_subdivided}. These vertices are 
    \begin{align*}
        v_{x} = (1,0,0),\quad
        v_{y} = (0,1,0),\quad
        v_{z} = (0,0,1),\quad \text{and} \quad
        v_{w} = (-1,-1,-1).
    \end{align*}
    The initial forms at these vertices are
    \begin{align*}
        \mathrm{in}_{v_{x}}(f) &= xyz+y^4+z^4+1\\
        \mathrm{in}_{v_{y}}(f) &= xyz+x^4+z^4+1,\\
        \mathrm{in}_{v_{z}}(f) &= xyz+x^4+y^4+1,\\
        \mathrm{in}_{v_{w}}(f) &= xyz+x^4+y^4+z^4.
    \end{align*}
    These span a single simplex in the tropicalization. Write
\[
  p=(1,1,1),\qquad
  m_x=(4,0,0),\quad m_y=(0,4,0),\quad
  m_z=(0,0,4),\quad m_w=(0,0,0),
\]
and 
\[
  \Xi=\ker(\mathds{1}_4)
    =\{(u_0,u_1,u_2,u_3)\in\ZZ^4:u_0+u_1+u_2+u_3=0\},
  \qquad
  E:=\{e_j-e_i:0\le i<j\le3\}.
\]
The initial forms at the four vertices give rise to four natural affine simplices $\psi_a:\ZZ^4\to M$, with
\[
    \begin{array}{|c|cccc|}
     \hline
      a&\psi_a(e_0)&\psi_a(e_1)&\psi_a(e_2)&\psi_a(e_3)\\ 
      \hline
      x&p&m_y&m_z&m_w\\
      y&p&m_x&m_z&m_w\\
      z&p&m_x&m_y&m_w\\
      w&p&m_x&m_y&m_z\\
      \hline
    \end{array}
\]
All four restrictions $\psi_a|_\Xi$ have the same image lattice
\[
  \Lambda
    =\ZZ p+4\ZZ^3
    =\{(r,s,t)\in\ZZ^3:r\equiv s\equiv t\pmod4\}.
\]
We have 
\[
  [M:\Lambda] = 16, \qquad M/\Lambda\cong(\ZZ/4\ZZ)^2.
\]

\paragraph{\underline{Finding suitable cellulation data}}

We now explain how to find cellulations of the different angle spaces. We distinguish between three levels:
\begin{itemize}
    \item Level 0: All the cells are of the form $(\tau,0)$. The maps of angle sets are injections, so we have to find a suitable set of $Q_{v}$ that define compatible cellulation data as in Definition \ref{def:TriangulationData}. 
    \item Level 1: We then consider cells of the form $(\tau,\sigma)$, where $\dim(\sigma) = 1$. We take projections of cells from level $0$ and find suitable $Q_{I}$. 
    \item Level 2: For cells $(\tau,\sigma)$ with $\dim(\sigma) = 2$, the corresponding angle spaces are zero-dimensional, so there is no need to find cellulations. 
\end{itemize}

We emphasize here that the cellulations coming from the local Coxeter arrangements are not compatible, as in Example \ref{ex:refinedCoxeterArrangement1}.     

The computations amount to elementary polyhedral geometry and linear algebra, so we have omitted them. We calculate the relevant $f$-vectors using Theorem \ref{thm:refined_Coxter_F_vect}.

\begin{enumerate}
    \item [\underline{Level 0}:]

    For every affine simplex $\psi_{a}$, it suffices to take $Q_{a}=Q$, where 
    \[
    Q=\{q_{0},q_{1},q_{2},q_{3}\}
    \]
    and 
    \[
    \begin{aligned}
      q_0&=(3,-1,-1,-1),&
      q_1&=(4,-2,-1,-1),\\
      q_2&=(4,-1,-2,-1),&
      q_3&=(4,-1,-1,-2).
    \end{aligned}
    \] 

Calculating the $f$-vectors of $\Delta_{[3],S,H}(Q)$ for different $S$ (with the index of $\Lambda$ in $H$ being $16$), we find the data as in Table \ref{tab:dwork-level-0-f-vectors}. 
\begin{table}[htbp]
\centering
\begin{tabular}{|ccc|}
\hline
\(S\) & \(f\)-vector on \(V_{[3]}/2L_{[3]}\)
      & \(f\)-vector on \(V_{[3]}/H\) \\
\hline
\(\{0,a\}\)
  & \((24,72,48,0)\)
  & \((384,1152,768,0)\) \\
\(\{a,b\}\)
  & \((40,112,72,0)\)
  & \((640,1792,1152,0)\) \\
\(\{0,a,b\}\)
  & \((76,328,396,144)\)
  & \((1216,5248,6336,2304)\) \\
\(\{1,2,3\}\)
  & \((84,360,420,144)\)
  & \((1344,5760,6720,2304)\) \\
\(\{0,1,2,3\}\)
  & \((109,534,714,288)\)
  & \((1744,8544,11424,4608)\) \\
\hline
\end{tabular}
\vspace{2mm}
\caption{\(f\)-vectors of the Level 0 cellulations
\(\Delta_{[3],S,H}(Q)\). The second column gives the $f$-vector of the complex modulo the full lattice $2L_{[3]}$. The third column gives the $f$-vector of the complex modulo the smaller lattice $H\subset 2L_{[3]}$. It is obtained by multiplying the second column by $16$.}
\label{tab:dwork-level-0-f-vectors}
\end{table}

\item[\underline{Level 1}:]

As an example, we will take $\psi_{z}$. Here only $I=\{1,2,3\}$ is relevant for the subdivision at infinity. The corresponding subgroup is $H_{I}=8L_{I}$, which is of index $16$ in $2L_{I}$. It now suffices to take 
\begin{equation*}
    Q_{I}={\rm Perm}(2,-1,-1)\cup {\rm Perm}(3,-2,-1).
\end{equation*}
Here ${\rm Perm}(2,-1,-1)=\{(2,-1,-1),(-1,2,-1),(-1,-1,2)\}$ gives the vectors obtained by permuting the coordinates. We thus find $|Q_{I}|=9$. Calculating the $f$-vectors of $\Delta_{I,S,H}(Q_{I})$ for different $S\subset I$ (with the index in $H$ again being $16$), we find the data in Table \ref{tab:dwork-level-1-f-vectors}. 
\begin{table}[htbp]
\centering
\renewcommand{\arraystretch}{1.2}
\setlength{\tabcolsep}{8pt}
\begin{tabular}{|ccc|}
\hline
\(S \subset  I\) & \(f\)-vector on \(V_I/2L_I\)
                & \(f\)-vector on \(V_I/H_I\) \\
\hline
Any two-element subset
  & \((8,8,0)\)
  & \((128,128,0)\) \\
\(I=\{1,2,3\}\)
  & \((65,174,108)\)
  & \((1040,2784,1728)\) \\
\hline
\end{tabular}
\vspace{2mm}
\caption{\(f\)-vectors of the Level 1 cellulations
\(\Delta_{I,S,H_I}(Q_I)\). The last column is obtained by
multiplying the middle column by \(16\).}
\label{tab:dwork-level-1-f-vectors}
\end{table}

\end{enumerate}

We now collect all of these triangulations and compute the global $f$-vector. We start with $\mathfrak{f}_{0}$:
\[
\begin{aligned}
\mathfrak{f}_0={}&6\cdot(9\cdot(640,1792,1152,0,0)+6\cdot(8,8,0,0,0)+(4,0,0,0,0))\\
&+4\cdot((1040,2784,1728,0,0)+(1344,5760,6720,2304,0))\\
&+4\cdot(1744,8544,11424,4608,0)+6\cdot(1216,5248,6336,2304,0)+4\cdot(384,1152,768,0,0).
\end{aligned}
\]
The $6$ in the first line comes from the $6$ unbounded $2$-dimensional regions in the tropicalization, the $4$ in the second line from the four unbounded rays, and the last line comes from the polytope. Note that the coefficients in front of the vectors can be seen as describing the number of vertices that give rise to this vector.

For $\mathfrak{f}_{1}$ we find the following
\[
\begin{aligned}
\mathfrak{f}_{1}={\rm sh}^{1}({}&6\cdot(38\cdot(640,1792,1152,0,0)+8\cdot(8,8,0,0,0))\\
&+8\cdot(1344,5760,6720,2304,0)+24\cdot(384,1152,768,0,0)+12\cdot(1216,5248,6336,2304,0)).
\end{aligned}
\]
Finally, for $\mathfrak{f}_{2}$ we find 

\[
\mathfrak{f}_{2}={\rm sh}^{2}(180\cdot(640,1792,1152,0,0)+24\cdot(384,1152,768,0,0)).
\]

The result is \begin{align*}
\mathfrak{f}_0 &= (60216, 201504, 182784, 41472, 0), \\
\mathfrak{f}_{1} &= (0, 180864, 545664, 410880, 46080), \\
\mathfrak{f}_2 &=   (0,     0,      124416, 350208, 225792).
\end{align*}

Adding these together yields the following $f$-vector for the CW complex structure arising from the Dwork family: 
\begin{equation*}
    \mathfrak{f} = (60216, 382368, 852864, 802560, 271872).
\end{equation*}
We again check that $\chi(\Delta(\mathcal{F}^{\triangle}))=24$.

\end{example}

\subsection{Calabi-Yau monodromy}\label{sec:CYMonodromy}

In this section, we will show how our polar decomposition works in the case of another K3-surface. We will work with the Ap\'{e}ry-Fermi pencil, whose generic fiber is birational to a K3 surface. Using a suitable toric compactification, this in fact becomes a family of K3-surfaces.

\begin{figure}[htpb]
\centering
\includegraphics[scale=0.5]{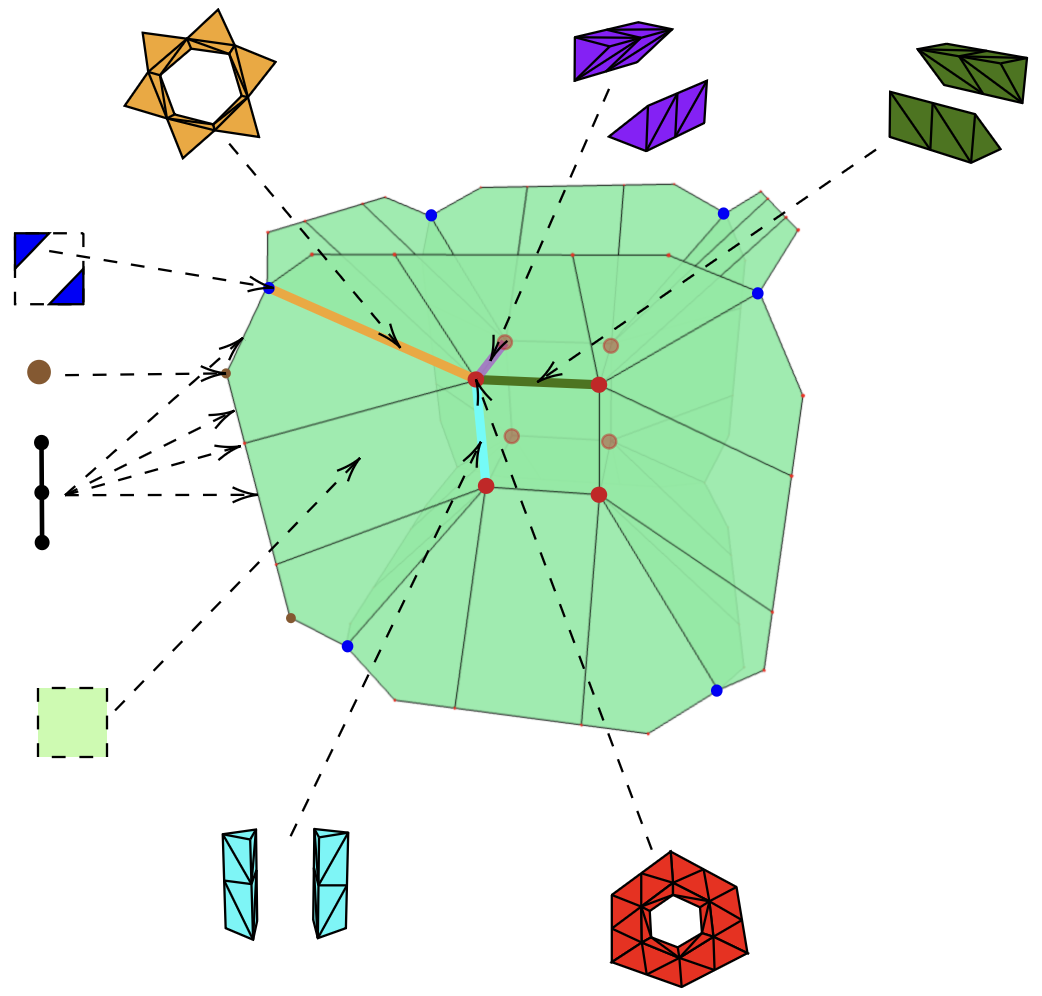}
\caption{\label{fig:AperyK32} The diagram of angle spaces associated to the Ap\'{e}ry-Fermi pencil. The $24$ singular points for the integral affine structure come from the fact that the angle spaces degenerate. For instance, next to the orange line we have three flaps where the diagram turns into a contractible space. The angle spaces here are drawn in a fundamental domain of $\RR^d / 2\ZZ^d$ with $d = 0, 1, 2, 3$.
}
\end{figure}

We start by calculating the $f$-vector using Theorem~\ref{thm:CWComplexStructureTriangulatedDiagram}. Consider the polynomial 
    \begin{equation*}
        f=1+t(x+y+z+x^{-1}+y^{-1}+z^{-1}).
    \end{equation*}
    Note that the subdivided Newton polygon consists of $8$ simplices, corresponding to $8$ vertices in the tropicalization. The ambient toric variety is not $\mathbb{P}^{3}$ in this case. It is a certain blow-up of $\PP^{3}$. Its fan can easily be deduced from Figure~\ref{fig:AperyK32}.

    We take the corresponding compactified tropicalization. Note that we have to add rays to obtain a smooth K3 surface. These can be found in Figure~\ref{fig:AperyK3}. As before, the $f$-vectors of the spaces are

    \[\begin{aligned}
        f(\Delta_{3,4}) &=& (7,42,60,24,0), \quad 
        f(\Delta_{3,3}) &=& (6,30,36,12,0), \quad
        f(\Delta_{3,2}) &=& (4,12,8,0,0),\\
        f(\Delta_{2,3}) &=& (3,6,2,0,0), \quad 
        f(\Delta_{2,2}) &=& (2,2,0,0,0), \quad
        f(\Delta_{1,1}) &=& (1,0,0,0,0).
    \end{aligned}
    \]
    
    We now calculate 
    {\footnotesize
    \begin{align*}
        \mathfrak{f}_0 &= 8f(\Delta_{3,4}) + 20f(\Delta_{3,3}) + {66}f(\Delta_{3,2}) + 8f(\Delta_{2,3}) + 84f(\Delta_{2,2}) + 24f(\Delta_{1,1})=(656, 1944, 1744, 432, 0),\\
\mathfrak{f}_{1}&={\rm sh}^{1}(40f(\Delta_{3,3}) + 384f(\Delta_{3,2}) + 0f(\Delta_{2,3}) + 120f(\Delta_{2,2}))=(0, 2016, 6048, 4512, 480),\\
        \mathfrak{f}_{2}&={\rm sh}^{2}(
        336f(\Delta_{3,2}))=(0, 0, 1344, 4032, 2688).
    \end{align*}
    }

This yields
\begin{equation*}
    \mathfrak{f}=\sum \mathfrak{f}_{k}=(656, 3960, 9136, 8976, 3168).
\end{equation*}
As before, we find that the Euler characteristic is $24$. From our homotopy colimit expression, we can also deduce $\pi_{1}(X)=(0)$. Indeed, we can obtain any loop by composing a loop in the tropicalization with a local loop in each angle space. Since the tropicalization is simply connected, it suffices to show that we can contract all the loops in the angle spaces. We start with a fixed loop in an angle space. We can now move this to the angle space at one of the flaps. By moving it to infinity, we then see that it is trivial. We conclude that $\pi_{1}(X)=(0)$.

\subsubsection{Monodromy of integral affine structures}

We now show how to interpret the monodromy of integral structures from the SYZ program in terms of angle spaces in the Ap\'{e}ry-Fermi example. This shows the monodromy arises canonically from infinite paths on the compactified tropicalization. We will follow \cite{Yamamoto22}, but with some adaptations for our affine approach. We also refer the reader to \cite{Gross01} and \cite{Goto23}.

\begin{definition}[Integral affine manifold structure]\label{def:IntegralAffineStructure}
    Consider the Ap\'{e}ry-Fermi tropicalization $\Sigma$, with the unimodular cellulation indicated in Figure \ref{fig:AperyK32} for the unbounded cones, and Figure \ref{fig:TriangulationOfK3Surface} for the sides of the cube $B$. We call the $8$ distinguished outer vertices of $B$ the corner vertices. 
    Let $\Pscr(B)$ be the induced poset and write $\Delta(B)$ for the corresponding order complex, which can be identified with the barycentric subdivision of $B$. 

    The \emph{singular points} of $\Delta(B)$ are the vertices on the $12$ edges of the cube that are adjacent to a corner vertex. 
    There are $24$ of these. We denote them by $\Gamma$.    

    We turn $\Delta(B)\backslash{\Gamma}$ into an integral affine manifold by providing certain charts, as in \cite[Section 3]{Yamamoto22}. We write $W_{v}$ for the star of $v$, with respect to the polyhedral structure $\Delta(B)$. We now assign the following two types of charts:
    \begin{enumerate}
        \item If $v$ is a vertex in the interior of a maximal $2$-dimensional cell $\tau$, then we assign the chart 
    \begin{equation*}
        W_{v}\to A_{\tau},
    \end{equation*}
    where $A_{\tau}$ is the affine linear subspace  spanned by $\tau$. 
    \item If $v$ is one of the $8$ corner vertices, then we consider the integral vector $\eta_{v}$ generating the ray in $\Sigma$ whose span with each adjacent plane in the cube is $\RR^{3}$. In particular, $\eta_{v}$ does not lie in any of the adjacent planes of the cube. We then assign $W_{v}\to N_{\RR}/{\rm Span}_{\RR}(\eta_{v})$.  
    \item Similarly, if $v$ is a vertex lying on an edge of the cube, then we take the vector $\eta_{v}$ perpendicular to the edge and in the direction of the outward $2$-dimensional cell. This then again defines a chart 
    \begin{equation*}
        W_{v} \to N_{\RR}/{\rm Span}_{\RR}(\eta_{v}).
    \end{equation*}
    \end{enumerate} 
    This endows  $B\backslash \Gamma$ with the structure of an affine integral manifold.  
\end{definition}

\begin{figure}
    \centering
    \scalebox{0.8}{

\tikzset{every picture/.style={line width=0.75pt}} %

\begin{tikzpicture}[x=0.75pt,y=0.75pt,yscale=-1,xscale=1]
\draw   (328.33,53.36) -- (398.33,96.21) -- (398.33,173.36) -- (328.33,216.21) -- (291.25,193.51) -- (258.33,173.36) -- (258.33,96.21) -- (292.26,75.44) -- cycle ;
\draw [line width=2.25]    (258.33,96.21) -- (293.71,117.87) -- (328.33,139.07) ;
\draw [line width=2.25]    (328.33,139.07) -- (398.33,96.21) ;
\draw [line width=2.25]    (328.33,216.21) -- (328.33,178.79) -- (328.33,139.07) ;
\draw    (328.33,53.36) -- (328.33,96.5) -- (328.33,139.07) ;
\draw  [color={rgb, 255:red, 74; green, 144; blue, 226 }  ,draw opacity=1 ][fill={rgb, 255:red, 74; green, 144; blue, 226 }  ,fill opacity=1 ] (254.42,96.21) .. controls (254.42,94.36) and (256.17,92.86) .. (258.33,92.86) .. controls (260.5,92.86) and (262.25,94.36) .. (262.25,96.21) .. controls (262.25,98.07) and (260.5,99.57) .. (258.33,99.57) .. controls (256.17,99.57) and (254.42,98.07) .. (254.42,96.21) -- cycle ;
\draw  [color={rgb, 255:red, 74; green, 144; blue, 226 }  ,draw opacity=1 ][fill={rgb, 255:red, 74; green, 144; blue, 226 }  ,fill opacity=1 ] (324.42,216.21) .. controls (324.42,214.36) and (326.17,212.86) .. (328.33,212.86) .. controls (330.5,212.86) and (332.25,214.36) .. (332.25,216.21) .. controls (332.25,218.07) and (330.5,219.57) .. (328.33,219.57) .. controls (326.17,219.57) and (324.42,218.07) .. (324.42,216.21) -- cycle ;
\draw  [color={rgb, 255:red, 74; green, 144; blue, 226 }  ,draw opacity=1 ][fill={rgb, 255:red, 74; green, 144; blue, 226 }  ,fill opacity=1 ] (394.42,96.21) .. controls (394.42,94.36) and (396.17,92.86) .. (398.33,92.86) .. controls (400.5,92.86) and (402.25,94.36) .. (402.25,96.21) .. controls (402.25,98.07) and (400.5,99.57) .. (398.33,99.57) .. controls (396.17,99.57) and (394.42,98.07) .. (394.42,96.21) -- cycle ;
\draw  [color={rgb, 255:red, 126; green, 211; blue, 33 }  ,draw opacity=1 ][fill={rgb, 255:red, 126; green, 211; blue, 33 }  ,fill opacity=1 ] (324.42,53.36) .. controls (324.42,51.5) and (326.17,50) .. (328.33,50) .. controls (330.5,50) and (332.25,51.5) .. (332.25,53.36) .. controls (332.25,55.21) and (330.5,56.71) .. (328.33,56.71) .. controls (326.17,56.71) and (324.42,55.21) .. (324.42,53.36) -- cycle ;
\draw   (568.33,53.36) -- (638.33,96.21) -- (638.33,173.36) -- (568.33,216.21) -- (531.25,193.51) -- (498.33,173.36) -- (498.33,96.21) -- (532.26,75.44) -- cycle ;
\draw [line width=2.25]    (498.33,96.21) -- (533.71,117.87) -- (568.33,139.07) ;
\draw [line width=2.25]    (568.33,139.07) -- (638.33,96.21) ;
\draw [line width=2.25]    (568.33,216.21) -- (568.33,178.79) -- (568.33,139.07) ;
\draw    (568.33,53.36) -- (568.33,96.5) -- (568.33,139.07) ;
\draw    (568.33,139.07) -- (528.38,158.64) -- (498.33,173.36) ;
\draw    (638.33,173.36) -- (568.33,139.07) ;
\draw  [dash pattern={on 4.5pt off 4.5pt}]  (518.33,139.07) -- (498.33,96.21) ;
\draw  [dash pattern={on 4.5pt off 4.5pt}]  (498.33,139.07) -- (518.33,139.07) ;
\draw  [dash pattern={on 4.5pt off 4.5pt}]  (498.33,173.36) -- (518.33,139.07) ;
\draw  [dash pattern={on 4.5pt off 4.5pt}]  (568.33,139.07) -- (518.33,139.07) ;
\draw  [dash pattern={on 4.5pt off 4.5pt}]  (533.71,117.87) -- (518.33,139.07) ;
\draw  [dash pattern={on 4.5pt off 4.5pt}]  (518.33,139.07) -- (528.38,158.64) ;
\draw  [dash pattern={on 4.5pt off 4.5pt}]  (528.38,158.64) -- (542.67,176.21) ;
\draw  [dash pattern={on 4.5pt off 4.5pt}]  (542.67,176.21) -- (568.33,139.07) ;
\draw  [dash pattern={on 4.5pt off 4.5pt}]  (498.33,173.36) -- (542.67,176.21) ;
\draw  [dash pattern={on 4.5pt off 4.5pt}]  (542.67,176.21) -- (568.33,216.21) ;
\draw  [dash pattern={on 4.5pt off 4.5pt}]  (542.67,176.21) -- (568.33,178.79) ;
\draw  [dash pattern={on 4.5pt off 4.5pt}]  (542.67,176.21) -- (531.25,193.51) ;
\draw  [dash pattern={on 4.5pt off 4.5pt}]  (568.33,53.36) -- (546.55,99.3) ;
\draw  [dash pattern={on 4.5pt off 4.5pt}]  (498.33,96.21) -- (546.55,99.3) ;
\draw  [dash pattern={on 4.5pt off 4.5pt}]  (546.55,99.3) -- (568.33,139.07) ;
\draw  [dash pattern={on 4.5pt off 4.5pt}]  (568.33,96.5) -- (546.55,99.3) ;
\draw  [dash pattern={on 4.5pt off 4.5pt}]  (532.26,75.44) -- (546.55,99.3) ;
\draw  [dash pattern={on 4.5pt off 4.5pt}]  (546.55,99.3) -- (533.71,117.87) ;
\draw  [dash pattern={on 4.5pt off 4.5pt}]  (638.33,96.21) -- (589.55,97.59) ;
\draw  [dash pattern={on 4.5pt off 4.5pt}]  (568.33,53.36) -- (589.55,97.59) ;
\draw  [dash pattern={on 4.5pt off 4.5pt}]  (589.55,97.59) -- (568.33,139.07) ;
\draw  [dash pattern={on 4.5pt off 4.5pt}]  (589.55,97.59) -- (568.33,96.5) ;
\draw  [dash pattern={on 4.5pt off 4.5pt}]  (603.33,117.64) -- (589.55,97.59) ;
\draw  [dash pattern={on 4.5pt off 4.5pt}]  (603.33,74.79) -- (589.55,97.59) ;
\draw  [dash pattern={on 4.5pt off 4.5pt}]  (616.67,140.21) -- (568.33,139.07) ;
\draw  [dash pattern={on 4.5pt off 4.5pt}]  (603.33,117.64) -- (616.67,140.21) ;
\draw  [dash pattern={on 4.5pt off 4.5pt}]  (638.33,96.21) -- (616.67,140.21) ;
\draw  [dash pattern={on 4.5pt off 4.5pt}]  (638.33,173.36) -- (616.67,140.21) ;
\draw  [dash pattern={on 4.5pt off 4.5pt}]  (638.33,134.79) -- (616.67,140.21) ;
\draw  [dash pattern={on 4.5pt off 4.5pt}]  (616.67,140.21) -- (603.33,156.21) ;
\draw  [dash pattern={on 4.5pt off 4.5pt}]  (591.67,177.07) -- (568.33,178.79) ;
\draw  [dash pattern={on 4.5pt off 4.5pt}]  (568.33,139.07) -- (591.67,177.07) ;
\draw  [dash pattern={on 4.5pt off 4.5pt}]  (603.33,156.21) -- (591.67,177.07) ;
\draw  [dash pattern={on 4.5pt off 4.5pt}]  (603.33,194.79) -- (591.67,177.07) ;
\draw  [dash pattern={on 4.5pt off 4.5pt}]  (638.33,173.36) -- (591.67,177.07) ;
\draw  [dash pattern={on 4.5pt off 4.5pt}]  (591.67,177.07) -- (568.33,216.21) ;
\draw  [color={rgb, 255:red, 74; green, 144; blue, 226 }  ,draw opacity=1 ][fill={rgb, 255:red, 74; green, 144; blue, 226 }  ,fill opacity=1 ] (494.42,96.21) .. controls (494.42,94.36) and (496.17,92.86) .. (498.33,92.86) .. controls (500.5,92.86) and (502.25,94.36) .. (502.25,96.21) .. controls (502.25,98.07) and (500.5,99.57) .. (498.33,99.57) .. controls (496.17,99.57) and (494.42,98.07) .. (494.42,96.21) -- cycle ;
\draw  [color={rgb, 255:red, 74; green, 144; blue, 226 }  ,draw opacity=1 ][fill={rgb, 255:red, 74; green, 144; blue, 226 }  ,fill opacity=1 ] (564.42,216.21) .. controls (564.42,214.36) and (566.17,212.86) .. (568.33,212.86) .. controls (570.5,212.86) and (572.25,214.36) .. (572.25,216.21) .. controls (572.25,218.07) and (570.5,219.57) .. (568.33,219.57) .. controls (566.17,219.57) and (564.42,218.07) .. (564.42,216.21) -- cycle ;
\draw  [color={rgb, 255:red, 126; green, 211; blue, 33 }  ,draw opacity=1 ][fill={rgb, 255:red, 126; green, 211; blue, 33 }  ,fill opacity=1 ] (494.42,173.36) .. controls (494.42,171.5) and (496.17,170) .. (498.33,170) .. controls (500.5,170) and (502.25,171.5) .. (502.25,173.36) .. controls (502.25,175.21) and (500.5,176.71) .. (498.33,176.71) .. controls (496.17,176.71) and (494.42,175.21) .. (494.42,173.36) -- cycle ;
\draw  [color={rgb, 255:red, 126; green, 211; blue, 33 }  ,draw opacity=1 ][fill={rgb, 255:red, 184; green, 233; blue, 134 }  ,fill opacity=1 ] (634.42,173.36) .. controls (634.42,171.5) and (636.17,170) .. (638.33,170) .. controls (640.5,170) and (642.25,171.5) .. (642.25,173.36) .. controls (642.25,175.21) and (640.5,176.71) .. (638.33,176.71) .. controls (636.17,176.71) and (634.42,175.21) .. (634.42,173.36) -- cycle ;
\draw  [color={rgb, 255:red, 139; green, 87; blue, 42 }  ,draw opacity=1 ][fill={rgb, 255:red, 139; green, 87; blue, 42 }  ,fill opacity=1 ] (529.8,117.87) .. controls (529.8,116.02) and (531.55,114.52) .. (533.71,114.52) .. controls (535.87,114.52) and (537.63,116.02) .. (537.63,117.87) .. controls (537.63,119.73) and (535.87,121.23) .. (533.71,121.23) .. controls (531.55,121.23) and (529.8,119.73) .. (529.8,117.87) -- cycle ;
\draw  [color={rgb, 255:red, 139; green, 87; blue, 42 }  ,draw opacity=1 ][fill={rgb, 255:red, 139; green, 87; blue, 42 }  ,fill opacity=1 ] (599.42,117.64) .. controls (599.42,115.79) and (601.17,114.29) .. (603.33,114.29) .. controls (605.5,114.29) and (607.25,115.79) .. (607.25,117.64) .. controls (607.25,119.5) and (605.5,121) .. (603.33,121) .. controls (601.17,121) and (599.42,119.5) .. (599.42,117.64) -- cycle ;
\draw  [color={rgb, 255:red, 139; green, 87; blue, 42 }  ,draw opacity=1 ][fill={rgb, 255:red, 139; green, 87; blue, 42 }  ,fill opacity=1 ] (564.42,178.79) .. controls (564.42,176.93) and (566.17,175.43) .. (568.33,175.43) .. controls (570.5,175.43) and (572.25,176.93) .. (572.25,178.79) .. controls (572.25,180.64) and (570.5,182.14) .. (568.33,182.14) .. controls (566.17,182.14) and (564.42,180.64) .. (564.42,178.79) -- cycle ;
\draw  [color={rgb, 255:red, 189; green, 16; blue, 224 }  ,draw opacity=1 ][fill={rgb, 255:red, 189; green, 16; blue, 224 }  ,fill opacity=1 ] (528.34,75.44) .. controls (528.34,73.59) and (530.1,72.08) .. (532.26,72.08) .. controls (534.42,72.08) and (536.18,73.59) .. (536.18,75.44) .. controls (536.18,77.3) and (534.42,78.8) .. (532.26,78.8) .. controls (530.1,78.8) and (528.34,77.3) .. (528.34,75.44) -- cycle ;
\draw  [color={rgb, 255:red, 189; green, 16; blue, 224 }  ,draw opacity=1 ][fill={rgb, 255:red, 189; green, 16; blue, 224 }  ,fill opacity=1 ] (524.46,158.64) .. controls (524.46,156.79) and (526.22,155.28) .. (528.38,155.28) .. controls (530.54,155.28) and (532.3,156.79) .. (532.3,158.64) .. controls (532.3,160.49) and (530.54,162) .. (528.38,162) .. controls (526.22,162) and (524.46,160.49) .. (524.46,158.64) -- cycle ;
\draw  [color={rgb, 255:red, 189; green, 16; blue, 224 }  ,draw opacity=1 ][fill={rgb, 255:red, 189; green, 16; blue, 224 }  ,fill opacity=1 ] (599.42,74.79) .. controls (599.42,72.93) and (601.17,71.43) .. (603.33,71.43) .. controls (605.5,71.43) and (607.25,72.93) .. (607.25,74.79) .. controls (607.25,76.64) and (605.5,78.14) .. (603.33,78.14) .. controls (601.17,78.14) and (599.42,76.64) .. (599.42,74.79) -- cycle ;
\draw  [color={rgb, 255:red, 189; green, 16; blue, 224 }  ,draw opacity=1 ][fill={rgb, 255:red, 189; green, 16; blue, 224 }  ,fill opacity=1 ] (494.42,139.07) .. controls (494.42,137.22) and (496.17,135.71) .. (498.33,135.71) .. controls (500.5,135.71) and (502.25,137.22) .. (502.25,139.07) .. controls (502.25,140.93) and (500.5,142.43) .. (498.33,142.43) .. controls (496.17,142.43) and (494.42,140.93) .. (494.42,139.07) -- cycle ;
\draw  [color={rgb, 255:red, 189; green, 16; blue, 224 }  ,draw opacity=1 ][fill={rgb, 255:red, 189; green, 16; blue, 224 }  ,fill opacity=1 ] (527.34,193.51) .. controls (527.34,191.66) and (529.09,190.15) .. (531.25,190.15) .. controls (533.42,190.15) and (535.17,191.66) .. (535.17,193.51) .. controls (535.17,195.37) and (533.42,196.87) .. (531.25,196.87) .. controls (529.09,196.87) and (527.34,195.37) .. (527.34,193.51) -- cycle ;
\draw  [color={rgb, 255:red, 189; green, 16; blue, 224 }  ,draw opacity=1 ][fill={rgb, 255:red, 189; green, 16; blue, 224 }  ,fill opacity=1 ] (599.42,194.79) .. controls (599.42,192.93) and (601.17,191.43) .. (603.33,191.43) .. controls (605.5,191.43) and (607.25,192.93) .. (607.25,194.79) .. controls (607.25,196.64) and (605.5,198.14) .. (603.33,198.14) .. controls (601.17,198.14) and (599.42,196.64) .. (599.42,194.79) -- cycle ;
\draw  [color={rgb, 255:red, 189; green, 16; blue, 224 }  ,draw opacity=1 ][fill={rgb, 255:red, 189; green, 16; blue, 224 }  ,fill opacity=1 ] (538.75,176.21) .. controls (538.75,174.36) and (540.5,172.86) .. (542.67,172.86) .. controls (544.83,172.86) and (546.58,174.36) .. (546.58,176.21) .. controls (546.58,178.07) and (544.83,179.57) .. (542.67,179.57) .. controls (540.5,179.57) and (538.75,178.07) .. (538.75,176.21) -- cycle ;
\draw  [color={rgb, 255:red, 189; green, 16; blue, 224 }  ,draw opacity=1 ][fill={rgb, 255:red, 189; green, 16; blue, 224 }  ,fill opacity=1 ] (514.42,139.07) .. controls (514.42,137.22) and (516.17,135.71) .. (518.33,135.71) .. controls (520.5,135.71) and (522.25,137.22) .. (522.25,139.07) .. controls (522.25,140.93) and (520.5,142.43) .. (518.33,142.43) .. controls (516.17,142.43) and (514.42,140.93) .. (514.42,139.07) -- cycle ;
\draw  [color={rgb, 255:red, 189; green, 16; blue, 224 }  ,draw opacity=1 ][fill={rgb, 255:red, 189; green, 16; blue, 224 }  ,fill opacity=1 ] (542.63,99.3) .. controls (542.63,97.45) and (544.38,95.95) .. (546.55,95.95) .. controls (548.71,95.95) and (550.46,97.45) .. (550.46,99.3) .. controls (550.46,101.16) and (548.71,102.66) .. (546.55,102.66) .. controls (544.38,102.66) and (542.63,101.16) .. (542.63,99.3) -- cycle ;
\draw  [color={rgb, 255:red, 189; green, 16; blue, 224 }  ,draw opacity=1 ][fill={rgb, 255:red, 189; green, 16; blue, 224 }  ,fill opacity=1 ] (585.63,97.59) .. controls (585.63,95.73) and (587.38,94.23) .. (589.55,94.23) .. controls (591.71,94.23) and (593.46,95.73) .. (593.46,97.59) .. controls (593.46,99.44) and (591.71,100.95) .. (589.55,100.95) .. controls (587.38,100.95) and (585.63,99.44) .. (585.63,97.59) -- cycle ;
\draw  [color={rgb, 255:red, 189; green, 16; blue, 224 }  ,draw opacity=1 ][fill={rgb, 255:red, 189; green, 16; blue, 224 }  ,fill opacity=1 ] (612.75,140.21) .. controls (612.75,138.36) and (614.5,136.86) .. (616.67,136.86) .. controls (618.83,136.86) and (620.58,138.36) .. (620.58,140.21) .. controls (620.58,142.07) and (618.83,143.57) .. (616.67,143.57) .. controls (614.5,143.57) and (612.75,142.07) .. (612.75,140.21) -- cycle ;
\draw  [color={rgb, 255:red, 189; green, 16; blue, 224 }  ,draw opacity=1 ][fill={rgb, 255:red, 189; green, 16; blue, 224 }  ,fill opacity=1 ] (587.75,177.07) .. controls (587.75,175.22) and (589.5,173.71) .. (591.67,173.71) .. controls (593.83,173.71) and (595.58,175.22) .. (595.58,177.07) .. controls (595.58,178.93) and (593.83,180.43) .. (591.67,180.43) .. controls (589.5,180.43) and (587.75,178.93) .. (587.75,177.07) -- cycle ;
\draw  [color={rgb, 255:red, 189; green, 16; blue, 224 }  ,draw opacity=1 ][fill={rgb, 255:red, 189; green, 16; blue, 224 }  ,fill opacity=1 ] (564.42,96.5) .. controls (564.42,94.65) and (566.17,93.14) .. (568.33,93.14) .. controls (570.5,93.14) and (572.25,94.65) .. (572.25,96.5) .. controls (572.25,98.35) and (570.5,99.86) .. (568.33,99.86) .. controls (566.17,99.86) and (564.42,98.35) .. (564.42,96.5) -- cycle ;
\draw  [color={rgb, 255:red, 189; green, 16; blue, 224 }  ,draw opacity=1 ][fill={rgb, 255:red, 189; green, 16; blue, 224 }  ,fill opacity=1 ] (634.42,134.79) .. controls (634.42,132.93) and (636.17,131.43) .. (638.33,131.43) .. controls (640.5,131.43) and (642.25,132.93) .. (642.25,134.79) .. controls (642.25,136.64) and (640.5,138.14) .. (638.33,138.14) .. controls (636.17,138.14) and (634.42,136.64) .. (634.42,134.79) -- cycle ;
\draw  [color={rgb, 255:red, 189; green, 16; blue, 224 }  ,draw opacity=1 ][fill={rgb, 255:red, 189; green, 16; blue, 224 }  ,fill opacity=1 ] (599.42,156.21) .. controls (599.42,154.36) and (601.17,152.86) .. (603.33,152.86) .. controls (605.5,152.86) and (607.25,154.36) .. (607.25,156.21) .. controls (607.25,158.07) and (605.5,159.57) .. (603.33,159.57) .. controls (601.17,159.57) and (599.42,158.07) .. (599.42,156.21) -- cycle ;
\draw    (258.33,173.36) -- (328.33,139.07) ;
\draw    (328.33,139.07) -- (398.33,173.36) ;
\draw  [color={rgb, 255:red, 126; green, 211; blue, 33 }  ,draw opacity=1 ][fill={rgb, 255:red, 126; green, 211; blue, 33 }  ,fill opacity=1 ] (394.42,173.36) .. controls (394.42,171.5) and (396.17,170) .. (398.33,170) .. controls (400.5,170) and (402.25,171.5) .. (402.25,173.36) .. controls (402.25,175.21) and (400.5,176.71) .. (398.33,176.71) .. controls (396.17,176.71) and (394.42,175.21) .. (394.42,173.36) -- cycle ;
\draw  [color={rgb, 255:red, 126; green, 211; blue, 33 }  ,draw opacity=1 ][fill={rgb, 255:red, 126; green, 211; blue, 33 }  ,fill opacity=1 ] (254.42,173.36) .. controls (254.42,171.5) and (256.17,170) .. (258.33,170) .. controls (260.5,170) and (262.25,171.5) .. (262.25,173.36) .. controls (262.25,175.21) and (260.5,176.71) .. (258.33,176.71) .. controls (256.17,176.71) and (254.42,175.21) .. (254.42,173.36) -- cycle ;
\draw  [color={rgb, 255:red, 255; green, 0; blue, 0 }  ,draw opacity=1 ][fill={rgb, 255:red, 255; green, 0; blue, 0 }  ,fill opacity=1 ] (322.17,139.07) .. controls (322.17,136.15) and (324.93,133.79) .. (328.33,133.79) .. controls (331.74,133.79) and (334.5,136.15) .. (334.5,139.07) .. controls (334.5,141.99) and (331.74,144.36) .. (328.33,144.36) .. controls (324.93,144.36) and (322.17,141.99) .. (322.17,139.07) -- cycle ;
\draw  [color={rgb, 255:red, 126; green, 211; blue, 33 }  ,draw opacity=1 ][fill={rgb, 255:red, 126; green, 211; blue, 33 }  ,fill opacity=1 ] (564.42,53.36) .. controls (564.42,51.5) and (566.17,50) .. (568.33,50) .. controls (570.5,50) and (572.25,51.5) .. (572.25,53.36) .. controls (572.25,55.21) and (570.5,56.71) .. (568.33,56.71) .. controls (566.17,56.71) and (564.42,55.21) .. (564.42,53.36) -- cycle ;
\draw [color={rgb, 255:red, 245; green, 166; blue, 35 }  ,draw opacity=1 ][line width=1.5]    (568.33,139.07) -- (588.18,100.26) ;
\draw [shift={(589.55,97.59)}, rotate = 117.08] [color={rgb, 255:red, 245; green, 166; blue, 35 }  ,draw opacity=1 ][line width=1.5]    (14.21,-4.28) .. controls (9.04,-1.82) and (4.3,-0.39) .. (0,0) .. controls (4.3,0.39) and (9.04,1.82) .. (14.21,4.28)   ;
\draw [color={rgb, 255:red, 245; green, 166; blue, 35 }  ,draw opacity=1 ][line width=1.5]    (589.55,97.59) -- (635.33,96.3) ;
\draw [shift={(638.33,96.21)}, rotate = 178.39] [color={rgb, 255:red, 245; green, 166; blue, 35 }  ,draw opacity=1 ][line width=1.5]    (14.21,-4.28) .. controls (9.04,-1.82) and (4.3,-0.39) .. (0,0) .. controls (4.3,0.39) and (9.04,1.82) .. (14.21,4.28)   ;
\draw [color={rgb, 255:red, 245; green, 166; blue, 35 }  ,draw opacity=1 ][line width=1.5]    (638.33,96.21) -- (617.99,137.52) ;
\draw [shift={(616.67,140.21)}, rotate = 296.22] [color={rgb, 255:red, 245; green, 166; blue, 35 }  ,draw opacity=1 ][line width=1.5]    (14.21,-4.28) .. controls (9.04,-1.82) and (4.3,-0.39) .. (0,0) .. controls (4.3,0.39) and (9.04,1.82) .. (14.21,4.28)   ;
\draw [color={rgb, 255:red, 245; green, 166; blue, 35 }  ,draw opacity=1 ][line width=1.5]    (616.67,140.21) -- (571.33,139.14) ;
\draw [shift={(568.33,139.07)}, rotate = 1.35] [color={rgb, 255:red, 245; green, 166; blue, 35 }  ,draw opacity=1 ][line width=1.5]    (14.21,-4.28) .. controls (9.04,-1.82) and (4.3,-0.39) .. (0,0) .. controls (4.3,0.39) and (9.04,1.82) .. (14.21,4.28)   ;
\draw  [color={rgb, 255:red, 255; green, 0; blue, 0 }  ,draw opacity=1 ][fill={rgb, 255:red, 255; green, 0; blue, 0 }  ,fill opacity=1 ] (562.17,139.07) .. controls (562.17,136.15) and (564.93,133.79) .. (568.33,133.79) .. controls (571.74,133.79) and (574.5,136.15) .. (574.5,139.07) .. controls (574.5,141.99) and (571.74,144.36) .. (568.33,144.36) .. controls (564.93,144.36) and (562.17,141.99) .. (562.17,139.07) -- cycle ;
\draw  [color={rgb, 255:red, 74; green, 144; blue, 226 }  ,draw opacity=1 ][fill={rgb, 255:red, 74; green, 144; blue, 226 }  ,fill opacity=1 ] (634.42,96.21) .. controls (634.42,94.36) and (636.17,92.86) .. (638.33,92.86) .. controls (640.5,92.86) and (642.25,94.36) .. (642.25,96.21) .. controls (642.25,98.07) and (640.5,99.57) .. (638.33,99.57) .. controls (636.17,99.57) and (634.42,98.07) .. (634.42,96.21) -- cycle ;
\draw [line width=2.25]    (16.6,94.69) -- (76.59,137.58) ;
\draw [line width=2.25]    (16.6,94.69) -- (95.66,54.91) ;
\draw [line width=2.25]    (76.59,137.58) -- (114.32,118.6) -- (155.65,97.8) ;
\draw [line width=2.25]    (95.66,54.91) -- (155.65,97.8) ;
\draw [line width=2.25]    (76.59,225.09) -- (76.59,137.58) ;
\draw [line width=2.25]    (16.6,182.2) -- (16.6,94.69) ;
\draw [line width=2.25]    (155.65,185.31) -- (155.65,97.8) ;
\draw [line width=2.25]    (16.6,182.2) -- (76.59,225.09) ;
\draw    (95.66,54.91) -- (76.59,137.58) ;
\draw    (16.6,94.69) -- (155.65,97.8) ;
\draw    (125.66,76.35) -- (46.59,116.13) ;
\draw    (115.22,118.14) -- (57.03,74.34) ;
\draw  [color={rgb, 255:red, 126; green, 211; blue, 33 }  ,draw opacity=1 ][fill={rgb, 255:red, 126; green, 211; blue, 33 }  ,fill opacity=1 ] (81.94,96.24) .. controls (81.94,94.52) and (83.81,93.13) .. (86.13,93.13) .. controls (88.44,93.13) and (90.31,94.52) .. (90.31,96.24) .. controls (90.31,97.96) and (88.44,99.36) .. (86.13,99.36) .. controls (83.81,99.36) and (81.94,97.96) .. (81.94,96.24) -- cycle ;
\draw [line width=2.25]    (76.59,225.09) -- (114.32,206.11) -- (155.65,185.31) ;
\draw    (116.12,205.2) -- (116.12,117.69) ;
\draw    (76.59,225.09) -- (155.65,97.8) ;
\draw    (155.65,185.31) -- (76.59,137.58) ;
\draw    (155.65,141.56) -- (76.59,181.34) ;
\draw    (76.59,181.34) -- (18.39,137.54) ;
\draw    (46.59,203.65) -- (46.59,116.13) ;
\draw    (76.59,225.09) -- (16.6,94.69) ;
\draw    (16.6,182.2) -- (76.59,137.58) ;
\draw  [color={rgb, 255:red, 126; green, 211; blue, 33 }  ,draw opacity=1 ][fill={rgb, 255:red, 126; green, 211; blue, 33 }  ,fill opacity=1 ] (111.93,161.45) .. controls (111.93,159.73) and (113.81,158.33) .. (116.12,158.33) .. controls (118.43,158.33) and (120.31,159.73) .. (120.31,161.45) .. controls (120.31,163.17) and (118.43,164.56) .. (116.12,164.56) .. controls (113.81,164.56) and (111.93,163.17) .. (111.93,161.45) -- cycle ;
\draw  [color={rgb, 255:red, 126; green, 211; blue, 33 }  ,draw opacity=1 ][fill={rgb, 255:red, 126; green, 211; blue, 33 }  ,fill opacity=1 ] (42.4,159.89) .. controls (42.4,158.17) and (44.28,156.77) .. (46.59,156.77) .. controls (48.91,156.77) and (50.78,158.17) .. (50.78,159.89) .. controls (50.78,161.61) and (48.91,163.01) .. (46.59,163.01) .. controls (44.28,163.01) and (42.4,161.61) .. (42.4,159.89) -- cycle ;
\draw  [color={rgb, 255:red, 74; green, 144; blue, 226 }  ,draw opacity=1 ][fill={rgb, 255:red, 74; green, 144; blue, 226 }  ,fill opacity=1 ] (110.13,118.6) .. controls (110.13,116.87) and (112.01,115.48) .. (114.32,115.48) .. controls (116.64,115.48) and (118.51,116.87) .. (118.51,118.6) .. controls (118.51,120.32) and (116.64,121.71) .. (114.32,121.71) .. controls (112.01,121.71) and (110.13,120.32) .. (110.13,118.6) -- cycle ;
\draw  [color={rgb, 255:red, 74; green, 144; blue, 226 }  ,draw opacity=1 ][fill={rgb, 255:red, 74; green, 144; blue, 226 }  ,fill opacity=1 ] (42.4,116.13) .. controls (42.4,114.41) and (44.28,113.02) .. (46.59,113.02) .. controls (48.91,113.02) and (50.78,114.41) .. (50.78,116.13) .. controls (50.78,117.85) and (48.91,119.25) .. (46.59,119.25) .. controls (44.28,119.25) and (42.4,117.85) .. (42.4,116.13) -- cycle ;
\draw  [color={rgb, 255:red, 74; green, 144; blue, 226 }  ,draw opacity=1 ][fill={rgb, 255:red, 74; green, 144; blue, 226 }  ,fill opacity=1 ] (111.93,205.2) .. controls (111.93,203.48) and (113.81,202.09) .. (116.12,202.09) .. controls (118.43,202.09) and (120.31,203.48) .. (120.31,205.2) .. controls (120.31,206.93) and (118.43,208.32) .. (116.12,208.32) .. controls (113.81,208.32) and (111.93,206.93) .. (111.93,205.2) -- cycle ;
\draw  [color={rgb, 255:red, 74; green, 144; blue, 226 }  ,draw opacity=1 ][fill={rgb, 255:red, 74; green, 144; blue, 226 }  ,fill opacity=1 ] (42.4,203.65) .. controls (42.4,201.93) and (44.28,200.53) .. (46.59,200.53) .. controls (48.91,200.53) and (50.78,201.93) .. (50.78,203.65) .. controls (50.78,205.37) and (48.91,206.76) .. (46.59,206.76) .. controls (44.28,206.76) and (42.4,205.37) .. (42.4,203.65) -- cycle ;
\draw  [color={rgb, 255:red, 74; green, 144; blue, 226 }  ,draw opacity=1 ][fill={rgb, 255:red, 74; green, 144; blue, 226 }  ,fill opacity=1 ] (12.41,138.44) .. controls (12.41,136.72) and (14.28,135.33) .. (16.6,135.33) .. controls (18.91,135.33) and (20.79,136.72) .. (20.79,138.44) .. controls (20.79,140.16) and (18.91,141.56) .. (16.6,141.56) .. controls (14.28,141.56) and (12.41,140.16) .. (12.41,138.44) -- cycle ;
\draw  [color={rgb, 255:red, 74; green, 144; blue, 226 }  ,draw opacity=1 ][fill={rgb, 255:red, 74; green, 144; blue, 226 }  ,fill opacity=1 ] (151.46,141.56) .. controls (151.46,139.84) and (153.34,138.44) .. (155.65,138.44) .. controls (157.97,138.44) and (159.84,139.84) .. (159.84,141.56) .. controls (159.84,143.28) and (157.97,144.67) .. (155.65,144.67) .. controls (153.34,144.67) and (151.46,143.28) .. (151.46,141.56) -- cycle ;
\draw  [color={rgb, 255:red, 74; green, 144; blue, 226 }  ,draw opacity=1 ][fill={rgb, 255:red, 74; green, 144; blue, 226 }  ,fill opacity=1 ] (121.47,76.35) .. controls (121.47,74.63) and (123.34,73.24) .. (125.66,73.24) .. controls (127.97,73.24) and (129.85,74.63) .. (129.85,76.35) .. controls (129.85,78.07) and (127.97,79.47) .. (125.66,79.47) .. controls (123.34,79.47) and (121.47,78.07) .. (121.47,76.35) -- cycle ;
\draw  [color={rgb, 255:red, 74; green, 144; blue, 226 }  ,draw opacity=1 ][fill={rgb, 255:red, 74; green, 144; blue, 226 }  ,fill opacity=1 ] (52.84,74.34) .. controls (52.84,72.62) and (54.71,71.23) .. (57.03,71.23) .. controls (59.34,71.23) and (61.22,72.62) .. (61.22,74.34) .. controls (61.22,76.06) and (59.34,77.46) .. (57.03,77.46) .. controls (54.71,77.46) and (52.84,76.06) .. (52.84,74.34) -- cycle ;
\draw  [color={rgb, 255:red, 255; green, 0; blue, 0 }  ,draw opacity=1 ][fill={rgb, 255:red, 255; green, 0; blue, 0 }  ,fill opacity=1 ] (149.06,97.8) .. controls (149.06,95.09) and (152.01,92.9) .. (155.65,92.9) .. controls (159.3,92.9) and (162.25,95.09) .. (162.25,97.8) .. controls (162.25,100.51) and (159.3,102.71) .. (155.65,102.71) .. controls (152.01,102.71) and (149.06,100.51) .. (149.06,97.8) -- cycle ;
\draw  [color={rgb, 255:red, 255; green, 0; blue, 0 }  ,draw opacity=1 ][fill={rgb, 255:red, 255; green, 0; blue, 0 }  ,fill opacity=1 ] (89.07,54.91) .. controls (89.07,52.2) and (92.02,50) .. (95.66,50) .. controls (99.31,50) and (102.26,52.2) .. (102.26,54.91) .. controls (102.26,57.62) and (99.31,59.81) .. (95.66,59.81) .. controls (92.02,59.81) and (89.07,57.62) .. (89.07,54.91) -- cycle ;
\draw  [color={rgb, 255:red, 255; green, 0; blue, 0 }  ,draw opacity=1 ][fill={rgb, 255:red, 255; green, 0; blue, 0 }  ,fill opacity=1 ] (10,94.69) .. controls (10,91.98) and (12.95,89.78) .. (16.6,89.78) .. controls (20.24,89.78) and (23.19,91.98) .. (23.19,94.69) .. controls (23.19,97.39) and (20.24,99.59) .. (16.6,99.59) .. controls (12.95,99.59) and (10,97.39) .. (10,94.69) -- cycle ;
\draw  [color={rgb, 255:red, 255; green, 0; blue, 0 }  ,draw opacity=1 ][fill={rgb, 255:red, 255; green, 0; blue, 0 }  ,fill opacity=1 ] (10,182.2) .. controls (10,179.49) and (12.95,177.29) .. (16.6,177.29) .. controls (20.24,177.29) and (23.19,179.49) .. (23.19,182.2) .. controls (23.19,184.91) and (20.24,187.1) .. (16.6,187.1) .. controls (12.95,187.1) and (10,184.91) .. (10,182.2) -- cycle ;
\draw  [color={rgb, 255:red, 255; green, 0; blue, 0 }  ,draw opacity=1 ][fill={rgb, 255:red, 255; green, 0; blue, 0 }  ,fill opacity=1 ] (69.99,137.58) .. controls (69.99,134.87) and (72.94,132.67) .. (76.59,132.67) .. controls (80.23,132.67) and (83.18,134.87) .. (83.18,137.58) .. controls (83.18,140.29) and (80.23,142.49) .. (76.59,142.49) .. controls (72.94,142.49) and (69.99,140.29) .. (69.99,137.58) -- cycle ;
\draw  [color={rgb, 255:red, 255; green, 0; blue, 0 }  ,draw opacity=1 ][fill={rgb, 255:red, 255; green, 0; blue, 0 }  ,fill opacity=1 ] (69.99,225.09) .. controls (69.99,222.38) and (72.94,220.19) .. (76.59,220.19) .. controls (80.23,220.19) and (83.18,222.38) .. (83.18,225.09) .. controls (83.18,227.8) and (80.23,230) .. (76.59,230) .. controls (72.94,230) and (69.99,227.8) .. (69.99,225.09) -- cycle ;
\draw  [color={rgb, 255:red, 255; green, 0; blue, 0 }  ,draw opacity=1 ][fill={rgb, 255:red, 255; green, 0; blue, 0 }  ,fill opacity=1 ] (149.06,185.31) .. controls (149.06,182.61) and (152.01,180.41) .. (155.65,180.41) .. controls (159.3,180.41) and (162.25,182.61) .. (162.25,185.31) .. controls (162.25,188.02) and (159.3,190.22) .. (155.65,190.22) .. controls (152.01,190.22) and (149.06,188.02) .. (149.06,185.31) -- cycle ;

\end{tikzpicture}
     }
    \caption{Left: the triangulation on every side of the cube $B$ inside $\Sigma=\trop(X)$. Every triangle has normalized volume $1$, so that the induced model $\mathcal{X}\to U\subset \Spec(\CC[t])$ is regular. The induced structure of an affine integral manifold comes from the order complex $\Delta(B)$ of $B$. This is topologically the same as the barycentric subdivision of $B$. 
    \\
    Middle: The local picture at one corner of the triangulated cube. 
    \\
    Right: the order complex $\Delta(B)$ around a corner vertex of $B$. The brown vertices give the singular points, the green lines give the edges of the cube. The shaded blue regions give the cells in $W_{v}$ for $v$ the corner vertex. The orange path indicates a nontrivial path in $\Delta(B)\backslash \Gamma$. The monodromy of integral affine structures is a unipotent matrix, as we will see in Example \ref{exa:CalculationMonodromyIntegralAffineStructures}  }
    \label{fig:TriangulationOfK3Surface}
\end{figure}

\begin{remark}
    The choice of singular points is not compatible with base changes of the form $\CC[t]\to \CC[t^{1/m}]$. Namely, in order to obtain a regular model, one has to subdivide the given polyhedral structure $\Sigma$ to a finer polyhedral structure $\Sigma'$ where the triangles all have normalized volume $1$. 
    Then one takes a vertex adjacent to a corner vertex. This lies closer to the corner vertex than the original vertex. In a sense, the singular points are supposed to be points that lie infinitesimally close to the corner vertex.  From a Berkovich-theoretic perspective, it would be better to assign the corresponding tangent point, which we can view as an element of the corresponding Zariski-Riemann space.     
\end{remark}

We now interpret the linear spaces in terms of maps of homology groups. To calculate with these explicitly, we use angle spaces. We view the vector spaces arising from the integral affine manifold interpretation as the \emph{universal covering spaces} of the corresponding angle spaces.  

To compare the two approaches, consider a morphism $(\tau_{1},\sigma_{1})\to (\tau_{2},\sigma_{2})$ arising from the tropical poset $\Pscr(\overline{\Sigma})$. This induces a map 
\begin{equation*}
    Z^{\sigma_{1},\log}_{\tau_{1}}\to Z^{\sigma_{2},\log}_{\tau_{2}}.
\end{equation*}
This induces 
\begin{equation*}
    H_{1}(Z^{\sigma_{1},\log}_{\tau_{1}},\RR)\to H_{1}(Z^{\sigma_{2},\log}_{\tau_{2}},\RR).  
\end{equation*}
If $\sigma_{i}=(0)$, then the corresponding cell lies in the interior of the tropicalization, and the comparison map is an injection. Moreover, the map of homology groups is an injection as well. 

To deal with the spaces at infinity, we will define certain local cofibrant replacements. Let $\tau$ be a vertex on the cube. We will assume that it does not lie in the interior of a $2$-plane of the cube. In other words, it is either one of the $8$ corner vertices, or one of the $24$ edge vertices. There is a unique unbounded ray $\tau'$ adjacent to $\tau$. This corresponds to the generator $\eta_{\tau}$ in Definition \ref{def:IntegralAffineStructure}. We write $\sigma=\sigma_{\tau'}$ for the corresponding ray in the recession fan of $\Sigma$. Note that we have an inclusion 
\begin{equation*}
    Z^{\log}_{\tau'}\to Z^{\log}_{\tau}.
\end{equation*}
Similarly, we have an inclusion 
\begin{equation*}
    Z^{\log}_{\tau'}\to Z^{\sigma,\log}_{\tau',c}.
\end{equation*}
Here we took the cofibrant replacement from Definition \ref{def:CofibrantReplacements2} to make  sure the comparison map is injective. We now define  $Z^{\log}_{\tau\to \tau'}$ to be the colimit of 
\begin{equation*}
    Z^{\log}_{\tau}\leftarrow Z^{\log}_{\tau'}\to Z^{\tau',\log}_{\tau,c}.
\end{equation*}
This can be seen as partially sealing the Kato-Nakayama space $Z^{\log}_{\tau}$. That is, in the homotopy colimit of the full Kato-Nakayama functor, we take colimits of this type over all cells, but now we are only taking the colimit over two maps. Note that we  automatically have a map
\begin{equation*}
    Z^{\log}_{\tau}\to Z^{\log}_{\tau\to \tau'},
\end{equation*}
and thus a map of homology groups
\begin{equation*}
    H_{1}(Z^{\log}_{\tau},\RR)\to H_{1}(Z^{\log}_{\tau\to \tau'},\RR).
\end{equation*}
We now describe all of the homology groups and maps explicitly in the following lemma. 
\begin{lemma}\label{lem:HomologyMapsFermiApery}
    Suppose that $\tau$ is a maximal $2$-dimensional cell. Then 
    \begin{equation*}
        H_{1}(Z^{\log}_{\tau},\RR)\simeq \RR^{2}.
    \end{equation*}
    Let $\tau$ be a vertex on the boundary of the cube with unbounded rays $\tau'\in\Sigma$ and $\sigma=\sigma_{\tau'}\in\Delta$. We then have
    \begin{align*}
        H_{1}(Z^{\log}_{\tau},\RR) \simeq \RR^{3} \quad \text{and} \quad
        H_{1}(Z^{\log}_{\tau\to \tau'},\RR) \simeq \RR^{2}.
    \end{align*}
    The maps
    \begin{equation*}
        \RR^{3}\simeq H_{1}(Z^{\log}_{\tau},\RR)\to H_{1}(Z^{\log}_{\tau\to \tau'},\RR)\simeq \RR^{2}
    \end{equation*}
    can be identified with $N_{\RR}
    \to N_{\RR} / {\rm Span}_{\RR}(\sigma)$.
\end{lemma}
\begin{proof}
    This immediately follows from our description of the corresponding Kato-Nakayama spaces and the quotient maps, see Section \ref{sec:InitialDegenerationsKNSpaces}. 
\end{proof}

\begin{definition}[Kato-Nakayama local systems on integral affine manifolds]

    We associate a two-dimensional vector space to every vertex $\tau \in \Delta(B)\backslash \Gamma$. Note that we can identify a vertex $\tau$ with an element of $\Pscr(B)$ and thus with an element of $\Pscr(\overline{\Sigma})$. Suppose that $\tau$ lies in the interior of a maximal-dimensional cell. Then we associate $H_{1}(Z^{\log}_{\tau},\RR)$. If $\tau$ is one of the vertices on the boundary of the cube, then we assign $H_{1}(Z^{\log}_{\tau\to \tau'},\RR)$.  

    We now define transport morphisms between different adjacent vertices. Here adjacent means that they are connected by an edge in $\Delta(B)$. We will use Lemma \ref{lem:HomologyMapsFermiApery} for this. If two vertices are connected by an edge that lies in the same affine linear space, then it is induced from the natural isomorphism between Kato-Nakayama spaces. Suppose now that $\tau_{1}$ is a vertex in the interior of a cell of the cube and let $\tau_{2}$ be a vertex on the boundary. We have an injection 
    \begin{equation*}
        H_{1}(Z^{\log}_{\tau_{1}},\RR)\to H_{1}(Z^{\log}_{\tau_{2}},\RR)
    \end{equation*}
    of a $2$-dimensional vector space into a $3$-dimensional vector space. We take the composite of the following two maps
    \begin{equation*}
         H_{1}(Z^{\log}_{\tau_{1}},\RR)\to H_{1}(Z^{\log}_{\tau_{2}},\RR)\to  H_{1}(Z^{\log}_{\tau_{2}\to \tau'_{2}},\RR).
    \end{equation*}
    Note that this is an isomorphism.
\end{definition}

\begin{corollary}
    The transport maps of vector spaces $$H_{1}(Z^{\log}_{\tau_{1}},\RR)\to H_{1}(Z^{\log}_{\tau_{2}},\RR)$$ for vertices in the interior of a $2$-plane, and the transport maps 
    $$H_{1}(Z^{\log}_{\tau_{1}},\RR)\to H_{1}(Z^{\log}_{\tau_{2}\to \tau'_{2}},\RR)$$
    for an interior vertex $\tau_{1}$ to a boundary vertex $\tau_{2}$
    are compatible with the comparison maps of integral affine manifolds.
\end{corollary}
\begin{proof}
    This follows from the fact that the quotient map of Kato-Nakayama spaces is induced from $N_{\RR}\to N_{\RR}/{\rm Span}(\sigma)$, see Section \ref{sec:InitialDegenerationsKNSpaces}. 
\end{proof}

From this we see that the monodromy of integral affine structures can be interpreted purely in terms of homology comparison maps. If we wish to calculate with these explicitly, note that we can use our angle maps.

\begin{example}
\label{exa:CalculationMonodromyIntegralAffineStructures}
We calculate the monodromy around a singular point explicitly here. Consider the set-up in Figure \ref{fig:TriangulationOfK3Surface}. We will take the orange path in that picture. We start at the vertex on top. To make things explicit, we will work at the vertex $v=(-1,-1,-1)$. The initial ideal at this point is 
\[
    {\rm in}_{w}(f) = 1 + x + y + z,
\]
so we can work with the corresponding fan. Throughout this explanation, we will often refer to Figure~\ref{fig:TriangulationOfK3Surface}. 

\begin{enumerate}
    \item \underline{From $1+y\rightarrow 1+z$ through the edge}.
    
    We start on the top, where the initial ideal is $1+y$. The local Kato-Nakayama space $Z^{\log}_{\tau}$ is then $(S^{1})^{2}$, with coordinates $\theta_{x}$ and $\theta_{z}$. We then move to the edge, where we have to change coordinates. Note that the initial ideal at the edge is $1+y+z$. The homology group $H_{1}(Z^{\log}_{\tau},\RR)$ at this edge is three-dimensional, see Figure \ref{fig:AngleSetInitialDegeneration}, where the angle set is given by a toblerone shape. Note that the corresponding Kato-Nakayama space is homotopy equivalent  
    to $\CC^{\times}\times V(1+y+z)$. 
    
    We now move from this set to $Z^{\log}_{\tau\to\tau'}$ using the vector $(0,-1,-1)$. Note that under the map $N_{\RR}\to N_{\RR}/{\rm Span}(\eta_{\tau})$, we have that $(0,0,1)$ is identified with $(0,-1,0)$. In terms of angles, we find that the cycle determined by $\theta_{z}$ is identified with minus the cycle determined by $\theta_{y}$. We will use these cycles on the maximal cell where $1+z$ is the initial ideal.

    \item \underline{From $1+z\rightarrow 1+y$ using the vertex}.
    
    We now complete the path using the corner vertex. We have coordinates $x$ and $y$ for the angles, as the initial ideal is generated by $1+z$. At the corner vertex, we have the initial ideal $1+x+y+z$, and the corresponding homology group $H_{1}(Z^{\log}_{\tau},\RR)$ is $3$-dimensional. In Figure \ref{fig:AngleSetInitialDegeneration}, this is the space corresponding to the $2$-skeleton of the cube. The angle set for the $2$-torus injects as a plane. We now consider the corresponding quotient map 
    \begin{equation*}
        H_{1}(Z^{\log}_{\tau},\RR)\to H_{1}(Z^{\log}_{\tau\to \tau'},\RR).
    \end{equation*}
    Here $\sigma$ corresponds to the vector $\eta_{\tau}=(-1,-1,-1)$. Note that we are in fact first moving to the locus in the tropicalization where the initial ideal is $V(x+y+z)$, and we are then taking the quotient by the corresponding torus. The corresponding quotient map identifies $(0,1,0)$ with $(-1,0,-1)$. In other words, the angle cycle generated by $\theta_{y}$ now becomes a sum of two cycles with respect to the coordinates $\theta_{x}$ and $\theta_{z}$. Namely, it becomes the cycle $-\gamma_{x}-\gamma_{z}$.

    \item \underline{Composing the two maps}. 
    
    We now compose these two maps and find that the cycle $\gamma_{z}$ generated by $z$ in $(S^{1})^{2}$ is sent to the cycle $\gamma_{z}-\gamma_{x}$ (note here that we first obtain a minus sign from the identification through $y$, which cancels with the minus sign of the second transformation to give $\gamma_{z}$). In other words, the monodromy matrix and its inverse are given by 
    \begin{equation*}
        A_{\gamma}=\begin{pmatrix}
            1 & -1\\
            0 & 1
        \end{pmatrix}
        \quad
        \text{and}
        \quad
        A_{\gamma^{-1}}=\begin{pmatrix}
            1 & 1\\
            0 & 1
        \end{pmatrix},
        \end{equation*}
        which is the standard monodromy matrix arising from a focus-focus singularity.
        \end{enumerate}
\end{example}

\begin{remark}
    We now note how the calculation in Example~\ref{exa:CalculationMonodromyIntegralAffineStructures} shows how this monodromy also occurs over the $p$-adic numbers. 
    We first replace $t$ with $p$ in the Ap\'{e}ry-Fermi pencil, and work over the field $\QQ_{p}$ of $p$-adic numbers. 
    Note that the tools of tropical geometry work equally well over this field, and we obtain a natural model $\mathcal{X}\to \Spec(\ZZ_{p})$ this time. 
    Rather than taking Kato-Nakayama spaces, one should then work with the shape of the corresponding logarithmic schemes with respect to the prime-to-$p$ \'{e}tale or Kummer site. 
    In particular, the $Z^{\sigma,\log}_{\tau}$ should be interpreted as shapes of logarithmic schemes with respect to these topologies. 
    We can also use the theory of rigid analytic spaces or Berkovich analytic spaces here. 
    From this point of view, each closed subset $Z_{\tau}\subset \mathcal{X}_{s}$ defines a natural open analytic space $Z^{\an}_{\tau}$ by anticontinuity of the reduction map $X^{\an}\to \mathcal{X}_{s}$. 
    One thus obtains a cover of $X^{\an}$, and one expects the homotopy colimit of this cover to again return the \'{e}tale shape of $X^{\an}$, as in Theorem \ref{thm:KNFunctorHomotopyColimit3}. 
    We will continue with the logarithmic spaces in this remark.

    Let $G=(\hat{\ZZ}')^{2}$. Here $\hat{\ZZ}'$ is the inverse limit of all $\ZZ/m\ZZ$ for $m$ coprime to $p$.
    The corresponding $Z^{\log}_{\tau}$ for $\tau$ maximal $2$-dimensional cells in $\trop(X)$ should be $K(G,1)$-spaces in the \'{e}tale topology on the corresponding logarithmic scheme. 
    Note that in the complex case, we have a $2$-dimensional torus, which is a $K(\ZZ^{2},1)$. 
    Unfortunately, circles and real tori do not occur naturally over the $p$-adic numbers, hence we work with $K(G,1)$-spaces. 
    These can be seen as the $p$-adic analogues of the $2$-dimensional real tori found over the complex numbers. 

    Now the monodromy in Example \ref{exa:CalculationMonodromyIntegralAffineStructures} can be calculated purely in terms of the framework of  \'{e}tale cohomology on logarithmic schemes. Indeed, we are simply taking the push-forwards of certain cycles under natural comparison morphisms. Moreover, these cycles can be reconstructed  exactly in terms of Kummer coverings. 
The monodromy of integral affine structures thus also naturally occurs over the $p$-adic numbers.
\end{remark}

\begin{example}
    To construct the torus fibration using our setup, one has to be slightly careful with the spaces at the corner vertices. Note that one cofibrant replacement of the diagram can be found in Figure \ref{fig:CofibrantReplacementK3}.
       \begin{figure}[ht]
\begin{center}
    \scalebox{0.3}{
\definecolor{xdxdff}{rgb}{0.49019607843137253,0.49019607843137253,1}
\definecolor{ududff}{rgb}{0.30196078431372547,0.30196078431372547,1}
\begin{tikzpicture}[line cap=round,line join=round,x=1cm,y=1cm]
\clip(-2,1) rectangle (23,18);
\draw [line width=2pt] (4,3)-- (4,15);
\draw [line width=2pt] (4,3)-- (16,3);
\draw [line width=2pt] (16,3)-- (16,15);
\draw [line width=2pt] (4,15)-- (16,15);
\draw [line width=2pt] (10,9)-- (4,3);
\draw [line width=2pt] (10,9)-- (4,15);
\draw [line width=2pt] (10,9)-- (4,9);
\draw [line width=2pt] (10,9)-- (10,15);
\draw [line width=2pt] (10,9)-- (16,15);
\draw [line width=2pt] (10,9)-- (16,9);
\draw [line width=2pt] (10,9)-- (10,3);
\draw [line width=2pt] (10,9)-- (16,3);
\draw (5.7,9) node[anchor=north west] [xscale=3, yscale=3] {$(S^{1})^{2}$};
\draw (16,17) node[anchor=north west] [xscale=3, yscale=3]   {{$\mathbb{C}^{2}$}};
\draw (8,17) node[anchor=north west] [xscale=3, yscale=3]   {{$\mathbb{C}\times S^{1}$}};
\draw (16.4,10) node[anchor=north west] [xscale=3, yscale=3]   {{{$S^{1}\times \mathbb{C}$}}};
\draw (16,3) node[anchor=north west] [xscale=3, yscale=3]   {{$\text{P}\times \mathbb{C}$}};
\draw (1,17) node[anchor=north west] [xscale=3, yscale=3]   {{$\mathbb{C}\times S^{1}$}};
\draw (0,10) node[anchor=north west] [xscale=3, yscale=3]   {{$(S^{1})^{2}$}};
\draw (1.6,3) node[anchor=north west] [xscale=3, yscale=3]   {{$\text{C}$}};
\draw (8,2.5) node[anchor=north west] [xscale=3, yscale=3]   {{$\text{P}\times S^{1}$}};
\begin{scriptsize}
\draw [fill=ududff] (4,3) circle (7pt);
\draw [fill=ududff] (4,15) circle (7pt);
\draw [fill=ududff] (16,3) circle (7pt);
\draw [fill=ududff] (16,15) circle (7pt);
\draw [fill=ududff] (10,9) circle (7pt);
\draw [fill=xdxdff] (4,9) circle (7pt);
\draw [fill=xdxdff] (10,15) circle (7pt);
\draw [fill=xdxdff] (16,9) circle (7pt);
\draw [fill=xdxdff] (10,3) circle (7pt);
\end{scriptsize}
\end{tikzpicture}     }
    \caption{\label{fig:CofibrantReplacementK3}A cofibrant replacement of the local diagram of spaces for the Ap\'{e}ry-Fermi family, at one of the pieces attached to $\trop(X)$. That is, for the given subposet, we have depicted the order complex, together with the space over it. Here ${\rm P}$ is the Kato-Nakayama space associated to the compactification $\mathbb{P}^{1}\backslash\{0,1,\infty\}\subset \mathbb{P}^{1}$, and ${\rm C}$ is the Kato-Nakayama space associated to $\mathbb{P}^{2}\backslash V\big(XYZ(-X-Y-Z)\big) \subset \mathbb{P}^{2}$. The arrows in the diagram are all pointed outward from the central vertex. Note that they are all cofibrations.}
    \end{center}       
\end{figure}
   Here we write ${\rm P}$ for the Kato-Nakayama space associated to $\mathbb{P}^{1}\backslash\{0,1,\infty\}$, and ${\rm C}$ for the Kato-Nakayama space associated to $\mathbb{P}^{2}\backslash V\big(XYZ(-X-Y-Z)\big)$.
    Note that the three pieces attached to the corner vertex give three copies of this diagram, which are glued along the bottom line of the diagram. It follows from this that the local spaces are in fact homotopy equivalent to a two-dimensional torus. This is then compatible with each of the induced comparison maps. 

   This can be interpreted as follows. One can avoid the spaces $\mathbb{P}^{2}\backslash D$ for $D = V\big(XYZ(-X-Y-Z)\big)$ at the corner vertices. That is, they can be replaced by $2$-dimensional real tori. To see this, one  starts with the diagram of semistable models $\mathcal{X}\to U\subset \Spec(\CC[t])$ associated to the rational tropical manifold. The local Kato-Nakayama spaces associated to all the components are then homotopy equivalent to tori, so we obtain that a K3-surface is a homotopy colimit of $2$-tori over a $2$-sphere. This interpretation makes it hard to work with the comparison maps however, which is why we worked with the extended tropical diagram here. In general, this seems to lead to a description of Calabi-Yau varieties as homotopy colimits of $k$-tori over a $k$-dimensional sphere $S^{k}$. To work explicitly with the corresponding topological model, one has to know the explicit comparison maps. This is where tropical geometry can help, as we saw above.

 \end{example}

\bibliographystyle{acm}
\bibliography{references}
\end{document}